\documentclass{amsart}
\usepackage{amsmath}          
\usepackage{amsthm}           
\usepackage{amssymb}          
\usepackage{graphicx}
\usepackage[section]{placeins}
\usepackage{comment}
\usepackage{listings}
\usepackage{xcolor}
\usepackage{silence}
\usepackage{longtable}
\usepackage[colorlinks=true,linkcolor=blue,citecolor=blue,urlcolor=blue,
            breaklinks=true]{hyperref}
\usepackage{xurl}
\usepackage{pifont}
\usepackage{xcolor}
\newcommand{\nomark}{\colorbox{red}{\textcolor{white}{\ding{55}}}}
\newcommand{\yesmark}{\colorbox{green!60!black}{\textcolor{white}{\ding{51}}}}

\newtheorem{theorem}{Theorem}
\newtheorem{lemma}{Lemma}
\newtheorem{corollary}{Corollary}

\theoremstyle{definition} 
\newtheorem{definition}{Definition}
\newtheorem*{remark}{Remark}
\newtheorem*{example}{Example}
\newtheorem*{examples}{Examples}
\newtheorem{convention}{Convention}

\def\Q{{\mathbb Q}}
\def\R{{\mathbb R}}
\def\G{{\mathbb G}}
\def\HH{{\mathbb H}}  
\def\Z{{\mathbb Z}}

\def\d{{\bf d}}

\def\mod{\mbox{\ mod\ }}

\def\lcm{{\rm lcm}}

\lstdefinelanguage{Sage}[]{Python}{
  morekeywords={sage, var, solve, factor, simplify, QQ, ZZ, RR, CC,
                EllipticCurve, PolynomialRing, matrix, vector, srange},
  sensitive=true
}

\def\vector#1#2#3{\left(
\begin{array}{l} 
#1\\ 
#2\\ 
#3 
\end{array}\right)}

\def\matrix#1#2#3#4#5#6#7#8#9{{
 \left( \begin{array}{ccc}
#1 & #2 & #3 \\
#4 & #5 & #6 \\
#7 & #8 & #9 
\end{array} \right)}}

\title{Triangle Tiling: The case $3\alpha +2\beta = \pi$}          \author{Michael Beeson}
\thanks{ Thanks to Vico Bonfioli for pointing out errors in the  previous version of this paper.  
This version corrects all known errors.
}%
\thanks{Claude Opus (Anthropic) was used to check proofs and prove lemmas.  All text has been
edited by the author, who remains responsible for the content.}
\date{\today}          

\begin{document}

\begin{abstract}  
An $N$-tiling of triangle $ABC$ by triangle $T$ (the ``tile'') 
is a way of writing $ABC$ as a union of $N$ copies of $T$ overlapping only at their boundaries. 
Let the tile $T$ have angles $(\alpha,\beta,\gamma)$, and 
sides $(a,b,c)$. 
This paper takes up the case when  $3\alpha + 2\beta = \pi$.
Then there are (as was already known)
 exactly five possible shapes of $ABC$ not similar to $T$:  either 
$ABC$ is isosceles with base angles $\alpha$, $\beta$, or $\alpha+\beta$,
or the angles of $ABC$ are  $(2\alpha,\beta,\alpha+\beta)$, or
the angles of $ABC$ are $(2\alpha, \alpha, 2\beta)$.  
We prove that the tile is rational (the ratios of its sides are rational).
We then exhibit several families of previously unknown tilings.
We also exhibit, for each shape, a Diophantine equation in $N$ and the 
parameters of the tile that must have a solution
whenever there is a tiling of some triangle of that shape. 
 By means of these Diophantine equations,
some conclusions about the possible values of $N$ are drawn. 
We prove, for example, that there is no $N$-tiling with $N$ prime 
when $3\alpha + 2\beta = \pi$, for three of the five shapes;
that question remains open for the shape $(2\alpha,\beta,\alpha+\beta)$ and 
for isosceles triangles with base angles $\beta$. 
These equations also imply that for each $N$, there is a finite
set of possibilities for the tile $(a,b,c)$ and the triangle $ABC$. 
(Usually, but not always, there is just one possible tile.) 
For each of the five shapes, we give a construction method that 
works under certain sufficient conditions, but there remain values of $N$ 
that satisfy our necessary conditions and for which we do not know 
whether a tiling exists.  All results have been formalized in Lean.
\end{abstract}

\maketitle

\section{Introduction}  This paper is part of a series of papers on triangle tiling, continuing a research program begun by Laczkovich \cite{laczkovich1995}.
 An $N$-tiling of triangle $ABC$ by triangle $T$ (the ``tile'') 
is a way of writing $ABC$ as a union of $N$ copies of $T$ overlapping only at their boundaries.
The general aim of this research program is to understand the nature of triangle tilings, which can be amazingly complex.    We can 
exhibit quite a few families of triangle tilings, some of which are very familiar, 
but at least six new families of tilings have been discovered in the course
of our work,  
five of which are described in this paper.
  It would be too much to ask that every triangle tiling belongs to one of these families, as sometimes it is possible to rearrange some of the tiles
within a given tiling, and also there are systematic ways to combine tilings.
Our aim instead is to completely classify
 the triples ($ABC, T, N)$  such that there
exists a tiling of $ABC$ by $T$ using $N$ tiles.   In particular, this
classification should 
enable us to answer more specific questions, such as,  for which $N$ does there exist
a triangle $ABC$ and a tile $T$ and an $N$-tiling of $ABC$ by $T$? 

Cases of specific $N$ can be  quite interesting.
The question that first sparked our interest was whether any triangle can be $7$-tiled.
We gave a long Euclidean-style proof that no 7-tiling exists, but it was clear that a similar 
proof for 11-tilings might be a thousand pages long.   Now we know that there is 
also no 11-tiling.  A self-contained proof of those two theorems is
presented in \cite{beeson-noseven}. 
   But there {\em does} exist a 28-tiling, which we think is a new discovery.   We generalize
 the 28-tiling, showing it to be the simplest member of a new 
 family of tilings of triangles of the shape $(2\alpha, \beta,\alpha+\beta)$.
   These are the ``triquadratic tilings'';  when a tiling of that shape triangle exists, 
 the ``tiling equation'' 
 $$M^2 + N = 2K^2$$
 has a solution in integers $(K,M)$.   Conversely,  if that equation has a solution in which
 $K$ divides $N$, or equivalently $K$ 
 divides $M^2$, then we can construct a tiling.  
 Each such solution determines the tile of the corresponding 
 triquadratic tiling: the tile must be similar to the triangle 
 with sides $a=M$, $c=K$, and $b = K-M^2/K$.  All three sides of
 that tile are integers, and the tile then satisfies the condition
 $3\alpha+2\beta = \pi$.
 
 In previous work of others (\cite{laczkovich1995,snover1991})
 and the present author (\cite{beeson-noseven}) the questions of triangle tiling
 have been successfully divided into a small finite number of cases
 according to the shapes of the tile $(a,b,c)$, whose angles
 are $(\alpha,\beta,\gamma)$, and the triangle $ABC$.  In this paper
 we take up the case $3\alpha + 2 \beta = \pi$ (although a few
lemmas are proved also for the case $\gamma = 2\pi/3$, when the 
same proof applies).   In case $3\alpha + 2\beta =\pi$, 
 according to \cite{laczkovich1995},  there are just five possible 
 shapes of $ABC$: 
 \begin{itemize} 
\item $ABC$ is isosceles with base angles $\alpha$, $\beta$, or $\alpha+\beta$,
\item or the angles of $ABC$ are  $(2\alpha,\beta,\alpha+\beta)$, or
\item the angles of $ABC$ are $(2\alpha, \alpha, 2\beta)$.
\end{itemize}
  In this paper, we give new tilings in each of these cases.
In each case, we give a necessary condition on $N$ depending on 
the angles of the tile, such that if there is an $N$-tiling, the
condition is satisfied.  We also give further conditions sufficient
for the construction of tilings.  In each case there are some unknown 
values of $N$ that fall between the necessary and the sufficient conditions.
 We now state our
results more precisely.
 
 Our results for the case when $ABC$ has angles $(2\alpha,\beta,\alpha+\beta)$
 are as follows:
\smallskip

(i)  If the tiling equation $M^2 + N = 2K^2$ has a solution $(K,M)$
in positive integers such that $K$ divides $M^2$ and $M^2 < N$, then 
there is a triquadratic $N$-tiling of $ABC$.  The exact description and 
pictures of these new tilings are given in the paper.  The smallest $N$
allowing a triquadratic tiling is $N=28$.  
\smallskip

(ii)  The tiling equation has, for each $N$, at most finitely
many solutions $(K,M)$.
  Each solution with $K$ dividing $N$ determines a possible tile shape,
similar to the triangle
with sides $a=M$, $c=K$, and $b = K-M^2/K$, of which all three sides are 
integers.   Each of these finitely many tiles may be used for 
a triquadratic tiling.
\smallskip

(iii) Except for the ``quadratic tilings'' (that always exist
when $N$ is a square),   
any $N$-tiling of $ABC$ must use a tile similar to one of 
the tiles in this finite set.  In particular, if the tiling equation
has no solution, and $N$ is not a square, then there do not exist
any $N$-tilings of $ABC$.
\smallskip

(iv) If the tiling equation has solutions, then 
$N$ is a square times a product of distinct primes 
of the form $8n \pm 1$, or 2.  
\smallskip

These results do not provide a complete classification of the 
triples $(ABC, N, T)$ such that $ABC$ can be $N$-tiled by 
the tile $T$, in case $ABC$ has angles $(2\alpha, \beta,\alpha+\beta)$.
In case the tiling equation has no solution, then unless $N$
is a square, there are no $N$-tilings of $ABC$; and if $N$ is a square, 
there are no $N$-tilings except the quadratic tilings by a 
tile similar to $ABC$, which always exist when $N$ is a square.
In case the tiling equation has a solution in which $K \mid M^2$, we can 
construct a tiling.  But in case it has a solution, and no solution 
satisfies $K \mid M^2$, then in general we get no information about the 
existence of a tiling.  That is what happens for  $N = 7$, $14$, $17$, $23$, which are
ruled out by separate arguments in this paper; for $N = 31$, $41$, which are prime
and hence ruled out by \cite{beeson-rosu2026};  and for $63$, which will be discussed
in Section~\ref{section:search}.
Note that we can always multiply the size of the triangle $ABC$ by $K$.  
That replaces $(K,M)$ by $(K^\prime, M^\prime) = (K^2, MK)$, which is again 
a solution, of 
$$ N^\prime \ = \ 2(K^\prime)^2 - (M^\prime)^2 \ = \ K^2(2K^2 - M^2) 
   \ = \ K^2 N, $$
and this time the condition $K^\prime \mid (M^\prime)^2$ holds 
automatically, since $(M^\prime)^2 = M^2K^2$.  So we do get a tiling of 
that shape, just for a larger $N$ --- that is how one gets a 28-tiling in 
place of a 7-tiling, since $28 = 2^2 \cdot 7$.

Now consider the case when $ABC$ has angles $(2\alpha,\alpha,2\beta)$.
Then there is another 
equation, the ``second tiling equation'', that governs the 
existence of tilings.  That equation is 
$$ N = M^2 \frac {(2-s^2)(3-s^2)}{(1-s)^2(2+s)^2} $$
where $s =2 \sin(\alpha/2)$ determines the angle $\alpha$ and 
the ratio of sides $a/c = s$, 
and the tile $(a,b,c)$ is given by $b = c-a^2/c$.   Here $M$
has to be a positive integer and $s$ a rational between 0 and 1.
 There is a 
tiling of $ABC$ for each solution of the second tiling equation, 
and if $N$ is not a perfect square, there are tilings of $ABC$ only 
if the second tiling equation has a solution.  Some of these tilings
belong to the new family we call the four-component tilings; pictures
are provided.  

The second tiling equation has infinitely many solutions, since we 
can start with the tile $(a,b,c)$ such that $a$ divides $c$ and 
$c$ divides $a^2$, and produce a tiling.  Conversely, we can determine
by an algorithm whether the tiling equation is solvable for a particular $N$,
although we do not have a concise number-theoretic characterization of
those $N$ for which it is solvable.  We used that algorithm to compute
a table of solutions for $N$ up to several thousand.

 We now sketch how the particular shapes $(2\alpha, \beta, \alpha+\beta)$
and $(2\alpha,\alpha,2\beta)$
arise naturally in this subject.  Let the sides of the tile 
opposite angles $\alpha$, $\beta$, and $\gamma$ respectively be 
$a$, $b$, and $c$. 
Between the three numbers $\alpha$, $\beta$, and $\gamma$,  we have two linear relations.
The first one is $\alpha + \beta + \gamma = \pi$, because they are the three angles of a 
triangle.  The second relation arises from the ``vertex splitting'' at the vertices
of the tiled triangle $ABC$.  Consider all the tiles that share any of the three vertices,
some (or none) have $\alpha$ angles at the vertex, some (or none) have $\beta$ angles, and some (or none) have 
$\gamma$ angles.   Thus for some nonnegative integers $P$, $Q$, and $R$, we have 
$$ P \alpha + Q \beta + R \gamma = \pi.$$
If $P=Q=R=1$ then $ABC$ is similar to $T$; that case was treated in \cite{golomb, snover1991, laczkovich1995}. In this paper we assume $T$ is not similar to $ABC$.  Hence the two 
linear relations are independent relations;  substituting $\gamma = \pi-\alpha-\beta$ we 
have a non-trivial linear relation between $\alpha$, $\beta$, and $\pi$.  
If we could get a third independent linear relation between the angles, we could solve 
for $\alpha$, $\beta$ and $\gamma$.  At each vertex of the tiling,  some angles 
add up to either $\pi$ (for a boundary vertex, or a ``non-strict'' interior vertex,
i.e. one lying on a side of another tile), or to $2\pi$, for a strict interior vertex.
If $\alpha$ is 
not a rational multiple of $\beta$, 
it will never be the case that the same angle can be composed of some $\alpha$ angles and 
also of some $\beta$ angles.   Hence the possibilities for angles meeting at interior
vertices can be controlled to some extent.   This analysis was carried out 
by Laczkovich \cite{laczkovich1995}, who showed in this way that when the tile $T$ 
is not a right triangle, nor similar to $ABC$,  there
are only a few possibilities:  either $ABC$ is equilateral and $\beta = \pi/3$,
or $\gamma = 2\alpha$ (which is another way of saying $3\alpha+\beta=\pi$)
 and $ABC$ is isosceles with base angles $\alpha$, or 
$\gamma = 2\pi/3$, or $3\alpha + 2\beta = \pi$.

It follows from work of Laczkovich \cite{laczkovich1995} (see 
Lemma~\ref{lemma:angles} below for details) 
 that if $3\alpha + 2\beta = \pi$,  then
 $\alpha$ and $\beta$ are not rational multiples of 
$\pi$ and either $ABC$ is isosceles (in which case its base angles are $\alpha$, $\beta$,
or $\alpha + \beta$), or the angles of $ABC$ are $(2\alpha,\alpha,2\beta)$, or 
the angles of $ABC$ are $(2\alpha, \beta, \alpha+\beta)$. All 
five cases  are taken up in this paper.  From $3\alpha + 2\beta = \pi$
we have $\gamma = 2\alpha + \beta$ and $\alpha < \pi/3$ and $\beta < \pi/2$.

We next explain the meaning of the integer $M$ in the first and second
tiling equations. 
Our work in this paper begins with the introduction of colored tilings,
in which the tiles are colored black and white alternately,
 and a ``coloring equation''
is derived that involves the difference $M$ between the numbers of black and white tiles. (Figure~\ref{figure:coloring} illustrates the kind 
of coloring in question.) 
The coloring equation is then combined with the ``area equation'',
that equates the area of triangle $ABC$, as computed from the sine of one 
angle and the lengths of two sides,  to $N$ times the area of the tile.
Eliminating one variable from these two equations, we obtain the
first and second tiling equations discussed above.

The next important step is to prove that the tile has to be 
 rational (that is, the ratios of its 
sides are rational), so in suitable units the sides can be taken to be 
integers.  For this, we use the coloring equation together with an 
analysis of the ``edge relations'' in a tiling.  An ``edge relation''
arises when a line (segment) interior to the tiling has different 
edges on one side than on the other.  For example, we might have 
15 $b$ edges on one side, and 18 $a$ edges on the other side, giving
rise to the edge relation $15b = 18a$.   The coloring equation and 
the analysis of edge relations work for all five possible shapes
of $ABC$ (subject to $3\alpha+2\beta=\pi$, including the three 
isosceles shapes).  After we prove that the tile has to be rational,
then the tiling equations can be used to analyze the possible 
values of $N$ and $(a,b,c)$ that can lead to tilings.  In short, when 
the tiling equation has no solutions, there are no tilings.  

When the tiling equation does have solutions, are there tilings? 
We answer this question in the affirmative, in case the solution satisfies $K^2 \mid M$,
 by explicitly constructing
the tilings in question.  None of these tilings were previously known.
Many illustrations of these tilings are given in the figures below.
In the $(2\alpha,\beta,\alpha+\beta)$ case we call these the 
 ``triquadratic tilings.''  The smallest of these is a 28-tiling.
  In the $(2\alpha,\alpha,2\beta)$ case, we call these
the ``four-component tilings''.  

For a given $N$,  if either of the tiling equations 
 has a solution, one can explicitly compute a finite set of possible tiles
(one for each solution of the tiling equation)
and for each tile, just one triangle $ABC$,  such that if any $N$-tiling exists,
it uses one of those tiles and tiles the corresponding $ABC$.

The first tiling equation does sometimes have more than one solution 
$(M,K)$ for a given $N$, but the smallest such example is $N=87800$.
The second tiling equation also has more than one solution for $N = 38962$.  

Finally,  we treat the three cases when $ABC$ is isosceles.  In 
each of these cases we follow the same method as for the cases 
already described.  Namely,  we write down the ``area equation'' 
that says the area of $ABC$ is equal to $N$ times the area of the tile.
Then we write down the ``coloring equation'' that relates the coloring 
number $M$ to the lengths $X$ and $Y$ of the sides of $ABC$.  Then we
eliminate $X$ and $Y$ from those two equations, yielding a Diophantine
equation in $N$, $M$, and the rational number $s = a/c$.  This is the 
``tiling equation'' for that shape of $ABC$.  We make a table of 
its solutions for $N$ up to some limit (such as 200 or 500 or 1000),
and then exhibit such tilings as we can,  and rule out by computational
search for ``boundary tilings'' such $N$ as we can, leaving question
marks in the table for other entries.  But for values of $N$ that do not 
occur in any of the five tables, we can definitely assert there is 
no $N$-tiling of any triangle not similar to $T$ by a tile satisfying $3\alpha + 2\beta = \pi$.

By an analysis of each of the five tiling equations, we were able to 
prove for three of the five shapes of $ABC$, that
if $3\alpha+2\beta = \pi$ and $ABC$ is $N$-tiled
by a tile with angles $(\alpha,\beta,\gamma)$, then $N$ cannot be prime.  
Two  cases remains open:  the triquadratic shape $(2\alpha, \beta, \alpha+\beta)$,
and the case when $ABC$ is isosceles with base angles $\beta$.

\section{Some basic lemmas}

A basic fact is that in a tiling with $3\alpha + 2\beta = \pi$, 
 the angles $\alpha$ and $\beta$ are not rational 
multiples of $\pi$.  This is a consequence of more general theorems 
of Laczkovich \cite{laczkovich1995}.  We here give the details as to 
how this fact follows from the theorems stated by Laczkovich, which are 
 \cite[Theorems~5.1 and 5.3]{laczkovich1995}. 

\begin{lemma} \label{lemma:angles}
Let $3\alpha+2\beta = \pi$. Suppose there is an $N$-tiling of triangle $ABC$
by tile $T$ with angles
$(\alpha,\beta,\gamma)$.  Suppose also that $ABC$ is not similar to $T$.
Then
  $\alpha$ and $\beta$ are not rational 
multiples of $\pi$, and every rational linear relation between $\pi$, $\alpha$,
and $\beta$ is a rational multiple of $3\alpha + 2\beta = \pi$.  
\end{lemma}

\begin{proof}  The last statement follows from the first, 
since a linear relation independent of  $3\alpha+2\beta = \pi$ would 
permit solving for $\alpha$ and $\beta$ in terms of $\pi$.  

Theorem~5.1 of \cite{laczkovich1995} does not mention the relation
$3\alpha+2\beta = \pi$.  Instead, it assumes that triangle $ABC$ 
can be dissected into similar (not necessarily congruent) triangles 
with angles $(\alpha,\beta,\gamma)$ that are rational multiples of $\pi$.
The conclusion is that either $ABC$  itself has angles $(\alpha,\beta,\gamma)$,
or $(\alpha,\beta,\gamma)$ belongs to 
a certain finite list of possible values. 
The first alternative is excluded by hypothesis, since it would make $ABC$ similar to $T$. 
So we have to check whether  
any of the triples in that list satisfy $3\alpha+2 \beta = \pi$.
We do not reproduce the entire list here, but the reader may easily
check that the only triple that does satisfy $3\alpha + 2\beta = \pi$,
under any labeling of its angles, 
is $(\pi/4,\pi/8,5\pi/8)$.  But Laczkovich's Theorem~5.3 says
that that triple is impossible for dissections into {\em congruent}
triangles.%
\footnote{Laczkovich proves several cases of Theorem~5.3
explicitly and remarks that similar arguments apply in the remaining
cases, which include the one needed here.}
\end{proof}

\begin{convention} \label{convention:tile}
Throughout this paper, whenever $T$ is a triangle (the ``tile'') with angles
$(\alpha,\beta,\gamma)$ and opposite sides $(a,b,c)$ satisfying
$$3\alpha + 2\beta = \pi,$$
  we write
$$s \ = \ 2 \sin \frac \alpha 2.$$
Note that $\gamma = 2\alpha+\beta$, and that $\alpha < \pi/3$ and
$\beta < \pi/2$; hence $0 < s < 1$.
\end{convention}
This convention is useful because 
the ratios $a/c$ and $b/c$ can be expressed simply in terms of $s$, as shown in the following
lemma.

\begin{lemma} \label{lemma:singamma}
Let $T$ have angles $(\alpha,\beta,\gamma)$ with $3\alpha + 2\beta = \pi$.   Then
$\sin\gamma = \cos(\alpha/2)$, $\ a/c = s$, and $\ b/c = 1-s^2$.
\end{lemma}

\begin{proof}   
Since $\gamma = \pi-(\alpha + \beta)$, we have 
\begin{eqnarray*}
\sin \gamma &=& \sin(\pi-(\alpha + \beta))  \\
 &=& \sin(\alpha + \beta) \nonumber \\
 &=& \cos(\pi/2- (\alpha + \beta))   \\
 &=& \cos \frac \alpha 2  \mbox{\qquad since $\pi/2 - \beta = 3\alpha/2$} 
 \end{eqnarray*}
Then $c = \sin \gamma = \cos \alpha/2$, and $a = \sin \alpha = 2 \sin(\alpha/2) \cos(\alpha/2)$.
Hence 
$$\frac a c = 2\sin \alpha/2.$$
Since $3 \alpha + 2\beta = \pi$, we have 
\begin{eqnarray*}
\sin \beta &=& \sin (\pi/2 - 3 \alpha/2) \\
&=& \cos(3 \alpha/2) \\
&=& 4 \cos^3 \frac \alpha 2 - 3 \cos \frac \alpha 2 
\end{eqnarray*}
Hence 
\begin{eqnarray*}
b/c &=& 4 \cos^2(\alpha/2) - 3 \\
&=& 4(1-\sin^2 \alpha/2) - 3 \\
&=& 1-4\sin^2 \alpha/2
\end{eqnarray*}
 Then we have 
\begin{eqnarray*}
\frac a c &=& s \\
\frac b c &=& 1-s^2
\end{eqnarray*}
establishing the second equation of the lemma.   
\end{proof}
\smallskip

\begin{lemma} \label{lemma:anglesOK}
Let $T$ be a  triangle with sides $a$, $b$, and $c$, 
 and let $\alpha$ and $\beta$ be the angles opposite $a$ and $b$ 
respectively.  Then 
    $3\alpha + 2\beta = \pi$ if and only if  
$b = c-a^2/c$. 
\end{lemma}

\begin{proof}  Suppose $3\alpha + 2\beta = \pi$.  By
Lemma~\ref{lemma:singamma}, $b/c = 1-s^2 = 1-(a/c)^2$.  Multiplying by $c$
we have the desired conclusion, $b = c-a^2/c$.

It remains to prove the right-to-left implication.  Assume $b = c-a^2/c$.
Since that condition is invariant if the triangle is re-scaled, we may
as well assume $c=1$.  Then $b = 1-a^2$, and since $b > 0$ we have
$0 < a < 1$.  By the law of cosines,
\begin{eqnarray*}
a^2 &=& 1 + (1-a^2)^2 - 2(1-a^2)\cos \alpha
\end{eqnarray*}
so that, since $1 + (1-a^2)^2 - a^2 = (a^2-1)(a^2-2)$,
\begin{eqnarray*}
\cos \alpha &=& \frac {1 + (1-a^2)^2 - a^2}{2(1-a^2)}
 \ = \ \frac {(a^2-1)(a^2-2)}{-2(a^2-1)} \ = \ 1 - \frac {a^2} 2
\end{eqnarray*}
Hence $a^2 = 2(1-\cos\alpha) = 4\sin^2(\alpha/2)$, and since $a>0$,
$$ a \ = \ 2 \sin \frac \alpha 2. $$
By the law of sines, $\sin \gamma = c (\sin\alpha) / a = (\sin \alpha) / a$, so
\begin{eqnarray*}
\sin \gamma &=& \frac {2 \sin \frac \alpha 2 \cos \frac \alpha 2}{2\sin \frac \alpha 2 }
 \ = \ \cos \frac \alpha 2 \ = \ \sin \Big( \frac \pi 2 - \frac \alpha 2 \Big)
\end{eqnarray*}
Therefore $\gamma = \pi/2 - \alpha/2$ or $\gamma = \pi/2 + \alpha/2$.
In the first case
$$ \beta \ = \ \pi - \alpha - \gamma \ = \ \frac \pi 2 - \frac \alpha 2
 \ = \ \gamma, $$
so $b = c$; but $b = 1-a^2 < 1 = c$, a contradiction.  Hence
$\gamma = \pi/2+\alpha/2$, and
$$ \beta \ = \ \pi - \alpha - \gamma \ = \ \frac \pi 2 - \frac{3\alpha} 2, $$
which is to say $3\alpha + 2\beta = \pi$, as required.
\end{proof}

\begin{lemma} \label{lemma:gsquare}
Suppose $3\alpha + 2 \beta = \pi$, and $(a,b,c)$ are the sides
of a triangle with angles $(\alpha,\beta,\gamma)$. Suppose $(a,b,c)$
are integers with no common factor. Let $a/c = m/k$  in lowest terms.
Then 
$$(a,b,c) = (mk,k^2-m^2,k^2),$$
and $k = \gcd(a,c)$ and $m = a/k$ and $m < k$.
\end{lemma}

\begin{proof}  Since $a$ and $c$ are positive we may take $k > 0$, and
then $m > 0$.  By Lemma~\ref{lemma:anglesOK},
$$\frac b c \ = \ 1 - \Big(\frac a c\Big)^2 \ = \ \frac{k^2-m^2}{k^2},$$
and therefore
$$(a,b,c) \ = \ \frac c {k^2}\,(mk,\ k^2-m^2,\ k^2).$$
The triple on the right has no common factor: a prime dividing both $k^2$
and $k^2-m^2$ would divide $k$ and $m$, contradicting $\gcd(m,k)=1$.
Since $(a,b,c)$ is likewise a triple of integers with no common factor,
and the two triples are positively proportional, they are equal.
Finally $\gcd(a,c) = \gcd(mk,k^2) = k \gcd(m,k) = k$, and $m = a/k$.
We have $m <k$ since $b > 0$ and $b = k^2-m^2$.
\end{proof}

\begin{lemma}\label{lemma:gammagreaterpiover2}
Suppose $3\alpha + 2\beta = \pi$.  Then $\gamma > \pi/2$.
\end{lemma}

\begin{proof}
\begin{eqnarray*}
\pi &=& 3\alpha + 2\beta \\
    &=& \alpha + 2(\alpha + \beta) \\
    &=& \alpha + 2 (\pi - \gamma)
\end{eqnarray*}
Hence $2\gamma = \pi + \alpha > \pi$, since $\alpha > 0$, so $\gamma > \pi/2$.
\end{proof}

\begin{lemma} \label{lemma:nogamma}
Let triangle $ABC$ be $N$-tiled by a tile $T$ with angles
$(\alpha,\beta,\gamma)$, where $3\alpha+2\beta = \pi$.  Suppose $ABC$
is not similar to $T$.  Then no tile has its $\gamma$ angle at a vertex of $ABC$,
and there are exactly five tile angles with vertices at the vertices of $ABC$.
\end{lemma}

\begin{proof}
By Lemma~\ref{lemma:angles}, $\alpha$ and $\beta$ are linearly
independent over the rationals.  Each angle of $ABC$ is the sum of the
angles at that vertex of the tiles having a vertex there, so it is a
nonnegative integral combination of $\alpha$, $\beta$, and $\gamma$;
and since $\gamma = 2\alpha+\beta$, it is a nonnegative integral
combination $p\alpha + q\beta$.  Writing the three angles of $ABC$ as
$p_i\alpha + q_i\beta$ for $i=1,2,3$, and summing, we have
$$ (p_1+p_2+p_3)\alpha + (q_1+q_2+q_3)\beta \ = \ \pi
 \ = \ 3\alpha + 2\beta, $$
so by linear independence $\sum p_i = 3$ and $\sum q_i = 2$.

Now suppose some tile has its $\gamma$ angle at a vertex of $ABC$, say
at the vertex with angle $p_1\alpha + q_1\beta$.  Then $p_1 \ge 2$ and
$q_1 \ge 1$, so $p_2 + p_3 \le 1$ and $q_2 + q_3 \le 1$.  Since the
angles of $ABC$ are positive, neither $p_2\alpha + q_2\beta$ nor
$p_3\alpha + q_3\beta$ is zero, so one of them is $\alpha$ and the
other is $\beta$.  Hence $p_1 = 2$ and $q_1 = 1$, and the angles of
$ABC$ are $\alpha$, $\beta$, and $2\alpha+\beta = \gamma$; that is,
$ABC$ is similar to $T$, contrary to hypothesis.  Hence no tile
has a $\gamma$ angle at a vertex of $ABC$.  Since $\sum p_i = 3$ and $\sum q_i = 2$,
and each tile contributes just one of those five angles,  there are exactly
five tiles at the vertices.
\end{proof}


\begin{lemma} \label{lemma:pigeonhole} Suppose triangle $ABC$ is $N$-tiled by a 
tile with angles $(\alpha,\beta,\gamma)$ and $\gamma > \pi/2$. 
 Suppose all the tiles along one side of $ABC$ do not have their $c$ sides
along that side of $ABC$.  Then there is a tile with a $\gamma$ angle at one of the endpoints of that side of $ABC$.
\end{lemma}

\begin{proof}
Let $PQ$ be the side of $ABC$ with no $c$ sides of tiles along it, and let
$n$ be the number of tiles with a side on $PQ$.  Each such tile has its $a$
or $b$ side on $PQ$, so the angle opposite that side is $\alpha$ or $\beta$;
hence the tile's $\gamma$ angle is at one of the two endpoints of its own
side on $PQ$, and in particular at a vertex of the tiling lying on $PQ$.

The edges of those $n$ tiles divide $PQ$ at $n-1$ interior points. Since
$\gamma > \pi/2$, two $\gamma$ angles at a single vertex on $PQ$ would
exceed the straight angle $\pi$; so each of those $n-1$ vertices carries at
most one of our $n$ angles.  By the pigeonhole principle at least one of
them is at $P$ or $Q$.
\end{proof}

\begin{lemma} \label{lemma:twocs}
Suppose triangle $ABC$  is 
$N$-tiled by a tile $T$ with angles $(\alpha,\beta,\gamma)$.  Suppose 
\smallskip

(i) $\gamma > \pi/2$, and
\smallskip
 
(ii) $\alpha/\beta$ is not rational, and 
\smallskip

(iii) Every angle of triangle $ABC$ is less than $\gamma$, and 
\smallskip

(iv) One of the following two conditions holds:  Either 
neither of $a,b$ is an integer multiple of the other, or else
neither of the vertices $A,C$ has two tiles with equal angles there
(that is, two $\alpha$ angles or two $\beta$ angles).
\smallskip

(v) $AC$ is not a single $c$ edge.
\smallskip

Then there 
are at least two $c$ edges of tiles on side $AC$.  
\end{lemma}

\begin{remark}  This lemma is proved in \cite{beeson-noseven}; 
we give the proof here in a corrected and expanded form.  One can prove by the same method that the 
$c$ edges must occur in adjacent blocks of at least two edges, but
we found no use for that result. 
\end{remark}

\begin{proof}[Proof of Lemma~\ref{lemma:twocs}]   
We first record a fact about boundary vertices.  Let $P$ be a vertex of the
tiling lying in the interior of a side of $ABC$, and suppose some tile has
its $\gamma$ angle at $P$.  The angles of the tiles at $P$ sum to $\pi$, so
they are $p\alpha+q\beta+r\gamma = \pi$ with $r \ge 1$.  Since
$\gamma > \pi/2$ by hypothesis~(i), we have $r = 1$; then
$p\alpha+q\beta = \pi-\gamma = \alpha+\beta$, and since $\alpha/\beta$ is
irrational by hypothesis~(ii), $p = q = 1$.  So exactly three tiles touch
$P$, contributing angles $\alpha$, $\beta$, and $\gamma$.

By hypothesis~(iii), no tile has a $\gamma$ angle at a vertex of $ABC$;
hence by Lemma~\ref{lemma:pigeonhole}, each side of $ABC$ has at least one
$c$ edge.  The present lemma, however, claims more: there must be at least
two $c$ edges.  Suppose, to the contrary, that there is just one tile with
a $c$ edge on the side $AC$ of triangle $ABC$; call it Tile~1.  Then every
other tile with an edge on $AC$ has its $a$ or $b$ edge there, and so has a
$\gamma$ angle at one of the two endpoints of that edge.
We visualize $AC$ as horizontal with triangle $ABC$ above, 
and use the words ``north'' and ``northwest'' accordingly. 
See Figure~\ref{figure:twocs}.

\begin{figure}[ht]
\begin{center}
\includegraphics[width=0.8\textwidth]{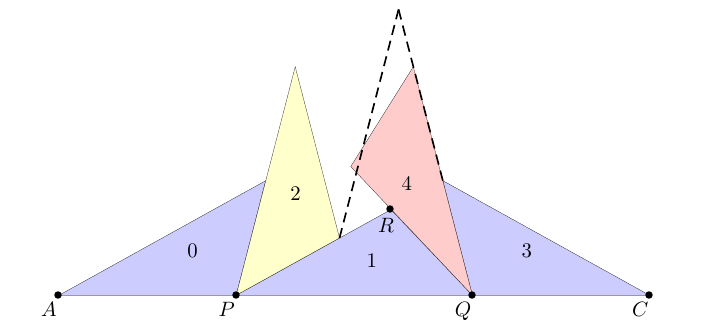}
\end{center}
\caption{Proof of Lemma~\ref{lemma:twocs}: another tile won't fit next to Tile~2}
\label{figure:twocs}
\end{figure}

\noindent
Since the hypotheses of the lemma remain true if (the names of)
$\alpha$ and $\beta$ are interchanged, we may assume without loss of
generality that $\alpha < \beta$.  Then by the law of sines, $a < b$.
Since $\gamma > \pi/2$ we also have $a < b < c$ (by the law of sines,
since $\beta < \pi - \gamma < \gamma$).

Let $PQ$ be the $c$ edge of Tile~1 lying on $AC$, and let $R$ be the
northern vertex of Tile~1.  Interchanging $P$ and $Q$ if necessary, suppose
that Tile~1 has its $\beta$ angle at $Q$.  Then Tile~1 has its $\alpha$
angle at $P$ and its $\gamma$ angle at $R$, the side $PR$ opposite $Q$
has length $b$, and the side $QR$ opposite $P$ has length $a$.
If $P$ is not a vertex of $ABC$, let Tile~0 be the tile on $AC$ immediately
west of $P$; if $Q$ is not a vertex of $ABC$, let Tile~3 be the tile on
$AC$ immediately east of $Q$.

We record some facts that will be used repeatedly.
\smallskip

(A) Tile~0 has its $\gamma$ angle at $P$, and Tile~3 has its $\gamma$ angle
at $Q$.  Indeed, let $n$ be the number of tiles with an edge on the segment
$AP$.  Their edges divide $AP$ at $n-1$ interior points.  Each of these
tiles has its $a$ or $b$ edge on $AP$, and so has a $\gamma$ angle at an
endpoint of its own edge; none has it at $A$, by hypothesis~(iii); and
since $2\gamma > \pi$, no point carries two $\gamma$ angles.  So $n$
angles are distributed among the $n-1$ interior points together with $P$,
one to a point, and in particular one lies at $P$.  The only one of these
tiles having $P$ as an endpoint of its edge is Tile~0.  The argument for
Tile~3 is the same.
\smallskip

(B) No tile has its $\gamma$ angle at $P$ except a tile with an edge on
$AC$.  Indeed, if $P$ is a vertex of $ABC$ this is hypothesis~(iii);
and if $P$ is not a vertex of $ABC$, then Tile~0 has its $\gamma$ angle at
$P$ by (A), and since $2\gamma > \pi$ there can be no second $\gamma$ angle
at $P$.
\smallskip

(C) The angle of $ABC$ at $P$ is greater than $\alpha$; consequently there
is a tile adjacent to $PR$ at $P$, which we call Tile~2.  Indeed, if $P$ is
not a vertex of $ABC$ the angle there is $\pi$.  Suppose then that $P$ is a
vertex of $ABC$ and its angle is $\alpha$.  Tile~1 fills that corner, so
$PR$ lies along a side of $ABC$, and $R$ lies on that side.  Then $Q$ and
$R$ both lie on the boundary of $ABC$; since $ABC$ is convex, no tile edge
along the line $QR$ extends beyond $Q$ or beyond $R$.  As $|QR| = a$ and
$a$ is the shortest side of the tile, exactly one tile lies on the far side
of $QR$ from Tile~1, and its $a$ edge is $QR$.  The two angles adjacent to
the side $a$ are $\beta$ and $\gamma$, so that tile has its $\gamma$ angle
at $Q$ or at $R$.  Not at $R$: the total angle at $R$ is $\pi$, since $R$
lies on a side of $ABC$, while Tile~1 contributes $\gamma$ there and
$2\gamma > \pi$.  Not at $Q$: if $Q$ is a vertex of $ABC$ that contradicts
hypothesis~(iii), and otherwise Tile~3 has its $\gamma$ angle at $Q$ by
(A), and again $2\gamma > \pi$.
\smallskip

(D) Suppose no tile edge lying along the line $PR$ contains $R$ in its
interior.  Then Tile~2 has its $a$ edge on $PR$, and so does every other
tile northwest of $PR$.  Indeed, under that supposition every tile edge
along $PR$ is contained in $PR$, so none is a $c$ edge, as $c > b = |PR|$.
If Tile~2 had its $b$ edge on $PR$, that edge would be all of $PR$, so its
$\gamma$ angle would be at $P$ or at $R$.  Not at $P$, by (B).  Not at $R$,
provided $R$ lies on a side of $ABC$ or in the interior of a tile edge:
then the angles at $R$ on the side of that line where Tile~1 lies sum to
$\pi$, while Tile~1 and Tile~2 would contribute $2\gamma > \pi$.  (The
remaining case allowed by (E), where $R$ is an interior vertex at which
another tile has its $\gamma$ angle, occurs only in Case~1c, and is
treated there.) So
Tile~2 has its $a$ edge on $PR$, and each remaining edge along $PR$ is
shorter than $b$ and hence is an $a$ edge.
\smallskip

(E) The hypothesis of (D) holds as soon as $R$ lies on a side of $ABC$, or
some tile other than Tile~1 has an angle of $\gamma$ at $R$, or some tile
contains $R$ in the interior of one of its edges.  For a tile edge through
$R$ along $PR$ contributes a straight angle at $R$ on the northwest side,
and Tile~1 contributes $\gamma$; the total angle at $R$ is $\pi$ in the
first case and $2\pi$ otherwise, and $\pi + \gamma > \pi$, while
$\pi + \gamma + \pi > 2\pi$ and $\pi + \gamma + \gamma > 2\pi$, the last
because $2\gamma > \pi$.
\smallskip

We also record the pigeon-hole argument that will be used three times.
Suppose every tile northwest of $PR$ has its $a$ edge on $PR$, and let $n$
be the number of such tiles; since $|PR| = b > a$ we have $n \ge 2$, and
the edges divide $PR$ at $n-1$ interior points.  Each of these tiles has
its $a$ edge on $PR$, so each has a $\gamma$ angle at an endpoint of its
own edge.  None is at $P$, by (B); and none is at $R$, since Tile~1 has its
$\gamma$ angle at $R$ and $2\gamma > \pi$.  Since $\gamma > \pi/2$, no
interior point of $PR$ carries two $\gamma$ angles.  So $n$ angles must
occupy $n-1$ points, which is impossible.  We refer to this as {\em the
pigeon-hole argument}.
\smallskip

We first show that one of $P$, $Q$ is a vertex of $ABC$.  Assume, for
proof by contradiction, that neither is.  Then Tile~0 and Tile~3 both
exist.  By (A), Tile~3 has its $\gamma$ angle at $Q$; so by the fact about boundary
vertices recorded at the start of the proof, exactly three tiles touch $Q$,
contributing $\alpha$, $\beta$, and $\gamma$.  Tile~1 contributes the
$\beta$ and Tile~3 the $\gamma$, so a single tile, which we call Tile~4,
fills the remaining angle $\alpha$ between them.
 Then Tile~4 has its $b$ or $c$ edge along $QR$, since
its angle at $Q$ is $\alpha$.  As $|QR| = a$ and $a < b < c$, the edge of
Tile~4 along $QR$ extends past $R$, so $R$ lies in the interior of an edge
of Tile~4.  By (E) the hypothesis of (D) holds, so by (D) every tile
northwest of $PR$ has its $a$ edge on $PR$, and the pigeon-hole argument
gives a contradiction.  This proves that one of $P$, $Q$ is a vertex of
$ABC$.
\smallskip

Now we argue by cases. 
\smallskip

Case 1: $Q$ is a vertex of $ABC$, i.e.\ $Q = C$.  If the angle of $ABC$
at $Q$ equals $\beta$, then $QR$ lies along the side $CB$ of $ABC$, and
there is no tile adjacent to Tile~1 at $Q$; call this Case~1a.
Otherwise let Tile~4 be the tile adjacent to Tile~1 at $Q$, sharing the
edge $QR$ with Tile~1.  Tile~4 cannot have its $\gamma$ angle at $Q$,
by hypothesis~(iii).  So either Tile~4 has its $\alpha$ angle at $Q$
(Case~1b) or its $\beta$ angle at $Q$ (Case~1c).  These three subcases
are exhaustive.
\smallskip

Case 1a: the angle of $ABC$ at $Q$ is $\beta$.  Then $QR$ lies along the
side $CB$, so $R$ lies on a side of $ABC$.  By (E) the hypothesis of (D)
holds, so by (D) every tile northwest of $PR$ has its $a$ edge on $PR$, and
the pigeon-hole argument gives a contradiction.  That completes Case~1a.
\smallskip

Case 1b: Tile~4 has its $\alpha$ angle at $Q$.  Then the two edges of
Tile~4 at $Q$ are its $b$ edge and its $c$ edge, and since
$|QR| = a < b < c$, the edge of Tile~4 along $QR$ extends past $R$, so $R$
lies in the interior of an edge of Tile~4.  By (E) the hypothesis of (D)
holds, and the pigeon-hole argument gives a contradiction as in Case~1a.
That completes Case~1b.
\smallskip

Case 1c: Tile~4 has its $\beta$ angle at $Q$.  Then Tile~1 and Tile~4
both have a $\beta$ angle at $Q$, so there are two equal angles at the
vertex $C$ of $ABC$; hence by hypothesis~(iv), $b$ is not an integer
multiple of $a$.  The two edges of Tile~4 at $Q$ are its $a$ edge and its
$c$ edge.  If the one along $QR$ is the $c$ edge, then since
$|QR| = a < c$ it extends past $R$, so $R$ lies in the interior of an edge
of Tile~4, and the argument of Case~1b gives a contradiction.  So Tile~4
has its $a$ edge on $QR$; that edge is $QR$ itself, and since the angles
adjacent to the $a$ edge are $\beta$ and $\gamma$, Tile~4 has its $\gamma$
angle at $R$.  By (E) the hypothesis of (D) holds.  Here $R$ is an interior
vertex, so we must check directly that Tile~2 does not have its $b$ edge on
$PR$.  Suppose it does.  That edge is all of $PR$; the angles adjacent to
the $b$ edge are $\alpha$ and $\gamma$, and by (B) the $\gamma$ is not at
$P$, so Tile~2 has its $\alpha$ angle at $P$.  If $P$ is not a vertex of
$ABC$, the angles at $P$ sum to $\pi$; Tile~0 contributes $\gamma$, by (A),
and Tiles~1 and~2 contribute $\alpha$ each, leaving
$\pi - \gamma - 2\alpha = \beta - \alpha > 0$.  But no sum of tile angles
equals $\beta - \alpha$: a $\gamma$ is too large, and
$p\alpha + q\beta = \beta - \alpha$ would give $(p+1)\alpha = (1-q)\beta$,
which is impossible since $\alpha/\beta$ is irrational and $p + 1 > 0$.
Hence $P$ is a vertex of $ABC$, so $P = A$; since also $Q = C$, the side
$AC$ is the single $c$ edge $PQ$, contrary to hypothesis~(v).  So Tile~2
has its $a$ edge on $PR$, and as in (D), every tile northwest of $PR$ has
its $a$ edge on $PR$.  Since $|PR| = b$, this makes $b$ an integer
multiple of $a$, a contradiction.  That completes Case~1c, and with it
Case~1.
\smallskip

Case 2: $P$ is a vertex of $ABC$, and $Q$ is not.  Then Tile~3 exists, and
as in the argument above, Tile~4 fills the open $\alpha$ angle between
Tile~1 and Tile~3, and its edge along $QR$ extends past $R$, so $R$ lies in
the interior of an edge of Tile~4.  Tile~2 exists by (C).  By (E) the
hypothesis of (D) holds, and the pigeon-hole argument gives a
contradiction.  That completes Case~2, and the proof of the lemma.
\end{proof}

\begin{lemma} \label{lemma:cminusb}
Suppose $3\alpha + 2\beta = \pi$ and $T$ is a triangle with angles $(\alpha,\beta,\gamma)$
and sides $(a,b,c)$.  Then $c-b$ is not an integer multiple of $a$.
\end{lemma}

\begin{proof} By Lemma~\ref{lemma:anglesOK}, we have $b = c-a^2/c$.  Then $c-b = a^2/c$, so 
$$
 \frac{c-b} a = \frac a c = s,
$$
which lies between 0 and 1 and so cannot be an integer.
\end{proof}

\section{Types of vertices}

  Suppose triangle $ABC$ is tiled by a tile $T$ with angles $(\alpha,\beta,\gamma)$.
    The {\em type} of a vertex  $V$ 
 is a triple $(n,k,j)$, meaning 
that there are  $n$ tiles with an $\alpha$ angle at $V$, $k$ tiles with a $\beta$ angle at $V$, and $j$
tiles with a $\gamma$ angle at $V$.    Vertices of types $(1,1,1)$ and 
$(2,2,2)$ are called ``standard vertices''.  A vertex of type $(0,1,3)$ is called a ``center''.   Vertices of types
$(3,2,0)$, $(6,4,0)$, and $(4,3,1)$ are called ``sporadic vertices''. 
(They do occur in some tilings.)

\begin{lemma} \label{lemma:standardvertices}
Suppose $3\alpha + 2\beta = \pi$.  Suppose triangle $ABC$
is $N$-tiled by a tile $T$ with angles $(\alpha,\beta,\gamma)$, and $ABC$ is not similar to $T$.
  Then 
  
 (i) the number of centers is $N_C = 1 + N_1 + 2N_2+ N_3$, 
where $N_1$ is the number of vertices of type $(3,2,0)$, $N_2$ is the number of vertices of type $(6,4,0)$,
and $N_3$ is the number of $(4,3,1)$ vertices, and 

(ii) every vertex other than $A$, $B$, and $C$ has one of the types  $(3,2,0)$ or $(1,1,1)$ (which have angle sum $\pi$),
or $(2,2,2)$, $(6,4,0)$, $(4,3,1)$, or $(0,1,3)$ (which have angle sum $2\pi$).
\end{lemma}

\begin{proof} 

First we prove (ii).  
By Lemma~\ref{lemma:angles}, $\alpha$ and $\beta$ are linearly independent over $\Q$. 
 Each vertex other than $A$, $B$, or $C$ has its tile angles summing to
$\pi$, if it lies on a side of $ABC$ or on a side of another tile, or
to $2\pi$, if it is a strict interior vertex.  Since
$\gamma = 2\alpha+\beta$, a vertex of type $(n,k,j)$ has angle sum
$(n+2j)\alpha + (k+j)\beta$; so by linear independence, the sum is
$\pi = 3\alpha+2\beta$ exactly when $n+2j = 3$ and $k+j = 2$, and is
$2\pi$ exactly when $n+2j = 6$ and $k+j = 4$.  The first pair of
equations gives the types $(3,2,0)$ and $(1,1,1)$; the second gives
$(6,4,0)$, $(4,3,1)$, $(2,2,2)$, and $(0,1,3)$.  The angle sums are as
in the statement of the lemma.

Now to prove (i):
  By Lemma~\ref{lemma:nogamma}, 
there are exactly five tiles with vertices at the vertices of $ABC$.  
In an $N$-tiling there are $N$ triangles.  Each has one $\alpha$, one $\beta$, and one $\gamma$ angle, 
so the total number of each kind of angle is $N$.  At the vertices of $ABC$ we have three $\alpha$ angles 
and two $\beta$ angles,  making an excess of three $\alpha$ angles over $\gamma$ angles at the 
vertices of $ABC$.  Similarly, at the sporadic vertices we get an excess of $3N_1 + 6 N_2+3N_3$.  At 
the centers,  the number of $\gamma$ angles exceeds the number of $\alpha$ 
 angles by 3 per center, for a total of $3N_C$.  Standard vertices contribute no excess, since $n=j$ there.
 Therefore the total excess of $\alpha$ 
angles over $\gamma$ angles is 
\begin{eqnarray*}
3 + 3N_1 + 6 N_2 +3N_3 - 3N_C &=& 0
\end{eqnarray*}
Solving for $N_C$ we have 
\begin{eqnarray*}
N_C &=& 1 + N_1 + 2N_2 + N_3
\end{eqnarray*}
as claimed.  
\end{proof}

\begin{lemma} \label{lemma:shapes}
Suppose $3\alpha + 2\beta = \pi$ and triangle $ABC$ is $N$-tiled by a
tile $T$ with angles $(\alpha,\beta,\gamma)$, with $ABC$ not similar to
$T$.  Then the angles of $ABC$ are one of the following five triples:
\begin{eqnarray*}
&& (\alpha,\ \alpha+2\beta,\ \alpha) \qquad\quad
   \mbox{($ABC$ isosceles with base angles $\alpha$)} \\
&& (\beta,\ 3\alpha,\ \beta) \qquad\qquad\ \
   \mbox{($ABC$ isosceles with base angles $\beta$)} \\
&& (\alpha+\beta,\ \alpha,\ \alpha+\beta) \quad\ \,
   \mbox{($ABC$ isosceles with base angles $\alpha+\beta$)} \\
&& (2\alpha,\ \beta,\ \alpha+\beta) \mbox{ \qquad\ \,(the triquadratic shape)} \\
&& (2\alpha,\ \alpha,\ 2\beta)
\end{eqnarray*}
and there are exactly three $\alpha$ angles and two $\beta$ angles of tiles
at the vertices of $ABC$.
\end{lemma}
 \begin{proof}
By Lemma~\ref{lemma:nogamma}, no tile has its $\gamma$ angle at a
vertex of $ABC$.  Hence, summing the angles of all the tiles having a
vertex at a vertex of $ABC$, we have $P\alpha + Q\beta = \pi$ for some
nonnegative integers $P$ and $Q$.  Since $\pi = 3\alpha+2\beta$ and
$\alpha$, $\beta$ are linearly independent over $\Q$ by
Lemma~\ref{lemma:angles}, we have $P = 3$ and $Q = 2$; so there are
three $\alpha$ angles and two $\beta$ angles in all, distributed among
the three vertices of $ABC$, with at least one at each vertex.

There are exactly six such distributions, up to the order of the
vertices; writing each vertex as the pair (number of $\alpha$ angles,
number of $\beta$ angles), they are
$$ \{(0,1),(0,1),(3,0)\},\quad \{(1,0),(1,0),(1,2)\},\quad
   \{(1,0),(1,1),(1,1)\}, $$
$$ \{(0,1),(1,1),(2,0)\},\quad \{(0,2),(1,0),(2,0)\},\quad
   \{(0,1),(1,0),(2,1)\}. $$
The last of these gives $ABC$ the angles $\beta$, $\alpha$, and
$2\alpha+\beta = \gamma$, so that $ABC$ is similar to $T$, contrary to
hypothesis.  The remaining five give the angles
$(\beta,\beta,3\alpha)$, 
$(\alpha,\alpha,\alpha+2\beta)$,
$(\alpha,\alpha+\beta,\alpha+\beta)$, 
$(\beta,\alpha+\beta,2\alpha)$,
and $(2\beta,\alpha,2\alpha)$, 
which are the five triples listed, and
the tile counts at the vertices are as stated,  three $\alpha$ and two $\beta$.
\end{proof}

\section{Coloring tilings black and white}  \label{section:coloring}
\begin{figure} [ht]
\begin{center} \includegraphics[width=0.8\textwidth]{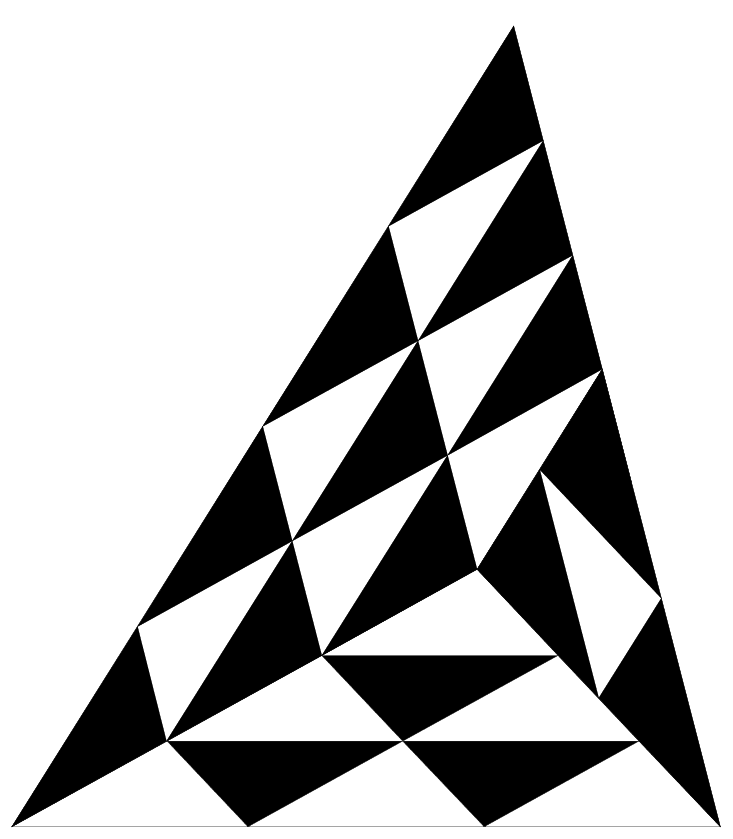}
\end{center}
\caption{A tiling colored so that touching tiles have different colors.}
\label{figure:coloring}
\end{figure}

 We think of coloring each tile black or white, in such a way 
that tiles touching along a line segment have different colors.  Technically we can 
represent the two colors as $\pm 1$ and speak of the ``sign'' of a triangle.  We show 
that the equation $3\alpha +2\beta = \pi$ implies that tilings can be 2-colored in this way,
and we derive from this a ``coloring equation'' relating the lengths of the sides
of $ABC$ to the difference $M$ between the numbers of black and white tiles.
See Figure~\ref{figure:coloring}.

\begin{lemma}\label{lemma:signedtiles}
Suppose triangle $ABC$ is $N$-tiled by
a tile with angles $(\alpha,\beta,\gamma)$,
and 
\smallskip

(i) $3\alpha + 2\beta = \pi$, and 
\smallskip

(ii) $ABC$ is not similar to the tile.
\smallskip

\noindent
Then it is possible to 
assign a color ``black'' or ``white'' 
 to each of the $N$ tiles in such a way that along each interior edge of the tiling, the triangles
on opposite sides of the edge receive opposite colors, 
and the color of one chosen tile is specified to be black. 
\end{lemma}

\begin{proof} 
We first prove that an even number of (segments of) edges meet at 
each interior vertex.  By Lemma~\ref{lemma:standardvertices}, every interior
vertex with angle sum $2\pi$  has type   $(2,2,2)$,  $(6,4,0)$, $(4,3,1)$, or $(0,1,3)$.  
At an interior vertex, the number of edges is the same as the number of tiles,
and in each of those four cases, the sum is even.  Now consider a vertex in the interior of $ABC$ 
with angle sum $\pi$, i.e.,  it occurs on the interior of some tile edge.  
Then the number of edges is one more than the number of tiles, so we need the
number of tiles to be odd.  By Lemma~\ref{lemma:standardvertices}, such 
vertices have type $(3,2,0)$ or $(1,1,1)$, and indeed the number of tiles is odd.
Then, as claimed, the number of edges is even.

Let $T_k$ be any one of the $N$ tiles, and let $P$ be a point in $T_k$.  Let $\sigma$ be a path 
from vertex $B$ to $P$ that does not pass through any vertex of the tiling, and passes transversally through 
each edge it crosses.  (Transversally means it is not tangent to the edge.)   Then the sign we wish to assign to $T_k$
is the parity of the number of edges crossed by $\sigma$.  That is, $1$ for even and $-1$ for odd. 
We claim this sign does not depend on which path $\sigma$ is chosen, 
but only on the tile $T_k$.  To show that, it suffices to show that the number is invariant under homotopies of $\sigma$
fixing the two endpoints.  That in turn follows from the fact that an even number of (segments of) edges meet at 
each interior vertex,  which was proved above, since the two ways of passing around a vertex cross complementary sets of the edge-segments at that vertex. 
\end{proof}

\begin{theorem}[Coloring theorem] \label{theorem:coloring}
 Let $3\alpha + 2\beta = \pi$, and assume triangle $ABC$ 
  is $N$-tiled by a tile with angles $(\alpha,\beta,\gamma)$, and 
  $ABC$ is not similar to the tile.  Label the vertices so that $A$ and
  $C$ have the same number of tiles at them (this is always possible),
  and let the tiles be assigned signs in accordance with
  Lemma~\ref{lemma:signedtiles}, with a tile at $B$ with an edge on $AB$ colored black.  Let
  $(X,Y,Z)$ be the lengths of $(BC, AC, AB)$ respectively.
   Then we have the ``coloring
  equation''
$$ M(a+b+c) = X \pm Y+Z$$
where $M$ is the number of black tiles minus the number of white tiles,
and the sign is determined by the shape of $ABC$:  the plus sign if
$ABC$ is isosceles with base angles $\alpha$ or $\beta$ (at $A$ and
$C$), and the minus sign in the three remaining cases, namely if the
angles of $ABC$ are $(2\alpha,\alpha,2\beta)$ or
$(2\alpha,\beta,\alpha+\beta)$, or $ABC$ is isosceles with base angles
$\alpha+\beta$.
\end{theorem}
  
\begin{proof}    By Lemma~\ref{lemma:angles},
which summarizes work of Laczkovich \cite{laczkovich1995},
$\alpha$ and $\beta$ are linearly independent over $\Q$.  
Since $ABC$ is not similar to the tile, no vertex of $ABC$ has a $\gamma$ angle, 
by Lemma~\ref{lemma:nogamma}, and there are five tiles at the vertices of $ABC$.
By Lemma~\ref{lemma:shapes}, $ABC$ has one of the five shapes listed in the
statement of the theorem, and the numbers of tiles at the three vertices are
$1$, $1$, $3$ if $ABC$ is isosceles with base angles $\alpha$ or $\beta$, and
$1$, $2$, $2$ in the three remaining cases; in either case two of the three
vertices have the same number of tiles at them, so the labelling assumed in
the statement of the theorem is possible.
 We think of $ABC$ with 
 vertex $B$ at the north.  
By Lemma~\ref{lemma:signedtiles}, all the tiles can be colored such that 
adjacent tiles have opposite colors. 

We first record a fact used repeatedly below.  Let $\ell$ be a line, and
let $T$ and $T'$ be two tiles on the same side of $\ell$, each with an
edge on $\ell$, whose edges on $\ell$ share an endpoint $V$.  Then $T$
and $T'$ have the same color. Indeed, the tiles having a vertex at $V$ on that side of $\ell$ form a fan
across the straight angle at $V$, beginning with $T$ and ending with $T'$.
Their angles sum to $\pi$, so by the computation in the proof of
Lemma~\ref{lemma:standardvertices} they comprise either three $\alpha$
angles and two $\beta$ angles, or one of each of $\alpha$, $\beta$,
$\gamma$; in either case the fan consists of an odd number of tiles, five
or three.  Consecutive tiles of the fan are adjacent along an edge and so
have opposite colors; since their number is odd, the first and last have
the same color.  We refer to this as the {\em fan argument}.
Consider a {\em maximal segment} $PQ$ in the tiling (if any exist), i.e.\ a
part of a straight line consisting of interior edges of the tiling, which
cannot be extended to a longer such segment.  
By the fan argument, all the tiles on one side of
$PQ$ with an edge on $PQ$ have the same color, and all those on the
other side have the same color; and since the two tiles meeting $PQ$
from opposite sides at any point are adjacent along an edge, those two
colors are opposite.  Hence the signed sum of the lengths of the edges
that these tiles share with $PQ$ is zero.  Every interior edge of the
tiling lies on exactly one maximal segment, so summing over all maximal
segments, the total length of the positive interior edges equals the
total length of the negative interior edges.

At vertex $B$ there is either exactly one tile or three, and in either
case the tiles at $B$ with an edge on $AB$ or on $BC$ are black:  when
there is one tile it is black by hypothesis, and when there are three
they alternate in color, so the two outer ones are both black.  By the
fan argument, consecutive tiles with an edge on $AB$ have the same
color; hence every tile with an edge on $AB$ is black, and likewise
every tile with an edge on $BC$.

The color of the tiles on side $AC$ is determined by the number of tiles
at $A$.  If there is only one tile at $A$, it has an edge on $AB$ and an
edge on $AC$, and it is black, being a tile with an edge on $AB$; if
there are two tiles at $A$, the one with an edge on $AB$ is black, and
the one with an edge on $AC$ is adjacent to it along an edge and so is
white.  In either case the fan argument shows that all the tiles with an
edge on $AC$ have the same color as the tile at $A$ with an edge on
$AC$.  So $AC$ is black if there is one tile at $A$, and white if there
are two.

Let $M$ be the number of positive tiles minus the number of negative
tiles.  Then the difference between the total length of the positive
edges and the total length of the negative edges is $M$ times the
perimeter $a+b+c$ of a tile.  Since the positive and negative interior
lengths are equal, we have
\begin{eqnarray*}
M(a+b+c) &=& X\pm Y+Z
\end{eqnarray*}
with the plus sign if $AC$ is black and the minus sign if $AC$ is white.
By the tile counts computed at the beginning of the proof, and the
labelling of the vertices, there is one tile at $A$ exactly when $ABC$
is isosceles with base angles $\alpha$ or $\beta$.  In that case $AC$ is
black and the plus sign applies; in the three remaining cases there are
two tiles at $A$, so $AC$ is white and the minus sign applies.
\end{proof}

\begin{lemma} \label{lemma:Mnotzero} Let $3\alpha + 2\beta = \pi$ and assume triangle $ABC$ is $N$-tiled by
a tile $T$ with angles $(\alpha,\beta,\gamma)$, and $T$ is not similar to $ABC$.   Let the tiles be 
assigned signs (or colors black and white) in accordance with Lemma~\ref{lemma:signedtiles}, and let
 $M$ be the number of positive (black) tiles minus the number of negative (white) tiles.  Then $M$ is not zero.
 \end{lemma}
 
\begin{proof}  By Theorem~\ref{theorem:coloring}, with the vertices
labeled as there, we have $M(a+b+c) = X \pm Y + Z$, where
$(X,Y,Z)$ are the lengths of $(BC, AC, AB)$.  Assume, for proof by
contradiction, that $M=0$.  Then either $X+Y+Z=0$, which is absurd since
all the sides are positive, or $X-Y+Z=0$, i.e.\ $Y = X+Z$, contradicting
Euclid I.20, that any two sides of a triangle are together greater than
the third.
\end{proof}

\section{Some definitions}
The {\em angle sum} of a vertex is the sum of angles of tiles at that vertex.  Except for the vertices
$A$, $B$, and $C$, it will be $2\pi$ or $\pi$.
 The latter occurs when the vertex is on the boundary of $ABC$, 
but it can also occur in the interior of $ABC$.

\begin{definition} \label{definition:supports}
A line (segment) $pq$ {\bf supports} a tile $T$ if a side, or a 
positive-length part of a side, of $T$ lies on the segment $pq$.
\end{definition}  

\begin{definition}  \label{definition:interiorsegment}
 An {\em \bf interior segment} is a line segment formed of boundaries of tiles 
in a tiling, such that its endpoints are vertices of tiles, and it 
does not lie on the boundary of the tiled triangle $ABC$.   
\end{definition}

An interior
segment $pq$ is {\em terminated at $q$} if on each side of $pq$ 
there is a tile supported by $pq$ with a vertex at $q$.
 An interior segment terminated at both endpoints
is called {\em doubly terminated}. 
\smallskip

A {\em  maximal segment} is an interior segment that does not lie on 
any longer interior segment, i.e. cannot be extended (in either direction).

\begin{definition} \label{definition:essentialedgerelation}
Given a tiling of triangle $ABC$ by triangle $T$, 

(i) A {\em  $b$-relation} is a relation of the form 
$jb = ua + vc$ with nonnegative integers $j$, $u$, and $v$
and $j > 0$. 
Similarly an $a$-relation is $ja = pb + qc$, and a $c$-relation is $jc = pa + qb$.

(ii) An interior segment
 $PQ$ is said to {\em  correspond to},
or to {\em  witness},  
a relation $jb = ua + vc$ (with $j>0$ and $u,v \ge 0$) 
if there are tiles on both sides of $PQ$ with vertices at $P$, and tiles on both 
sides with vertices at $Q$, and on one side of $PQ$ there are $j$ more $b$ edges than on the other,
and on the other side, there are $u$ more $a$ edges and $v$ more $c$ edges than on the 
first side.   Similarly for a relation $ja = ub + vc$ or $jc = ua + vb$.

(iii) An {\em essential segment} $PQ$ of a tiling is a segment supporting tiles on both
sides,  one of which has a vertex at $P$ and one has a vertex at $Q$, and does not have
the same number of tiles of each edge length on both sides.
\end{definition}

\begin{remark} 
Equating the sum of the edges on one side of a doubly
terminated interior segment to the sum of the edges on the other side
gives rise to an integral relation (linear relation with integral coefficients) between $a$, $b$, and $c$.
That relation will be nontrivial for an essential interior segment.
\end{remark}

\begin{examples}  The 28-tiling illustrated in Figure~\ref{figure:28}  contains two
essential maximal segments.
The one at the lower right in the figure 
 has three $a$ edges on one side and two $b$ edges on the other.  Thus for this tiling we have 
 $3a = 2b$.   
 Each biquadratic
tiling contains a maximal segment on the altitude of $ABC$ connecting the right angle to the hypotenuse.
\end{examples}

\begin{remark}
Note that the definition of $b$-relation does not require that there be a tiling in which
the relation is witnessed.  A $b$-relation is just a number-theoretic fact.
An essential segment gives rise to a relation $ja + kb + \ell c = pa + qb + rc$, where
it is not the case that $(j,k,\ell) = (p,q,r)$. Then by subtraction, we can get an $a$-relation, a $b$-relation, or a
$c$-relation.  If all three of $a$, $b$, $c$ occur with different counts on
the two sides, exactly one of the three relations results; if one of them
occurs equally often on both sides, then two of the three result at once,
as with $3a = 2b$ in the 28-tiling.  The essential segment witnesses the
relations so obtained.

One way in which a $b$-relation could be witnessed is if $PQ$ supports on one side only
tiles with $b$~edges on $PQ$.  Although that will be of interest in many proofs, it is not
the only way to witness a $b$-relation.  
\end{remark}

The following definition is equivalent to the one given 
in \cite[p.~346]{laczkovich2012}, except that there condition (iv) is
automatic because of an additional assumption about the tile. 

\begin{definition}
 [The directed graph $\Gamma_a$]
  \label{definition:Gamma} Given a tiling
of some triangle, 
the nodes of the graph $\Gamma_a$ are the endpoints of links as defined below. 
A link of $\Gamma_a$ connects
vertices $X$ and $Y$ if each (all) of the following holds:
 
(i) the segment $XY$ is covered by $a$ edges of tiles on one side
(say Side~1) of $XY$, and no $a$ edge on Side~1 extends $XY$ past $X$;
\smallskip

(ii) at vertex $Y$, there is
another tile supported by $XY$ extended past $Y$ on Side~1 of $XY$ with a vertex at $Y$,
 that does not have its
$a$ edge on Side~1, and
\smallskip

(iii) $Y$ lies in the interior of an edge of a tile on Side~2 of $XY$.

   The directed graphs $\Gamma_b$ and $\Gamma_c$
are defined
similarly. 
\end{definition}

\begin{remark}  We use {\em link} instead of {\em edge} with 
these graphs, to avoid confusion with tile ``edges.''
If $XY$ is a link in $\Gamma_a$, then $XY$ does not terminate at $Y$,
because there is another tile past $Y$ whose $a$ side is not on $XY$;
and also because $Y$ lies on the interior of an edge of a tile on 
the other side of $XY$.  In graph theory, sometimes the word ``vertex''
is used.  We use ``node'' here because ``vertex'' is already in use for 
triangle tilings.  Every node of these graphs is a vertex of the tiling,
but not every vertex is a node.  
\end{remark}

\begin{remark}
Laczkovich's definition of $\Gamma_a$ in \cite{laczkovich2012} differs from
ours in three ways.  He requires $X$ to be an endpoint of a maximal internal
segment, so that the tiles on both sides of $XY$ have a vertex at $X$; and he
requires that the first tile on Side~2 (the one with a vertex at $X$) does not
have its $a$ edge on $XY$.  On the other hand, he does not assume condition~(iii);
it follows from his hypothesis that $a$ is not a linear combination of $b$
and $c$ with nonnegative rational coefficients, a hypothesis not available
here.  Our definition is the one used in the Lean formalization
(Appendix~\ref{section:lean}).  The links produced by
Lemma~\ref{lemma:apropagation} satisfy Laczkovich's two extra conditions, and
the proofs need only that lemma, Lemma~\ref{lemma:noboundarylink}, and the
fact that at most one link ends at any vertex.  If there are no witnessed
$a$-relations, both graphs are empty.
\end{remark}

\section{There is an $a$-relation}

%
%

We will use our version of Laczkovich's graph $\Gamma_a$, as given in Definition~\ref{definition:Gamma}.

\begin{lemma} \label{lemma:noboundarylink}
No link of $\Gamma_a$ terminates at a point of the boundary of $ABC$.
\end{lemma}

\begin{proof}
Let $XY$ be a link and suppose $Y$ lies on the boundary of $ABC$.  Since
$XY$ is an internal segment it meets the interior of $ABC$, and $ABC$ is
convex, so the extension of $XY$ past $Y$ lies outside $ABC$.  Hence no
tile is supported by it, contrary to condition~(iii).
\end{proof}

\begin{lemma}[Propagation in $\Gamma_a$]\label{lemma:apropagation}  
Let triangle $ABC$ be $N$-tiled by a tile with sides $(a,b,c)$,
whose angles are  $(\alpha,\beta,\gamma)$ with $3\alpha + 2\beta = \pi$. 
 Suppose 
$ABC$ is not similar to the tile, and that the tiling witnesses no $a$-relation.
 Then  \begin{itemize}
\item The graph $\Gamma_a$ propagates links,   which means that if $PQ$ is a link in $\Gamma_a$,
then there is a link $QR$ beginning at $Q$, and 
\item The in-degree and out-degree of every node is $1$.
\item No link can start on the boundary of $ABC$.
\end{itemize}
\end{lemma}

\begin{proof}
Let $PQ$ be a link in $\Gamma_a$.  We must show there is 
a link $QR$ for some $R$.  Since $PQ$ is a link, $Q$ is a vertex of angle sum $\pi$,
since $Q$ is not a vertex on the non-$a$ side of $PQ$.  The last tile on the 
$a$ side has an $a$ edge on $PQ$ and the next tile does not.  We claim
there is an outgoing $a/c$ or $a/b$ edge from $Q$.    The type of 
the vertex $Q$ is $(1,1,1)$ or $(3,2,0)$. We argue by cases.

{\em Case~1}, the type of $Q$  is $(1,1,1)$. Then tiles with $\beta$ and $\gamma$ angles at $Q$
each have an $a$ edge at $Q$, but the tile with an $\alpha$ angle does not.  One of those $a$ edges is 
on $PQ$ before $Q$, since if not, $PQ$ would witness an $a$-relation.
Neither of those two $a$ edges is on $PQ$ after $Q$.  The other one must form an $a/c$ or $a/b$
edge outgoing from $Q$.  

{\em Case~2}, the type of $Q$ is $(3,2,0)$. Then the tiles with $\beta$ angles each have an $a$ edge, but the
$\alpha$ angles do not. So there are two $a$ edges at $Q$, just as in Case~1, so again there is
an outgoing $a/c$ or $a/b$ edge.

Now let $QR$ be the line containing that $a/c$ or $a/b$ edge, with $R$ chosen so that 
all the tiles on one side of $QR$ have their $a$ edges on $QR$, but the next tile after $R$ does not,
or else $QR$ terminates at $R$.   $QR$ does not terminate at $R$, since that would make $QR$ witness
an $a$-relation.  Therefore $QR$ is the desired link in $\Gamma_a$. 

That completes the proof of propagation.
The in-degree of any node of $\Gamma_a$ is at most 1, since two different lines of the tiling
cannot pass through a vertex of angle sum $\pi$. 
We have just proved that the out-degree is at least 1. Therefore the out-degree is always 
greater than or equal to the in-degree.  Since the total out-degree equals the total in-degree 
(every link has a head and a tail),  the in-degree and out-degree of every node is exactly 1.
By Lemma~\ref{lemma:noboundarylink}, no link terminates at a point of the
boundary of $ABC$.  Hence such a point has in-degree $0$, so by the
preceding paragraph it has out-degree $0$; that is, no link starts on the
boundary of $ABC$.
\end{proof}

\FloatBarrier

\section{A new invariant}

In this section we introduce a new invariant, applicable to the case
$3\alpha + 2\beta = \pi$, but similar to the invariant introduced by
Laczkovich on p.~351 of \cite{laczkovich2012} for the case $\gamma = 2\pi/3$,
and also similar to the variations of Laczkovich's invariant that were used in
\cite{beeson-zhang2026} and \cite{bonfioli2026}.  What is new here is that a
direction occurring in such a tiling requires two indices rather than one, and
that there are then exactly two invariants of this kind
(Lemma~\ref{lemma:Phiproperties}): one is the black-and-white coloring of
Section~\ref{section:coloring}, and the other is our new invariant $\Phi$.

Throughout this section and the next, $T$ is a triangle with angles $\alpha$,
$\beta$, $\gamma$ opposite the sides $a$, $b$, $c$, satisfying
$3\alpha+2\beta = \pi$; $ABC$ is a triangle tiled by $N$ copies of $T$; and
$ABC$ is not similar to $T$.  We write
$$\delta = \alpha + \beta .$$

\begin{lemma} \label{lemma:incommensurable} 
$\delta/\pi$ is irrational.
\end{lemma}

\begin{proof}
By Lemma~\ref{lemma:angles}, $\alpha$ is not a rational multiple of $\pi$.
Since $\alpha = \pi - 2\delta$, if $\delta$ were a rational multiple of $\pi$
then so would $\alpha$ be.
\end{proof}

\begin{lemma} \label{lemma:deltacoords}
$\alpha = \pi - 2\delta$, $\beta = 3\delta - \pi$ and $\gamma = \pi - \delta$,
and $\pi/3 < \delta < \pi/2$.  Normalizing $c=1$ and putting $s = a/c$ we have
$$s = 2\cos\delta, \qquad 0 < s < 1, \qquad b = 1-s^2, \qquad {\mbox and\ }
  \cos\alpha = 1 - \tfrac{1}{2}s^2 .$$
\end{lemma}

\begin{proof}
From $\delta = \alpha+\beta$ and $3\alpha+2\beta = \pi$ we get
$2\delta = 2\alpha+2\beta = \pi - \alpha$, so $\alpha = \pi-2\delta$ and
$\beta = \delta - \alpha = 3\delta-\pi$, and
$\gamma = \pi-\alpha-\beta = \pi-\delta$.  Then $\alpha>0$ gives
$\delta < \pi/2$ and $\beta>0$ gives $\delta > \pi/3$.  By the law of sines,
\begin{eqnarray*}
s &=& \frac{\sin\alpha}{\sin\gamma} = \frac{\sin 2\delta}{\sin\delta}
    = 2\cos\delta, \\
  b  &=& \frac{\sin\beta}{\sin\gamma} \ = \ \frac{-\sin 3\delta}{\sin\delta} \\
  &=&  4\sin^2\delta - 3 = 1 - 4\cos^2\delta \\
     &=& 1-s^2, \mbox{\qquad and\ } \\
\cos\alpha  &=& -\cos 2\delta = 1-2\cos^2\delta = 1-s^2/2.
\end{eqnarray*}
  Finally
$\pi/3<\delta<\pi/2$ gives $0<2\cos\delta<1$.
\end{proof}

\begin{lemma} \label{lemma:deltadirections}
 All differences of directions of edges of
the tiling lie in the subgroup of $\R/2\pi\Z$ generated by $\delta$ and $\pi$,
and the map
$$\Z \times \Z/2\Z \longrightarrow \R/2\pi\Z, \qquad
  (n,\varepsilon) \longmapsto n\delta + \varepsilon\pi$$
is injective.  Consequently, if we fix $AB$ as the base direction
$0$, every edge direction of the tiling is $n\delta+\varepsilon\pi$ for a unique
pair $(n,\varepsilon)$ with $\varepsilon \in \{0,1\}$.
\end{lemma}

\begin{proof}
Traversing a tile counterclockwise, consecutive edge directions differ by the
exterior angles $\pi-\gamma = \delta$, $\pi-\alpha = 2\delta$ and
$\pi-\beta = 2\pi-3\delta$, all in the subgroup generated by $\delta$; and the
two sides of any segment of the tiling give directions differing by $\pi$.  For
injectivity, suppose $n\delta+\varepsilon\pi \equiv 0 \pmod{2\pi}$.  If $n\ne0$
then $\delta/\pi \in \Q$, contrary to Lemma~\ref{lemma:incommensurable}; and if $n=0$ then
$\varepsilon\pi \equiv 0 \pmod{2\pi}$ forces $\varepsilon = 0$.
\end{proof}

\begin{definition} \label{definition:Phi}
For a direction $\theta = n\delta+\varepsilon\pi$ put
$$\Phi_1(\theta) = (-1)^{\varepsilon}, \qquad
  \Phi_2(\theta) = (-1)^{n+\varepsilon} .$$
For a region $R$ whose boundary edges have directions $\theta_1,\dots,\theta_r$
and lengths $\ell_1,\dots,\ell_r$ when $\partial R$ is traversed
counterclockwise, put $\Phi_i(R) = \sum_j \ell_j\,\Phi_i(\theta_j)$.
\end{definition}

Changing the base direction of Lemma~\ref{lemma:deltadirections} by $\delta$
shifts every $n$ by $1$ and so replaces $\Phi_2$ by $-\Phi_2$, leaving $\Phi_1$
unchanged; changing it by $\pi$ replaces both by their negatives.  Thus $\Phi_1$
and $\Phi_2$ are each defined only up to an overall sign, independently of one
another.  We fix this in Lemma~\ref{lemma:Phisigns}.

\begin{lemma} \label{lemma:Phitile}
Let $t$ be a tile of the tiling.  Then the three edges of $t$ have the same
$\varepsilon$, and their indices are $n$, $n+1$, $n+3$ for the edges $a$, $b$,
$c$ respectively, for some integer $n$ depending on the placement (or $n$,
$n-1$, $n-3$ if $t$ is a reflected copy).  Consequently
$$\Phi_1(t) = (-1)^{\varepsilon}(a+b+c), \qquad
  \Phi_2(t) = (-1)^{n+\varepsilon}(a-b-c),$$
so $|\Phi_1(t)| = a+b+c$ and $|\Phi_2(t)| = b+c-a$ for every tile.
\end{lemma}

\begin{proof}
Traverse $t$ counterclockwise starting along the edge $a$.  The exterior angles
at $C$, $A$, $B$ are $\pi-\gamma = \delta$, $\pi-\alpha = 2\delta$ and
$\pi-\beta = 2\pi-3\delta$, with coordinates $(1,0)$, $(2,0)$ and $(-3,0)$
respectively; these sum to $(0,0)$, as they must.  Since each has
$\varepsilon$-part $0$, $\varepsilon$ is constant on the three edges, and the
indices are $n$, $n+1$, $n+1+2 = n+3$, returning to $n+3-3 = n$.  Then
$\Phi_2$ takes the values $(-1)^{n+\varepsilon}$, $-(-1)^{n+\varepsilon}$,
$-(-1)^{n+\varepsilon}$ on the edges $a$, $b$, $c$.  Finally $b+c-a>0$ by the
triangle inequality.
\end{proof}

\begin{lemma} \label{lemma:Phiproperties}
$\Phi_1$ and $\Phi_2$ take the values $\pm1$ and satisfy
$$\Phi(\theta+\pi) = -\Phi(\theta), \qquad \Phi(-\theta) = \Phi(\theta).$$
\end{lemma}

\begin{proof}
Adding $\pi$ changes $\varepsilon$ by $1$ and leaves $n$ fixed, which negates
both; and 
$$-(n\delta+\varepsilon\pi) = (-n)\delta + \varepsilon\pi \pmod {2\pi},$$
 so evenness holds because $(-1)^{-n} = (-1)^n$.  For the last statement,
put $\chi(n) = \Phi(n\delta)$; then $\Phi(\theta+\pi) = -\Phi(\theta)$ forces
$\Phi(n\delta+\varepsilon\pi) = (-1)^{\varepsilon}\chi(n)$, and {\em a priori}
$\chi\colon \Z\to\{\pm1\}$ is arbitrary.  By Lemma~\ref{lemma:Phitile} the edges
$a$, $b$, $c$ of a placed tile lie in directions with indices $n$, $n+1$, $n+3$
(or $n$, $n-1$, $n-3$ for a reflected placement), all with the same
$\varepsilon$, so its tile value is
$$(-1)^{\varepsilon}\bigl[\pm a\,\chi(n) \pm b\,\chi(n+1) \pm c\,\chi(n+3)\bigr]$$
with signs fixed by the orientation convention.  The tile is scalene, so $a$,
$b$, $c$ are distinct, and $a+b>c$; hence the four numbers $a+b+c$, $-a+b+c$,
$a-b+c$, $a+b-c$ are positive and pairwise distinct, and the absolute value of
the tile value determines the triple
$(\chi(n),\chi(n+1),\chi(n+3))$ up to a global sign.  If that absolute value is
the same for every placement, then $\chi(n)\chi(n+1)$ is a constant $c_1$, so
$\chi(n) = \chi(0)c_1^{\,n}$; normalizing $\chi(0)=1$, $\chi$ is a homomorphism
$\Z\to\{\pm1\}$.  There are exactly two, given by $c_1 = \pm1$, namely the
trivial one and $(-1)^n$, yielding $\Phi_1$ and $\Phi_2$.
\end{proof}

\begin{lemma} \label{lemma:Phicancellation}
Let $ABC$ be tiled by the tiles $t_1,\dots,t_N$.  Then for $i = 1,2$,
$$\Phi_i(\partial ABC) \;=\; \sum_{j=1}^{N} \Phi_i(t_j) .$$
\end{lemma}

\begin{proof}
Each tile edge lies in a unique maximal segment of the tiling, and each maximal
segment either lies in a side of $ABC$ or is interior.  Group the tile edges by
the segment containing them.  On an interior segment, the tile edges on the two
sides have total length equal to the length of the segment, and their directions
differ by $\pi$; by Lemma~\ref{lemma:Phiproperties} the two contributions
cancel.  The remaining contributions are exactly those of $\partial ABC$.
\end{proof}

\begin{theorem} \label{theorem:Mdef}
There are integers $M_1$, $M_2$
with $|M_1|, |M_2| \le N$ and $M_i \equiv N \pmod 2$ such that
$$\Phi_1(\partial ABC) = M_1\,(a+b+c), \qquad
  \Phi_2(\partial ABC) = M_2\,(a-b-c).$$
Indeed $M_1 = \sum_j (-1)^{\varepsilon_j}$ and
$M_2 = \sum_j (-1)^{n_j+\varepsilon_j}$, summed over the tiles.
\end{theorem}

\begin{proof}
Immediate from Lemmas~\ref{lemma:Phitile} and~\ref{lemma:Phicancellation}: each
tile contributes $\pm(a+b+c)$ to the first sum and $\pm(a-b-c)$ to the second.
\end{proof}

\begin{remark} \label{remark:Phi1iscoloring}
$\Phi_1$ is the black-and-white coloring of Section~\ref{section:coloring}:
a tile is black or white according to the value of $\varepsilon$ on its edges,
and $M_1$ is the difference of the numbers of black and white tiles.
\end{remark}

\begin{lemma} \label{lemma:Phisigns}
Label $ABC$ counterclockwise, and normalize so that $\Phi_i = +1$ on
$AB$.  Then, with the angles of $ABC$ listed at $A$, $B$, $C$:
$$
\begin{array}{l| ccc|ccc}
 & \multicolumn{3}{c}{\Phi_1} & \multicolumn{3}{c}{\Phi_2} \\
\mbox{angles of } ABC & AB & BC & CA & AB & BC & CA \\
\hline
(2\alpha,\ \beta,\ \alpha+\beta)       & + & + & - & + & - & - \\
(\beta,\ 3\alpha,\ \beta)              & + & + & + & + & + & - \\
(\alpha+\beta,\ \alpha,\ \alpha+\beta) & + & + & - & + & + & + \\
(\alpha,\ \alpha+2\beta,\ \alpha)      & + & + & + & + & + & + \\
(2\alpha,\ \alpha,\ 2\beta)            & + & + & - & + & + & -
\end{array}
$$
\end{lemma}

\begin{proof}
Give $AB$ the $\delta$-coordinates $(0,0)$.  Turning left at $B$ by $\pi$
minus the angle of $ABC$ there gives the coordinates of $BC$, and turning
left at $C$ those of $CA$.  By Lemma~\ref{lemma:singamma} the turn $\pi-A$
has coordinates $(0,1)$ minus those of $A$, so its $\varepsilon$-part is
$1$ plus that of $A$; these need not vanish, which is why $\varepsilon$ is
not constant on $\partial ABC$ as it is on a tile.  For
instance in the case $(\beta,3\alpha,\beta)$ the turn at $B$ is
$\pi-3\alpha$, of coordinates $(6,-2)\equiv(6,0)$, so $BC$ has $(6,0)$; the
turn at $C$ is $\pi-\beta$, of coordinates $(-3,0)$, so $CA$ has $(3,0)$;
and the turn at $A$ returns the direction to $(0,0)$, as it must.  The
other rows are the same computation.
\end{proof}

Note that for the last two shapes the two rows agree, so
$\Phi_1(\partial ABC) = \Phi_2(\partial ABC)$; for the first three they
differ.  That is the whole source of what follows.

To compute the boundary values we need the sides of $ABC$ in terms of the
tile.  For that we need to compute the sines of the angles occurring in the five shapes of
Lemma~\ref{lemma:shapes}.

\begin{lemma} \label{lemma:sines}
With the tile normalized so that $c = 1$, so that $a = s$, $b = 1-s^2$ and
$\cos\alpha = 1 - s^2/2$ by Lemma~\ref{lemma:singamma}, we have
$$\frac{\sin\alpha}{\sin\delta} = s, \qquad
  \frac{\sin\beta}{\sin\delta} = 1-s^2, \qquad
  \frac{\sin\gamma}{\sin\delta} = \frac{\sin(\alpha+\beta)}{\sin\delta} = 1,$$
$$\frac{\sin 2\alpha}{\sin\delta} = \frac{\sin(\alpha+2\beta)}{\sin\delta}
   = s\,(2-s^2), \quad\mbox{\ and\quad }
  \frac{\sin 3\alpha}{\sin\delta} = \frac{\sin 2\beta}{\sin\delta}
   = s\,(1-s^2)(3-s^2).$$
\end{lemma}

\begin{proof}
Since $\alpha = \pi-2\delta$ we have
$\sin\alpha = \sin 2\delta = 2\sin\delta\cos\delta$, and $s = 2\cos\delta$,
so $\sin\alpha/\sin\delta = s$.  Since $\beta = 3\delta-\pi$,
$$\sin\beta = -\sin 3\delta = -\sin\delta\,(4\cos^2\delta - 1)
            = \sin\delta\,(1-s^2).$$
Since $\gamma = \pi-\delta$ and $\alpha+\beta = \delta$, $\sin \gamma = \sin(\alpha+\beta) = \sin\delta$.

Because $3\alpha+2\beta = \pi$, the angles $2\alpha$ and $\alpha+2\beta$ are
supplementary, and so are $3\alpha$ and $2\beta$; hence
$\sin 2\alpha = \sin(\alpha+2\beta)$ and $\sin 3\alpha = \sin 2\beta$.  Now
$$\frac{\sin 2\alpha}{\sin\delta} = 2\,\frac{\sin\alpha}{\sin\delta}\cos\alpha
  = 2s\Bigl(1-\frac{s^2}{2}\Bigr) = s\,(2-s^2),$$
and, since $4\cos^2\alpha - 1 = 4\bigl(1-\tfrac{1}{2}s^2\bigr)^2 - 1
= 3 - 4s^2 + s^4 = (1-s^2)(3-s^2)$,
$$\frac{\sin 3\alpha}{\sin\delta}
  = \frac{\sin\alpha}{\sin\delta}\,(4\cos^2\alpha - 1)
  = s\,(1-s^2)(3-s^2). \qedhere$$
\end{proof}

\begin{theorem} \label{theorem:M1M2}
Let $ABC$ be one of the five shapes of Lemma~\ref{lemma:shapes}, tiled by
copies of $T$, and let $\Lambda$ be the ratio of the circumdiameter of $ABC$ to
that of $T$.  Normalize the tile by $c=1$ as in Lemma~\ref{lemma:singamma} and
the signs of $\Phi_1$, $\Phi_2$ as in Lemma~\ref{lemma:Phisigns}.  Then, with 
the angles of $ABC$ listed in counterclockwise order, 
$$\begin{array}{lll}
\mbox{angles of } ABC & M_1 & M_2 \\[2pt]
\hline
(2\alpha,\ \beta,\ \alpha+\beta)       & \Lambda s & \Lambda s \\
(\beta,\ 3\alpha,\ \beta)              & \Lambda(1-s)(1+s)^2 & -\Lambda(1-s)^2(1+s) \\
(\alpha+\beta,\ \alpha,\ \alpha+\beta) & \Lambda/(1+s) & -\Lambda/(1-s) \\
(\alpha,\ \alpha+2\beta,\ \alpha)      & \Lambda s(2+s)/(1+s) & -\Lambda s(2-s)/(1-s) \\
(2\alpha,\ \alpha,\ 2\beta)            & \Lambda s(1-s)(2+s) & -\Lambda s(2-s)(1+s)
\end{array}$$
In particular $M_1 \ne 0$ and $M_2 \ne 0$, and $\Lambda = M_1/s$ in the first
row.
\end{theorem}

\begin{proof}
By the law of sines each side of $ABC$ is $\Lambda \sin A'/\sin\delta$, where
$A'$ is the opposite angle, so Lemma~\ref{lemma:sines} gives each side as
$\Lambda$ times a polynomial in $s$.  Combining with the signs of
Lemma~\ref{lemma:Phisigns},
$$\Phi_1(\partial ABC) = M_1\,(a+b+c), \qquad
  \Phi_2(\partial ABC) = M_2\,(a-b-c),$$
where, by Lemma~\ref{lemma:singamma},
$$a+b+c = 2+s-s^2 = (2-s)(1+s), \qquad
  a-b-c = s^2+s-2 = -(2+s)(1-s).$$
We compute the five rows in turn, writing each side of $ABC$ as $\Lambda$ times
its value from Lemma~\ref{lemma:sines}.

{\em Case $(2\alpha,\ \beta,\ \alpha+\beta)$.}  Here
$|AB| = \Lambda$, $|BC| = \Lambda s(2-s^2)$, $|CA| = \Lambda(1-s^2)$, with signs
$+,+,-$ for $\Phi_1$ and $+,-,-$ for $\Phi_2$.  So
\begin{eqnarray*}
\Phi_1(\partial ABC) &=& \Lambda\bigl[1 + s(2-s^2) - (1-s^2)\bigr]
   \ = \ \Lambda s\,(2+s-s^2), \\
\Phi_2(\partial ABC) &=& \Lambda\bigl[1 - s(2-s^2) - (1-s^2)\bigr]
   \ = \ \Lambda s\,(s^2+s-2),
\end{eqnarray*}
giving $M_1 = M_2 = \Lambda s$.

{\em Case $(\beta,\ 3\alpha,\ \beta)$.}  Here
$|AB| = |BC| = \Lambda(1-s^2)$ and $|CA| = \Lambda s(1-s^2)(3-s^2)$, with signs
$+,+,+$ for $\Phi_1$ and $+,+,-$ for $\Phi_2$.  Since
$2+3s-s^3 = (2-s)(1+s)^2$ and $2-3s+s^3 = (2+s)(1-s)^2$,
\begin{eqnarray*}
\Phi_1(\partial ABC) &=& \Lambda(1-s^2)\bigl[2 + s(3-s^2)\bigr]
   \ = \ \Lambda(1-s^2)(2-s)(1+s)^2, \\
\Phi_2(\partial ABC) &=& \Lambda(1-s^2)\bigl[2 - s(3-s^2)\bigr]
   \ = \ \Lambda(1-s^2)(2+s)(1-s)^2,
\end{eqnarray*}
so 
$$M_1 = \Lambda(1-s^2)(1+s) = \Lambda(1-s)(1+s)^2 \mbox{\ and} $$
$$M_2 = -\Lambda(1-s^2)(1-s) = -\Lambda(1-s)^2(1+s).$$

{\em Case $(\alpha+\beta,\ \alpha,\ \alpha+\beta)$.}  Here
$|AB| = |BC| = \Lambda$ and $|CA| = \Lambda s$, with signs $+,+,-$ for $\Phi_1$
and $+,+,+$ for $\Phi_2$.  So
$$\Phi_1(\partial ABC) = \Lambda\,(2-s), \qquad
  \Phi_2(\partial ABC) = \Lambda\,(2+s),$$
giving $M_1 = \Lambda/(1+s)$ and $M_2 = -\Lambda/(1-s)$.

{\em Case $(\alpha,\ \alpha+2\beta,\ \alpha)$.}  Here
$|AB| = |BC| = \Lambda s$ and $|CA| = \Lambda s(2-s^2)$, with signs $+,+,+$ for
both $\Phi_1$ and $\Phi_2$.  So both boundary values equal
$$\Lambda s\bigl[2 + (2-s^2)\bigr] \ = \ \Lambda s\,(4-s^2)
   \ = \ \Lambda s\,(2-s)(2+s),$$
giving $M_1 = \Lambda s(2+s)/(1+s)$ and $M_2 = -\Lambda s(2-s)/(1-s)$.

{\em Case $(2\alpha,\ \alpha,\ 2\beta)$.}  Here
$|AB| = \Lambda s(1-s^2)(3-s^2)$, $|BC| = \Lambda s(2-s^2)$ and
$|CA| = \Lambda s$, with signs $+,+,-$ for both $\Phi_1$ and $\Phi_2$.  So both
boundary values equal
$$\Lambda s\bigl[(1-s^2)(3-s^2) + (2-s^2) - 1\bigr]
  \ = \ \Lambda s\,(s^4-5s^2+4)
  \ = \ \Lambda s\,(1-s^2)(4-s^2),$$
giving $M_1 = \Lambda s(1-s)(2+s)$ and $M_2 = -\Lambda s(2-s)(1+s)$.

Each entry is nonzero because $0<s<1$.
\end{proof}

\begin{definition} \label{definition:relationlattice}
The \emph{relation lattice} of the tile $T$ with sides $a$, $b$, $c$ is
$$L \;=\; \{(P,Q,R) \in \Z^3 : Pa+Qb+Rc = 0\},$$
the kernel of the homomorphism $\Z^3 \to \R$ given by
$(P,Q,R) \mapsto Pa+Qb+Rc$.  A \emph{relation} is a nonzero element of $L$.
\end{definition}

\begin{lemma} \label{lemma:relationlattice}
$L$ is a subgroup of $\Z^3$ such that if $k\nu \in L$ for some nonzero integer
$k$ then $\nu \in L$.
 When $3\alpha+2\beta=\pi$,
so that $a/c = s$ and $b/c = 1-s^2$, a vector $(P,Q,R)$ lies in $L$ if and only
if $Qs^2 - Ps - (Q+R) = 0$.
\end{lemma}

\begin{proof}
If $k\nu \in L$ with $k \ne 0$, then $k\,(\nu\cdot(a,b,c)) = 0$ in $\R$, so
$\nu\cdot(a,b,c) = 0$ and $\nu \in L$.

Let $V \subseteq \Q^3$ be the set of rational $(P,Q,R)$ with $Pa+Qb+Rc = 0$; it
is a $\Q$-subspace of dimension $3 - \dim_\Q\operatorname{span}_\Q\{a,b,c\}$,
being the kernel of the linear map $(P,Q,R) \mapsto Pa+Qb+Rc$.  Clearing
denominators carries each element of $V$ to an integer multiple of itself lying
in $L$, so $L$ spans $V$ over $\Q$. 
 The last statement follows by dividing
$Pa+Qb+Rc = 0$ by $c$ and substituting $a/c = s$, $b/c = 1-s^2$, since
$\operatorname{span}_\Q\{a,b,c\} = c\cdot\operatorname{span}_\Q\{1,s,s^2\}$.
\end{proof}

\subsection{Geometric interpretation of $\Phi$}
This subsection is for readers wishing to see pictures of $\Phi$ and develop a geometric
intuition for it.  Nothing below depends on it. 
Since the coloring interpretation $\Phi_1$ takes two values, and $\Phi$ takes
two values, together they take four values.  If $n\delta + \varepsilon\pi$ is
the direction of the $a$ edge of a tile, directed so that the tile is described
counterclockwise, then $\Phi_1$ is determined by $\varepsilon$ and $\Phi$ by
the parity of $n+\varepsilon$.  Therefore if we want to color a tiling showing
the values of both interpretations, we will need four colors.  The black tiles
we used for $\Phi_1$ will become black or red, and the white tiles will become
white or yellow.

A placed tile has its $a$, $b$, $c$ edges at indices $n$, $n+1$, $n+3$ with a
common $\varepsilon$, and by Lemma~\ref{lemma:Phitile} its tile values are
$\Phi_1(t) = (-1)^{\varepsilon}(a+b+c)$ and
$\Phi_2(t) = (-1)^{n+\varepsilon}(a-b-c)$.  So $\Phi_1$ depends only on
$\varepsilon$ and $\Phi_2$ only on the parity of $n+\varepsilon$.  An isometry
carrying one tile to another acts as follows, the last two columns giving the
factor by which the tile value is multiplied; the mirror in the last row has
direction $p\delta + \eta\pi$.
$$
\begin{array}{lcccc}
\mbox{isometry} & \varepsilon \mapsto & n \mapsto & \Phi_1 & \Phi_2 \\
\hline
\mbox{translation}                  & \varepsilon   & n    & + & + \\
\mbox{rotation by } \delta          & \varepsilon   & n+1  & + & - \\
\mbox{rotation by } \pi             & \varepsilon+1 & n    & - & - \\
\mbox{reflection in a lattice line} & \varepsilon+1 & 2p-n & - & -  \\
\mbox{rotation by } \delta+\pi      & \varepsilon+1 & n+1  & - & + \
\end{array}
$$
In the last row $2p-n \equiv n \pmod 2$, so $n+\varepsilon$ changes parity and
$\Phi_2$ is negated, whatever the mirror line.

The colors are the four values of the pair $(\Phi_1,\Phi)$:
$$R = (+,+), \qquad B = (+,-), \qquad Y = (-,+), \qquad W = (-,-),$$
so that the black class of $\Phi_1$ is $\{B,R\}$ and the white class is
$\{W,Y\}$.  An isometry carrying one tile of the tiling to another changes the
color as follows:
$$
\begin{array}{lcccc}
\mbox{isometry} & B & R & W & Y \\
\hline
\mbox{translation}                 & B & R & W & Y \\
\mbox{rotation by } \delta         & R & B & Y & W \\
\mbox{rotation by } \pi            & Y & W & R & B \\
\mbox{reflection in a lattice line}& Y & W & R & B \\
\mbox{rotation by } \delta+\pi     & W & Y & B & R
\end{array}
$$

\begin{figure} [ht]
\begin{center} \includegraphics[width=0.8\textwidth]{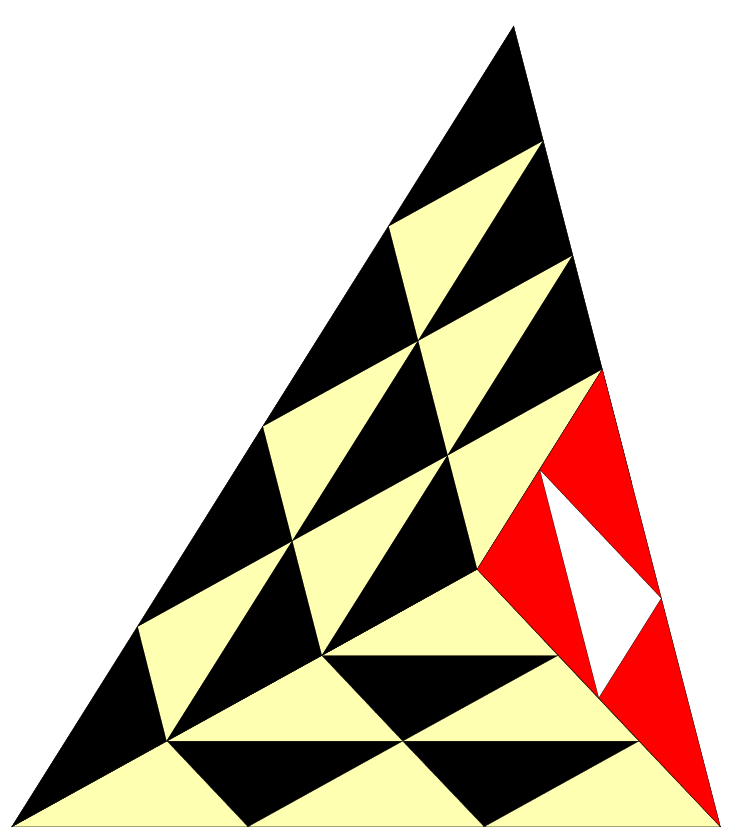}
\end{center}
\caption{A tiling colored to show values of $\Phi_1$ and $\Phi$}
\label{figure:PhiColoring}
\end{figure}

Such a coloring is illustrated in Figure~\ref{figure:PhiColoring}.

Write $R$, $B$, $Y$, $W$ also for the numbers of tiles of each color.  Since
$\Phi_1 = +1$ on $R$ and $B$ and $\Phi = +1$ on $R$ and $Y$,
$$N = R+B+Y+W, \qquad M_1 = R+B-Y-W, \qquad M_2 = R+Y-B-W,$$
and therefore
$$R-W = \frac{M_1+M_2}{2}, \qquad B-Y = \frac{M_1-M_2}{2},$$
both integers because $M_1 \equiv M_2 \equiv N \pmod 2$.  For a triquadratic
$ABC$ we have $M_1 = M_2 = M$ by Theorem~\ref{theorem:M1M2}, so in that case
$$R - W = M \qquad\mbox{and}\qquad B = Y.$$
In Figure~\ref{figure:coloring}, $N = 28$ and $M = 2$, and the counts are
$R = 3$, $B = 12$, $Y = 12$, $W = 1$, so indeed $R-W = 2$ and $B = Y$.

The remaining combination $(R+W)-(B+Y)$ is not determined by $N$, $M_1$ and
$M_2$.  It is the sum over tiles of the product $\Phi_1(t)\Phi(t)$, since that
product is $+1$ on $R$ and $W$ and $-1$ on $B$ and $Y$; and as a function of
direction, $\Phi_1\Phi$ sends $n\delta+\varepsilon\pi$ to $(-1)^n$.  That
function assigns each tile a well-defined value $(-1)^n(a-b-c)$, but it takes
the \emph{same} value at $\theta$ and at $\theta+\pi$, so the two sides of an
interior segment reinforce instead of cancelling and
Lemma~\ref{lemma:Phicancellation} fails for it.  There is therefore no boundary
formula for this combination, and the four counts are not determined by the
invariants.

\section{Rationality of the tile}
The argument in this section is uniform.  Both $M_1$ and $M_2$ are proportional
to $\Lambda$ and $N$ is proportional to $\Lambda^2$, so the quantities
$$\kappa = \frac{M_1+M_2}{M_2-M_1} \qquad\mbox{and}\qquad \frac{N}{M_1M_2}$$
do not involve $\Lambda$ at all: each is a rational number and, by
Theorem~\ref{theorem:M1M2} and Lemma~\ref{lemma:shapeareas}, an explicit
rational function of the single unknown $s$.  For two of the five shapes the
first of these already gives $s$ outright.  For two others it gives only a
quadratic relation, and the second quantity is needed to finish.  For the
triquadratic shape $M_1 = M_2$, so $\kappa$ does not exist and only the area
equation is available; it gives $s^2$, and one further ingredient --- a
relation with a nonzero $a$-coefficient --- gives $s$.

We keep the notation of the previous section: the tile is normalized so that
$c=1$, $a=s$ and $b=1-s^2$ with $0<s<1$ (Lemma~\ref{lemma:singamma}),
$\Lambda$ denotes the ratio of the circumdiameter of $ABC$ to that of $T$,
$\delta/\pi$ is irrational, and the signs of $\Phi_1$ and $\Phi_2$ are
normalized as in Lemma~\ref{lemma:Phisigns}, so that $\Phi_i = +1$ on $AB$.
 $T$ is \emph{rational} when $a/c$, $b/c$ and $a/b$ are all
rational, which by Lemma~\ref{lemma:singamma} happens exactly when
$s \in \Q$.

\begin{lemma} \label{lemma:shapeareas}
Let $ABC$ be tiled by $N$ copies of $T$.  Then, with the angles of $ABC$
listed at $A$, $B$, $C$,
$$
\begin{array}{ll}
\mbox{angles of } ABC & N/\Lambda^2 \\
\hline
(2\alpha,\ \beta,\ \alpha+\beta)       & 2-s^2 \\[2pt]
(\beta,\ 3\alpha,\ \beta)              & (1-s^2)^2(3-s^2) \\[2pt]
(\alpha+\beta,\ \alpha,\ \alpha+\beta) & \dfrac{1}{1-s^2} \\[6pt]
(\alpha,\ \alpha+2\beta,\ \alpha)      & \dfrac{s^2(2-s^2)}{1-s^2} \\[6pt]
(2\alpha,\ \alpha,\ 2\beta)            & s^2(2-s^2)(3-s^2)
\end{array}
$$
By Lemma~\ref{lemma:shapes} these are the only possible shapes.
\end{lemma}

\begin{proof}
A triangle with angles $A$, $B$, $C$ and circumradius $R$ has area
$2R^2\sin A\sin B\sin C$, so
$$\frac{N}{\Lambda^2} \;=\;
  \frac{\sin A\,\sin B\,\sin C}{\sin\alpha\,\sin\beta\,\sin\gamma},$$
where $A$, $B$, $C$ are the angles of $ABC$.  Divide all six sines by
$\sin\delta$; by Lemma~\ref{lemma:sines} the values of these ratios are
\begin{eqnarray*}
  \frac{\sin\alpha}{\sin\delta} &=& s, \\
  \frac{\sin\beta}{\sin\delta} &=& 1-s^2, \\
  \frac{\sin\gamma}{\sin\delta} &=& \frac{\sin(\alpha+\beta)}{\sin\delta}
     \ =\ 1, \\
  \frac{\sin 2\alpha}{\sin\delta} &=&
     \frac{\sin(\alpha+2\beta)}{\sin\delta} \ =\ s\,(2-s^2), \\
  \frac{\sin 2\beta}{\sin\delta} &=& \ \frac{\sin 3\alpha}{\sin\delta}
    \ =\ s\,(1-s^2)(3-s^2).
\end{eqnarray*}
In particular the denominator is the same in all five cases, namely
$$\frac{\sin\alpha\,\sin\beta\,\sin\gamma}{\sin^3\delta}
   \;=\; s\,(1-s^2)\cdot 1 \;=\; s\,(1-s^2).$$
It remains to compute the numerator in each case and divide.

{\em Case $(2\alpha,\ \beta,\ \alpha+\beta)$.}
$$\frac{N}{\Lambda^2} \;=\; \frac{s(2-s^2)\cdot(1-s^2)\cdot 1}{s(1-s^2)}
   \;=\; 2-s^2 .$$

{\em Case $(\beta,\ 3\alpha,\ \beta)$.}
$$\frac{N}{\Lambda^2}
   \;=\; \frac{(1-s^2)\cdot s(1-s^2)(3-s^2)\cdot(1-s^2)}{s(1-s^2)}
   \;=\; (1-s^2)^2(3-s^2) .$$

{\em Case $(\alpha+\beta,\ \alpha,\ \alpha+\beta)$.}
$$\frac{N}{\Lambda^2} \;=\; \frac{1\cdot s\cdot 1}{s(1-s^2)}
   \;=\; \frac{1}{1-s^2} .$$

{\em Case $(\alpha,\ \alpha+2\beta,\ \alpha)$.}
$$\frac{N}{\Lambda^2} \;=\; \frac{s\cdot s(2-s^2)\cdot s}{s(1-s^2)}
   \;=\; \frac{s^2(2-s^2)}{1-s^2} .$$

{\em Case $(2\alpha,\ \alpha,\ 2\beta)$.}
$$\frac{N}{\Lambda^2}
   \;=\; \frac{s(2-s^2)\cdot s\cdot s(1-s^2)(3-s^2)}{s(1-s^2)}
   \;=\; s^2(2-s^2)(3-s^2) .$$
\end{proof}

\begin{remark}
For the tile $(a,b,c) = (2,3,4)$, so $s = 1/2$, the five shapes give
$$
\begin{array}{lrrrr}
\mbox{angles of } ABC & \Lambda & N & M_1 & M_2 \\
\hline
(2\alpha,\ \beta,\ \alpha+\beta)       & 2   & 7  & 1 & 1 \\
(\beta,\ 3\alpha,\ \beta)              & 8/3 & 11 & 3 & -1 \\
(\alpha+\beta,\ \alpha,\ \alpha+\beta) & 3/2 & 3  & 1 & -3 \\
(\alpha,\ \alpha+2\beta,\ \alpha)      & 6   & 21 & 5 & -9 \\
(2\alpha,\ \alpha,\ 2\beta)            & 8   & 77 & 5 & -9
\end{array}
$$
The reader may check that $\kappa = -1/2$, $1/2$, $2/7$, $2/7$ in the last four
rows, recovering $s = 1/2$ in each case by
Theorem~\ref{theorem:rationaliso} or Theorem~\ref{theorem:shape45}, and that
$2M^2/(N+M^2) = 2/8 = 1/4 = s^2$ in the first.
\end{remark}

\begin{definition} \label{definition:kappa}
When $M_1 \ne M_2$ we put
$$\kappa \;=\; \frac{M_1+M_2}{M_2-M_1} \;\in\; \Q .$$
\end{definition}

\begin{lemma} \label{lemma:kappa}
Let $ABC$ be tiled by copies of $T$.  If $ABC$ is triquadratic then
$M_1 = M_2$ and $\kappa$ is undefined.  Otherwise $\kappa$ is defined and
$$\kappa = -s, \qquad \kappa = s, \qquad \kappa = \frac{s}{2-s^2}, \qquad
  \kappa = \frac{s}{2-s^2}$$
for the shapes $(\beta,3\alpha,\beta)$,
$(\alpha+\beta,\alpha,\alpha+\beta)$, $(\alpha,\alpha+2\beta,\alpha)$ and
$(2\alpha,\alpha,2\beta)$ respectively.  In the last two cases
$0 < \kappa < 1$, and
$$s = \kappa\,(2-s^2), \qquad\text{equivalently}\qquad
  \kappa^2 s^4 - (4\kappa^2+1)s^2 + 4\kappa^2 = 0 .$$
\end{lemma}

\begin{proof}
For the last two shapes the two rows of Lemma~\ref{lemma:Phisigns} agree, so
$\Phi_1(\partial ABC) = \Phi_2(\partial ABC)$, that is,
$M_1(a+b+c) = M_2(a-b-c)$, and therefore
$$\frac{M_2}{M_1} = \frac{a+b+c}{a-b-c}
  = \frac{2+s-s^2}{-(2-s-s^2)} = -\frac{(2-s)(1+s)}{(2+s)(1-s)} .$$
Hence $M_1+M_2$ and $M_2-M_1$ are $M_1$ times
\begin{eqnarray*}
\frac{(2+s)(1-s)-(2-s)(1+s)}{(2+s)(1-s)} &=& \frac{-2s}{(2+s)(1-s)}
   \qquad\mbox{and} \\
\frac{-(2-s)(1+s)-(2+s)(1-s)}{(2+s)(1-s)} &=& \frac{-2(2-s^2)}{(2+s)(1-s)}
\end{eqnarray*}
respectively, so $\kappa = s/(2-s^2)$; in particular $M_1 \ne M_2$, and
$0<s<1$ gives $1 < 2-s^2 < 2$ and hence $0<\kappa<1$.  Consequently
$\kappa(2-s^2)+s > 0$, so
$$s = \kappa(2-s^2) \iff \kappa^2(2-s^2)^2 - s^2 = 0
  \iff \kappa^2 s^4 - (4\kappa^2+1)s^2 + 4\kappa^2 = 0,$$
since the middle expression is
$\bigl[\kappa(2-s^2)-s\bigr]\bigl[\kappa(2-s^2)+s\bigr]$.

For the first three shapes read $M_1$ and $M_2$ off
Theorem~\ref{theorem:M1M2}: in the triquadratic case $M_1 = M_2 = \Lambda s$;
for $(\beta,3\alpha,\beta)$,
\begin{eqnarray*}
M_1+M_2 &=& \Lambda(1-s^2)\bigl[(1+s)-(1-s)\bigr] = 2\Lambda s(1-s^2), \\
  M_2-M_1 &=& -\Lambda(1-s^2)\bigl[(1-s)+(1+s)\bigr] = -2\Lambda(1-s^2),
\end{eqnarray*}
so $\kappa = -s$; and for $(\alpha+\beta,\alpha,\alpha+\beta)$, over the common
denominator $1-s^2$,
$$M_1+M_2 = \frac{-2\Lambda s}{1-s^2}, \qquad
  M_2-M_1 = \frac{-2\Lambda}{1-s^2},$$
so $\kappa = s$.
\end{proof}

\begin{theorem} \label{theorem:rationaliso}
If a triangle of shape $(\beta,3\alpha,\beta)$ or
$(\alpha+\beta,\alpha,\alpha+\beta)$ is tiled by copies of $T$, then $T$ is
rational, and 
$$s = \frac{M_1+M_2}{M_1-M_2} \qquad\text{and}\qquad
  s = \frac{M_1+M_2}{M_2-M_1}$$
respectively.
\end{theorem}

\begin{proof}
By Lemma~\ref{lemma:kappa}, $s = -\kappa$ in the first case and $s = \kappa$ in
the second, and $\kappa$ is a ratio of integers.
\end{proof}

\begin{theorem} \label{theorem:shape45}
Let $ABC$ have shape $(\alpha,\alpha+2\beta,\alpha)$ or
$(2\alpha,\alpha,2\beta)$ and suppose $ABC$ is tiled by $N$ copies of $T$.
Then $N+M_1M_2 \ne 0$, and
$$s^2 = \frac{2\,(2N+M_1M_2)}{N+M_1M_2} \quad\text{in the first case}, \qquad
  s^2 = \frac{(M_1+M_2)^2}{2\,(N+M_1M_2)} \quad\text{in the second}.$$
In either case $s^2 \in \Q$, hence $s = \kappa(2-s^2) \in \Q$ by
Lemma~\ref{lemma:kappa}, and $T$ is rational.
\end{theorem}

\begin{proof}
Write $w = s^2$ and $t = N/(M_1M_2)$.  Both $\kappa$ and $t$ are independent
of $\Lambda$, since $M_1$ and $M_2$ are proportional to $\Lambda$ and $N$ to
$\Lambda^2$, and both are rational.  By Theorem~\ref{theorem:M1M2} and
Lemma~\ref{lemma:shapeareas},
$$t = -\frac{2-w}{4-w} \quad\text{for } (\alpha,\alpha+2\beta,\alpha),
\qquad
t = -\frac{(2-w)(3-w)}{(1-w)(4-w)} \quad\text{for } (2\alpha,\alpha,2\beta),$$
and in both cases $t \ne -1$, since $t=-1$ would give $2=4$ in the first case
and $6=4$ in the second; equivalently $N+M_1M_2 \ne 0$.

In the first case, solving $t = -(2-w)/(4-w)$ for $w$ gives
$w = (4t+2)/(t+1)$, which on substituting $t = N/(M_1M_2)$ is the first
displayed formula.

In the second case put $y = w^2-5w$, so that $t = -(y+6)/(y+4)$ and hence
$y \in \Q$ with $y+4 = -2/(t+1)$.  Substituting $w^2 = y+5w$ into the quartic
$\kappa^2w^2-(4\kappa^2+1)w+4\kappa^2 = 0$ of Lemma~\ref{lemma:kappa} and
collecting terms gives
$$(\kappa^2-1)\,w + \kappa^2(y+4) = 0 ,$$
and $\kappa^2 \ne 1$ because $0<\kappa<1$, so
$w = -2\kappa^2/\bigl((1-\kappa^2)(t+1)\bigr) \in \Q$.  Finally
$$1-\kappa^2 = \frac{(M_2-M_1)^2-(M_1+M_2)^2}{(M_2-M_1)^2}
             = \frac{-4M_1M_2}{(M_2-M_1)^2}, \qquad
  t+1 = \frac{N+M_1M_2}{M_1M_2},$$
so $(1-\kappa^2)(t+1) = -4(N+M_1M_2)/(M_2-M_1)^2$, and substituting gives the
second displayed formula.  In both cases $w \in \Q$ and $s = \kappa(2-w)$ with
$\kappa \in \Q$, so $s \in \Q$ and $b/c = 1-s^2 \in \Q$.
\end{proof}

\begin{theorem}\label{theorem:arelation}  
Let triangle $ABC$ be $N$-tiled by a tile with sides $(a,b,c)$, whose 
 angles are $(\alpha,\beta,\gamma)$ with $3\alpha + 2\beta = \pi$.
Suppose $ABC$ is not similar to the tile. 
 Then there is an $a$-relation $ja = pb + qc$,  where $p,q,j$ 
 are nonnegative integers with  $j > 0$.
\end{theorem}

\begin{proof}
We argue by cases on the possible shapes of $ABC$ as given in Lemma~\ref{lemma:shapes}.
For the three isosceles shapea and the shape $(2\alpha, \alpha,2\beta)$,  
Theorem~\ref{theorem:rationaliso} and Theorem~\ref{theorem:shape45}
tell us that the tile is rational, and a rational tile has an $a$-relation $ka=mc$ with $a/c = m/k$.
That leaves only the triquadratic shape $(2\alpha,\beta,\alpha+\beta)$ to take up here.

Suppose, for proof by contradiction, that there is no $a$-relation.
Suppose there is an $a$ edge on a side $AB$ of $ABC$.  Either $AB$ is entirely composed of $a$ 
edges, or has no $a$ edges,  or there is a vertex $Q$ with an $a$ edge on one side and a non-$a$ edge on the other.
Then, by the same reasoning as for a link $PQ$, there is an outgoing link from $Q$ into the interior
of $ABC$, which would be a contradiction.  Therefore,  every side of triangle  $ABC$ is either 
composed entirely of $a$ edges or has no $a$ edges.  

Suppose a side $PQ$ of $ABC$ has only $a$ edges.
By Lemma~\ref{lemma:angles}, no tile at a vertex of $ABC$ has its $\gamma$ angle there. 
By Lemma~\ref{lemma:gammagreaterpiover2}, $\gamma > \pi/2$.  Hence there is exactly one $\gamma$ angle
at each interior vertex on $PQ$ and none at the endpoints.  Since there is one more tile than 
interior vertices, that violates the pigeonhole principle.  That contradiction shows that no side 
of $ABC$ is composed only of $a$ edges.  Therefore,  there is no $a$ edge anywhere on the boundary of $ABC$.

The triquadratic shape has a standalone $\beta$ vertex.  By Lemma~\ref{lemma:shapes},
there is just one tile at that vertex.
That tile must have its $a$ side on one of the adjacent sides
of $ABC$.  But we have proved that there is no $a$ edge on the boundary of $ABC$, so 
that is a contradiction.
\end{proof}

\subsection{The triquadratic shape}

Here $ABC$ has angles $(2\alpha,\beta,\alpha+\beta)$, with $2\alpha$ at $A$,
$\beta$ at $B$ and $\alpha+\beta$ at $C$.  Since
$\sin\gamma = \sin(\alpha+\beta)$, the law of sines gives
$$|BC| = 2\Lambda a\cos\alpha = \Lambda s(b+c), \qquad
  |CA| = \Lambda b, \qquad |AB| = \Lambda c ,$$
and by Theorem~\ref{theorem:M1M2}, $M_1 = M_2 = \Lambda s$; we write $M$ for
this common value.  Thus the two invariants coincide here, and the arguments of
Theorem~\ref{theorem:shape45} are unavailable.

\begin{theorem} \label{theorem:tqsquare}
If a triquadratic triangle is tiled by $N$ copies of $T$, then
$$s^2 = \frac{2M^2}{N+M^2} .$$
In particular $s^2$ and $b/c = 1-s^2$ are rational,  and $N > M^2$.
\end{theorem}

\begin{proof}
$M = \Lambda s$ and $0<s<1$, so $M \ne 0$ and $\Lambda = M/s$.  Substituting
into Lemma~\ref{lemma:shapeareas} gives $Ns^2 = M^2(2-s^2)$, that is,
$s^2(N+M^2) = 2M^2$.  Since $M$ and $N$ are integers, $s^2 \in \Q$, and
$s^2<1$ gives $N>M^2$.  Finally $1-s^2 \in \Q$ provides a nonzero element
$(0,Q,R)$ of $L$.
\end{proof}

\begin{theorem} \label{theorem:tqrational}
Suppose a triquadratic triangle is tiled by $N$ copies of $T$.  Then $T$ is
rational, and $2(N+M^2)$ is a perfect square.
\end{theorem}

\begin{proof}
By Theorem~\ref{theorem:arelation} the tiling witnesses an $a$-relation, so some
segment carries only $a$ edges on one side, giving $ja = Qb+Rc$ with $j \ge 1$;
that is, $L$ contains a vector with nonzero first coordinate.  By
Lemma~\ref{lemma:relationlattice} such a vector $(P,Q,R)$ satisfies
$Qs^2 - Ps - (Q+R) = 0$, and $s^2 \in \Q$ by Theorem~\ref{theorem:tqsquare}, so
$Ps = Qs^2 - (Q+R) \in \Q$; since $P \ne 0$, $s \in \Q$ and $T$ is rational.
For the last statement put $D = N+M^2$, so that $s = M\sqrt{2/D}$ by
Theorem~\ref{theorem:tqsquare}; as $s \in \Q$, $2/D = 2D/D^2$ is the square of
a rational, so $2D$ is a perfect square.
\end{proof}

\begin{corollary} \label{corollary:tqparameters}
If a triquadratic triangle is tiled by $N$ copies of tile, written
$(a,b,c) = \tau(mk,\,k^2-m^2,\,k^2)$ as in Lemma~\ref{lemma:gsquare} with
$\gcd(m,k)=1$ and $m<k$, then $m \mid M$, and with $j = M/m$,
$$N = j^2(2k^2-m^2), \qquad M = jm, \qquad \Lambda = jk,$$
$$|BC| = jm(2k^2-m^2)\,\tau, \qquad |CA| = jk(k^2-m^2)\,\tau, \qquad
  |AB| = jk^3\,\tau .$$
\end{corollary}

\begin{proof}
Here $s = m/k$, so Theorem~\ref{theorem:tqsquare} gives
$N = M^2(2k^2-m^2)/m^2$, and in particular $m^2 \mid M^2(2k^2-m^2)$.  Since
$\gcd(m,k)=1$ we have
$\gcd(m^2,\,2k^2-m^2) = \gcd(m^2,\,2k^2) = \gcd(m^2,2)$.  If $m$ is odd this
is $1$, so $m^2 \mid M^2$ and $m \mid M$.  If $m$ is even, write $m = 2m_1$;
then $2m_1^2 \mid M^2$, which gives $v_p(m_1) \le v_p(M)$ for every odd prime
$p$ and $2v_2(m_1)+1 \le 2v_2(M)$, hence
$v_2(M) \ge v_2(m_1)+1 = v_2(m)$; so again $m \mid M$.  Putting $j = M/m$ gives
$N = j^2(2k^2-m^2)$ and $\Lambda = M/s = jk$, and the side lengths follow from
the formulas preceding Theorem~\ref{theorem:tqsquare}.
\end{proof}

\begin{theorem} \label{theorem:rationaltile}
Let $ABC$ be a triangle, not similar to $T$, tiled by copies of a triangle $T$
with $3\alpha+2\beta = \pi$.  Then $T$ is rational.
\end{theorem}

\begin{proof}
By Lemma~\ref{lemma:shapes}, $ABC$ has one of five shapes.  
$(\beta,3\alpha,\beta)$ and $(\alpha+\beta,\alpha,\alpha+\beta)$ apply
Theorem~\ref{theorem:rationaliso}; for $(\alpha,\alpha+2\beta,\alpha)$ and
$(2\alpha,\alpha,2\beta)$ apply Theorem~\ref{theorem:shape45}; for
$(2\alpha,\beta,\alpha+\beta)$ apply Theorem~\ref{theorem:tqrational}.
\end{proof}

\section{The connected components of a tiling}

Suppose given an $N$-tiling of triangle $ABC$, by the tile with angles $\alpha$, $\beta$, and $\gamma$,
where $3\alpha + 2\beta = \pi$ and $\alpha$ is not a rational multiple of $\pi$.  We define a graph $\HH$
whose nodes are the tiles of the tiling.  Since the word ``edges'' is already in use for the sides of a tile,
we shall refer to the edges of this graph as ``connections'' instead, and speak of one tile being connected 
to another.  

\begin{definition}  Given a tiling of $ABC$, the graph $\HH$ has for its nodes the tiles used in the tiling,
and its edges (connections)  are defined as follows:  Tiles $T$ and $S$ are {\em \bf connected} if   $T$ and $S$ share two vertices (and hence the edge between those vertices), and  $S$ 
and $T$ do not have the same angles at the vertices of their shared side. 
A {\em \bf connected component}, or just a {\em \bf component},
 of a tiling is a maximal connected set in the graph $\HH$.
\end{definition}

\begin{figure}[ht]
\begin{center}
\includegraphics[width=0.8\textwidth]{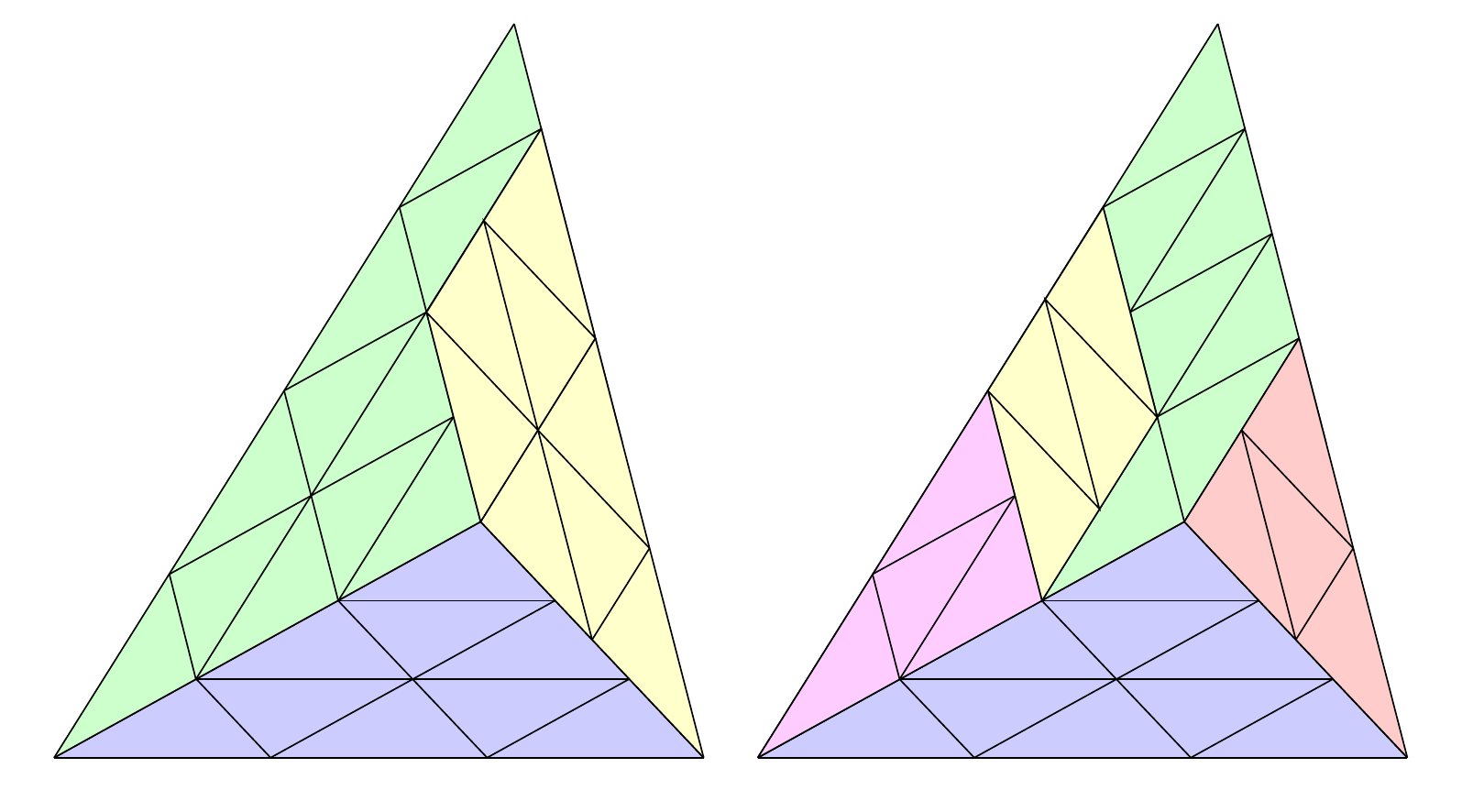}  
\end{center}
\caption{Components of two 28-tilings}
\label{figure:two28}
\end{figure}

\begin{examples}  Any quadratic tiling has only one connected component.   Any biquadratic 
tiling has two connected components, one for each of the two quadratic tilings that it contains.  
Fig. \ref{figure:two28} illustrates the connected components of two 
different 28-tilings, both different from the one in Figure~\ref{figure:28}.
\end{examples}

\begin{figure}[ht]
\begin{center}
\includegraphics[width=0.4\textwidth]{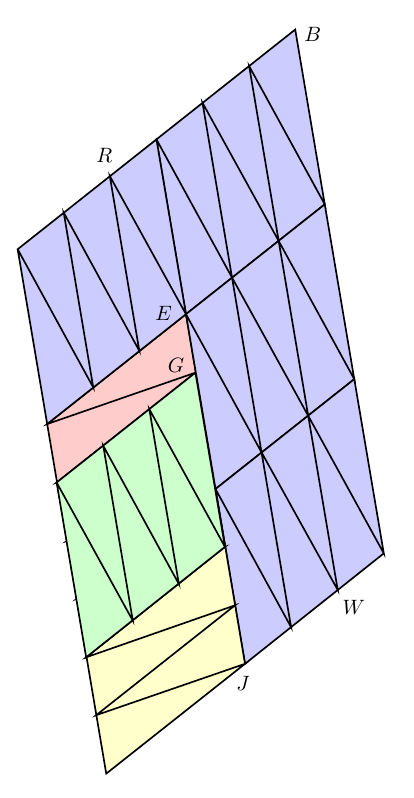} 
\end{center}
\caption{There are four connected components in this partial tiling, not three.}
\label{figure:105}
\end{figure}

\begin{remark}  We could have also allowed two tiles to be connected if 
they lie on opposite sides of an interior segment $L$ and have edges of the 
same length lying on $L$, i.e. both have their $a$ edges on $L$, or both have their $b$ edges on $L$,
or both have their $c$ edges on $L$.  This less stringent definition results in a graph 
with more connected components.  For example, consider the partial tiling illustrated
in Figure~\ref{figure:105}.  There are four components of the graph $\HH$; but with the 
less restrictive definition of connection, there would be only three.  The segment $FJ$
is a component boundary in $\HH$,  but with the less restrictive definition of connection,
the tiles with $c$ edges on  $FJ$ would be connected.  It will turn out to be 
important that they {\em not} be connected. 
\end{remark}

Another variation on the definition of $\HH$ is the graph $\G$, defined by removing the 
requirement that connected tiles have different angles at their common vertices. 
In this graph,  lines that have all the same length edges on both sides, with matching 
vertices but the tiles have equal angles at each vertex,  are not component boundaries 
as they are in $\HH$.  For example, in the 28-tiling shown in  Figure~\ref{figure:28},
  such a line occurs as the angle
bisector of angle $A$, and the two components of $\HH$ above and below that line will 
join into one $\G$-component.   We will make no use of these variations on the definition
of $\HH$, and the word ``component'' will be used according to the definition.

\begin{lemma} \label{lemma:connectedcomponents} 
Assume $a \neq b$, i.e. the tile is not isosceles.
In any connected set $\mathcal C$ in a tiling,
any two edges of the same length are parallel or on the same line.  Two tile edges of 
different lengths lying on the same line do not belong to the same connected component.
\end{lemma}

\begin{remark}  Of course $a$  and $b$ are always unequal if $3\alpha + 2\beta = \pi$
and $\alpha$ is not a rational multiple of $\pi$,
as we generally assume in this paper; but the lemma needs only $a \neq b$.  We
formalized only the $3\alpha + 2\beta = \pi$ case.
\end{remark}

\begin{proof}[Proof of Lemma~\ref{lemma:connectedcomponents}]
 The second sentence of the lemma follows from the first.  We prove 
the first by induction on the number $n$ of tiles in $\mathcal C$.  In the base case, 
when $n=1$, if we have two edges of the same length they are the same edge, so they lie
on the same line.  That completes the base case.  Now for the induction step.  Suppose $\mathcal C$
has $n+1$ tiles.  Since $\mathcal C$ is connected and finite, we may choose $T$ so that
$\mathcal C \setminus \{T\}$ is still connected --- for instance, a tile
at maximum distance from a fixed tile of $\mathcal C$ in the graph $\HH$.
Let $S$ be a tile in $\mathcal C$ that is 
directly connected to $T$; then by the definition of the graph $\HH$,  the two tiles $S$ and $T$
have different angles at the endpoints of the shared edge $UV$. 

Since $S$ and $T$ are congruent and $UV$ is the same side of each, the
two angles at the ends of $UV$ are the same pair for both tiles.  Since
$S$ and $T$ are connected, they have different angles at $U$; hence the
angle of $T$ at $U$ equals the angle of $S$ at $V$, and the angle of $T$
at $V$ equals the angle of $S$ at $U$.  These are alternate interior
angles for the transversal $UV$, so the corresponding sides of $S$ and
$T$ are parallel, and $S \cup T$ is a parallelogram.

  By hypothesis, $a \neq b$;
by the law of sines $\sin \alpha \neq \sin \beta$; hence $\alpha \neq \beta$.
 Since $\alpha \neq \beta$, the sides of 
$S$ and $T$ of the same lengths are parallel (or on the same line, in the case of the shared side).
We can apply the induction hypothesis to the set obtained by deleting $T$ from $\mathcal C$; then 
all the other tiles in $\mathcal C$ have their corresponding edges parallel to or collinear with 
the edges of $S$, and hence to those of $T$.  \end{proof}

\begin{definition} A {\em \bf lattice tiling} (of some region, not necessarily a triangle)
is a tiling whose vertices lie on the lattice points of 
some lattice.
\end{definition}

For example, a quadratic tiling, or any subtiling of a quadratic tiling, is a lattice tiling.

\begin{lemma} \label{lemma:Hcomponents}
Suppose  there is an $N$-tiling of $ABC$ by the tile with angles $\alpha$, $\beta$, and $\gamma$, 
where  $3\alpha+2\beta = \pi$ and $\alpha \neq \beta$.
Let $\mathcal E$ be a connected component of the tiling, or more generally, the union of any connected set of tiles.
Then the tiling (restricted to $\mathcal E$) is a lattice tiling of $\mathcal E$.  
\end{lemma}

\begin{remark} The 28-tiling shows that the lemma fails if the definition of $\HH$ is changed
to permit tiles with the same angles at their shared vertices to be connected.  The lemma actually 
needs only $\alpha \neq \beta$, but we formalized only the $3\alpha + 2\beta = \pi$ case.  
\end{remark}

\begin{proof} 
 We proceed by induction on the number of tiles in the connected set $\mathcal E$.
When  this number is $1$, the result is trivially true.    Suppose $\mathcal E$ is a connected set of 
$n+1 \ge 2$ tiles.  Since $\mathcal E$ is connected and finite, we may choose $T$ so that
$\mathcal E \setminus \{T\}$ is still connected --- for instance, a tile
at maximum distance from a fixed tile of $\mathcal E$ in the graph $\HH$.
 Remove $T$ from $\mathcal E$, and call what is left 
$\mathcal E^\prime$.  By the induction hypothesis, the conclusion of the lemma holds for $\mathcal E^\prime$. 
  Let $S$ be a tile in $\mathcal E^\prime$ that shares an edge and its 
vertices with $T$.  Then by the induction hypothesis, the vertices of $S$ lie in a lattice that contains 
all the vertices of $\mathcal E^\prime$.  Since $S$ and $T$ are connected,
they have different angles at the endpoints of the shared edge.  Therefore $S$ and $T$
form a parallelogram, as in the proof of Lemma~\ref{lemma:connectedcomponents}.  
Therefore the lattice also contains
the third vertex of $T$.  
\end{proof}
\medskip

The ``type'' of a component is specified by giving the directions 
of its $a$, $b$, and $c$ sides.  If a particular triangle $ABC$
is understood in the context, and a tile with angles $(\alpha,\beta,\gamma)$
is also understood, with $\alpha$ and $\beta$ not rational multiples of
$\pi$, and with a tile at $B$ having a side on $AB$ and its $\beta$
angle at $B$, then the following terminology makes sense.  A
``direction'' is an equivalence class of parallel lines; when we say a
line segment has direction $AB$, we mean that it belongs to the class of
lines parallel to $AB$, and similarly for the other three directions
defined below.

\begin{definition}
\label{definition:directions}
Consider a tile placed at vertex $B$ with its $\beta$ angle at $B$ and a
side along $AB$.  There are two such placements.  If its $a$ side lies
along $AB$, the side opposite $B$ defines {\em  Direction $A$}; if its
$c$ side lies along $AB$, the side opposite $B$ defines
{\em  Direction $C$}.  In a given tiling we write Tile~1 for the tile
at $B$ having a side along $AB$, and $E$ for the vertex of Tile~1 not on
$AB$; {\em Direction $BE$} is the direction of the side $BE$ and
{\em Direction $AB$} is the direction of  $AB$.
\end{definition}

So if $B$ is to the north,  Direction $A$ is northwest and Direction~$C$ is northeast.

\begin{definition} \label{definition:types}
We give names to the following types of components:
\smallskip

Type I:  $c$ edges have direction $BE$,  $a$ edges have direction $AB$, and $b$ edges have direction $A$.
\smallskip

Type II: $c$ edges have direction $AB$, $a$ edges have direction $BE$, and $b$ edges have direction $C$.
\smallskip

Type III:  $c$ edges have direction $BE$,  $a$ and $b$ edges have directions to make kites with Type~I.
\smallskip

Type IV: $c$ edges have direction $AB$,   $a$ and $b$ edges have directions to make kites with Type~II.
\smallskip

\end{definition}

These are the four types of tiles that have their $a$ and $c$ edges in directions $BE$ or $AB$.
Types~I and IV can make parallelograms, as can Types~II and III.  Types~I and III can make kites,
as can Types~II and IV.  Type~III and Type~IV tiles are not of Types~I or~II.   

\begin{remark}  $ABC$ may have just one tile at $B$, or it may have more than one.  Tile~1 is 
the one with a side on $AB$.  Remember that in this definition, we have not assumed $3\alpha + 2 \beta = \pi$,
so we are free to apply this definition when the vertex angle of Tile~1 is not $\beta$ but $\alpha$.
\end{remark}

These types of components are illustrated in the 28-tiling and the other triquadratic tilings.
Further types  also arise.  For example, along the base $AC$ of the 
28-tiling are tiles of Types~III and~IV.  
There are in principle infinitely many types, but they do not actually arise in the known tilings, so there is 
not much use in classifying them.  

Technically, components are of Type I or Type II, not tiles.  But we say that 
a tile (considered as part of a tiling)  is ``of Type I''
if it belongs to a component of Type I, and similarly for other types. 

Please have another look at Figure~\ref{figure:105}.  Note that there are two different Type I components 
in that figure, that are not connected even though some of their tiles share a common boundary with tile
edges of the same length on that boundary.  They 
still do not connect, since the tiles on that boundary do not share their vertices.  This motivates
the following definition:

\begin{definition}\label{definition:sync}  Two components $\mathcal C$ and $\mathcal D$ of the same 
type are said to be {\bf \em out of sync} if they share a common boundary segment with edges of the 
same length on that boundary segment that do not share vertices.
\end{definition}

\section{The case when $ABC$ has angles $(2\alpha,\beta,\alpha+\beta)$}
In this section, we derive the ``tiling equation'' for $ABC$ of 
the shape mentioned, and prove that it is a necessary 
condition for the existence of an $N$-tiling of some triangle $ABC$
of the specified shape.   

Our convention throughout this section is that the angle 
at $A$ is $2\alpha$, the angle at $B$ is $\beta$, and the angle 
at $C$ is $\alpha+\beta$.  The lengths of the sides opposite
vertices $(A,B,C)$ are $X$, $Y$, and $Z$ respectively.  In pictures,
we draw the vertex $B$ at the top.   

The $\d$ matrix (which occurs in \cite{snover1991} and the 
first edition of \cite{soifer}) 
is defined by 
$$ \d \vector a b c = \vector X Y Z$$
where $X$ is the length of $BC$ (opposite the $2\alpha$ angle), $Y$ is the length of $AC$ (opposite
the angle at $B$) and $Z$ is the length of $AB$ (opposite the angle at $C$).  We always draw
our pictures with $AC$ horizontal and $A$ at the lower left, or southwest, and $B$ at the 
top (or north) of the picture.  The $\d$ matrix 
tells how the sides of $ABC$ are composed of edges of tiles (although it contains no information 
about the order of those edges).
We use the following letters for elements of the  $\d$ matrix:
$$ \d = \matrix p d e  g m f  h {\ell} r $$

   In case the tile is similar to $ABC$,  then there is an 
eigenvalue equation involving the $\d$ matrix:
$$ \d \vector a b c = \lambda \vector a b c$$
since $(X,Y,Z) = \lambda (a,b,c)$ in that case.  In the present case we can make the equation look 
something like that:

\begin{lemma}\label{lemma:dmatrix}
Let the angles of $ABC$ be $(2\alpha, \beta, \alpha + \beta)$. 
Let $X$, $Y$, and $Z$ be the lengths of $BC$, $AC$, and $AB$ respectively.     
Let $s = 2\sin(\alpha/2)$.  Then for a suitable ``scaling factor'' $\mu$, 
$$ \vector X Y Z = \d \vector a b c = \mu \vector {(2-s^2)a} b c$$
\end{lemma}

\begin{remark} Because of the factor $2-s^2$, this is not quite an eigenvalue problem.  When triangle $ABC$ is similar to the 
tile, the factor $2-s^2$ does not appear, but still we can make this problem look {\em almost} like an 
eigenvalue problem.
\end{remark}

\begin{proof}
  According to the law of 
sines, for a certain positive $\lambda$ we have
\begin{eqnarray*} 
\vector X Y Z &=& \lambda \vector {\sin A} {\sin B} {\sin C} 
= \lambda\vector { \sin 2 \alpha} {\sin \beta} { \sin(\alpha + \beta)} 
\ = \  \lambda\vector { 2 \cos \alpha \sin \alpha} {\sin \beta} { \sin \gamma}
\end{eqnarray*}
Also by the law of sines, for some positive $\kappa$ we have 
$(a,b,c) = \kappa (\sin \alpha, \sin \beta, \sin \gamma)$. 
Define $\mu = \lambda/\kappa$.  Then we have
\begin{eqnarray*}
\vector X Y Z &=&  \mu\vector { 2 (\cos \alpha) a} b c 
\ = \   \mu \vector { 2 (1-2\sin^2 (\frac \alpha 2)) a} b c 
\ = \  \mu \vector { (2-s^2)a} b c 
\end{eqnarray*}
\end{proof}

\subsection{The tiling equation for the shape $(2\alpha,\beta,\alpha+\beta)$}
The ``tiling equation''  for a particular shape of $ABC$ is a Diophantine
equation that is solvable whenever there exists an $N$-tiling of a triangle of that shape.  

\begin{lemma} \label{lemma:52}
 Let $3\alpha + 2\beta = \pi$ and assume triangle $ABC$ has angles 
 $(2\alpha,\beta,\alpha+\beta)$ and is $N$-tiled by
a tile with angles $\alpha$ and $\beta$.
   Let the tiles be 
assigned signs (or colors) in accordance with Lemma~\ref{lemma:signedtiles}.  Let $s = 2 \sin(\alpha/2)$ and let 
$M$ be the number of black (positive) tiles minus the number of white (negative) tiles.   Then $M > 0$, 
$\mu = M/s$,  and 
$$  s^2 = \frac {2M^2} {M^2 + N}$$
\end{lemma}

\begin{proof}  By Theorem~\ref{theorem:coloring}, we have 
$M(a+b+c) = X - Y + Z$.  By 
Lemma~\ref{lemma:dmatrix},  we have 
\begin{eqnarray}
\vector X Y Z &=& \mu \vector {a(2-s^2)} b c  \label{eq:4322}
\end{eqnarray}
Hence $M(a+b+c) = \mu(a(2-s^2)-b+c)$.  Dividing by $c$  we have 
\begin{eqnarray*}
M(s+(1-s^2) + 1) &=& \mu (s(2-s^2) - (1-s^2) + 1) \\
M(s+2-s^2) &=& \mu(2s-s^3+s^2) \\
   &=& \mu s(2-s^2+s)
\end{eqnarray*}
Dividing both sides by $2-s^2+s$ (which cannot be zero since $0 < s < 1$) we have 
\begin{eqnarray}
M &=& \mu s   \label{eq:mu=M/s}
\end{eqnarray}
Since $\mu$ and $s$ are both positive, this implies $M > 0$. 
We claim $N = \mu^2 (2-s^2)$.  To prove this we use the 
area equation  
$$Nbc \sin \alpha = YZ \sin 2\alpha,$$ 
whose left side is twice the area of $N$ tiles, and whose right side
is twice the area of triangle $ABC$.
According to Lemma~\ref{lemma:singamma2}, 
$$\sin 2\alpha = (2-s^2)\sin \alpha.$$
Putting that into the area equation just above, we have
$$Nbc \sin \alpha = YZ (2-s^2) \sin \alpha.$$
Cancelling $\sin \alpha$ we have 
$$Nbc = YZ(2-s^2).$$
Replacing 
$Y$ and $Z$ by their values from (\ref{eq:4322}), we have 
\begin{eqnarray*}
Nbc &=& \mu^2 bc (2-s^2)
\end{eqnarray*}
Cancelling $bc$ we have the claimed equation $N = \mu^2(2-s^2)$.

Then squaring (\ref{eq:mu=M/s}) we have 
\begin{eqnarray*}
M^2 &=& \mu^2 s^2  \\
   &=& s^2 \frac N {2-s^2} \\
   &=& N \frac {s^2}{2-s^2}  \\
\frac N {M^2} &=& \frac {2-s^2} {s^2} \\
   &=& \frac 2 {s^2} - 1
\end{eqnarray*}
Solving for $s^2$ we have 
\begin{equation}
 s^2 = \frac {2M^2} {M^2 + N}  \label{eq:sMN}
 \end{equation}
\end{proof}

\begin{theorem}[Necessary condition for a tiling] \label{theorem:tilingequation}
 Let $3\alpha + 2\beta = \pi$ and assume triangle $ABC$ has angles $(2\alpha, \beta, \alpha+\beta)$ and is $N$-tiled by
a  tile with angles $\alpha$ and $\beta$. Then 
\smallskip

(i) $N$ satisfies the equation 
$$M^2 + N = 2K^2$$
for some positive integers $M$ and $K$, with $M^2 < N$, 
where 
$$s = 2 \sin(\alpha/2) = M/K,$$ 
and
\smallskip

(ii) The squarefree part of $N$ is a product of distinct primes of the form $8n \pm 1$, possibly times $2$.   
(The case when the product is empty and $N$ is a square is not ruled out.)
\end{theorem}

\begin{proof}  Recall that $s$ is rational if and only if the tile is rational, since $a/c = s$ and $b/c = 1-s^2$.  
If there is an $N$-tiling as in the theorem, then by Lemma~\ref{lemma:52}, 
there exists an integer $M$ such that $M^2+N = 2M^2/s^2$.
Setting $K=M/s$ we have $M^2 + N = 2K^2$.  Then $K$ is the square root of a half-integer (not yet known to be an integer).

Now to prove  $M^2<N$.  Recall from (\ref{eq:sMN}) that
$$ s^2 = \frac {2M^2} {M^2 + N} $$
But $s = 2 \sin(\alpha/2)$.  Both $\alpha $ and $\beta$
are less than $\gamma$, so both $\alpha$ and $\beta$ are less than $\pi/3$.  Hence $\alpha/2 < \pi/6$,
so $\sin(\alpha/2) < 1/2$, and $s < 1$.  Hence $2M^2/(M^2+N) < 1$.  Hence $2M^2 < M^2 + N$; 
hence $M^2 < N$ as claimed in the statement of the lemma.   Since $M^2 + N = 2K^2$ we have
$$s^2 = \frac {2M^2}{M^2+N} = \frac {M^2}{K^2}.$$  Taking the square root of both sides we have $s = M/K$.
  That completes the proof of the first part of the lemma, but we
still have to prove that $K$ is an integer.

Let $K = M/s$.  Then $M^2 + N = 2K^2$;  then we have $s = M/K$,
and $a/c = s$, $b/c = 1-s^2$,  so  $a$,
$b$, and $c$ are proportional to $(M,(N/K)-K,K)$.  
We have 
\begin{eqnarray}
 2-s^2 &=& (2K^2-M^2)/K^2 = N/K^2  \nonumber \\
 &=& s^2 \frac N {M^2} \mbox{\qquad since $s = M/K$}. \label{eq:4388}
\end{eqnarray}
By Theorem~\ref{theorem:rationaltile}, $s$ is rational.
Then $K = M/s$ is rational.  
Then $2(M^2+N) = 4K^2$ is a rational square,  and hence it is an
integer square.  Hence $2K$ is an integer.  If $K$ is not an integer,
then $2K^2$ is also not an integer.  But $2K^2 = M^2 + N$ is in 
fact an integer.  Hence $K$ is an integer.

Since $M^2 + N = 2K^2$ is solvable for  $K$ with $K > 0$,
 2 is a square mod $p$ for each odd prime $p$ dividing $N$ but not $M$.   But 2 is square mod $p$ just in case 
$p$ is congruent to $\pm 1$ mod 8. Hence, odd primes dividing $N$ but not $M$ are congruent to $\pm 1$ mod 8.
If $p$ is a prime that divides both $N$ and $M$, but divides $N$ to a power $p^j$ that does not divide $M^2$, 
then similarly 2 is a square mod $p^j$, and hence mod $p$,  so again $p$ is  congruent to $\pm 1$ mod 8.
Thus any odd primes $p$ dividing $N$ that are not congruent to $\pm 1$ mod 8, must occur to even powers
and divide $M^2$ to at least the power they divide $N$.  Therefore the squarefree part of 
$N$ is a square times some primes
congruent to $\pm 1$ mod 8, possibly times a power of 2.  
\end{proof}

\medskip
\noindent{\em Question}\/:  If I divide $M$ and $K$ both by 2 (assuming they are both even), then $s=M/K$ does not 
change, and since $s = 2\sin(\alpha/2)$ determines the shape of the tile, the new tile is similar to the old one,
just half the linear dimensions, so one-quarter the area.    
But $N$ is divided by 4 according 
to the tiling equation, whereas,  it should take four times as many tiles to tile $ABC$ with these smaller tiles.
What is going on here?   

{\em Answer}\/:  When you divide $M$ and $K$ by 2, the size of the tiled
triangle $ABC$ also changes.  With $\lambda$ held fixed, the formula
$(a,b,c) = \lambda(M,\,N/K-K,\,K)$ shows that halving $M$ and $K$
halves each of $a$, $b$, and $c$; so the tiles have half the linear
dimensions and one-quarter the area, as expected.  But the length $Y$ of
side $AC$ is $\mu b = Kb$ by Lemma~\ref{lemma:dmatrix} and
(\ref{eq:mu=K}), and halving both $K$ and $b$ divides $Y$ by $4$.  The
same holds for the other two sides, so triangle $ABC$ has one-sixteenth
its former area.  One-sixteenth the area, tiled by tiles of one-quarter
the area, requires one-fourth as many tiles --- which is just what the
tiling equation gives, since halving $M$ and $K$ replaces $N = 2K^2-M^2$
by $N/4$.
 
 \begin{lemma} \label{lemma:XYZ} Assume triangle $ABC$ has angles 
 $(2\alpha,\beta,\alpha+\beta)$,
   and there is an $N$-tiling of $ABC$ 
   by the tile with angles $(\alpha,\beta,\gamma)$,
 and $N +M^2 = 2K^2$ where $K = M/s$.
Then for some $\lambda > 0$ we have 
 $(a,b,c) = \lambda (M,N/K-K,K)$. Then 
 $$(X,Y,Z) = \lambda (MN/K, N-K^2, K^2)$$
 and the $\d$-matrix equation takes the form
 \begin{eqnarray*}
 MN/K &=& pM + d(N/K-K) + eK \\
 N-K^2 &=& gM + n(N/K-K) + fK \\
 K^2 &=& hM + \ell(N/K-K) + rK
 \end{eqnarray*}
 \end{lemma}
 
 \begin{proof}  
Since $K = M/s$, we have $s = M/K$.  Since
 $a/c = s$ and $b/c = 1-s^2$,   $a$,
$b$, and $c$ are proportional to $(M,N/K-K,K)$.
Thus, as claimed, we can choose 
\begin{eqnarray}
(a,b,c) &=&  \lambda (M,N/K-K,K). \label{eq:4467}
\end{eqnarray}
According to (\ref{eq:mu=M/s}), we have $\mu = M/s$.
Since $s = M/K$ that implies 
\begin{eqnarray}
\mu &=& K \label{eq:mu=K}
\end{eqnarray} 
Since $s = M/K$ we have
\begin{eqnarray}
 2-s^2 &=& \frac{2K^2-M^2}{K^2} \ = \ \frac N {K^2}\label{eq:4474}
 \end{eqnarray}
By Lemma~\ref{lemma:dmatrix} we have
\begin{eqnarray*}
 \vector X Y Z &=& \d \vector a b c \ = \ \mu \vector {(2-s^2)a} b c \\
 &=& K \vector {(2-s^2)a} b c  \mbox{\qquad since $\mu = K$ }\\
 &=& K \vector { (N/K^2) a} b c \mbox{\qquad by (\ref{eq:4474})}\\
 &=& \lambda K \vector { (N/K^2) M} {N/K-K} K \mbox{\qquad by (\ref{eq:4467})} \\
 &=& \lambda \vector{ MN/K} {N-K^2} {K^2}.
 \end{eqnarray*}
This gives the left-hand side of the equation in the theorem.
The right-hand side is obtained by substituting the values of
$(a,b,c)$ from (\ref{eq:4467}) into the $\d$-matrix.
\end{proof}

\begin{lemma} \label{lemma:n<K}
Suppose given an $N$-tiling of a triangle $ABC$ with angles 
 $(2\alpha,\beta,\alpha+\beta)$, and integers $M$ and $K$
satisfying the tiling equation with $s = M/K$.  Then $n < K$,
where $n$ is the number of $b$ edges on the side $AC$ of $ABC$ opposite
the $\beta$ angle.
\end{lemma}

\begin{proof} 
From the $\d$-matrix equation we have
$$ Y = ga + nb + fc, $$
and $Y = \mu b$ by Lemma~\ref{lemma:dmatrix}.
Therefore 
$$ \mu b = ga + nb + fc.$$
The first and third terms on the 
 right are nonnegative; and the third one is not just
 nonnegative, it is positive, since $f > 0$ by Lemma \ref{lemma:twocsapplies}.
Therefore $n < \mu$.   Since $\mu = K$ by Lemma~\ref{lemma:52},
we have $n < K$. 
 \end{proof}
 
  \begin{lemma}\label{lemma:X=M(b+c)}
  Assume triangle $ABC$ has angles 
 $(2\alpha,\beta,\alpha+\beta)$,  and there is an $N$-tiling,
 and $N + M^2 = 2K^2$.  
 Then the length of $BC$ is given by $X = M(b+c)$.
 \end{lemma}
 
 \begin{proof}  We give two proofs. 
 The length of side $BC$ (opposite vertex $A$) is $X=\mu a (2-s^2)$ according to the $\d$ matrix equation, 
 and according to Lemma \ref{lemma:XYZ}, $X = \lambda MN/K$, where  
 \begin{eqnarray*}
 a &=&  \lambda M,  \\
 c &=& \lambda K , \\ 
 b &=&c(1-s^2) =  \lambda K(1-(M/K)^2)  \\
 &=&   \lambda(K^2-M^2)/K  \\
 &=&  \lambda (N-K^2)/K \\
 &=& \lambda  (N/K-K)
 \end{eqnarray*}
 Hence $b+c =  \lambda(N/K-K)+  \lambda K = \lambda N/K$.  Hence 
 $$BC = X = \lambda MN/K = M(\lambda N/K) =  M(b+c).$$   That completes the first proof.
 
 The same result can also be obtained directly from the relation 
 $X-Y+Z = M(a+b+c)$, since $Y = \mu b = Kb$ and $Z = \mu c = Kc$ (by Lemma~\ref{lemma:dmatrix}).  So
 $X = \vert BC \vert = M(a+b+c)-Kc+Kb$. Now put $a=\lambda M$,
 $c=\lambda K$, and $b = \lambda(K-M^2/K)$, 
  we have  
 \begin{eqnarray*}
 X &=& M(b+c) + Ma-Kc+Kb \\
 &=& M(b+c) + \lambda M^2-\lambda K^2 + \lambda (K^2-M^2) \\
 &=& M(b+c)
 \end{eqnarray*}
That completes the second proof of the lemma.
\end{proof} 

\subsection{For each $N$, there are only a few possibilities for the tile and for $ABC$}

In this section, we assume that triangle $ABC$ has angles 
 $(2\alpha,\beta,\alpha+\beta)$.  Then each tiling corresponds
 to a solution of the tiling equation
$M^2+N = 2K^2$ with $M^2 < N$, 
according to Theorem~\ref{theorem:tilingequation}. 
Hence, for a given $N$, there are finitely many solutions.  
In practice, there are comparatively few such $M$
and $K$, often just one for a given $N$.   
  Since $s = 2 \sin \tfrac \alpha 2 = M/K$,
each solution of the tiling equation 
determines one possible value of $\alpha$ (which may or may not 
actually correspond to a tiling).   Then 
the numbers $\mu$ is also determined as $\mu = M/s$, by (\ref{eq:mu=M/s}).

We note that the equation does sometimes have more than one solution.  The least such $N$ 
is $N=119$, where we have $(M,K) = (3,8)$ or $(9,10)$.  These yield tiles 
$(24,55,64)$ and $(90,19,100)$.  In that case $K$ does not divide $N$,
and as we will see below, no tiling exists in that case.  The least $N$
such that the equation has more than one solution in which $K$ divides $N$ is 
$ N = 87808$.  In that case there are two solutions with $K$ dividing $N$, namely
$ M= 112, K= 224 $ and 
 $ M= 208, K= 256$.  These correspond to the tiles 
 $(2,3,4)$ and $(208,87,256)$
 respectively.  The corresponding triangles $ABC$ are quite large; it is not practical
 to make a picture that fits on one page.  There are many more $N$ with 
 two solutions, but none with three solutions less than ten million.

\subsection{No $N$-tilings for $N < 28$}
\begin{lemma} \label{lemma:jmk}
Assume triangle $ABC$ has angles $(2\alpha,\beta,\alpha+\beta)$ and there
is an $N$-tiling, and let $M$ and $K$ be as in
Theorem~\ref{theorem:tilingequation}, so that $M^2+N = 2K^2$ and
$s = M/K$.  Then the tiling may be scaled so that $(a,b,c)$ are integers
with no common factor; with $m$ and $k$ as in
Lemma~\ref{lemma:gsquare}, so that
$$ (a,b,c) \ = \ (mk,\ k^2-m^2,\ k^2) $$
and $a/c = m/k$ in lowest terms, put $j = \gcd(M,K)$.  Then
$$ M = jm, \qquad K = jk, \qquad N = j^2(2k^2-m^2), $$
and
$$ X = M(b+c), \qquad Y = Kb, \qquad Z = Kc . $$
Moreover $K \mid M^2$ if and only if $k \mid j$.
\end{lemma}

\begin{proof}
By Theorem~\ref{theorem:rationaltile} the ratios of the sides of the
tile are rational, so the tiling may be scaled so that $(a,b,c)$ are
integers with no common factor; then Lemma~\ref{lemma:gsquare} applies
and gives $(a,b,c) = (mk,\,k^2-m^2,\,k^2)$ with $a/c = m/k$ in lowest
terms.

By Lemma~\ref{lemma:singamma}, $s = a/c$, and by
Theorem~\ref{theorem:tilingequation}, $s = M/K$.  Hence $M/K = m/k$, and
since $m/k$ is in lowest terms we have $M = jm$ and $K = jk$ with
$j = \gcd(M,K)$.  Then
$$ N \ = \ 2K^2 - M^2 \ = \ j^2(2k^2-m^2). $$
The formula $X = M(b+c)$ is Lemma~\ref{lemma:X=M(b+c)}, and
$Y = \mu b = Kb$ and $Z = \mu c = Kc$ by Lemma~\ref{lemma:dmatrix}
together with $\mu = K$.

Finally, $K \mid M^2$ says $jk \mid j^2m^2$, that is, $k \mid jm^2$;
since $\gcd(m,k) = 1$, that holds if and only if $k \mid j$.
\end{proof}

\begin{lemma} \label{lemma:congruences}
Under the hypotheses and notation of Lemma~\ref{lemma:jmk}, the entries
of the $\d$ matrix satisfy
$$ k \mid n, \qquad k \mid \ell, \qquad d \equiv jm \pmod k, $$
where $n$ is the number of $b$ edges on $AC$, $\ell$ the number of $b$
edges on $AB$, and $d$ the number of $b$ edges on $BC$.
\end{lemma}

\begin{proof}
Modulo $k$ we have $a = km \equiv 0$, $b = k^2-m^2 \equiv -m^2$, and
$c = k^2 \equiv 0$.  The second row of the $\d$ matrix says
$Y = ga + nb + fc$, so
$$ -n m^2 \ \equiv \ Y \ = \ Kb \ = \ jk(k^2-m^2) \ \equiv \ 0
   \pmod k, $$
and since $\gcd(m,k)=1$ we may cancel $m^2$ to get $k \mid n$.  The
third row gives, in the same way,
$-\ell m^2 \equiv Z = Kc = jk^3 \equiv 0$, so $k \mid \ell$.  The first
row gives
$$ -d m^2 \ \equiv \ X \ = \ M(b+c) \ = \ jm(2k^2-m^2) \ \equiv \
   -jm^3 \pmod k, $$
and canceling $m^2$ gives $d \equiv jm \pmod k$.
\end{proof}

\begin{corollary} \label{corollary:nokt}
Under the hypotheses of Lemma~\ref{lemma:jmk}, we have $n = kt$ for some
integer $t$ with $0 \le t < j$.  In particular, if $\gcd(M,K) = 1$ then
$n = 0$; that is, there are no $b$ edges at all on the side $AC$.
\end{corollary}

\begin{proof}
By Lemma~\ref{lemma:congruences}, $k \mid n$, say $n = kt$ with
$t \ge 0$.  By Lemma~\ref{lemma:n<K}, $n < K = jk$, so $kt < jk$ and
$t < j$.  If $j = \gcd(M,K) = 1$ then $t < 1$, so $t = 0$ and $n = 0$.
\end{proof}
\begin{lemma} \label{lemma:twocsapplies}
Assume triangle $ABC$ has angles $(2\alpha,\beta,\alpha+\beta)$ with
$3\alpha+2\beta = \pi$, that $\alpha$ is not a rational multiple of
$\pi$, and that there is an $N$-tiling.  Suppose that neither of $a$,
$b$ is an integer multiple of the other.  Then $f \ge 2$, where $f$ is
the number of $c$ edges on the side $AC$.
\end{lemma}

\begin{proof}
We verify the hypotheses of Lemma~\ref{lemma:twocs}.  By
Lemma~\ref{lemma:gammagreaterpiover2}, $\gamma > \pi/2$, which is (i).
For (ii): $\alpha$ is not a rational multiple of $\pi$ by hypothesis,
and then neither is $\beta = (\pi-3\alpha)/2$; and $\alpha/\beta$ is
irrational, for if $\alpha/\beta = u/v$ with $u,v$ integers then
$\pi = 3\alpha+2\beta = (3 + 2v/u)\alpha$ would make $\alpha$ a rational
multiple of $\pi$.  For (iii): the angles of $ABC$ are $2\alpha$,
$\beta$, and $\alpha+\beta$, and $\gamma = 2\alpha+\beta$ exceeds each
of them, since $\gamma - 2\alpha = \beta > 0$,
$\gamma - \beta = 2\alpha > 0$, and $\gamma - (\alpha+\beta) = \alpha >
0$.  Next, (iv) holds by hypothesis.  Finally, for (v): by
Lemma~\ref{lemma:jmk}, $|AC| = Y = Kb$.  If $AC$ were a single $c$ edge,
then $Kb = c$, that is, $K(1-s^2) = 1$, since $b/c = 1 - s^2$.  With
$s = M/K$ this gives $M^2 = K(K-1)$.  But $0 < M < K$, so $K \ge 2$ and
$(K-1)^2 < K(K-1) < K^2$; hence $K(K-1)$ is not a square, a contradiction.
Hence Lemma~\ref{lemma:twocs} applies and gives at least two $c$ edges on
$AC$, that is, $f \ge 2$.
\end{proof}

\begin{lemma} \label{lemma:No14}
There is no $14$-tiling of a triangle with angles
$(2\alpha,\beta,\alpha+\beta)$, where $3\alpha+2\beta = \pi$ and
$\alpha$ is not a rational multiple of $\pi$.
\end{lemma}

\begin{proof}
Suppose there were.  The only solution of $M^2 + 14 = 2K^2$ with
$M < K$ is $(M,K) = (2,3)$, so by Lemma~\ref{lemma:jmk} we have
$j = \gcd(2,3) = 1$, $m = 2$, $k = 3$, and the tile is
$$ (a,b,c) \ = \ (6,\,5,\,9), $$
while
$$ X = M(b+c) = 28, \qquad Y = Kb = 15, \qquad Z = Kc = 27 . $$
Since $j = 1$, Corollary~\ref{corollary:nokt} gives $n = 0$: there are
no $b$ edges on $AC$.  Hence the second row of the $\d$ matrix
expresses $15$ using only $6$ and $9$, and the only way to do that is
$15 = 6+9$; that is, $(g,n,f) = (1,0,1)$, so $f = 1$.  But neither of
$a = 6$, $b = 5$ is a multiple of the other, so
Lemma~\ref{lemma:twocsapplies} gives $f \ge 2$, a contradiction.
\end{proof}

\begin{lemma} \label{lemma:No17}
There is no $17$-tiling of a triangle with angles
$(2\alpha,\beta,\alpha+\beta)$, where $3\alpha+2\beta = \pi$ and
$\alpha$ is not a rational multiple of $\pi$.
\end{lemma}

\begin{proof}
Suppose there were.  The only solution of $M^2 + 17 = 2K^2$ with
$M < K$ is $(M,K) = (1,3)$, so by Lemma~\ref{lemma:jmk} we have
$j = 1$, $m = 1$, $k = 3$, and the tile is
$$ (a,b,c) \ = \ (3,\,8,\,9), $$
while
$$ X = M(b+c) = 17, \qquad Y = Kb = 24, \qquad Z = Kc = 27 . $$
Since $j = 1$, Corollary~\ref{corollary:nokt} gives $n = 0$.  Hence the
second row of the $\d$ matrix expresses $24$ using only $3$ and $9$;
the possibilities are
$$ (g,n,f) \ = \ (8,0,0), \quad (5,0,1), \quad (2,0,2). $$
Neither of $a = 3$, $b = 8$ is a multiple of the other, so
Lemma~\ref{lemma:twocsapplies} gives $f \ge 2$, which rules out the
first two.  Hence $(g,n,f) = (2,0,2)$: the side $AC$ consists of two
$a$ edges and two $c$ edges, four tiles in all.

We have $\vert AB \vert = Kc = 27$, so its edge decomposition by 
the $\d$-matrix is $(3,0,2)$ or $(0,0,3)$ or $(6,0,1)$.  We will
rule out $(6,0,1)$:  then there would be 7 tiles with an edge on $AB$,
and at least six ``in-betweens'' with a vertex on the boundary, making
13 tiles touching $AB$, plus four on $AC$, with no double-counting at 
vertex $A$ as the tiles touching the two tiles at $A$ are white and black,
and so not identical.  At $C$ there is an $\alpha+\beta$ angle, hence
at least one more tile, making 18.  That rules out $(6,0,1)$ for $AB$.

We have $\vert BC \vert = 17 = M(b+c)$, so its edge decomposition is $(0,1,1)$
or $(3,1,0)$.  We can rule out $(3,1,0)$ by Lemma~\ref{lemma:pigeonhole}.
Then there are no $a$ edges on $BC$.  The tile at $B$ has a $\beta$ angle there,
so it must have its $c$ edge on $BC$ and its $a$ edge on $AB$.  That rules 
out $(0,0,3)$ for $AB$.

Therefore the edge decomposition
of $AB$ is exactly $(3,0,2)$, as the other possibilities have been 
ruled out.  Therefore the tile at $B$ has its $a$ edge on $AB$.
 So now we have the  edge
decompositions of all three sides:

\begin{center}
\begin{tabular}{rl}
$AB$ & $(3,0,2)$ \\
$BC$ & $(0,1,1)$ \\
$AC$ & $(2,0,2)$
\end{tabular}
\end{center}

That accounts for a total of 10 tiles with an edge on the boundary, since the 
tile at $B$ has edges on both adjacent sides.  Next we count
the ``in-between'' tiles with a vertex on the boundary.  There are at least four on $AB$,
one on $BC$, and three on $AC$, but there is double-counting at $B$, where the tile 
adjacent to the tile with a vertex at $B$ has an edge on two sides of $ABC$.  There is 
no double-counting at the other two corners, as the first in-between at either end of $AC$
is black, while the first in-betweens on $AB$ and $BC$ are white.  So the number
of in-betweens comes to $4+1+3-1 = 7$.  Together with 10 tiles on the boundary, that 
accounts for all 17 tiles.  It follows that every vertex on the boundary except $A$, $B$, and $C$
must be of type $(1,1,1)$, as $(3,2,0)$ would require more tiles. 
Therefore, on $AC$ every vertex has exactly one $\beta$ angle, since there are no $b$ edges there.
Therefore, the $\beta$ angles on $AC$ are all to the east.  Similarly, since there are no $b$
edges on $AB$, all the $\beta$ angles on $AB$ are to the north.
 
Since there are no $b$~edges on either  $AB$ or $AC$, the two tiles with an $\alpha$ angle at $A$
must share their $b$ edges.  Then those two tiles use $2\gamma$ at the interior vertex $P$ of 
the two corner tiles, so the adjacent in-betweens must share their $a$ edges with the corner tiles.
Therefore they have their $\beta$ or $\gamma$ angles at $P$.  Call them Tile~1 and Tile~2, with
Tile~2 having an edge on $AC$, 
and let Tile~3 be the in-between tile next to Tile~2, and Tile~4 the next tile on $AC$.
Tile~3 shares an $a$ edge with Tile~2, since Tile~1 blocks any longer edge against Tile~2.
Tile~2 has its $\beta$ angle at its eastern vertex $Q$, so Tile~3 cannot have its $\beta$
angle there, since the vertex has type $(1,1,1)$.  Therefore Tile~3 has its $\gamma$ angle
there (next to its $a$ edge).  Tile~4, next to Tile 3, has its $\alpha$ angle at $Q$,
so it must have its $c$ edge on $AC$, and Tile~3 has its $\beta$ angle at $P$.

Let Tile~5 be the tile next to Tile~1 on the north.  By the same reasoning as we applied
to $AC$,  which depended only on there being no $b$ edges on $AC$,  we find that Tile~5,
like Tile~2,  has its $\beta$ angle at $P$.  Therefore $P$ is not a center, but a 
vertex of type $(2,2,2)$, so there are two more tiles with a vertex at $P$, each with
an $\alpha$ angle there.  Call them Tile~6 and Tile~7.  Tile~6 has one side parallel to $AB$,
since the east edge of Tile~1 is a transversal making alternate interior angles both equal to $\gamma$.
Similarly Tile~7 has one side parallel to $AC$. The two sides of Tile~6 at $P$ lie along the rays from $P$ making angles
$\gamma+\beta$ and $\gamma+\beta+\alpha$ with the ray $PA$, and those of
Tile~7 along the rays making angles $\gamma+\beta+\alpha$ and
$\gamma+\beta+2\alpha$.  Since $AP$ bisects the angle $2\alpha$ at $A$,
and $\gamma+\beta = \pi-\alpha$, the first of these rays is parallel to
$AB$ and the last is parallel to $AC$; the middle one is the extension
of $AP$ beyond $P$.  Hence none of the three rays meets $AB$ or $AC$.
Taking $\vert AB \vert = 27$ and $\vert AC \vert = 24$, so that
$\vert AP \vert = b = 8$, one computes that the three rays meet $BC$ at
distances $18$, $16$, and $16$ from $P$ respectively.  But no edge of
the tile is longer than $c = 9$.  Therefore every vertex of Tile~6 and
of Tile~7 other than $P$ lies in the interior of $ABC$; in particular
neither tile has a vertex on the boundary, so neither is a boundary tile
nor an in-between.  That makes $17+2 = 19$ tiles, a contradiction.
\end{proof}

\begin{lemma} \label{lemma:No7}
There is no $7$-tiling of a triangle with angles
$(2\alpha,\beta,\alpha+\beta)$, where $3\alpha+2\beta = \pi$ and
$\alpha$ is not a rational multiple of $\pi$.
\end{lemma}

\begin{proof}
The only solution of $M^2+7 = 2K^2$ with $M<K$ is $(M,K) = (1,2)$, so by
Lemma~\ref{lemma:jmk} we have $j = 1$, $m = 1$, $k = 2$, and the tile is
$(a,b,c) = (2,3,4)$, while
$$ X = M(b+c) = 7, \qquad Y = Kb = 6, \qquad Z = Kc = 8. $$
Since $j = 1$, Corollary~\ref{corollary:nokt} gives $n = 0$, so the
second row of the $\d$ matrix expresses $6$ using only $2$ and $4$: the
possibilities are $(3,0,0)$ and $(1,0,1)$.  Neither of $a = 2$, $b = 3$
is a multiple of the other, so Lemma~\ref{lemma:twocsapplies} gives
$f \ge 2$, and both are ruled out.
\end{proof}

\begin{lemma} \label{lemma:No23}
There is no $23$-tiling of a triangle with angles
$(2\alpha,\beta,\alpha+\beta)$, where $3\alpha+2\beta = \pi$ and
$\alpha$ is not a rational multiple of $\pi$.
\end{lemma}

\begin{proof}
The only solution of $M^2 + 23 = 2K^2$ with $M < K$ is $(M,K) = (3,4)$,
so by Lemma~\ref{lemma:jmk} we have $j = \gcd(3,4) = 1$, $m = 3$,
$k = 4$, and the tile is
$$ (a,b,c) \ = \ (12,\,7,\,16), $$
while
$$ X = M(b+c) = 69, \qquad Y = Kb = 28, \qquad Z = Kc = 64 . $$
Since $j = 1$, Corollary~\ref{corollary:nokt} gives $n = 0$: there are
no $b$ edges on $AC$.  Hence the second row of the $\d$ matrix expresses
$28$ using only $12$ and $16$, and the only way to do that is
$28 = 12+16$; that is, $(g,n,f) = (1,0,1)$, so $f = 1$.  But neither of
$a = 12$, $b = 7$ is a multiple of the other, so
Lemma~\ref{lemma:twocsapplies} gives $f \ge 2$, a contradiction.
\end{proof}

\begin{lemma} \label{lemma:Nolessthan28}
Let $3\alpha + 2\beta = \pi$ and suppose triangle $ABC$ has angles $(A,B,C) = (2\alpha, \beta, \alpha + \beta)$, and $\alpha$ is not 
a rational multiple of $\pi$.  Then there is no $N$-tiling
of $ABC$ for $N < 28$.   
\end{lemma}

\begin{proof}  For $N < 28$, the equation $M^2 + N = 2K^2$ has a solution
with $0 < M < K$ (which is required, since $s = M/K$ and $0 < s < 1$)
only for $N = 7$, $14$, $17$, and $23$, as one verifies by trying each
$N$ in turn.  We have ruled each of these out in the previous four
lemmas.
\end{proof}

\subsection{The case $N=28$}
\begin{figure}[ht]
\begin{center}
\includegraphics[width=0.4\textwidth]{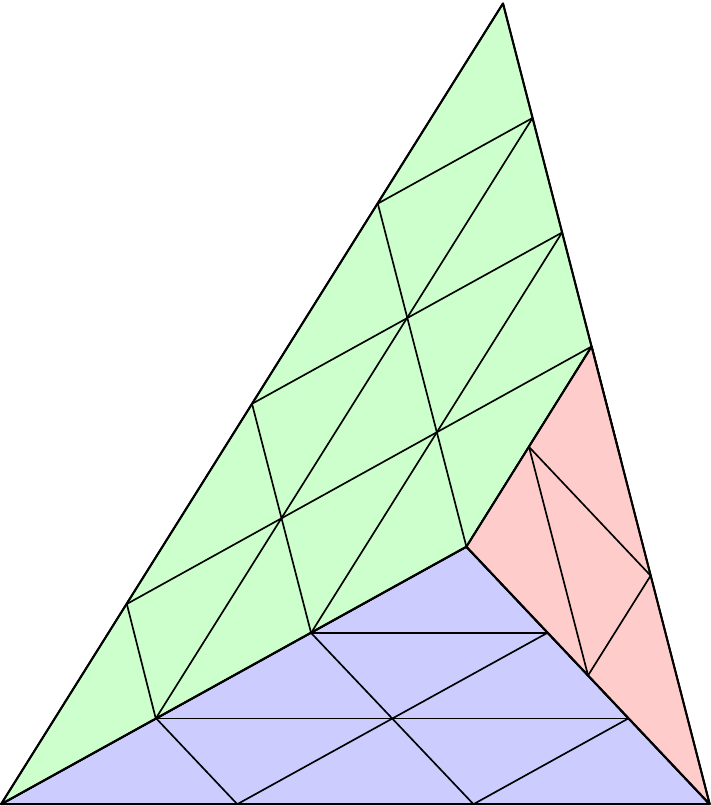} 
\end{center}
\caption{A 28-tiling}
\label{figure:28}
\end{figure}

The equation $M^2 + 28 = 2K^2$, with $M^2 < N$,
 has only one solution,
namely $M=2$ and $K=4$.  That yields the tile 
$(2,3,4)$ and the triangle $ABC$ is 
(14,12,16).  This is the same shape tile and triangle as in the case 
$N=7$, except now the triangle is twice as big in linear dimensions.
Its angles can be worked out from Lemma~\ref{lemma:singamma},
since $s= 1/2 = a/c$ and $b/c = 1-s^2$. One finds
\begin{eqnarray*}
\alpha &=& 28.955024^{\circ}\\
\beta &=&  46.567463^{\circ}
\end{eqnarray*}
For some time I thought that there was no 28-tiling.  I could not find one by hand searching, 
even with the aid of paper tiles.  I then wrote a computer program to search for one,  
expecting to show that none exist.  But on October 7, 2011, the program was producing 
``boundary tilings'', placing triangles in many possibly ways around the boundary of $ABC$, 
and I saw that one of these can be filled in.  See Figure~\ref{figure:28}.

\subsection{Triquadratic tilings}
In this section we generalize the 28-tiling, showing that it is just the simplest
member of a new family of tilings.

\begin{theorem}[Existence of tilings]
\label{theorem:triquadraticexistence}
Let $N$ be given.  Suppose $M$ and $K$ are positive integers such that
$M^2 + N  = 2K^2$  and $M^2 < N$ and $K  \mid N$ 
(or equivalently, $K \mid M^2$).  Let $ABC$ be the triangle with sides
$BC = MN/K$, $AC = N-K^2$, and $AB = K^2$, and let the triangle $T$ have 
sides $a=M$, $b = N/K-K$ (which must be a positive integer),
and $c = K$.   Then there is an $N$-tiling of triangle 
$ABC$ by tile $T$.  This tile satisfies $3\alpha + 2\beta = \pi$, and
 $ABC$ has angles $(2\alpha,\beta,\alpha+\beta)$.
\end{theorem}

\begin{proof}  First note that since $M^2 < N$, we have $2K^2 = M^2 + N < 2N$,
so $K^2 < N$.  Hence $K < N/K$.  Hence $b$, which is defined to be $N/K-K$, is 
positive. Since $K$ divides $N$, $b$ is a positive integer.

Define $J= c-b = K-(N/K-K) = 2K-N/K = (2K^2-N)/K = M^2/K$.  Since $K$ divides $N$,
$J$ is an integer.  We then have $c = K = M^2/J$ since $J = M^2/K$.  Since $a= M$ 
and $b = c-J$ we have
\begin{eqnarray*}
a &=& M \\
c &=& \frac {a^2} J \\
b &=& \frac {a^2} J - J
\end{eqnarray*}

Let $\alpha$, $\beta$, and $\gamma$ be the angles of the tile opposite sides
$a$, $b$, and $c$. Then we have $b = c-a^2/c$ since $b = c-J$ and $J= M^2/K = a^2/c$.
By Lemma \ref{lemma:anglesOK} we have $3\alpha + 2\beta = \pi$.

We can  construct a tiling as follows.  Fix the ``center point'' Q.  Construct
three quadratic tilings whose vertex angles meet at $Q$,  one with $a^2$ tiles, 
and two with $b^2$ tiles.  Let the two $b^2$ tilings share a common side $AQ$.
Let the $a^2$ tiling share a common side $CQ$ with one of the $b^2$ tilings.
This is possible because along $CQ$ there are, on the side with the $a^2$ tiling,
$a$ tiles, each with its $b$ edge on that side, so $CQ$ has length $ab$;  on the 
side of $CQ$ with the $b^2$ tiling, there are $b$ tiles, but each tile has its $a$
edge along $QC$, so the length on that side is also $ab$.  Hence the corners of 
the $a^2$ tiling and the $b^2$ tiling occur at the same point $C$.   There are 
 two $\alpha$ angles at $A$ and at $C$ there are an $\alpha$ and a $\beta$ angle.
 Let $D$ be the other vertex of the $b^2$ tiling that does not share side $QC$.
 Let $E$ be the other vertex of the $a^2$ tiling,  where the corner tile has a $\beta$
 angle.  Now construct point $B$ as the intersection of the line containing $AD$
 and the line containing $CE$.   Since angle $DAC$ is $2\alpha$ and angle $ACE$
 is $\alpha + \beta$, the sum of these two angles is $3\alpha + \beta < \pi$.
 Hence, by Euclid's fifth postulate, point $B$ does exist and lies on the same 
 side of $AC$ as the tilings.  
 
 Now consider the quadrilateral $BDQE$.   The interior angle at $B$ is $\beta$,
 since the angles of triangle $ABC$ must sum to $\pi$, and the angles at $A$
 and $C$ sum to $3\alpha + \beta = \pi-\beta$.   The exterior angle 
 $QEC$ is also $\beta$, so $AB$ is parallel to $QE$.   The exterior angle $ADQ$
 is also $\beta$, so $BC$ is parallel to $DQ$.   Hence quadrilateral $BDQE$ is a 
 parallelogram.   Its side $QD$ is composed 
 of $b$ edges of tiles, each of length $a$.   We can therefore divide quadrilateral
 $BDQE$ into $b$ parallegrams by drawing $b-1$ equally spaced lines parallel to $EB$. 
 Each of these parallelograms has a side equal to $QE$.   $QE$ is composed of $a$ 
 edges of tiles, with each edge 
 of length $a$, so the length of $QE$ is $a^2$.  But $a^2 = Jc$.  It is 
 therefore possible to break each of the $b$ small parallelgrams into $J$ yet 
 smaller parallelograms, with sides $c$ and $a$, and one angle $\gamma$.  
 Each of these smaller parallelograms can be cut into two copies of the tile $T$.
 We can thus tile quadilateral $BDQE$ by $2bJ$ tiles.  (Note that these tiles 
 extend the $b^2$ quadratic tiling of $QAD$, so they would be part of a $c^2$ 
 quadratic tiling,  but the vertex of that larger quadratic tiling would overlap 
 the $a^2$ tiling.)
 
 It remains to  count the tiles and verify that there are $N$ of them.  The 
 total number of tiles, which we temporarily call $Z$ until we prove it is equal to $N$, is
 \begin{eqnarray*}
  Z &=& 2b^2 + a^2 + 2bJ 
  \end{eqnarray*}
  We substitute $a=M$ and $b = N/K-K$, obtaining
  \begin{eqnarray*}
 Z &=& 2(N/K-K)^2 + M^2 + 2( N/K-K) J 
 \end{eqnarray*}
 Now substitute $N = 2K^2-M^2$;  then $N/K-K = 2K-M^2/K -K = K-J$, since $J=M^2/K$. 
 We obtain 
 \begin{eqnarray*}
  Z &=& 2(K-J)^2 + M^2 + 2(K-J)J \\
  &=&  2K^2 - 4KJ + 2J^2 + M^2 + 2KJ - 2J^2 \\
  &=& 2K^2 -2KJ + M^2 \\
  &=& 2K^2 -2K(M^2/K) + M^2 \mbox{\qquad since $J = M^2/K$}\\
  &=& 2K^2 -2M^2 + M^2 \\
  &=& 2K^2 - M^2 \\
  &=& N
 \end{eqnarray*}
 as desired.  The number of tiles is $N$.  \end{proof}
 
\begin{figure}[ht]
  \begin{center}
 \includegraphics[width=0.8\textwidth]{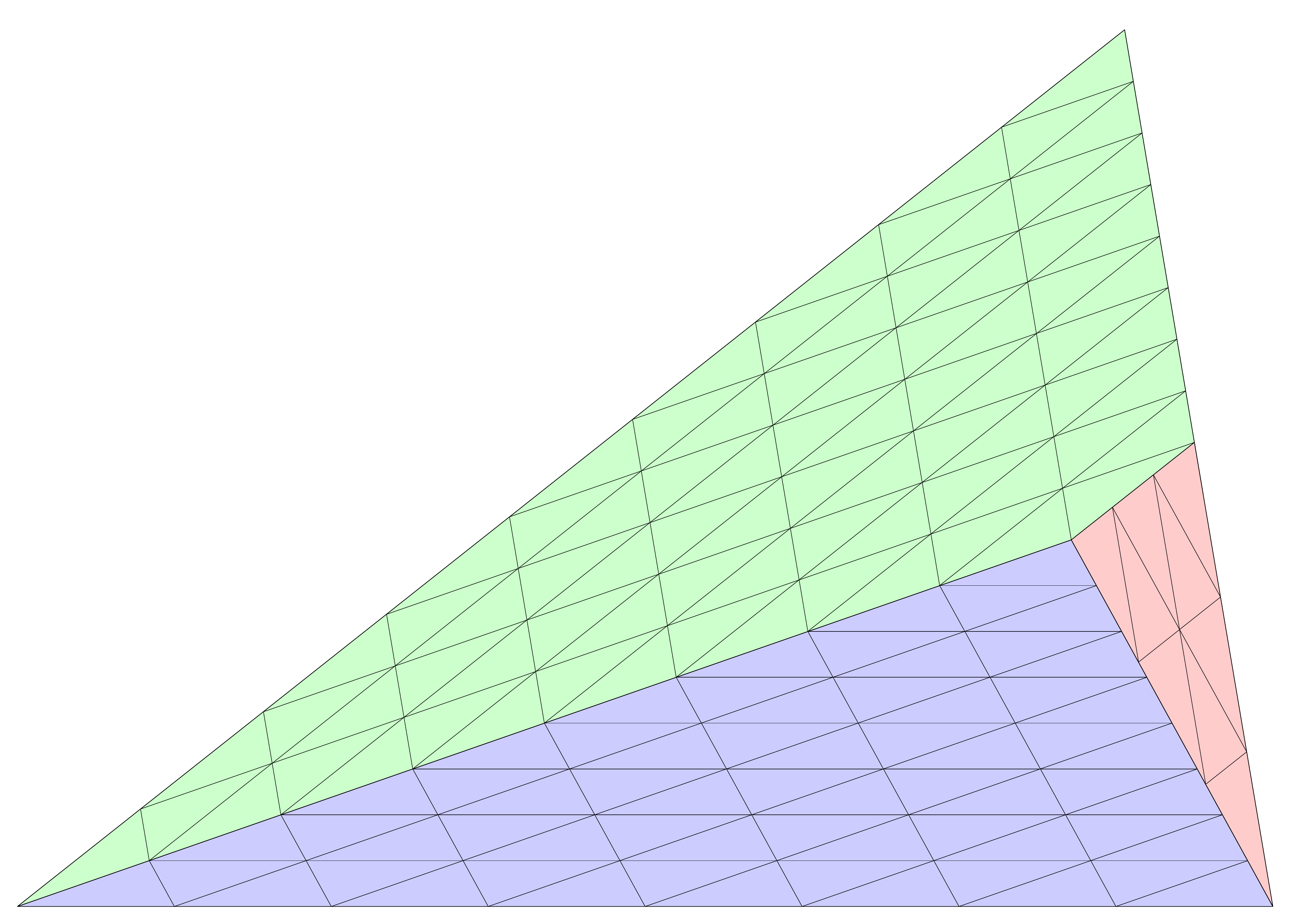} 
 \end{center}
 \caption{A triquadratic tiling with $N=153=9\cdot17, M=3, K=9$}
 \label{figure:triquadratic153}
 \end{figure}

\begin{figure}[ht]
  \begin{center}
 \includegraphics[width=0.4\textwidth]{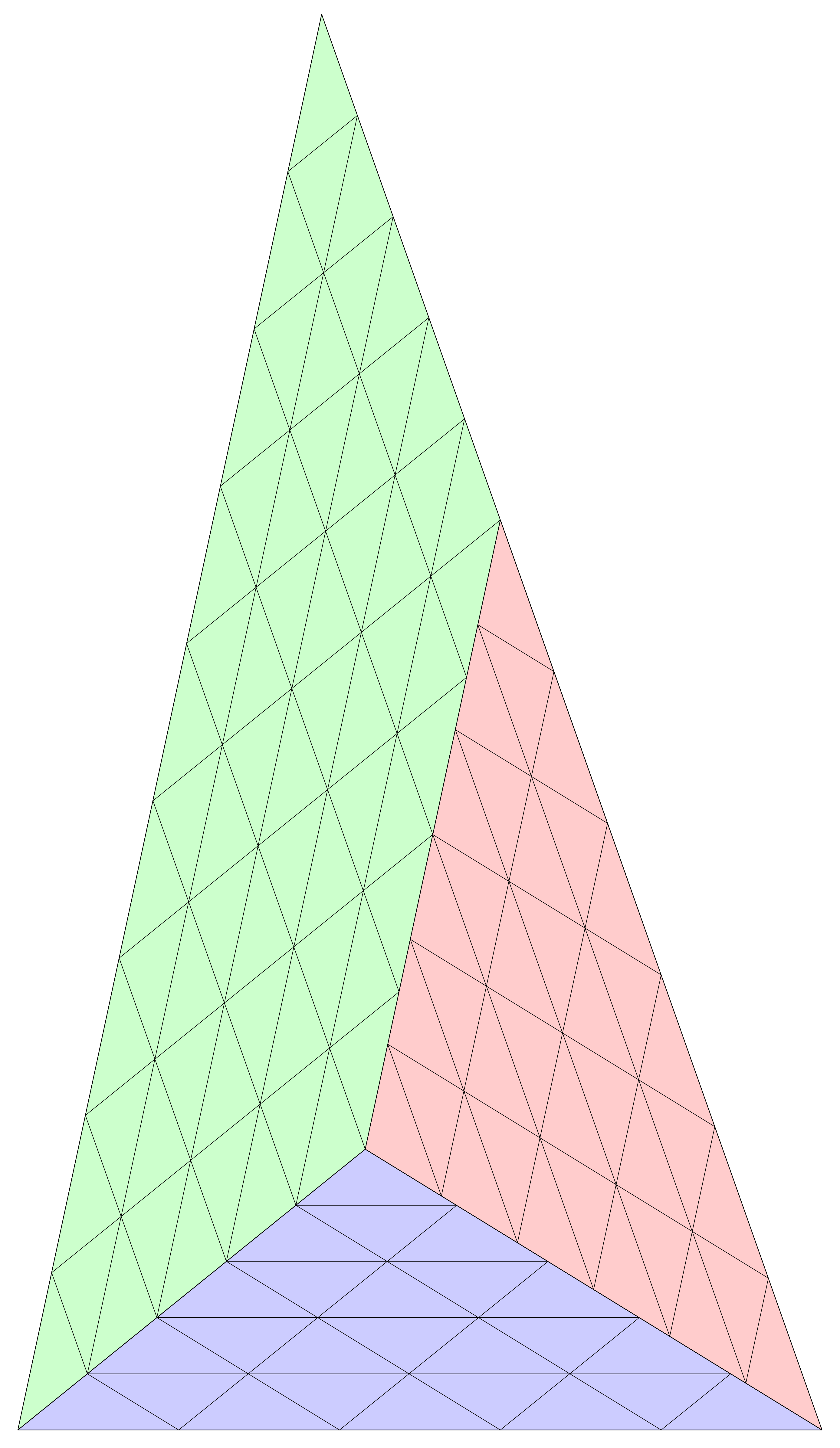}   
 \end{center}
 \caption{A triquadratic tiling with $N=126=9\cdot14, M=6, K=9$.}
 \label{figure:triquadratic126}
 \end{figure}

 Figures \ref{figure:triquadratic153}, \ref{figure:triquadratic126}, and \ref{figure:triquadratic612}
  illustrate triquadratic tilings.  Note that in Figure~\ref{figure:triquadratic126}, we have 
  $\alpha > \beta$ since $M > K/2$.  (One can show using $s = 2\sin(\alpha/2) = M/K$,  that $M > K/2$
  is equivalent to $\alpha > \beta$.)

 \begin{figure}[ht]
  \begin{center}
 \includegraphics[width=0.8\textwidth]{Figures/tiling612} 
 \end{center}
 \caption{A triquadratic tiling with $N=612=17\cdot6^2, M=6, K=18$}
 \label{figure:triquadratic612}
 \end{figure}
  
 \FloatBarrier
 
\subsection{Given $ABC$ and the tiling equation, a tiling exists for some $N$}

The existence theorem,  Theorem~\ref{theorem:triquadraticexistence}, 
has the number of tiles $N$  and the solution of $N = 2K^2- M^2$  
as given data;  then the tile must be $(a,b,c)$ with $s = a/c = M/K$.
Recall from Lemma~\ref{lemma:jmk} that there are coprime positive integers $m < k$
such that $$(a,b,c) = (mk,k^2-m^2,k^2).$$
Thus $ s= a/c = m/k=  M/K$. 

The existence theorem requires the additional hypothesis that $M \mid K^2$.
But the tiling equation is often solvable when that hypothesis is not satisfied.
In those cases we cannot construct an $N$-tiling.  But we {\em can} still 
construct a $j^2 N$ tiling, for some $j$.   We will prove that now.

\begin{theorem}\label{theorem:jshape}
Assume $3\alpha + 2\beta = \pi$.   
Let triangle $ABC$ have the angles $(2\alpha, \beta, \alpha + \beta)$,
and let $s = 2\sin \tfrac \alpha 2$.   Suppose $M$ and $K$ are positive
integers such that $s = M/K$.  (It is not assumed that $M/K$ is in lowest terms.)
 Then there is a tiling
of $ABC$ using a tile similar to $(\alpha,\beta,\gamma)$, 
and the number of tiles is a square times $2K^2-M^2$.
\end{theorem}

\begin{remark}
 It is not assumed that
$K \mid M^2$,  which is a hypothesis of Theorem~\ref{theorem:triquadraticexistence}.
That additional hypothesis guarantees that the number of tiles 
is {\em exactly} $N = 2K^2-M^2$.   But for example, with $M=2$ and $K=3$,
we have $14 =2K^2-M^2$,  but no 14-tiling.  If, however, we change $(M,K)$ to $(6,9)$,
we have  
$$3^2 \cdot 14 = 2K^2-M^2 = 126 = 2 \cdot 9^2  - 6^2,$$
and the existence of a 126-tiling is guaranteed since with the new $M$ and $K$,
we do have $K \mid M^2$.
\end{remark}

\begin{proof}[Proof of Theorem~\ref{theorem:jshape}]
Let $ABC$, $K$, and $M$ be as in the theorem.   Let $N = 2K^2-M^2$.
Let $j$ be any positive integer.  Then we also have $s = jM / (jK)$,
$$j^2 N = 2 (jK)^2 - (jM)^2. $$
Setting $M^\prime = jM$ and $K^\prime = jK$,
we have $s = M^\prime/K^\prime$.
Now the condition $K^\prime \mid (M^\prime)^2$ says $jK \mid j^2M^2$,
that is, $K \mid jM^2$; so it holds whenever $K \mid j$.  Choose, then,
$j = K$, or any multiple of $K$.  Then the hypothesis
$K^\prime \mid (M^\prime)^2$ of
Theorem~\ref{theorem:triquadraticexistence} is satisfied; and its
hypothesis $(M^\prime)^2 < j^2N$ holds too, since
$(M^\prime)^2 = j^2M^2 < j^2N$.  Hence that theorem guarantees the
existence of a $j^2N$-tiling of $ABC$.
\end{proof}

Until 2026, we did not know any example to refute the conjecture that 
if $N = 2K^2-M^2$ is solvable with $M<N$,  and there is an 
$N$-tiling, then $K \mid M^2$.  But Bonfioli's 63-tiling settled
that question, as $63 = 2\cdot 6^2 - 3^2$, but $6 \nmid 3^2$.  That
is why the $63$-tiling was not constructed like the 28-tiling, but instead
discovered by computer search (trial-and-error).

Let $j = \gcd(M,K)$ and write $M = jm$, $K = jk$, so that
$s = m/k$ in lowest terms and, by Lemma~\ref{lemma:jmk}, the tile is
$(a,b,c) = (mk,\,k^2-m^2,\,k^2)$.  Since $\gcd(m,k) = 1$, the
condition $K \mid M^2$ is equivalent to $k \mid j$.  Because $m < k$
we have $k \ge 2$, so the condition always fails when $j = 1$.
Putting $N_0 = 2k^2 - m^2$, the integers $N$ for which the tile
$(mk,k^2-m^2,k^2)$ occurs at all are $N = j^2N_0$, and those
satisfying $K \mid M^2$ are exactly $N = k^2t^2N_0$ for
$t = 1,2,\ldots$; only one in every $k$ of them.  For $N \le 500$
there are $137$ solutions of the tiling equation, of which just eight
satisfy $K \mid M^2$, as shown in Table~\ref{table:triquadratics}.
The proportion decreases: among the $155{,}216$ solutions with
$N \le 500{,}000$, only $816$ satisfy $K \mid M^2$.  Thus
Theorem~\ref{theorem:triquadraticexistence} settles the existence
question for only a sparse set of $N$.  It is shown in \cite{beeson-rosu2026}
that no tilings exist with $j=1$; in particular that include $N$ prime.
There remain three values of $N$ which have $j=2$, so are not settled by \cite{beeson-rosu2026};
namely 56, 68, and 92.   
Bonfioli \cite{bonfioli2026} reports Lean proofs that no triangle has an
$N$-tiling for $N = 56$ and $N = 92$.%
\footnote{ Theorems
{\tt Frontier34Wrap.tileCount\_iff\_le\_104} and {\tt N92.not\_tile\_92} in
the repository \url{https://github.com/ElVec1o/erdos_634_proof}, commit
{\tt 119bec0}); we have not checked them in Lean.
}
For $N = 68$, the triquadratic triangle has no $6$8-tiling \cite[Theorem {\tt W2E1.not\_wTilable\_1\_f\_2}]{bonfioli2026}, 
which excludes every $(m,k) = (1,k)$ with $j = 2$, $k \ge 3$), although $68 = 2^2(1^2+4^2)$ tiles some other triangle.

\begin{table}[ht]
\caption{All $N\le 500$ with $N = 2K^2-M^2$, $M^2 < N$, and $K \mid M^2$.
Tilings exist for these $N$.}
\label{table:triquadratics}
\begin{center}
\begin{tabular}{rrrrrcc}
$N$  &  $M$  &  $K$  &  $M^2/K$  &  $j=\gcd(M,K)$  &  $(m,k)$  &  $(a,b,c)$  \\
  \hline
28  & 2  & 4  & 1 & 2 & $(1,2)$ & $(2, 3, 4)$   \\
112 & 4  & 8  & 2 & 4 & $(1,2)$ & $(2, 3, 4)$   \\
126 & 6  & 9  & 4 & 3 & $(2,3)$ & $(6, 5, 9)$   \\
153 & 3  & 9  & 1 & 3 & $(1,3)$ & $(3, 8, 9)$   \\
252 & 6  & 12 & 3 & 6 & $(1,2)$ & $(2, 3, 4)$   \\
368 & 12 & 16 & 9 & 4 & $(3,4)$ & $(12, 7, 16)$ \\
448 & 8  & 16 & 4 & 8 & $(1,2)$ & $(2, 3, 4)$   \\
496 & 4  & 16 & 1 & 4 & $(1,4)$ & $(4, 15, 16)$ 
\end{tabular}
\end{center}
\end{table}

\begin{table}[ht]
\caption{All $N \le 100$ for which the equation $N = 2K^2-M^2$, $M^2<N$,
is solvable but $K \nmid M^2$, and $N$ is not squarefree.}
\label{table:undecided}
\begin{center}
\begin{tabular}{rrrrccl}
$N$  &  $M$  &  $K$  &  $j=\gcd(M,K)$  &  $(m,k)$  &  $(a,b,c)$  &  exists? \\
  \hline
56 & 4 & 6 & 2 & $(2,3)$ & $(6, 5, 9)$ &\nomark \cite{bonfioli2026}\\
68 & 2 & 6 & 2 & $(1,3)$ & $(3, 8, 9)$ &\nomark \cite{bonfioli2026}\\
92 & 6 & 8 & 2 & $(3,4)$ & $(12, 7, 16)$ &\nomark \cite{bonfioli2026}\\
\end{tabular}
\end{center}
\end{table}

 \subsection{Number theory of the tiling equation}
Here we show that the tiling equation always has a solution when $N$ has the right divisibility 
properties.  Whether that solution corresponds to a tiling depends on whether $K$ divides $M^2$ or not.

\begin{lemma}
\label{lemma:6902} If the tiling equation $N = 2K^2-M^2$ is solvable in
integers, then $N$ is a square times a product of distinct primes, each of
which is either $2$ or is of the form $8n\pm 1$.
\end{lemma}

\begin{proof} Write $N = f^2N_0$ with $N_0$ squarefree.  Since $N_0$ is
squarefree it is a product of distinct primes, so it suffices to prove that
every odd prime $p$ dividing $N_0$ satisfies $p \equiv \pm 1$ mod~8.  Note
that $p \mid N_0$ says exactly that $p$ divides $N$ to an odd power.  (It is
not true that every odd prime dividing $N$ must be congruent to $\pm 1$;
for example $63 = 2\cdot 6^2 - 3^2$ is divisible by 3.)

Let $p$ be such a prime, and let $v = v_p(N)$, which is odd.  Since $p$ is
odd we have $v_p(2K^2) = 2v_p(K)$ and $v_p(M^2) = 2v_p(M)$.  If those two
were different, then $v_p(2K^2-M^2)$ would be their minimum, which is even,
contradicting the oddness of $v$.  Hence $v_p(M) = v_p(K) = t$ for some
$t \ge 0$.  Write $M = p^tM_1$ and $K = p^tK_1$, so that $p \nmid M_1$ and
$p \nmid K_1$, and
$$ N \ = \ p^{2t}\,(2K_1^2 - M_1^2).$$
Then $v_p(2K_1^2-M_1^2) = v-2t$, which is odd and nonnegative, hence at
least 1.  Therefore
$$ 2K_1^2 \ \equiv \ M_1^2 \pmod p. $$
Since $p \nmid K_1$ we may divide by $K_1^2$, obtaining
$2 \equiv (M_1K_1^{-1})^2$ mod $p$.  Thus 2 is a quadratic residue mod $p$,
and hence $p \equiv \pm 1$ mod~8.
\end{proof}

\begin{lemma} \label{lemma:normtheorem}
Suppose $N$ is a square times a product of distinct primes, each of which is
either 2 or is of the form $8n\pm 1$.  Then the tiling equation
$N+M^2 = 2K^2$ has a solution in integers.  If moreover $N$ is neither a
square nor twice a square, then it has a solution in positive integers $M$
and $K$ with $M < K$.
\end{lemma}
\begin{remark}  But a solution with $K \mid M^2$ is not guaranteed.
So we cannot conclude that some triangle can be $N$-tiled,  but only
that some triangle can be $j^2N$-tiled, for some $j$.  For example,
there is no 7-tiling,  but there is a 28-tiling.
\end{remark}

\begin{proof}[Proof of Lemma~\ref{lemma:normtheorem}]%
\footnote{The idea of this proof was given to me on MathOverflow by Noam
Elkies and Will Jagy.}
We first prove that there exists an integer solution $(M,K)$.  After that we
will prove that, under the extra hypothesis, there is one with $0 < M < K$.

The tiling equation asks for an integer in the ring $\Z[\sqrt 2]$ of
integers of the field $\Q(\sqrt 2)$ whose norm is $-N$.  In detail, those
integers have the form $M + K\sqrt 2$, and the norm of such an integer is
$M^2-2K^2$.  The norm is multiplicative, and $1+\sqrt 2$ has norm $-1$; so
an integer of norm $n$, multiplied by $1+\sqrt 2$, gives one of norm $-n$,
and it does not matter whether we seek norm $N$ or norm $-N$.  Again by
multiplicativity, a square factor $f^2$ of $N$ can be contributed by the
integer $f$, and a factor 2, if present, by $\sqrt 2$, whose norm is $-2$.
We are therefore reduced without loss of generality to the case when $N$ is
an odd prime congruent to $\pm 1$ mod~8.  The form $x^2-2y^2$ has
discriminant 8, so by Theorem 4.23 on p.~74 of \cite{buehl}, it integrally
represents every prime congruent to $\pm 1$ mod~8.  Hence the tiling
equation does have some integer solution $(M,K)$.

It remains to prove that there is a solution with $0 < M < K$.  (Now we no
longer assume $N$ is prime, and we assume $N$ is neither a square nor twice
a square.)  One can verify by computation that if $(M,K)$ is a solution of
the tiling equation, then so is
$$ (3M-4K,\ 3K-2M); $$
one can verify it without computation by observing that
$$ (M+K\sqrt 2)(3-2\sqrt 2) \ = \ (3M-4K) + (3K-2M)\sqrt 2, $$
and $3-2\sqrt 2$ is a unit, of norm 1.  Since only $M^2$ and $K^2$ occur in
the equation, we may also change the sign of either coordinate.  Note that
every solution satisfies $M^2 < 2K^2$, since $N > 0$; hence $K \ge 1$, and
$M < K\sqrt 2$, which gives both $M < 2K$ and $2M < 3K$.

This set is not empty.  The first part gives some integer solution $(M,K)$;
we have $K \ne 0$, since $K=0$ would give $N = -M^2 < 0$, and $M \ne 0$,
since $M=0$ would give $N = 2K^2$, twice a square, contrary to hypothesis.
Changing signs as above, we may take $M > 0$ and $K \ge 0$.
Now choose a solution with $M$ positive and as small as possible, and with
$K \ge 0$.  If $3M-4K$ were positive, then $(3M-4K,\,3K-2M)$ would be a
solution with positive first coordinate, and $3M-4K < M$ because $M < 2K$,
contradicting the minimality of $M$.  Hence $3M-4K \le 0$.  If $M < K$ we
are done, so we may assume $K \le M$.  Put
$$ M^\prime \ = \ 4K-3M, \qquad K^\prime \ = \ 3K-2M. $$
Then $(M^\prime,K^\prime)$ is also a solution, $M^\prime \ge 0$ since
$3M-4K \le 0$, and $K^\prime > 0$ since $2M < 3K$.  Moreover
\begin{eqnarray*}
K^\prime - M^\prime &=& (3K-2M)-(4K-3M) \\
    &=& M-K \\
    &\ge& 0 \mbox{\qquad since $K \le M$,}
\end{eqnarray*}
so $M^\prime \le K^\prime$.  But if $M^\prime = K^\prime$ then
$N = 2K^{\prime 2} - M^{\prime 2} = K^{\prime 2}$ is a square.  And if
$M^\prime = 0$ then $4K = 3M$, so $4 \mid M$, and
$$ N \ = \ 2K^2-M^2 \ = \ 2\Bigl(\frac{3}{4}\Bigr)^2M^2 - M^2
   \ = \ \frac{M^2}{8} \ = \ 2\Bigl(\frac{M}{4}\Bigr)^2 $$
is twice a square.  Both are excluded by hypothesis, so
$0 < M^\prime < K^\prime$, which completes the proof.
\end{proof}

\subsection{Can $N$ be prime?}
In the previous version of this paper,  it was claimed that $N$ cannot be prime
in a tiling of a triangle of the shape $(2\alpha, \beta, \alpha+\beta)$. 
That proof was incorrect.  However, the result is correctly proved
(and verified in Lean) in \cite{beeson-rosu2026}.  
As a partial result with the methods at hand, we have the following lemma.

\begin{lemma} \label{lemma:tqprime}
If a triquadratic triangle is tiled by $N$ copies of $T$ and $N$ is prime, then
$N \equiv \pm 1 \pmod 8$.  Moreover $m$ is odd, and $k$ is even or odd according as $N \equiv -1$ or
$N \equiv 1 \pmod 8$.
\end{lemma}

\begin{proof}
By Corollary~\ref{corollary:tqparameters}, $N = j^2(2k^2-m^2)$ with
$\gcd(m,k)=1$; since $N$ is prime, $j=1$ and $N = 2k^2-m^2$.  Then
$m^2 \equiv 2k^2 \pmod N$, and $N \nmid k$, since $N \mid k$ would give
$N \mid m^2$ and hence $N \mid m$, contradicting $\gcd(m,k)=1$.  So $2$ is a
quadratic residue modulo $N$, whence $N \equiv \pm1 \pmod 8$.

If $m$ were even then $k$ would be odd, so $2k^2 \equiv 2$ and
$m^2 \equiv 0$ or $4 \pmod 8$, making $N$ even; so $m$ is odd and
$m^2 \equiv 1 \pmod 8$.  If $k$ is odd then $2k^2 \equiv 2 \pmod 8$ and
$N \equiv 1$; if $k$ is even then $2k^2 \equiv 0 \pmod 8$ and $N \equiv -1$.
\end{proof}

Next we connect question whether $N$ is prime 
 to the related proposition $j=1$, where $j=\gcd(M,K)$.

As shown above, the condition $K \mid M^2$ is equivalent to $k \mid j$, where $j = \gcd(M,K)$,
so $j(m,k) = (M,K)$.  Many of the examples where the tiling equation is solvable and $K \nmid M^2$
have $j=1$,  as far as possible from $k \mid j$.  Moreover,  we are interested in the case 
of prime $N$, and $N$ prime implies $j=1$, as we now show.

\begin{lemma}\label{lemma:primeimpliesjone}  Let $j,m,k,M,K$, and $N$ be positive integers. 
Suppose $j(m,k) = (M,K)$, with $M > 0$  and $N= 2K^2-M^2 > 0$.   Assume $N$ is prime.  Then $j=1$.  
\end{lemma}

\begin{proof}
From $j(m,k) = (M,K)$ we get
$$ N \ = \ 2K^2-M^2 \ = \ j^2(2k^2-m^2). $$
Suppose $j > 1$.  Then $j^2$ is a divisor of $N$ with $j^2 > 1$, so since $N$
is prime we must have $j^2 = N$ and $2k^2-m^2 = 1$.  But a prime is not a
perfect square.  Hence $j = 1$.
\end{proof}

Thus one way to prove there is no tiling with $N$ prime would be to prove there 
is no tiling with $j=1$.  This is the method used in \cite{beeson-rosu2026}.
Indeed no tiling with $j=1$ can arise from
Theorem~\ref{theorem:triquadraticexistence}, since its hypothesis
$K \mid M^2$ is equivalent to $k \mid j$, and $m<k$ forces $k \ge 2$, so
$j \ge k \ge 2$.  Every tiling we can construct therefore has $j \ge 2$; the
case $j=1$ is exactly the case in which our existence theorem gives nothing.
The proof that no tiling has $j = 1$ settles more than the
primality question: it also rules out, for example, $N = 46$, $49$,
$119$, and $158$, which are composite but have $j=1$.

\section{The case when $ABC$ has angles $(2\alpha, \alpha, 2\beta)$}
 
\begin{figure}
\centering
\includegraphics[width=0.4\textwidth]{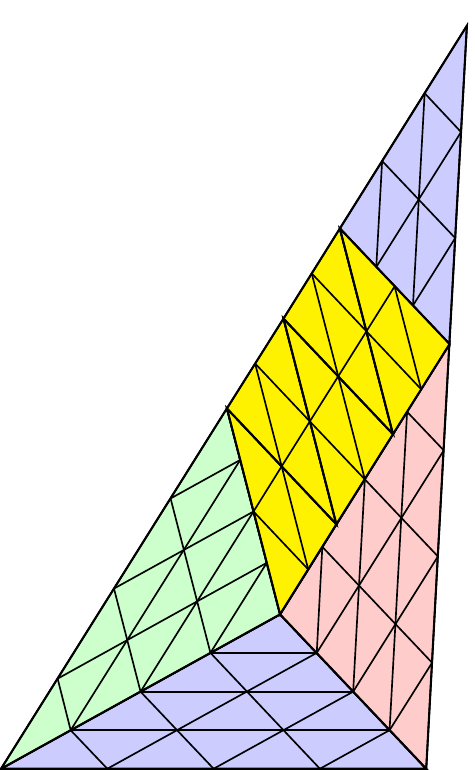}  
\caption{$N=77$.  The tile $(a,b,c)$ is $(2,3,4)$.}
\label{figure:77}
\end{figure}

\begin{figure}
\centering
\includegraphics[width=0.8\textwidth] {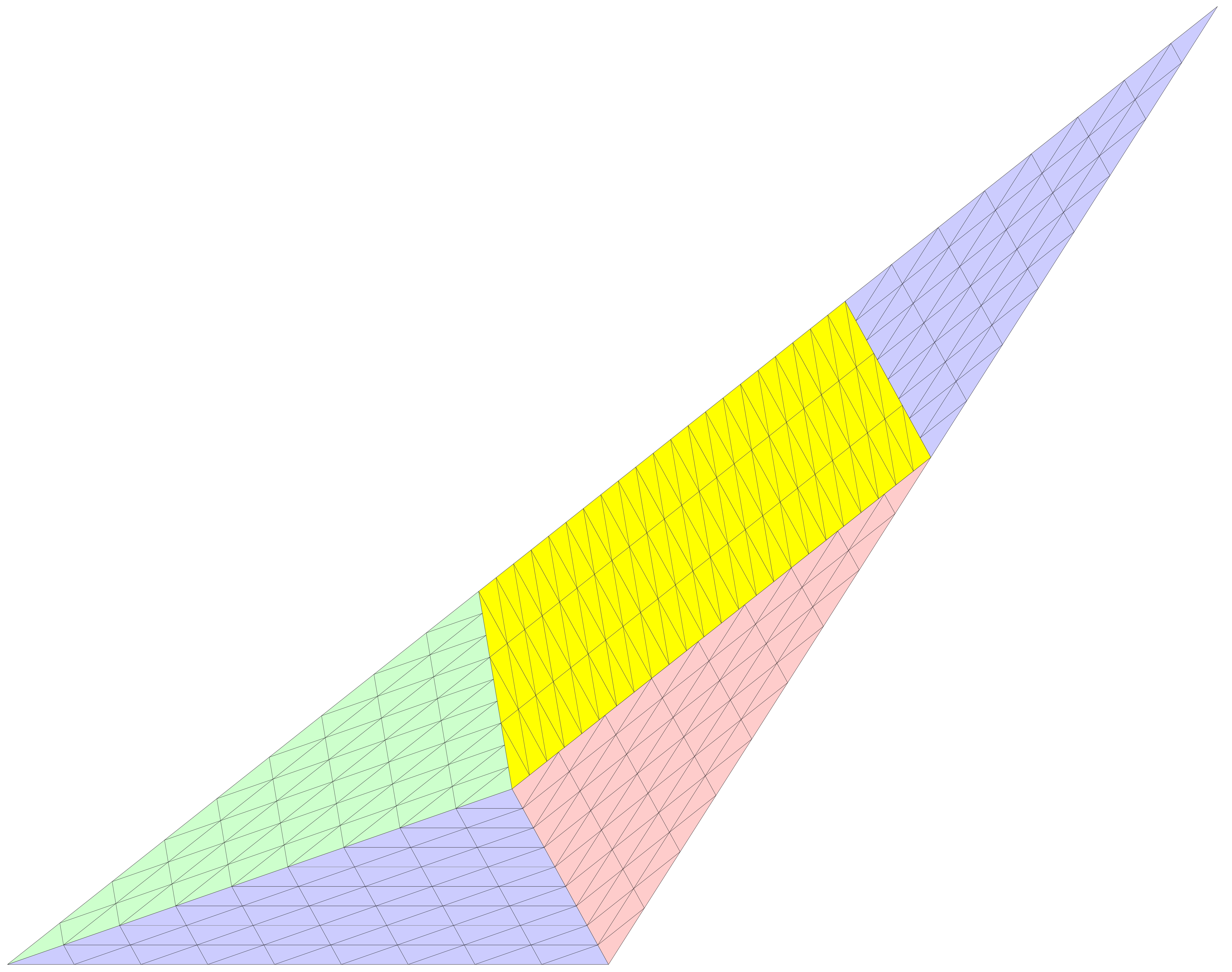} 
\caption{$N=442$. The tile $(a,b,c)$ is $(3,8,9).$}
\label{figure:442}
\end{figure}

\begin{figure}
\centering
\includegraphics[width=0.8\textwidth]{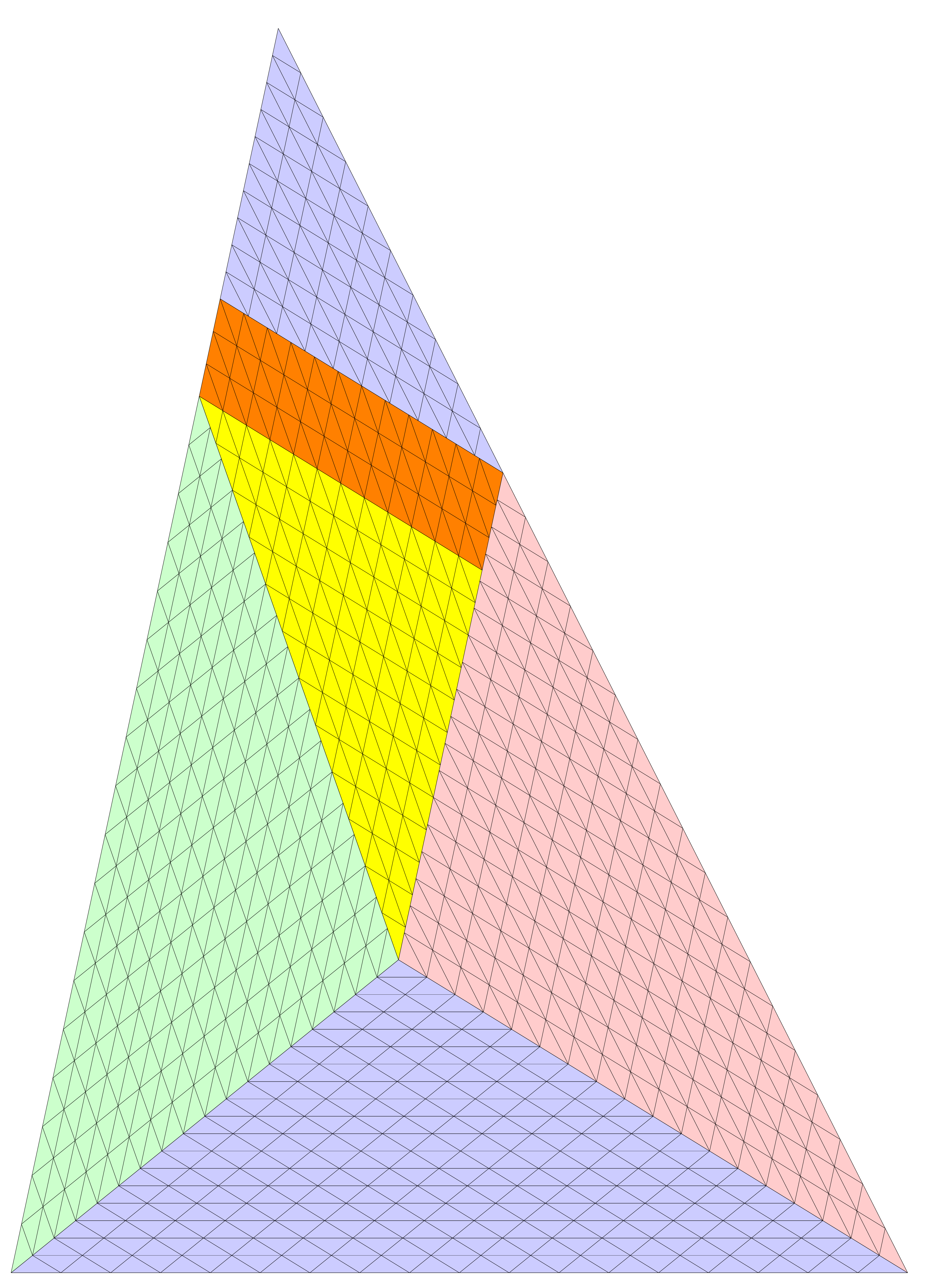} 
\caption{$N = 3\cdot 18^2 + 12^2 + 10^2 + 3\cdot2 \cdot 12 = 1288$.
The tile is $(a,b,c) = (6,5,9)$. Note, here $a > b$. The tile is 
not $(5,6,9)$ but $(6,5,9)$.
}
\label{figure:1288}
\end{figure}

 We begin by observing that for $ABC$ with angles
 $(2\alpha, \alpha, 2\beta)$,  there do exist some previously 
 unknown tilings, discovered in July, 2018.   Examples are
 given in Figs.~\ref{figure:77},\ref{figure:442}, and \ref{figure:1288}.
 It should be mentioned again that Laczkovich \cite{laczkovich1995} 
 studied decompositions of these and other triangles into similar 
 triangles, and knew that if these triangles were rational, then finer
 decompositions into some number of congruent triangles would be possible.
 Here we have proved that the tiles must be rational, and focused attention
 on $N$ as well as $ABC$; some of our tilings have much smaller $N$
 than would arise by Laczkovich's method.  
 
\subsection{The second tiling equation}
 In this section we derive the ``second tiling equation'',  which 
 provides a necessary condition for the 
 existence of an $N$-tiling of a triangle $ABC$ with angles $(2\alpha, \alpha, 2\beta)$.
 Let the sides opposite those angles be $(X,Y,Z)$ respectively.  If
 $ABC$ is $N$-tiled, let the sides of the tile be $a$, $b$, and $c$,
 which are proportional to $\sin \alpha$, $\sin \beta$,
 and $\sin \gamma$, respectively.  By Theorem~\ref{theorem:rationaltile}, the 
 tile is rational.

 \begin{lemma} \label{lemma:singamma2}
  Suppose $3\alpha + 2\beta = \pi$.  Let $s = 2\sin \tfrac \alpha 2$.  
Suppose the triangle with sides $(a,b,c)$ and angles $(\alpha,
\beta,\gamma)$ is scaled by $\kappa$, so that 
$$(a,b,c) = \kappa( \sin \alpha, \sin \beta, \sin \gamma).$$
Then 
\begin{eqnarray*}
\frac {\kappa \sin 2\alpha} c &=& s(2-s^2)\\ 
\sin 2 \alpha &=& (2-s^2) \sin \alpha \\
\frac {\kappa \sin 2\beta} c &=& \frac {\kappa \sin 3\alpha} c  \ = \ 
s(1+s)(1-s)(3-s^2)
\end{eqnarray*}
\end{lemma} 

\begin{proof} 
In Lemma~\ref{lemma:singamma} we proved
\begin{eqnarray*}
\sin \gamma &=& \cos \frac \alpha 2 \\
\frac a c &=& s \\
\frac b c &=& 1-s^2  
\end{eqnarray*}
 Using these formulas, we argue 
 as follows.
 
\begin{eqnarray*}
\frac {\sin 2\alpha} {\sin \gamma} &=& \frac {2 \sin \alpha \cos \alpha} {\sin \gamma} \\
&=& 2 \frac a c \cos \alpha \\
&=& 2s \cos \alpha \\
&=& 2s(1-2\sin^2 \frac \alpha 2) \\
&=& s(2-s^2 )
\end{eqnarray*}
Multiplying numerator and denominator on the left by $\kappa$, we obtain the first formula of the lemma.

To prove the second formula:
$$ \sin 2\alpha = s(2-s^2) \sin \gamma = s(2-s^2)c/\kappa = (2-s^2)a/\kappa
 = (2-s^2) \sin \alpha.$$

Next we work on the third formula.
\begin{eqnarray*}
\frac {\kappa\sin 2\beta} c &=& \frac{\kappa\sin 3\alpha} c \mbox{\qquad since $3\alpha + 2\beta = \pi$} \\
&=& \frac{\kappa\sin \alpha \cos 2\alpha + \kappa\sin 2\alpha \cos \alpha} c \\
&=& \frac a c \cos 2 \alpha + \frac{\kappa\sin 2\alpha}c \cos \alpha \\
&=& s \cos 2\alpha + s(2-s^2) \cos \alpha \\
&=& \ s (2\cos^2 \alpha-1)+ s(2-s^2) (1-2\sin^2 \frac \alpha 2) \\
&=& s(2(1-2\sin^2\frac \alpha 2)^2-1) + s(2-s^2)(1-\frac 1 2 s^2) \\
&=& s(2(1-\frac 1 2 s^2)^2 -1) + s(2-s^2)(1-\frac 1 2 s^2)\\
&=& s(1+s)(1-s)(3-s^2)
\end{eqnarray*}
\end{proof}
\smallskip

\begin{theorem}
\label{theorem:secondtilingequation}
Let $3\alpha+2\beta=\pi$, let $ABC$ have angles $(2\alpha,\alpha,2\beta)$,
and suppose $ABC$ is $N$-tiled by a tile with angles $(\alpha,\beta,\gamma)$
and $ABC$ is not similar to $T$.
Let the tiles be colored alternately black
and white with a black tile at $B$, and let $M$ be the number of black
tiles minus the number of white tiles. Let $s = a/c$.  Then $s$ is 
rational and the ``second tiling equation'' is satisfied:
\begin{eqnarray*}
\frac  N {M^2}  &=&  \frac{(2-s^2)(3-s^2)}{(1-s)^2(2+s)^2}
\end{eqnarray*}
\end{theorem}

\begin{example} Consider the 77-tiling from Figure~\ref{figure:77}. Then $N=77$, $M=5$, $s=\frac 1 2 $, $(a,b,c) = (2,3,4)$. 
 Both sides evaluate to $77/25$.
 \end{example}

\begin{example}
Consider the 1288-tiling from 
Figure~\ref{figure:1288}. Then $M=16$, as is verified using the 
observation that in a quadratic tiling with $m^2$ tiles, there is 
an excess of $m$ black over white tiles, if the top tile is black.
Therefore $M$ for this example is given by
 $M = 10-12+0+18+18-18 = 16$.  So now, 
let's check that the second tiling equation is satisfied in this example:
\begin{eqnarray*}
\frac{N}{M^2}   &=& \frac {(2-s^2)(3-s^2)}{(1-s)^2(2+s)^2}\\
N  &=& 16^2\frac {(14/9)(23/9)}{(1/3)^2(8/3)^2} \\
N &=& 16^2 \frac {7\cdot 23}{32} \\
N &=&  8 \cdot 7 \cdot 23 \\
N &=& 1288
\end{eqnarray*}
So the second tiling equation correctly predicts the number of tiles 
from $M$ and $s$.  
\end{example}

\begin{proof}[Proof of Theorem~\ref{theorem:secondtilingequation}] 
By Theorem~\ref{theorem:rationaltile}, $s$ 
is rational.  By the law of sines we have, for some $\lambda > 0$,  
\begin{eqnarray*}
 \vector X Y Z &=& \lambda \vector { \sin 2\alpha} {\sin \alpha} {\sin 2 \beta}  \end{eqnarray*} 
By Theorem~\ref{theorem:coloring} we have $M(a+b+c) = X-Y+Z$.
Hence
\begin{eqnarray*}
M(a+b+c) &=& \lambda(\sin 2\alpha - \sin \alpha + \sin 2\beta) \\
\end{eqnarray*}
Now we consider the size and scaling of the tile.  For some $\kappa > 0$,
we have $$(a,b,c) = \kappa (\sin \alpha, \sin \beta, \sin \gamma).$$
Dividing and multiplying the right side of the previous equation by $\kappa$,
we have
\begin{eqnarray*}
M(a+b+c) &=& \frac {\lambda}{\kappa}
           (\kappa\sin 2\alpha - \kappa\sin \alpha + \kappa\sin 2\beta) \\
\end{eqnarray*}
Dividing by $c$ and then using Lemmas~\ref{lemma:singamma}
and \ref{lemma:singamma2}, and $\kappa \sin \alpha = a$,  to 
express everything in terms of $s = 2 \sin(\alpha/2)$, we have
\begin{eqnarray}
M(\frac a c + \frac b c + 1) &=& \frac {\lambda} {\kappa}( s(2-s^2)-s + s(1+s)\nonumber(1-s)(3-s^2)) \nonumber\\
M(s + (1-s^2) + 1)&=& \frac {\lambda} {\kappa} s(2+s)(1+s)(1-s)(2-s)\nonumber\\
M(2-s)(1+s) &=&\frac {\lambda} {\kappa} s(2+s)(1+s)(1-s)(2-s) \nonumber \\
M &=& \frac {\lambda} {\kappa} s(1-s)(2+s) \label{eq:M}
\end{eqnarray}

Twice the area of one copy of the tile is  $bc \sin \alpha $.
Equating twice the area of $N$ tiles to twice the area of $ABC$, we have
$$Nbc\sin \alpha = XZ \sin \alpha$$
 since the angle 
opposite $Y$ is $\alpha$. 
Dividing by $c^2 \sin \alpha$ we have $Nb/c = XZ/c^2$.  Expressing this in terms of 
$s$ we have 
\begin{eqnarray}
N(1-s^2) &=&  \lambda^2\frac{ \sin 2\alpha} c \frac{\sin 2\beta} c \nonumber\\
        &=& \bigg(\frac \lambda \kappa\bigg)^2\frac{ \kappa\sin 2\alpha} c \frac{\kappa\sin 2\beta} c \nonumber\\
         &=& \bigg(\frac \lambda \kappa\bigg)^2 s(2-s^2)  s(1+s)(1-s)(3-s^2) \nonumber \\
 N    &=& \bigg(\frac \lambda \kappa\bigg)^2 s^2(2-s^2)(3-s^2) \label{eq:3895}
\end{eqnarray}
Solving for $(\lambda/\kappa)^2$ we have
\begin{eqnarray}
\bigg(\frac \lambda \kappa\bigg)^2 &=& \frac N {s^2(2-s^2)(3-s^2)}  \label{eq:lambdasquared}
\end{eqnarray}
Squaring (\ref{eq:M}) we have
\begin{eqnarray*}
M^2  &=&\bigg(\frac {\lambda} {\kappa}\bigg)^2 s^2(1-s)^2(2+s)^2 \\
\bigg(\frac \lambda \kappa\bigg)^2 &=& \frac { M^2}{s^2(1-s)^2(2+s)^2}
\end{eqnarray*}
Putting this expression into (\ref{eq:3895}),
we have
\begin{eqnarray*}
\frac N { M^2} &=& \frac {s^2(2-s^2)(3-s^2)}{s^2(1-s)^2(2+s)^2} \\
&=& \frac {(2-s^2)(3-s^2)}{(1-s)^2(2+s)^2} 
\end{eqnarray*}
That is the 
``second tiling equation''.  
\end{proof}

\begin{lemma} \label{lemma:Msquared}
If the second tiling equation has a solution with $0<s<1$, then $M^2 < N$. 
\end{lemma}

\begin{proof} It suffices to show that the right side of 
the second tiling equation is greater than 1.  That is, it suffices to show 
$$(1-s)^2(2+s)^2 < (2-s^2)(3-s^2).$$
To prove that, consider the difference:
\begin{eqnarray*}
 (2-s^2)(3-s^2)-(1-s)^2(2+s)^2 &=& -2s^3 - 2s^2 + 4s + 2
 \end{eqnarray*}
 which is positive at 0 and 1; its derivative is $-6s^2 -4s + 4$,
 which has one zero $s = 0.548\ldots$ in the interval $[0,1]$
 and is positive when $s=0$ and negative when $s=1$.  Hence 
 the displayed expression is positive on $[0,1]$.  
\end{proof}

\begin{lemma} \label{lemma:scaling}
Suppose  triangle $ABC$ is $N$-tiled, and 
the tile has angles $(\alpha,\beta,\gamma)$ satisfying
$3\alpha + 2\beta = \pi$.   
If we multiply each side of 
the tile $(a,b,c)$ by a positive number $\kappa$, and rescale 
triangle $ABC$ by the same factor $\kappa$, then we still have
a tiling, with the same $N$, $M$, and $s = a/c$.  If $s$ 
is rational, after a suitable scaling we can assume the tile 
has integer sides.
\end{lemma}

\begin{proof}
Let $s$ be rational, say $s=a/c$ with
$a$ and $c$ integers, not necessarily in lowest terms. 
 Then, 
by Lemma~\ref{lemma:anglesOK}, $b = c-a^2/c$.   
Perhaps $b$ is not an integer. 
For example with $s=2/3$ we could have $(a,b,c) = (2,5/3,3)$.
In that example, we could consider instead $(6,5,9)$.
In general, if $b$ is not an integer, we could  scale $(a,b,c)$
by multiplying by $c$.  Then the new $(a,b,c)$ is 
$(ac,c^2-a^2,c^2)$, and all sides of this triangle are integers.
\end{proof}

\begin{remark}
This scaling does not change $s$; it only changes the 
marks on the ruler that we use to measure lengths.  But therefore,
it does change not only $(a,b,c)$, but also $(X,Y,Z)$, as the 
lengths of the sides of $ABC$ remain the same when measured in 
terms of $a$, $b$, and $c$, but they change in absolute units.
  Thus, for example,
if the sides $a$, $b$, and $c$ are tripled, the sides of $ABC$ also
would triple, but $N$ would remain the same.  The expansion factor
$\lambda$ such that $Y = \lambda \sin \alpha$,  would triple.
But $\lambda/\kappa$ would not change, and $M$, which is a function 
of $s$ times $\lambda/\kappa$, would also not change.  That is good,
since $M$ is the excess number of black tiles over white tiles, and 
the tiles do not change when we change the scale.
\end{remark}

\begin{example}
With
$s=2/3$, $N=1288$, $M=16$, the second tiling equation is solved.
Here $s = 2/3$;  if we write $s = 6/9$ then 
$c$ divides $a^2$, and the tile $(6,5,9)$ has integer sides,
and there is an $N$-tiling of triangle $ABC$ as shown in
 Figure~\ref{figure:1288}.  But if we   
write $s$ as $2/3$ then $c$ does not divide $a^2$ 
and $b = c-a^2/c = 5/3$ is not an integer.  Nevertheless
there is an $N$-tiling of a triangle similar to $ABC$, but one-third
the size, by the tile $(2,5/3,3)$, which is one-third the size of $(6,5,9)$.
 \end{example}

\subsection{Tiling construction for $(2\alpha, \alpha,2\beta)$, first method.}
\begin{theorem}[Tiling construction]   \label{theorem:existence2}
Suppose $c$ divides $a^2$, where $a$ and $c$ are integers with $a < c$,
let $b = c-a^2/c$ (so $3\alpha+2\beta = \pi$ and $b$ is an integer), and
suppose $c < 2b$.  Let $s = a/c$, $\ell = \lcm(a,c)$, and $j = \ell/c$.  Then some
triangle $ABC$ with angles $(2\alpha, \alpha, 2\beta)$ can be $N$-tiled by
the tile $(a,b,c)$, with
\begin{eqnarray*}
 M &=& \ell(2-s^2-s) = j(b+c-a),  \\
 N &=& \ell^2(2-s^2)(3-s^2) = j^2(b+c)(b+2c).
 \end{eqnarray*}
 \end{theorem}

\begin{remark} \label{remark:existence2}
Theorem~\ref{theorem:existence2} does not produce a tiling for every $N$
for which the second tiling equation is solvable with $s = a/c$.  Since
$(1-s)(2+s) = (b+c-a)/c$ and $(2-s^2)(3-s^2) = (b+c)(b+2c)/c^2$,
the second tiling equation can be written as
$$ N \ = \ (b+c)(b+2c)\Bigl(\frac{M}{b+c-a}\Bigr)^2 , $$
which has the integer solution $M = b+c-a$, $N = (b+c)(b+2c)$: the
number of tiles required, by area, for the triangle $ABC$ whose side
opposite the angle $\alpha$ has length $c^2$, which is $1/j$ times the size of
the triangle constructed in the theorem.  When $j = 1$, that is, when $a$
divides $c$, the theorem and its quadratic refinements account for every
solution of the second tiling equation with this $s$.  When $j > 1$ they do not.
For example, with $(a,b,c) = (6,5,9)$ we have $j = 2$; the theorem gives
$M = 16$ and $N = 1288$, but $M = 8$, $N = 322$ also solves the second tiling
equation, and this method does not construct a $322$-tiling.
\end{remark}

\begin{remark} \label{remark:case2}
For every tile with $c \mid a^2$, including the case $c > 2b$ not covered
by Theorem~\ref{theorem:existence2}, some triangle with angles
$(2\alpha,\alpha,2\beta)$ can be tiled, by reducing to the isosceles
triangle with base angles $\alpha+\beta$; see
Corollary~\ref{corollary:existence2case2} below.  In general the
resulting $N$ is much larger than that of Theorem~\ref{theorem:existence2}.
\end{remark}

\begin{proof}[Proof of Theorem~\ref{theorem:existence2}]
Let $(a,b,c)$ be as in the theorem, so all three are integers.
We imitate and generalize the example
of the 1288-tiling shown in Figure~\ref{figure:1288}.
At one point in the construction we will need the hypothesis $c < 2b$.
 Let $\ell = \lcm(a,c)$.  Then $u = \ell /a$ is an integer, and
$j = \ell/c$ is an integer, and we have $ua = jc = \ell$.

 In the example we have
 \begin{eqnarray*}
(a,b,c) &=& (6,5,9)  \\
s &=& 6/9 \ = \ 2/3 \\
N &=& 1288 \\
\ell &=& 18 \\
u &=& 18/6\ = 3\ \\
j &=& 18/9\ =\ 2 \\
b &=& c - a^2/c \ = \ 5 \\
3a &=& 2c \\
18a &=& 12 c \mbox{\qquad at the yellow-green boundary}\\
18b &=& 15 a \mbox{\qquad at the yellow-pink and red-pink boundary} \\
10a &=& 12 b  \mbox{\qquad at the blue-red boundary}
\end{eqnarray*}
We have generalized $3a = 2c$ to $ua = jc$.  The
equation $18a = 12c$ becomes $\ell a = ua^2 = (ja)c$.
 We now want to
generalize $18b = 15a$.   That becomes $\ell b = ha$.
Here $h = \ell b/a$ is defined to make that equation true, and
the point to be proved is that $h$ is an integer.
 To prove that,
\begin{eqnarray*}
h &=& \frac {\ell b} a \\
&=& ub  \mbox{\qquad since $ u = \ell/a$}
\end{eqnarray*}
This is an integer, since $b = c - a^2/c$ and $c$ divides $a^2$ by
hypothesis.  In the example, $h=15$.

Finally we want to generalize $10a = 12b$.  The $12b$ becomes $(ja)b$.
The $10a$ becomes $(jb)a$.
Then, as in
the example, we can construct a tiling, following the pattern
of Figure~\ref{figure:1288}.  Let us count the number of tiles
required.  There are three quadratic tilings, each $\ell$ by $\ell$,
touching points $A$ and $C$ and
contributing $3\ell^2$.  There is another quadratic tiling touching
point $B$ at the top of the figure. That one contributes $(jb)^2$ tiles.
The yellow quadratic tiling in the middle of the figure
contributes $(ja)^2$ tiles.   Finally the red parallelogram consists
of $h-ja$ rows ($15-12 = 3$ in the example) and $ja$ columns (12 in the
example), each place containing two tiles,  so the total number of
tiles in the red parallelogram is $2(h-ja)(ja)$.  Here is where we need
the hypothesis $c < 2b$: the number of rows must be positive, and
$$ h-ja \ = \ \frac{\ell b}{a} - \frac{\ell a}{c}
       \ = \ \frac{\ell(bc-a^2)}{ac}
       \ = \ \frac{\ell(c^2-2a^2)}{ac}
       \ = \ \frac{\ell(2b-c)}{a}, $$
since $bc = c^2-a^2$ and $a^2 = c(c-b)$.

 Adding the numbers calculated for each color,
we compute the number
of tiles in this tiling.  This should come out to be $N$,  but
that is yet to be proved, so for now we use another letter:

\begin{eqnarray}
T &=& 3\ell^2 + (jb)^2 +(ja)^2 + 2(h-ja)(ja) \nonumber \\
  &=& 3 \ell^2 + (\ell b/c)^2 + (\ell a/c)^2 + 2(\ell b/a - ja)(ja)
           \mbox{\qquad since $h = \ell b/a$} \nonumber\\
  &=& 3 \ell^2 + \ell^2 (1-s^2)^2 + \ell^2 s^2 + 2(\ell b/a - \ell a/c)(\ell a/c) \mbox{\qquad since $ja = \ell a/c$} \nonumber\\
  &=& 3 \ell^2 + \ell^2 (1-s^2)^2 + \ell^2 s^2 +
           2\bigg( \frac {\ell (1-s^2)}{s} - \ell s\bigg)(\ell s)\nonumber\\
  &=& \ell^2 \bigg(3 + (1-s^2)^2 + s^2 + 2\big((1-s^2)-s^2\big)\bigg) \nonumber \\
  &=&\ell^2 ((2-s^2)^2 + 2-s^2) \nonumber\\
  &=& \ell^2 (2-s^2)(3-s^2)  \label{eq:3950}
\end{eqnarray}
As a check, in the example we have $\ell = 18$ and $s = 2/3$, yielding
$T = 1288$ as expected.

Reading off the figure, $M$ will be $jb-ja+\ell$, which is
an integer.  We have $M = jb-ja+\ell = \ell b/c - \ell a/c +\ell =
\ell(2-s^2-s)$, as claimed in the theorem.
Putting this value into the second tiling equation we find a value for
the number of tiles $N$:
\begin{eqnarray*}
N &=& M^2 \frac {(2-s^2)(3-s^2)}{(1-s)^2(2+s)^2}  \\
  &=& \ell^2 (2-s^2-s)^2 \frac {(2-s^2)(3-s^2)}{(1-s)^2(2+s)^2} \\
  &=& \ell^2(2+s)^2(1-s)^2\frac {(2-s^2)(3-s^2)}{(1-s)^2(2+s)^2} \\
  &=& \ell^2 (2-s^2)(3-s^2)
\end{eqnarray*}
which is exactly the value obtained for the number of tiles in
(\ref{eq:3950}).  That is, the constructed tiling does indeed
have $N$ tiles.
\end{proof}

\subsection{Paliy's $322$-tiling}
\begin{figure}[ht]
\centering
\includegraphics[width=0.5\textwidth]{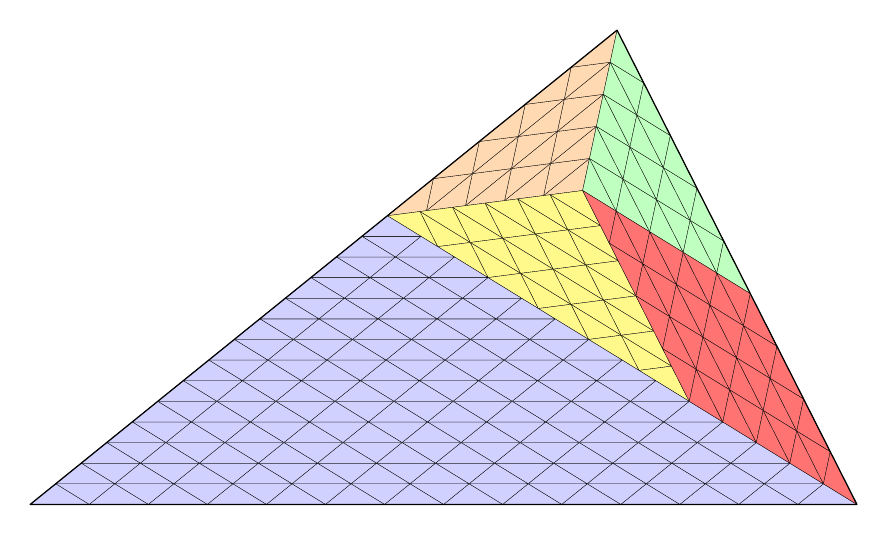}
\caption{$N=322$, $(a,b,c) = (6,5,9)$, $ABC$ of shape $(2\alpha, \alpha, 2\beta)$.}
\label{figure:322Paliy}
\end{figure}

\begin{figure}[ht]
\centering
\includegraphics[width=0.8\textwidth]{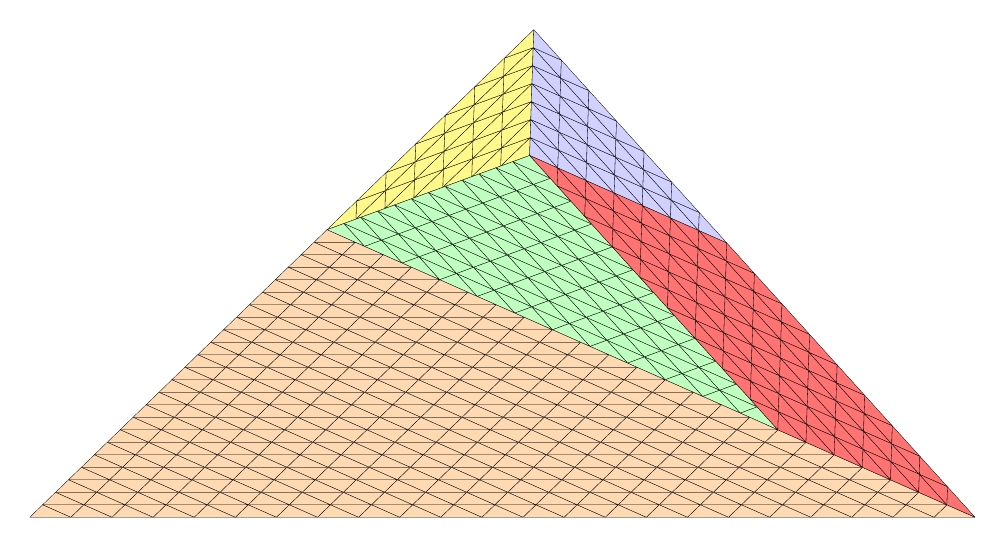}
\caption{$N=897$, $(a,b,c) = (12,7,16)$, $ABC$ of shape $(2\alpha, \alpha, 2\beta)$.}
\label{figure:897}
\end{figure}

On 29 September 2026, D. Paliy found a $322$-tiling and an $897$-tiling (of a triangle of shape $(2\alpha,\alpha,2\beta)$).
These tilings were constructed by a simple pasting operation:  start with a triquadratic tiling, with angles
$(2\alpha, \beta, \alpha + \beta)$.   Paste a quadratic tiling along the side opposite the $2\alpha$ angle, so
that the $\gamma$ angle of the quadratic tiling combines with the $\alpha+\beta$ angle to make $\pi$, a straight line.
The result is a tiling of a $(2\alpha, \alpha, 2\beta )$ triangle.   See Figure~\ref{figure:322Paliy} 
and Figure~\ref{figure:897}.

Paliy's construction leads to the following theorem, which
removes the hypothesis $c < 2b$ from Theorem~\ref{theorem:existence2}.
Not only does it produce tilings for more values of $N$, it generally produces
{\em different} tilings than Theorem~\ref{theorem:existence2}.
\begin{theorem}[Paliy]
\label{theorem:existencePaliy}
Let $3\alpha+2\beta=\pi$. 
Suppose positive integers $a$, $b$, and $c$ with no common factor satisfy
$$ (a,b,c) = (mk, k^2-m^2,k^2).$$
Let a triangle $R$ (the tile) have sides $(a,b,c)$ and angles $(\alpha,\beta,\gamma)$.
Let $j$ be any positive integer and let
$$ N = j^2(b+c)(b+2c).$$
Then some triangle $ABC$ with
angles $(2\alpha,\alpha,2\beta)$  can be $N$-tiled by $R$.
\end{theorem}

\begin{proof}  
Without loss of generality we may suppose $j=1$,
since after constructing a tiling for $j=1$, we can subdivide it quadratically.
Suppose $N$ satisfies the equation in the theorem with $j=1$.
Define $M:= mk$ and $K := k^2$.  Then $K$ divides $M^2$, and 
$k^2 N = k^2(2k^2-m^2)(3k^2-m^2) = (2K^2 - M^2)(3k^2-m^2)$.  
Since $m<k$, we have $m^2k^2 < k^2(2k^2-m^2)$.  That is,
$M^2 < 2K^2 - M^2$, which is the hypothesis of Theorem~\ref{theorem:triquadraticexistence}.
According to 
Theorem~\ref{theorem:triquadraticexistence}, we can construct a triangle $AEB$ with 
angles $(2\alpha, \alpha+\beta, \beta)$ at $(A,E,B)$ respectively, and
tiled by $R$ using $2K^2-M^2$ tiles.  Let $Q:=2k^2-m^2$. 
According to Theorem~\ref{theorem:triquadraticexistence},
the length of $EB$ (the side opposite $2\alpha$) is 
\begin{eqnarray*}
(M/K)(2K^2-M^2)= 2MK - M^3/K = 2mk^3-km^3  = mk(2k^2-m^2) = aQ.
 \end{eqnarray*}

 Let $C$ be the point
such that $E$ is between $A$ and $C$ and $\vert EC \vert = bQ$.  Since angle $AEB$ is
$\alpha+\beta$, angle $BEC$ is its supplement, $\gamma$.  Let $R'$ be the tile
scaled by $Q$, with sides $(aQ, bQ, cQ)$; in $R'$ the sides $aQ$ and $bQ$ include the
angle $\gamma$.  Triangle $BEC$ has $\vert EB \vert = aQ$, $\vert EC \vert = bQ$, and
included angle $\gamma$.  Hence, by SAS, triangle $BEC$ is congruent to $R'$, so it
can be quadratically tiled into $Q^2$ copies of $R$.  In $ABC$ the angle at $A$ is $2\alpha$, the angle at $B$ is
$\beta+\beta = 2\beta$, and the angle at $C$ is $\alpha$, since in $R'$ the angle
opposite $aQ$ is $\alpha$.  Then since $ABC$ is
the union of the triangles $ABE$ and $EBC$, which overlap only at their boundaries,
$ABC$ is tiled by $2K^2-M^2 + Q^2$ tiles.  Counting the tiles we have
\begin{eqnarray*}
 2K^2-M^2 + Q^2 &=& 2k^4 - k^2m^2 + Q^2 \\
 &=&2k^4 - k^2m^2 + Q(2k^2-m^2) \\
 &=& k^2Q + Q^2 \\
 &=& (k^2 + Q)Q \\
&=& (3k^2-m^2)Q  \\
&=&(2k^2-m^2)(3k^2-m^2)\\
&=& (b+c)(b+2c)
 \end{eqnarray*}
But that is $N$, since we have assumed $j=1$.
\end{proof}

\subsection{Another form of the second tiling equation}

\begin{theorem} \label{theorem:secondcasesolution}
Let the tile have angles $(\alpha,\beta,\gamma)$ with $3\alpha+2\beta=\pi$,
and suppose some triangle with angles $(2\alpha,\alpha,2\beta)$ is
$N$-tiled and not similar to the tile.
  Then $s = 2\sin\frac{\alpha}{2}$ is rational.  Scale the tile
so that $(a,b,c)$ are integers with no common divisor; then $s = a/c$,
$c$ divides $a^2$, and
$$ N \ = \  (b+c)(b+2c) \Bigl( \frac{M}{b+c-a} \Bigr)^2.$$
\end{theorem}

\begin{proof}[Proof of Theorem~\ref{theorem:secondcasesolution}]
By Theorem~\ref{theorem:secondtilingequation}, $s$ is rational and the
second tiling equation holds:
\begin{eqnarray*}
\frac  N {M^2}  &=&  \frac{(2-s^2)(3-s^2)}{(1-s)^2(2+s)^2}
\end{eqnarray*}
We have, since $s=a/c$ and $b = c-a^2/c$, 
$(1-s)(2+s) = (b+c-a)/c$ and 
$$(2-s^2)(3-s^2) = (b+c)(b+2c)/c^2.$$
So the equation takes the form 
$$ N \ = \  (b+c)(b+2c) \Bigl( \frac{M}{b+c-a} \Bigr)^2$$
  Since
$b = c - a^2/c$ is an integer, $c$ divides $a^2$. 
\end{proof}

\subsection{$N$ is not prime when $ABC$ has angles $(2\alpha,\alpha,2\beta)$}
\begin{theorem}
\label{theorem:notprime2}  Suppose $3\alpha + 2\beta = \pi$. Let 
triangle $ABC$ with angles $(2\alpha,\alpha,2\beta)$
 be $N$-tiled by a tile with 
angles $(\alpha,\beta,\gamma)$.  Suppose $ABC$ is not similar to the tile.
Then $N$ is not prime.
\end{theorem}

\begin{proof}  Suppose $ABC$ is $N$-tiled as in the
statement of the theorem, and suppose, for proof by contradiction, that $N$ is prime.
Let the tile have sides $(a,b,c)$ and define $s= a/c$.  We may assume that $(a,b,c)$
are integers with no common factor and that $c$ divides $a^2$.  By
Lemma~\ref{lemma:gsquare}, $(a,b,c) = (mk,\,k^2-m^2,\,k^2)$ with $\gcd(m,k)=1$
and $0 < m < k$.  According to Theorem~\ref{theorem:secondcasesolution},
$$
 N \ = \  (b+c)(b+2c) \Bigl( \frac{M}{b+c-a} \Bigr)^2.
$$
Clearing the denominator,
\begin{equation}
 N(b+c-a)^2 = (b+c)(b+2c)M^2.   \label{eq:4050}
 \end{equation}
Here $b+c-a = (k-m)(2k+m) > 0$, so $M \neq 0$.

Every positive integer $n$ can be written uniquely as $n = r u^2$ with $r$ squarefree;
we call $r$ the squarefree part of $n$.  The squarefree part of $nu^2$ is the
squarefree part of $n$, and if $\gcd(x,y) = 1$, the squarefree part of $xy$ is the
product of the squarefree parts of $x$ and $y$.

The factors $b+c = 2k^2-m^2$ and $b+2c = 3k^2-m^2$ are relatively prime: a prime
dividing both divides their difference $c = k^2$, hence divides $k$, hence divides
$m^2 = 2k^2 - (b+c)$, contradicting $\gcd(m,k) = 1$.
Let $p$ and $q$ be the squarefree parts of $b+c$ and $b+2c$ respectively.
Since $N$ is prime, the squarefree part of the left side of \eqref{eq:4050} is $N$,
and the squarefree part of the right side is $pq$.  Hence $N = pq$.

Now $b+2c$ is not a square, so $q > 1$.  Indeed, suppose $b+2c = 3k^2-m^2 = z^2$.
If $k$ is even and $m$ is odd, then $z^2 \equiv 3 \pmod 4$, which is impossible.
If $k$ is odd and $m$ is even, then $z^2 \equiv 3 \pmod 4$, which is impossible.
If $k$ and $m$ are both odd, then $z^2 \equiv 2 \pmod 4$, which is impossible.
Since $\gcd(m,k) = 1$, they cannot both be even.

Since $N = pq$ is prime and $q > 1$, we have $p = 1$ and $N = q$.  Then $b+c$ is a
square, $2k^2 - m^2 = x^2$.  Then $k$ and $m$ are both odd, as otherwise $x^2$ would
be $2$ or $3$ modulo $4$.  Then $b+2c = 3k^2-m^2 \equiv 2 \pmod 4$, so its squarefree
part $q$ is even.  Since $N = q$ is prime, $N = 2$.

But the angles of $ABC$ are $2\alpha$, $\alpha$, and $2\beta$, so at least three
$\alpha$ angles of tiles occur at the vertices of $ABC$; since each tile has only
one $\alpha$ angle, $N \ge 3$.  This contradiction shows that $N$ is not prime.
\end{proof}
 
\subsection{Complete solution of the tiling problem for $(2\alpha, \alpha, 2\beta)$}

In 2026, D. Paliy observed that the fraction $M/(b+c-a)$ in 
Theorem~\ref{theorem:secondcasesolution} is an integer.  This is 
a key observation, leading to a complete characterization of the possible
tilings for the shape  $(2\alpha, \alpha, 2\beta)$.

\begin{theorem} [Paliy, \cite{paliy2026}] \label{theorem:Paliy} 
Suppose $3\alpha + 2\beta = \pi$,
and let $(a,b,c)$ be positive integers with no common factor such that 
the triangle with sides $(a,b,c)$ has angles $(\alpha,\beta,\gamma)$ opposite
$(a,b,c)$ respectively.
Then the following are equivalent:

(i) $N = j^2(b+c)(b+2c)$ for some positive integer $j$.

(ii) There is an $N$-tiling of some triangle $ABC$ with angles $(2\alpha, \alpha, 2\beta)$
by the tile $(a,b,c)$.

Moreover,  when the equation
$$ N \ = \  (b+c)(b+2c) \Bigl( \frac{M}{b+c-a} \Bigr)^2$$
has a solution, $M/(b+c-a)$ is an integer;  and 
when such an $N$-tiling exists, it has coloring number $M = j(b+c-a)$. 
\end{theorem}

\begin{proof}  (i) implies (ii) is Theorem~\ref{theorem:existencePaliy}.
It remains to show (ii) implies (i).  Suppose there is an $N$-tiling of $ABC$
by tile $(a,b,c)$.
The triangle $ABC$ is not similar to the tile: otherwise $2\alpha$, $\alpha$, $2\beta$
would be a permutation of $\alpha$, $\beta$, $\gamma = 2\alpha+\beta$, and since
$3\alpha + 2\beta = \pi$ the only possibility is $\beta = 2\alpha$, $\alpha = \pi/7$;
but then $a/c = 2\sin(\pi/14)$ is irrational (it is a root of $x^3 - x^2 - 2x + 1$,
which has no rational root), contradicting that $a$ and $c$ are integers.
According to Theorem~\ref{theorem:secondcasesolution},
$$ N \ = \  (b+c)(b+2c) \Bigl( \frac{M}{b+c-a} \Bigr)^2.$$

We will prove that $j = M /(b+c-a)$ is an integer, just on the basis that 
that equation is satisfied (without needing a tiling further).
By Lemma~\ref{lemma:gsquare}, there are positive coprime integers $m$ and $k$
such that 
$$ (a,b,c) = (mk, k^2-m^2, k^2).$$
Let $p$ be an odd prime, and suppose that $p$ divides both $b+c$ and $b+c-a$.
 Then $a \equiv 0 \pmod p$, so $p \mid m$ or $p \mid k$.
 If $p \mid m$ then $b+c = 2k^2-m^2 \equiv 2k^2 \pmod p$.
 Then since $p \mid b+c$ and $p$ is odd, $p \mid k$, contradicting
 $\gcd(m,k) = 1$.  On the other hand, if $p \mid k$ then $b +c = 2k^2-m^2 \equiv -m^2 \pmod p$.
 Since $p \mid b+c$, $m \equiv 0 \pmod p$, contradicting $\gcd(m,k) = 1$.  
 Therefore $p$ cannot divide both $b+c$ and $b+c-a$.  
 
 Next we claim
 $p$ cannot divide both $b+c-a$ and $b+2c = 3k^2-m^2$. If $p$ divides
 both, then it divides their difference $c+a = k(k+m)$.
 Then $p \mid k$ or $p \mid k+m$.  If $p \mid k$ then $b \equiv -m^2 \pmod p$ and $c \equiv 0 \pmod p$,
 so $b+c-a \equiv -m^2 \pmod p$;  since $p \mid b+c-a$, we have $p \mid m$, contradicting $\gcd(m,k)= 1$.
If $p \mid k+m$ then $m \equiv -k \pmod p$, so $b+c-a = 2k^2-m^2-mk \equiv 2k^2 \pmod p$; since $p$ is odd, $p \mid k$, and then $p \mid m$, contradicting $\gcd(m,k)=1$.  Therefore, as claimed, $p$
cannot divide both $b+c-a$ and $b+2c = 3k^2-m^2$.

Therefore, if $p$ is an odd prime and $p^e$ divides $b+c-a$, then $p^e$ divides $M$,
 since $N(b+c-a)^2 = (b+c)(b+2c)M^2$ and $p \nmid (b+c)(b+2c)$. 
 
 Therefore, if $j$ is not an integer, it is because the denominator $b+c-a$ is even; more explicitly,
 the prime 2 divides the denominator  more times than it divides the numerator $M$. That is,
 $t = v_2(b+c-a) - v_2(M) > 0$, where $v_2$ is the exponent of $2$.
 Then because $M/(b+c-a)$ is squared, $2^{2t}$ must divide $(b+c)(b+2c)$.  We will show that is not possible,
 as $(b+c)(b+2c)$ is not divisible by 4.
 
Since $\gcd(m,k)=1$, $m$ and $k$ are not both even.  If $m$ is even, then $k$ is odd,
so $b+c = 2k^2-m^2 \equiv 2 \pmod 4$ and $b+2c = 3k^2-m^2$ is odd.  If $m$ is odd and
$k$ is even, then $b+c$ and $b+2c$ are both odd.  If $m$ and $k$ are both odd, then
$k^2 \equiv m^2 \equiv 1 \pmod 8$, so $b+c = 2k^2-m^2$ is odd and
$b+2c = 3k^2-m^2 \equiv 2 \pmod 8$.  In every case $(b+c)(b+2c)$ is not divisible by $4$.
\end{proof}

\begin{corollary} \label{corollary:everysolution}  Every solution $(N,M,s)$ of the second tiling 
equation, with $N$ and $M$ positive integers and $s$ rational with $0<s<1$,
 corresponds to an $N$-tiling of some triangle of the shape $(2\alpha, \alpha,2\beta)$ with
 coloring number $M$.
\end{corollary}

\begin{proof}  Suppose the equation is satisfied:
$$ \frac N { M^2} \ = \ \frac {(2-s^2)(3-s^2)}{(1-s)^2(2+s)^2}, $$
where $0 < s < 1$.  Let $\alpha$ be such that $2 \sin(\tfrac \alpha 2) = s$.
Let the tile have angles $(\alpha,\beta,\gamma)$ and sides $(a,b,c)$,
integers with no common factor.  
We have, since $s=a/c$ and $b = c-a^2/c$, 
$(1-s)(2+s) = (b+c-a)/c$ and 
$$(2-s^2)(3-s^2) = (b+c)(b+2c)/c^2.$$
Therefore
$$ N \ = \  (b+c)(b+2c) \Bigl( \frac{M}{b+c-a} \Bigr)^2,$$
now proved without assuming there is a tiling.
Then by Theorem~\ref{theorem:Paliy}, $M/(b+c-a)$ is an integer.
Therefore Part~(i) of the theorem holds.  Then since (i) implies (ii),
 there is an $N$-tiling as desired.
\end{proof}

\begin{remark}  Theorem~\ref{theorem:Paliy} makes the non-primality of $N$ immediate, 
since the factors $(b+c)$ and $(b+2c)$ are integers greater than 1.  The proof
of the theorem, however, is somewhat longer than the proof we gave for Theorem~\ref{theorem:notprime2}.
\end{remark}

\subsection{Solving the second tiling equation}

\begin{theorem} \label{theorem:solvingsecond} There is an 
algorithm that proceeds from input $N$ to termination, and 
at termination either produces a solution $(M,s)$ of the second tiling 
equation with $s$ rational, or reports that no such solution exists.
\end{theorem}

\begin{proof}  Recall the second tiling equation
$$ \frac N { M^2} \ = \ \frac {(2-s^2)(3-s^2)}{(1-s)^2(2+s)^2}. $$
By Lemma~\ref{lemma:Msquared}, any solution has $0 < M^2 < N$, so there
are only finitely many $M$ to consider.  For each such $M$, the equation is
equivalent to $p(s) = 0$, where
\begin{eqnarray*}
p(s) &=& N(1-s)^2(2+s)^2 - M^2  (2-s^2)(3-s^2)\\
&=& (N-M^2)s^4 + 2Ns^3 + (5M^2-3N)s^2 -4Ns +(4N- 6M^2),
\end{eqnarray*}
a polynomial with integer coefficients and nonzero leading coefficient
$N-M^2$.  By the rational root theorem, any rational root $s = e/f$ in
lowest terms has $f \mid N - M^2$ and $e \mid 4N - 6M^2$, and we need only
consider $0 < e < f$.  These are finitely many candidates, each of which
can be tested exactly.
\end{proof}

 Table~\ref{table:SageTest} lists all solutions of the second tiling
equation with $N \le 2000$, obtained by carrying out that algorithm.
The table also shows the smallest integer-sided tile $(a,b,c)$
with $a/c = s$ and $b = c-a^2/c$.

\begin{table}[ht]
\caption{All solutions of the second tiling equation for $N\le 2000$.}
\label{table:SageTest}
\begin{center}
\begin{tabular}{rrrr}
$N$  &       $M$   &     $s$   &  $(a,b,c)$ \\
  \hline
77  &      5   &     1/2 &  (2,3,4) \\
308   &   10    &    1/2 &  (2,3,4)\\
322  &     8    &    2/3  & (6,5,9)\\
442  &    14    &    1/3  & (3,8,9) \\
693  &   15    &    1/2 &  (2,3,4)\\
897  &    11   &     3/4 &  (12,7,16)\\
1232  &   20  &      1/2  & (2,3,4)\\
1288   &  16    &    2/3  & (6,5,9) \\
1457   &  27    &    1/4  & (4,15,16)\\
1768   &  28    &    1/3  & (3,8,9)\\
1925   &  25    &    1/2  & (2,3,4)
\end{tabular}
\end{center}
\end{table}

If the second tiling equation is solvable for $N$, then it is
solvable for every square multiple of $N$, since if $N$ is
multiplied by $\lambda^2$ and $M$ by $\lambda$, the equation is
still satisfied (with the same $s$).  Thus some of the
entries in Table~\ref{table:SageTest} were predictable, e.g.,
$308 = 4\cdot 77$ and $693 = 9\cdot 77$.  In fact, only the first
occurrence of each value of $s$ in the table has a square-free $N$.

In Table~\ref{table:SageTest}, for each $N$ there is either no solution,
or exactly one.  That is an accident of the small data set: $N = 38962$
has the two solutions $(M,s) = (88, 2/3)$ and $(152, 1/9)$, the first
from $38962 = 11^2\cdot 322$ and the second from the tile $(9,80,81)$,
for which $(b+c)(b+2c) = 161\cdot 242 = 38962$.  This is the smallest such
$N$; up to $10^6$ the only others are its square multiples.

Since Corollary~\ref{corollary:everysolution} shows that solutions of the equation
always correspond to tilings, there is no need to mark lines of the table that
correspond to tilings: they all do.

 \section{The case when $ABC$ is isosceles with base angles $\beta$}

\begin{figure}
\centering
\includegraphics[width=0.8\textwidth]{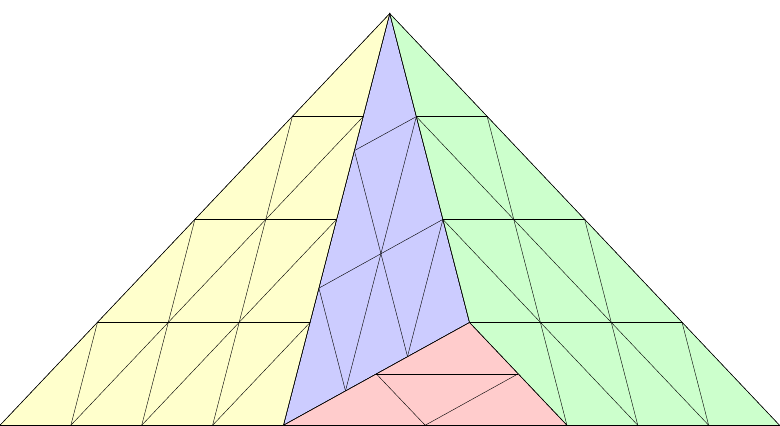}
\caption{$N=44$, $(a,b,c) = (2,3,4)$, $ABC$ isosceles with base angles $\beta$.}
\label{figure:44}
\end{figure}

\begin{figure}
\includegraphics[width=0.8\textwidth]{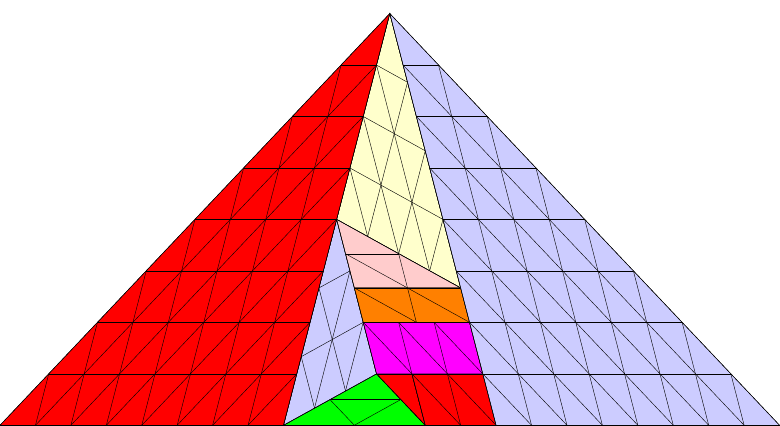}
\caption{$N=176$, $(a,b,c) = (2,3,4)$, $ABC$ isosceles with base angles $\beta$.}
\label{figure:176}
\end{figure}

In Figs.~\ref{figure:44} and \ref{figure:176}, we give examples
of tilings of triangles of this form.

Let $ABC$ have angles $\beta$ at $A$ and $C$. 
Then the vertex angle at $B$ is $3\alpha$, since $3\alpha +2\beta = \pi$.
Let $X$ be the length of $AB$ and $BC$; let $Y$ be the length of $AC$.
Let $\mu$ be the ``scale factor'', defined by  $X = \mu b$. 
As for other shapes of $ABC$,  we have the ``area equation'' 
expressing that the area of $ABC$ is equal to the area of $N$ tiles, 
and the ``coloring equation'' from Theorem~\ref{theorem:coloring}, 
namely 
\begin{eqnarray*}
M(a+b+c) &=& 2X+Y 
\end{eqnarray*}
Combining these, as we did for other shapes of $ABC$, will result 
in a ``tiling equation'' giving a necessary condition for 
the existence of an $N$-tiling of $ABC$.  Then the question 
is whether the solvability of that equation is also a sufficient
condition.  For the shapes considered so far it was not, at least not
for the smallest $N$ satisfying the equation.  Below we show that in
this case it is sufficient, provided $a+c$ divides $M$.

\subsection{The tiling equation for $ABC$ isosceles with base angles $\beta$}
In this section we show how, given $N$, to determine a finite
set of possible tiles $(a,b,c)$ 
(each of which then determines the triangle $ABC$)  such that,  if any triangle 
$ABC$ of the form considered here can be $N$-tiled, then one of these
specific triangles
$ABC$ can be $N$-tiled with the corresponding specific tile.  That is implied
by the ``tiling equation'' given in the following theorem.
\FloatBarrier

\begin{theorem}[Tiling equation] \label{theorem:tilingequation3}
Suppose $3\alpha + 2\beta = \pi$, 
and triangle $ABC$ is isosceles with base angles $\beta$,
and $ABC$ is $N$-tiled by a tile with angles $(\alpha,\beta,\gamma)$
and sides $(a,b,c)$ of integer lengths with no common factor.
Suppose $ABC$ is not similar to the tile.
Let $M$ be the coloring number of the tiling and $s = a/c$.
  Then  $0 < N/3 < M^2 < 2N$, and 
\begin{eqnarray*}
\frac N  {M^2} &=&  \frac {(3-s^2)}{(1+s)^2},
\end{eqnarray*}
and consequently 
$$s = \frac {\sqrt{3M^2 + 2N}M - N}{M^2 + N}. $$
\end{theorem}

\begin{remark}  Then $3M^2 + 2N$ is a square, since 
$$ \sqrt{3M^2+2N} = \frac{s(M^2+N)+N} M,$$
which is rational, so the integer $3M^2 + 2N$ is a square.%
\footnote{
The condition  $ \gcd(a,c) \mid M$ was claimed in earlier versions 
of this paper, but it is not correct.}
\end{remark}

\begin{proof}  Let $\kappa$ be defined by 
$b = \kappa \sin \beta$.  Let $\mu$ be defined
by $X = \mu b$.
We have by the law of sines, 
\begin{eqnarray*}
 \frac Y X &=& \frac {\sin 3\alpha} {\sin \beta}
 \end{eqnarray*}
By Lemma~\ref{lemma:singamma2} we have 
$$ \frac {\kappa \sin 3\alpha} c  \ = \ 
s(1+s)(1-s)(3-s^2)$$
Since $b = \kappa \sin \beta$  and $b/c = 1-s^2$, we have 
\begin{eqnarray*}
\frac Y X &=& \frac { s(1+s)(1-s)(3-s^2)} {1-s^2}
\end{eqnarray*}
Simplifying,
\begin{eqnarray}
 \frac Y X  &=& s(3-s^2) \label{eq:Y/X}
\end{eqnarray}
By Theorem~\ref{theorem:coloring}, the coloring equation is 
\begin{eqnarray*}
M(a+b+c) &=& 2X+Y  \\  
&=& X(2+Y/X) \\
&=& \mu b (2+s(3-s^2))  
\end{eqnarray*}
Dividing by $c$ and expressing everything in terms of $s$, 
\begin{eqnarray*}
M(s+(1-s^2) + 1) &=& \mu (1-s^2)(2+3s-s^3) \\
\frac{\mu} M &=& \frac {2+s-s^2}{(1-s^2)(2+3s-s^3)} \\
   &=& \frac {(1+s)(2-s)}{(1-s^2)(2-s)(1+s)^2}\\
   &=& \frac {1} {(1-s)(1+s)^2} \\
   &=& \frac 1 {(1-s^2)(1+s)}
\end{eqnarray*}
Multiplying by $M$ and squaring, we have 
\begin{eqnarray}
\mu^2 &=& \frac {M^2 } {(1-s^2)^2(1+s)^2} \label{eq:7442}
\end{eqnarray}

The area equation is 
\begin{eqnarray*}
X^2 \sin 3\alpha &=& Nbc \sin \alpha 
\end{eqnarray*}
Multiplying by $\kappa/c$ and using Lemma~\ref{lemma:singamma2} again,
as well as $\kappa \sin \alpha = a$, we have 
\begin{eqnarray}
X^2 s(1+s)(1-s)(3-s^2) &=& Nbc \frac a c  \ = \ Nab \label{eq:7505}
\end{eqnarray}
Putting $X = \mu b$ we have 
\begin{eqnarray*}
\mu^2 b^2 s(1+s)(1-s)(3-s^2) &=& Nab  \\
\mu^2 b  s(1+s)(1-s)(3-s^2) &=& Na
\end{eqnarray*}
Expressing everything in terms of $s$, and dividing by $cs$, 
\begin{eqnarray*}
N &=& \mu^2 (1-s^2)^2(3-s^2)
\end{eqnarray*}
Substituting the value of $\mu^2$ from (\ref{eq:7442}),
\begin{eqnarray}
N &=& \frac {M^2 (1-s^2)^2(3-s^2)} {(1-s^2)^2(1+s)^2} \nonumber \\
N &=&  M^2 \frac {(3-s^2)}{(1+s)^2} \label{eq:9073}
\end{eqnarray}
This is the ``tiling equation'' of the theorem; the 
formula for $s$ follows by solving this equation, which 
is quadratic in $s$.  ( The other solution is negative and
is therefore irrelevant.) 

It remains to prove the inequality $0 < N/3 < M^2 < 2N$.
Since $3\alpha + 2\beta = \pi$ we have $\gamma = (\pi+\alpha)/2 > \pi/2$, so $\gamma$
is the largest angle of the tile and $c$ is its longest side; hence $0 < s < 1$.
As $s$ varies over $(0,1)$, the right side of the tiling equation decreases
strictly from $3$ to $1/2$, since
$$\frac{d}{ds}\,\frac{3-s^2}{(1+s)^2} \ = \ \frac{-2(s+3)}{(1+s)^3} \ < \ 0 .$$
Hence $1/2 < N/M^2 < 3$, which is the inequality of the theorem.%
\footnote{One can make a graph by the SageMath command
\texttt{plot((3-s\textasciicircum 2)/(1+s)\textasciicircum 2,0,1)}.}
Hence the inequality  
is a consequence of the tiling equation.
\end{proof}

\begin{theorem}[Tiling equation, second form] \label{theorem:notprime4}
Suppose $3\alpha + 2\beta = \pi$, 
and triangle $ABC$ is isosceles with base angles $\beta$,
and $ABC$ is $N$-tiled by a tile with angles $(\alpha,\beta,\gamma)$,
coloring number $M$, 
and sides $(a,b,c)$ of integer lengths with no common factor. 
Suppose $ABC$ is not similar to the tile.
 Then  $m+k \mid M$, and for some integer $j$,  namely $j = M/(m + k) > 0$, 
$$N = j^2(3k^2-m^2) = j^2\left(3c - \frac {a^2}{c}\right),$$
and if $N$ is squarefree, $j=1$.
\end{theorem}

\begin{remark} It does not follow that $N$ cannot be prime,
since we can have $j=1$ and $3k^2-m^2$ can be prime, e.g., with $(m,k) = (1,2)$.
But if we could exclude tilings with $j=1$, it would follow that $N$ cannot be prime,
or even squarefree. 
\end{remark}

\begin{proof}[Proof of Theorem~\ref{theorem:notprime4}] 
Let $M$ be the coloring number of the tiling,
and let $(a,b,c)$ be the sides of the tile, scaled so that they
are integers with no common factor.  That is possible, by 
Theorem~\ref{theorem:rationaltile}.  
According to Theorem~\ref{theorem:tilingequation3}, we have  
\begin{eqnarray*}
\frac N  {M^2} &=&  \frac {(3-s^2)}{(1+s)^2}
\end{eqnarray*}
Using $s = a/c$ and 
clearing denominators, we have 
\begin{eqnarray}
N(a+c)^2 &=& M^2(3c^2-a^2) \label{eq:9152}
\end{eqnarray}
Since $M$ is the excess of the number of black tiles over the number of white tiles,
we have $\vert M \vert < N$, and hence $N \nmid M$.  
Let $k = \gcd(a,c)$.  
By Lemma~\ref{lemma:gsquare}, $c=k^2$. Let $m = a/k$.  Then
$$ (a,b,c) = (km, k^2-m^2, k^2) \mbox{\ with $s = a/c = m/k$ and $\gcd(m,k) = 1$}$$
and dividing \eqref{eq:9152} by $k^2$, we find
\begin{eqnarray}
N (m+k)^2 = M^2 (3k^2-m^2) \label{eq:8021}
\end{eqnarray}
We intend to prove that $m+k \mid M$.  We will do that by proving that for every
prime $p$, the order to which $p$ divides $M$ is at least the order to which $p$ divides $m+k$.

First we consider the case $p=2$.  Suppose $2^e$ is the exact power of
$2$ dividing $m+k$, with $e \ge 1$.  Then $m$ and $k$ are both odd, since
they are not both even; so $m^2 \equiv k^2 \equiv 1 \pmod 8$ and
$3k^2-m^2 \equiv 2 \pmod 8$.  Thus exactly one power of $2$ divides
$3k^2-m^2$.  Comparing the powers of $2$ on the two sides of
(\ref{eq:8021}), we have $v_2(N) + 2e = 2v_2(M) + 1$, so
$2v_2(M) \ge 2e$, that is, $2^e \mid M$, as desired.

Next,  we show that $(m+k)$ and $(3k^2-m^2)$
have no odd prime factor $p$ in common.   Suppose $p \mid m+k$
and $p \mid 3k^2-m^2$.  Then, assuming $p \neq 2$, we have  
\begin{eqnarray*}
m \equiv - k \pmod p && \mbox{\qquad since $p \mid m+k$}\\
m^2 \equiv k^2 \pmod p &&\\
2k^2 \equiv m^2 \equiv k^2 \pmod p && \mbox{\qquad since $p \mid 3k^2-m^2$}\\
k^2 \equiv 0 \pmod p  && \mbox{\qquad since $p \neq 2$ }\\
m^2 \equiv 0  \mod p  && \mbox{\qquad since $m^2 \equiv k^2 \pmod p$} \\
k \equiv m\equiv 0 \pmod p.  && 
\end{eqnarray*}
But that contradicts the fact that $k$ and $m$ are coprime,
completing the proof that no odd prime divides both $m+k$ and $3k^2-m^2$.
From \eqref{eq:8021}, it then follows that if $p^e \mid (m+k)$ then $p^e \mid M$.
That completes the proof that $m+k \mid M$.

 Then we set
\begin{equation} \label{eq:jdef}
j := \frac M {m+k} \in \Z.
\end{equation}
Then \eqref{eq:8021} can be written
\begin{equation}
N = j^2 (3k^2-m^2)
\end{equation}
If  $N$ is squarefree,  we can then conclude that $j=1$, $M = m+k$, and 
$$N = 3k^2-m^2 = 3c - \frac {a^2}{c}.$$
\end{proof}

\subsection{Construction of tilings when $ABC$ is isosceles with base angles $\beta$}

Let $ABC$ be isosceles with base angles $\beta$.  Let the 
vertex angle at $B$ be trisected into three $\alpha$ angles,
and let the trisecting lines meet the base $AC$ at $D$ and $E$,
in the order $ADEC$.  Then triangle $BDE$ is isosceles and 
triangle $BDC$ has angles $(2\alpha, \alpha + \beta, \beta)$
(at $B$, $D$, and $C$ respectively).  Therefore there are two 
obvious approaches to constructing tilings:
\medskip

(i) Tile the triangle $BDC$ with a triquadratic tiling, and tack 
on a quadratic tiling of $ABD$.  This is illustrated in Figure~\ref{figure:44}.
\smallskip

(ii) Tile the isosceles triangle $BDE$, which has base angles 
$\alpha + \beta$, and tack on two quadratic tilings.  This is 
illustrated in Figure~\ref{figure:176}.
\smallskip

We were able to find a sufficient condition 
for the existence of tilings of isosceles $ABC$ with base angles $\beta$
via method (i).  We  attempt to tile $ABC$ by 
constructing a triquadratic tiling of $BDC$.  

\begin{theorem} \label{theorem:isoscelesbetaexistence}
Suppose the tiling equation in Theorem~\ref{theorem:tilingequation3} 
is solvable by $(M,s)$, with $0 < s < 1$.
Let $(a,b,c)$ be integers with no common 
factor such that $s = a/c$ and $c$ divides $a^2$ (so the angles 
of $(a,b,c)$ satisfy $3\alpha + 2\beta = \pi$), 
and let $ABC$ be an isosceles triangle 
with base angles $\beta$ and scaled so that its 
area is equal to $N$ times the area of 
the tile $(a,b,c)$.  Suppose in addition that $a+c$ divides $M$.
 Then there exists an $N$-tiling of $ABC$ by $(a,b,c)$.
\end{theorem}

\begin{proof}
What fraction of the area of $ABC$ is the area of $ABD$?  
Let $X$ be the length of $AB$ and $BC$, and let $Z$ be the 
length of $BD$ and $BE$.  Then $X/Z = c/b$.  Twice the area of $ABD$ 
is $XZ \sin \alpha$; that is also twice the area of $CBE$.  The
area of $DBE$ is $Z^2 \sin \alpha$.  So the ratio of the area of $ABD$ to 
the area of $ABC$ is 
\begin{eqnarray*}
 \frac {XZ \sin \alpha} { 2XZ \sin \alpha + Z^2 \sin \alpha} &=&
  \frac X {2X+Z} \\
  &=& \frac 1 { 2 + (b/c)} \\
  &=& \frac 1 {3-s^2}
\end{eqnarray*}
So the ratio of $DBC$ to $ABC$ is 
$$ 1-  \frac 1 {3-s^2} = \frac {2-s^2}{3-s^2} $$
Now suppose there is an $N$-tiling of $ABC$.  Then 
the number of tiles needed to tile $DBC$ would be 
\begin{eqnarray*}
N^\prime &:=&  N\frac {2-s^2}{3-s^2}
\end{eqnarray*}
By the tiling equation in Theorem~\ref{theorem:tilingequation3},
\begin{eqnarray*}
N &=& M^2 \frac {(3-s^2)}{(1+s)^2}
\end{eqnarray*}
Substituting this in the previous equation we have 
\begin{eqnarray*}
N^\prime &=& M^2 \bigg(\frac {(3-s^2)}{(1+s)^2}\bigg)\bigg( \frac {2-s^2}{3-s^2}\bigg)\\
 &=& M^2 \frac {2-s^2}{(1+s)^2} 
\end{eqnarray*}
Then 
\begin{eqnarray}
 N^\prime &=& 2 \bigg( \frac M {1+s} \bigg)^2 - \bigg(\frac {Ms} {1+s}\bigg)^2
 \label{eq:7578}
\end{eqnarray} 
Let $k = \gcd(a,c)$ and $m = a/k$.  Then $s = m/k$, and
Lemma~\ref{lemma:gsquare},
$$ (a,b,c) = (km, k^2-m^2,k^2).$$
We claim that the fractions $M/(1+s)$ and $Ms/(1+s)$ are 
integers.  Since $s = a/c = m/k$, that is the same as claiming 
that $Mk/(m+k)$ and $Mm/(m+k)$  are integers. 
It suffices to show $m+k$ divides $M$.  That is part of
Theorem~\ref{theorem:notprime4}.
That completes the proof that $M/(1+s)$ and $Ms/(1+s)$ are
integers.  Define
 \begin{eqnarray*}
 K &:=& \frac M {1+s} \\
 J &:=& \frac  {Ms} {1+s}
 \end{eqnarray*}
 Then $K$ and $J$ are integers, and by (\ref{eq:7578}) we have 
 \begin{eqnarray*}
 N^\prime &=& 2K^2 - J^2
 \end{eqnarray*}
 Hence $N^\prime$ is an integer, and the triquadratic tiling
 equation is satisfied.  
 
 In order to show that a triquadratic tiling exists, we must 
 also show that $K$ divides $N^\prime$.  That is not true 
 without an additional hypothesis, as the example $N=11$, $M=3$,
 $K = 2$, $J=1$,  $N^\prime = 7$ shows.   Our additional 
 hypothesis is that $a+c \mid M$.

We now show that $K$ divides $N^\prime$.
  Since $N^\prime = 2K^2-J^2$, it suffices to show 
$K$ divides $J^2$.  Consider the quotient
\begin{eqnarray*}
\frac {J^2} K &=& \frac {(Ms/(1+s))^2}{M/(1+s)} \\
&=& \frac{ Ms^2 }{1+s} \\
&=& \bigg(\frac{ M } {a+c}\bigg)\bigg( \frac {a^2} c \bigg)
\end{eqnarray*}
Now $a^2/c$ is an integer, as remarked above, and $M/(a+c)$ 
is an integer, by hypothesis. 
Hence a triquadratic 
tiling of $BCD$ by $(a,b,c)$ exists.

  It remains to show that $ABD$ can be 
quadratically tiled by $(a,b,c)$. 
The number of tiles needed to tile $ABD$ is $N-N^\prime$.
We claim that this is a square.  We have
\begin{eqnarray*}
N-N^\prime &=& M^2 \frac {(3-s^2)}{(1+s)^2} - M^2 \frac {2-s^2}{(1+s)^2} \\
   &=& \bigg( \frac M {1+s} \bigg)^2 (3-s^2 - (2-s^2)) \\
   &=& \bigg( \frac M {1+s} \bigg)^2
\end{eqnarray*}
and since we proved above that $M/(1+s)$ is an integer, 
this is a square, as claimed.  Since $ABD$ is similar to the tile,
and its area is equal to that of a square number of tiles, when 
we start a quadratic tiling at $B$, the last row of tiles that fits 
into $ABD$ will fit exactly, completing a quadratic tiling.  
\end{proof}
\medskip

\subsection{Computational results for the isosceles, base angles $\beta$ case}
\begin{figure}
\centering
\includegraphics[width=0.5\textwidth]{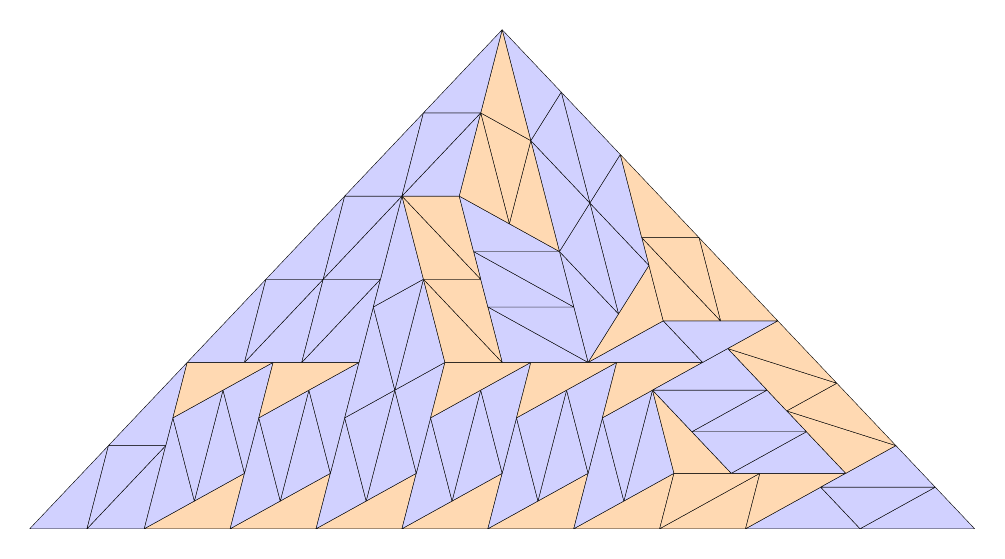}
\caption{$N=99$. $ABC$ is isosceles. The tile is $(2,3,4)$.}
\label{figure:99}
\end{figure}

Let us compute the values of $N$ up to 200 for which there is an $N$-tiling of 
an isosceles triangle with base angles $\beta$, by a tile with $3\alpha + 2\beta = \pi$
and $\alpha$ not a rational multiple of $\pi$.
A simple program computes the values  of $N$ and the possible tile $(a,b,c)$,
by checking for each $N$ in the specified range whether the conditions
in Theorem~\ref{theorem:notprime4} are satisfied.
There are $29$ such values, occupying $30$ rows if printed, since  $N = 143$ arises from two
different tiles. 
The hypothesis $a+c \mid M$ of
Theorem~\ref{theorem:isoscelesbetaexistence} is satisfied exactly by
 $N = 44$ and $N = 176$, yielding the tilings shown in Figure~\ref{figure:44} and
 Figure~\ref{figure:176}.  $N=7$ was ruled out in \cite{beeson-noseven}.
   Values of $N$ with $j = \gcd(N,M) = 1$ are
ruled out by \cite[Theorem~37]{beeson-rosu2026}. That includes all squarefree $N$;
in particular all prime $N$.   That eliminates 24 of the 29 values not
addressed by Theorem~\ref{theorem:isoscelesbetaexistence}, and leaves 
five values undecided.  Bonfioli \cite{bonfioli2026} found a $99$-tiling
by computer search. It is shown in Figure~\ref{figure:99}. 

That leaves four values of $N \le 200$.  Bonfioli ruled these four out by 
exhaustive search for a tiling.  We do not have to trust his searching 
program, because it produces a complete record of the search tree that can 
be checked by Lean, resulting in a Lean-certified proof that no tiling for 
that $N$ exists.  That theorem can be viewed at \cite{bonfioli2026};
in fact he proved there is no $N$-tiling of {\em any} triangle for $N=92, 156$, and $188$,
and no $N$-tiling with tile $(3,8,9)$ and $N = 104$.%
\footnote{
The statements of Bonfioli's theorems looks like this:
{\tt theorem not\_tile\_92 (D : CongruentDissection 92) : False}. 
  Here
{\tt CongruentDissection} is the formalization of tiling that we used in 
our own Lean formalizations.  Now the issue of trusting this result comes 
down to whether in fact Lean does check this theorem.  We downloaded
Bonfioli's repository and ran Lean on it ourselves,  which takes several
hours on a MacBook Air. It contained the proofs for $N = 92, 188$, and $104$,
but for $156$, it has instead a Rust program for generating the proof.  The 
proof itself will be about 130mb, about four times the size of the rest of the
repository, and Bonfolio chose not to store it.  So we did not check 156 
ourselves.   
}

Table~\ref{table:isoscelesbeta} displays the complete solution of 
the tiling problem for this shape and $N \le 200$.

\begin{table}[ht]
\caption{All solutions of the isosceles-$\beta$ tiling problem for $N\le 200$}
\label{table:isoscelesbeta}
\begin{center}
\begin{tabular}{rrrrrl}
$N$ & $M$ & $side$ & $base$ &$(a,b,c)$ & existence \\
\hline 
$44 = 2^2 \cdot 11$ & 6 & 16 & 22 & (2, 3, 4) & \yesmark Figure~\ref{figure:44}, Theorem~\ref{theorem:isoscelesbetaexistence} \\
$92 = 2^2 \cdot 23$ & 10 & 54 & 92 & (6, 5, 9) & \nomark \cite{bonfioli2026} \\ 
$99 = 3^2 \cdot 11$ & 9 & 24 & 33 & (2, 3, 4) & \yesmark Figure~\ref{figure:99}, \cite{bonfioli2026} \\
$104 = 2^2 \cdot 26$ & 8 & 54 & 52 & (3, 8, 9) & \nomark \cite{bonfioli2026} \\
$156  = 2^2 \cdot 39$ & 14 & 128 & 234 & (12, 7, 16) & \nomark \cite{bonfioli2026}\\
$176= 4^2 \cdot 11$ & 12 & 32 & 44 & (2, 3, 4) & \yesmark Figure~\ref{figure:176}, Theorem~\ref{theorem:isoscelesbetaexistence} \\
$188 = 4^2 \cdot 47$ & 10 & 128 & 94 & (4, 15, 16) & \nomark \cite{bonfioli2026}
\end{tabular}
\end{center}
\end{table}

\FloatBarrier

\section{The case when $ABC$ is isosceles with base angles $\alpha+\beta$}
Given a tile $(a,b,c)$ with angles $(\alpha,\beta,\gamma)$
(such that $3\alpha + 2\beta = \pi$) and $s = a/c > 0.656$, the diagrams of
Laczkovich \cite{laczkovich1995} can be used to construct
an $N$-tiling of some isosceles triangle with base angles $\alpha+\beta$,
for some, in general large, $N$.  For smaller $s$ his diagram does not
close up, as will be discussed below.  The facts that his method produces
only $N$ too large to draw,  and sometimes doesn't work, focus attention
on other methods of construction.

\begin{figure}
\centering
\includegraphics[width=0.4\textwidth]{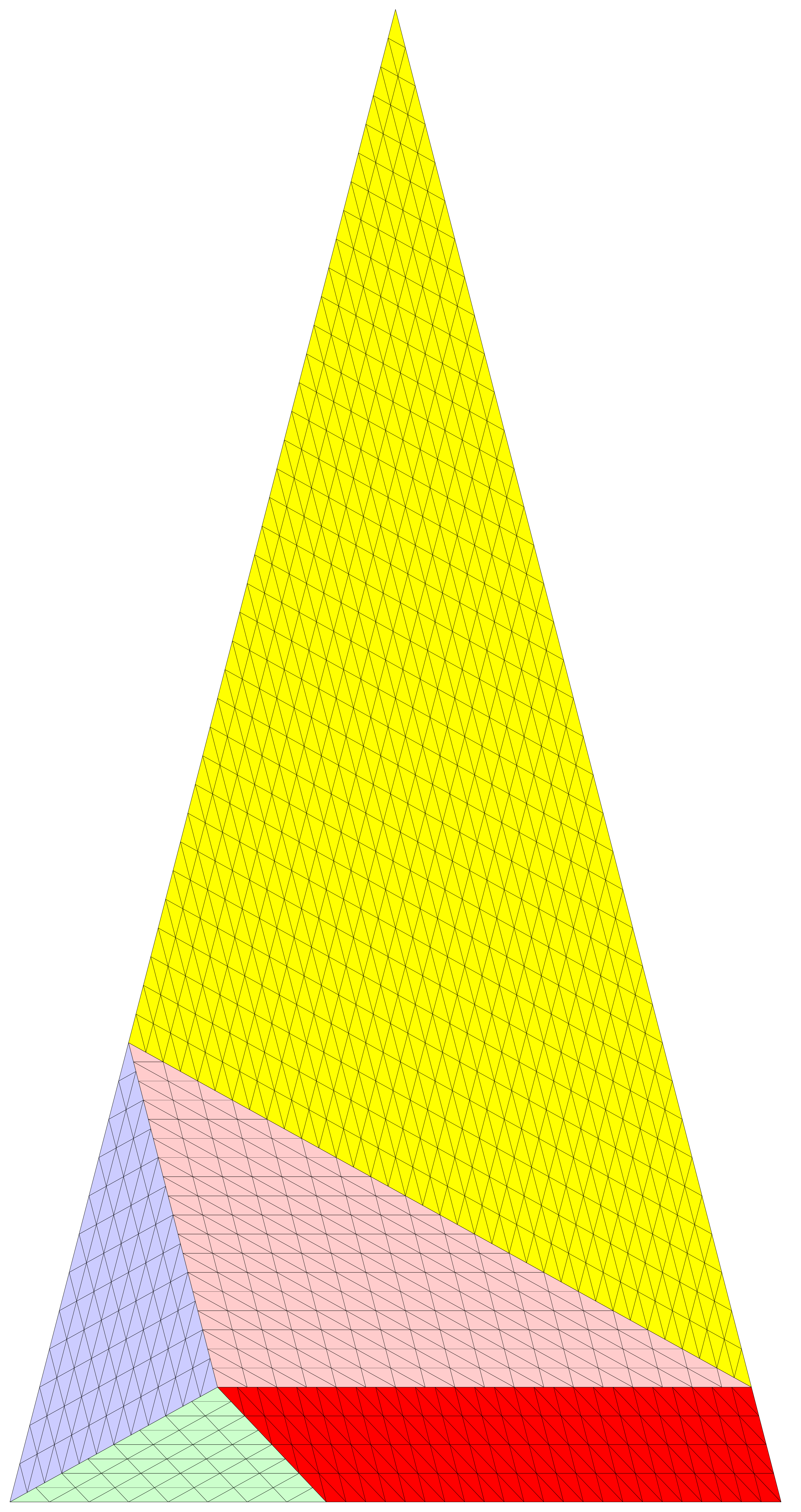} 
\quad 
\includegraphics[width=0.4\textwidth]{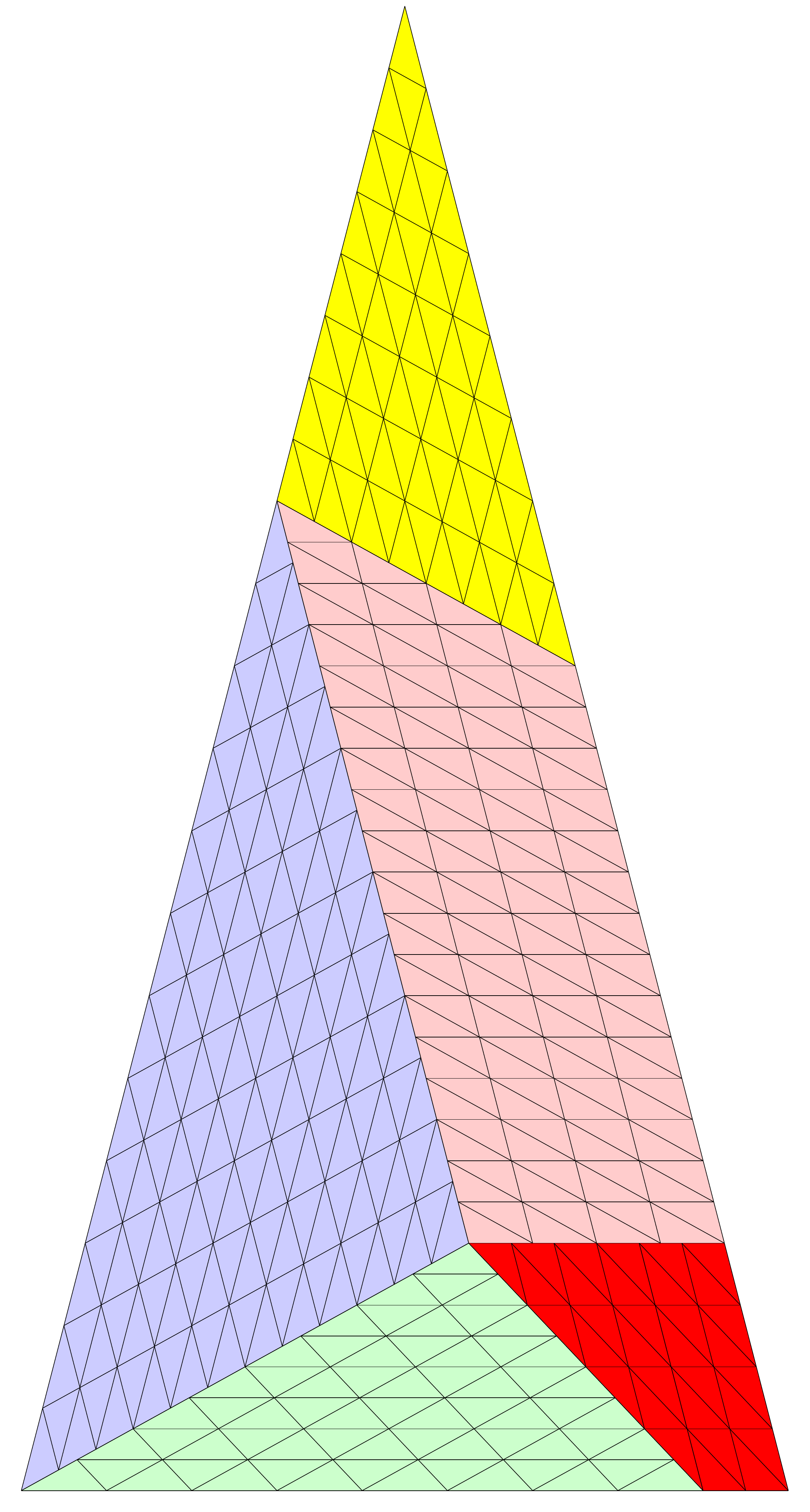}  
\caption{$N= 2028$ and $N=432$.  $ABC$ is isosceles.  The tile is $(2,3,4)$.}
\label{figure:2028}
\end{figure}

\begin{figure}
\centering
\includegraphics[width=0.34\textwidth]{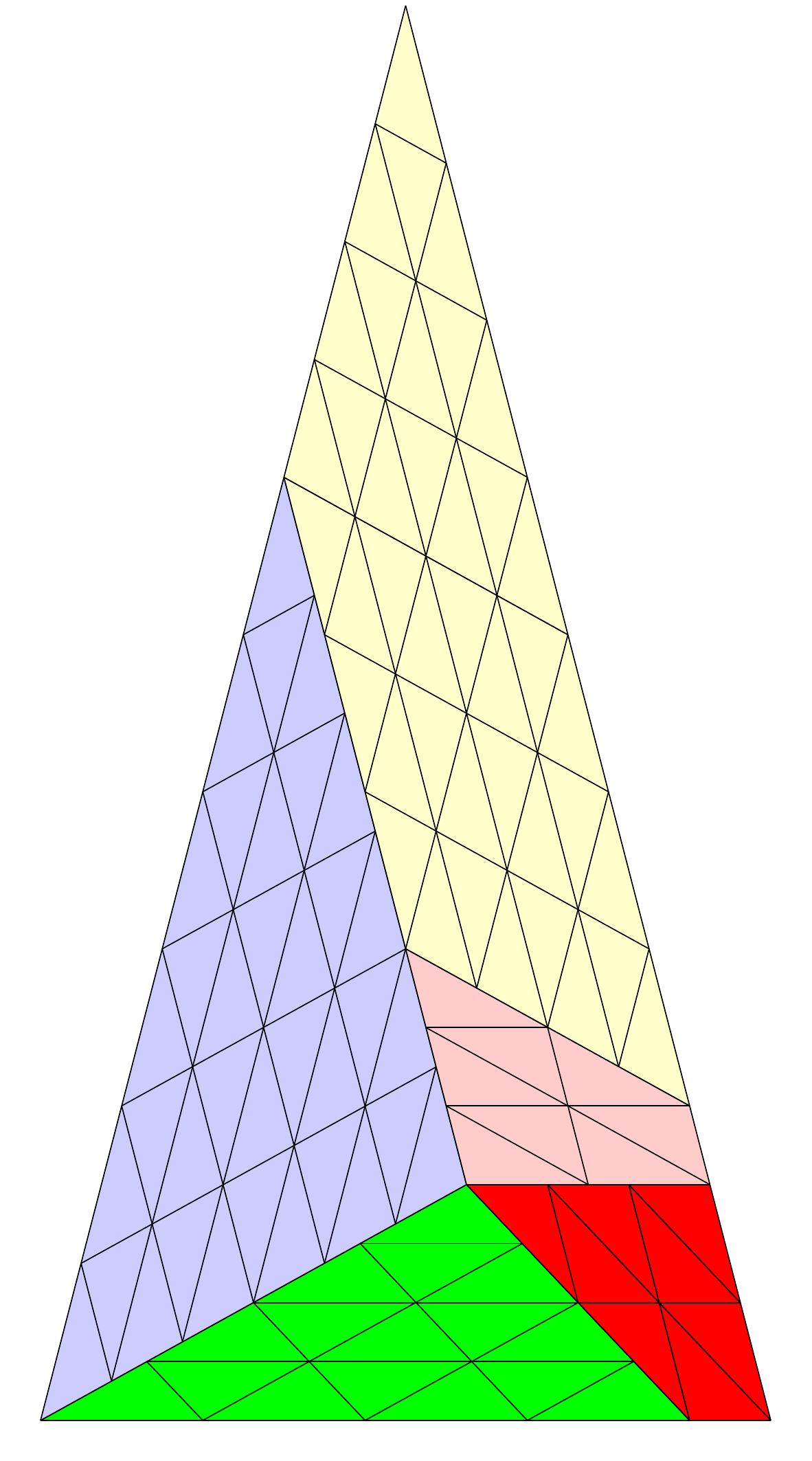}
\quad
\raisebox{2.3mm}{%
  \includegraphics[width=0.3\textwidth]{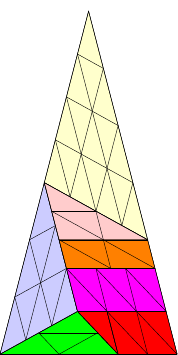}%
}
\caption{$N=108$ and $N=48$. $ABC$ is isosceles. The tile is $(2,3,4)$.}
\label{figure:48}
\end{figure}

\begin{figure}
\centering
\includegraphics[width=0.4\textwidth]{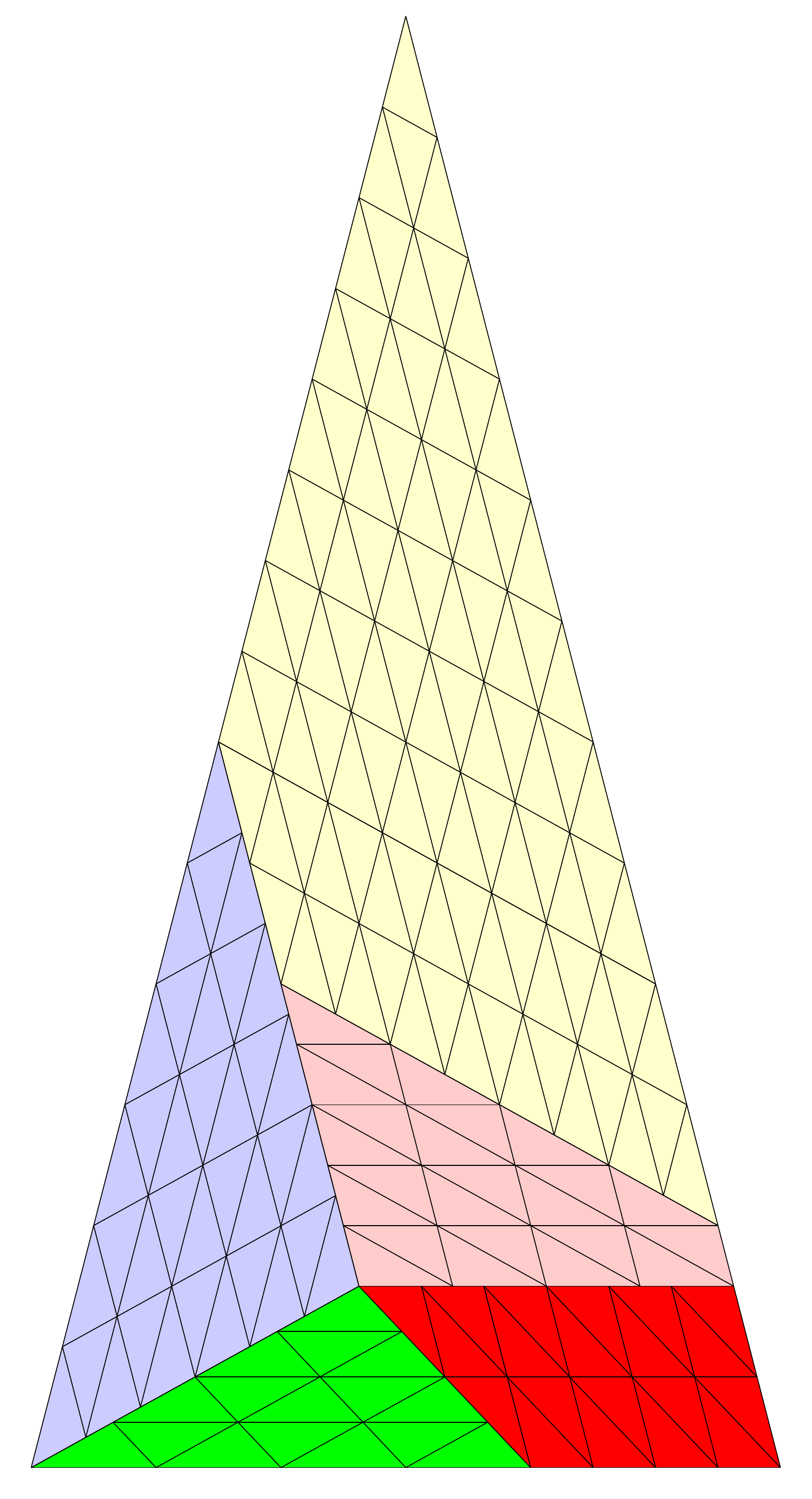}  
\quad
\includegraphics[width=0.4\textwidth]{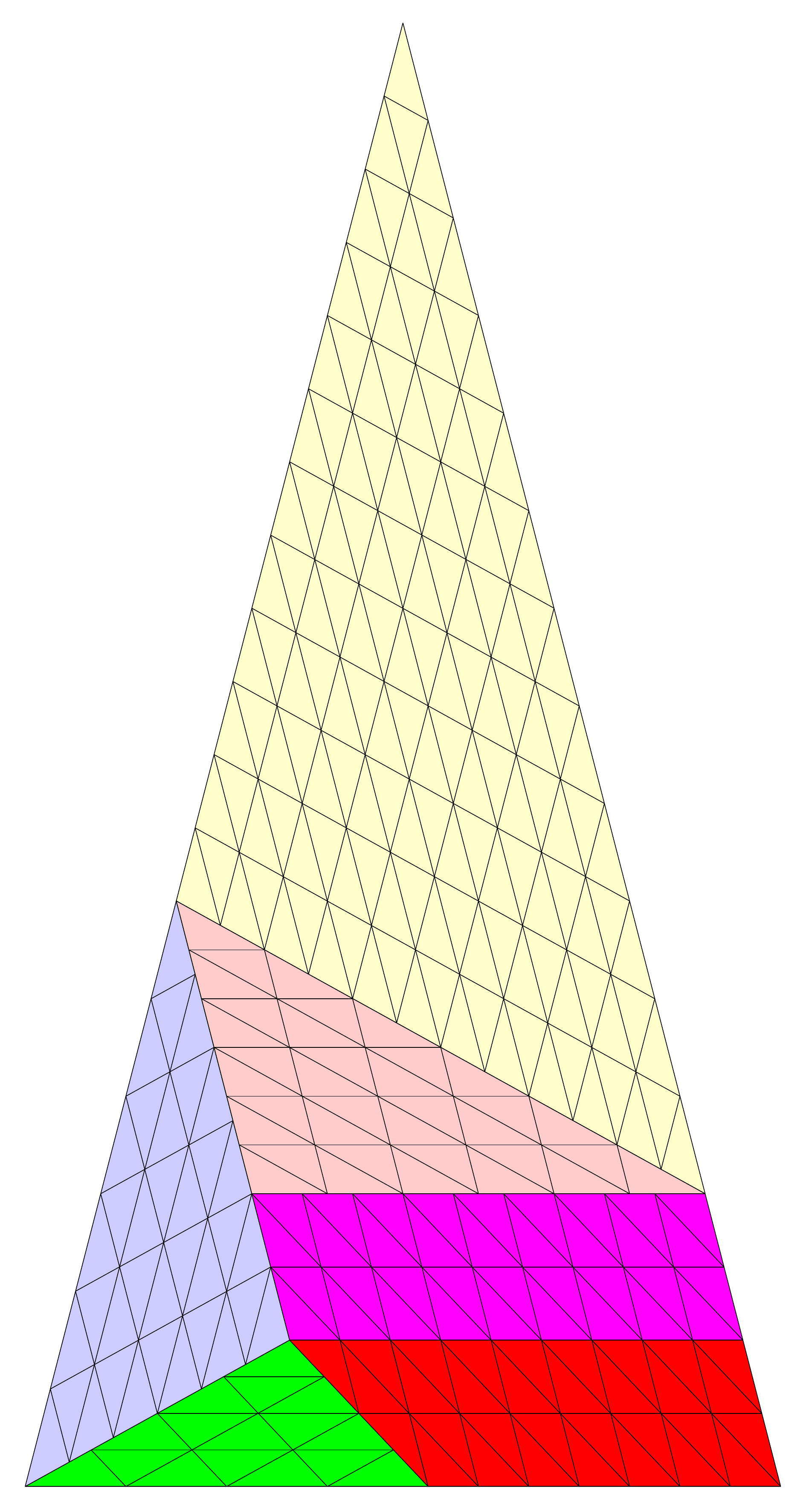} 
\caption{$N= 192$ and $N=300$.  $ABC$ is isosceles.  The tile is $(2,3,4)$.}
\label{figure:192}
\end{figure}

In July 2018, we discovered the tilings shown in Figures~\ref{figure:2028},
 \ref{figure:48}, and \ref{figure:192}.
 The particular values 
of $N$ were obtained from an equation to be explained below.  Like 
the first and second tiling equations, this equation does two things:
First, if it is not solvable, there is no $N$-tiling of any 
$ABC$ (that is isosceles of the form considered in this section).
Second, if the equation {\em is} solvable, each of the 
(finitely many) solutions determines the shape of the tile $(\alpha,\beta,\gamma)$.   That reduces the tiling problem to testing whether 
particular triangles $ABC$ can be tiled by particular tiles.
But if the equation is solvable, at present we have no efficient way
to determine whether there actually is a corresponding tiling.
That has to be checked on a case-by-case
basis.  In this section we will present these theoretical and empirical 
results.

\begin{figure}
\centering  
\includegraphics[width=0.45\textwidth]{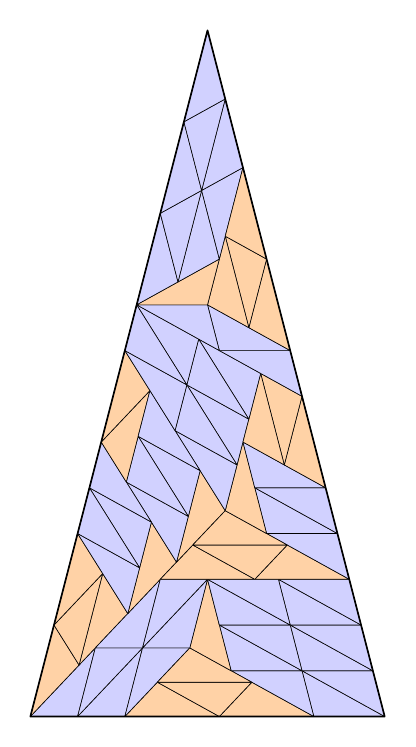} 
\quad
\includegraphics[width=0.35\textwidth]{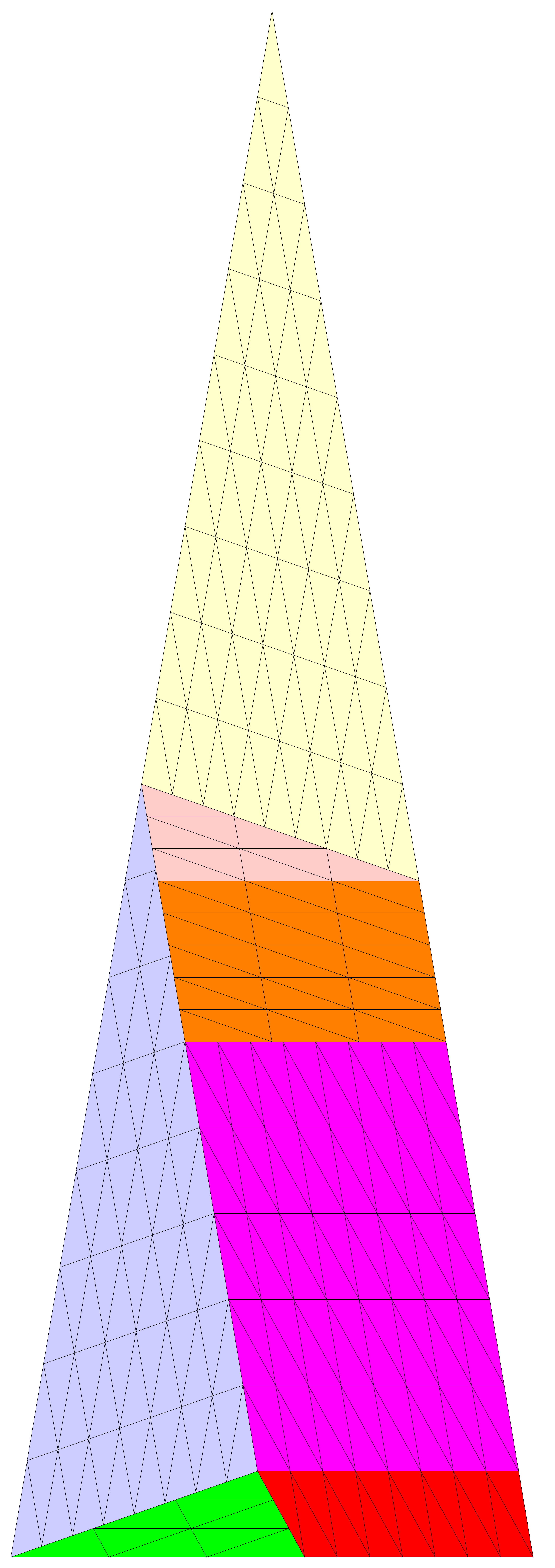} 
\caption{$N=75$, tile $(2,3,4)$, and $N=288$, tile $(3,8,9)$.}
\label{figure:288}
\end{figure}

\FloatBarrier

\subsection{A tiling equation for isosceles $ABC$ with base angles $\alpha + \beta$}

In this section we present necessary conditions for an $N$-tiling
of an isosceles triangle $ABC$ with base angles $\alpha+\beta$ 
to exist.  

\begin{definition} \label{definition:scaling}  We assume 
that the tile sides $(a,b,c)$ are integers with no common factor.
Define the ``scaling factors''
\begin{eqnarray*}
\lambda &=& Y / \sin \alpha \ = \ X/\sin(\alpha+\beta) \\
\kappa &=& a  / \sin \alpha \\
\mu &=& \lambda / \kappa
\end{eqnarray*}
\end{definition}

\begin{lemma} \label{lemma:XY}
Let triangle $ABC$ be isosceles with base angles $\alpha+\beta$ and vertex angle $\alpha$ at $B$, 
and  suppose $ABC$ is $N$-tiled by the integer-sided
tile $(a,b,c)$
with angles $(\alpha,\beta,\gamma)$ and $gcd(a,b,c) = 1$.  Let $X$
be the length of $AB$ and $Y$ the length of $AC$.  Then 
\begin{eqnarray*}
X &=& \mu c \\
Y &=& \mu a \\
\end{eqnarray*}
\end{lemma}

\begin{proof}  Let $\lambda$, $\kappa$, and $\mu$ be given 
by Definition~\ref{definition:scaling}.  
By the law of sines, 
\begin{eqnarray*}
X &=& \lambda \sin(\alpha+\beta) \\
Y &=& \lambda \sin \alpha\\
a &=& \kappa \sin \alpha \\
c &=& \kappa \sin \gamma
\end{eqnarray*}
 Since $\gamma = \pi - (\alpha+\beta)$,
the preceding equations imply
\begin{eqnarray*}
X &=& \bigg(\frac \lambda \kappa \bigg) c \ = \ \mu c  \\
Y &=& \bigg(\frac \lambda \kappa \bigg) a = \ \mu a 
\end{eqnarray*}
\end{proof}

\begin{theorem}[Necessary condition]\label{theorem:isosceles1}
Suppose triangle $ABC$ is $N$-tiled by
the tile $(a,b,c)$ with angles $(\alpha,\beta,\gamma)$, 
where  $3\alpha+2\beta=\pi$.  
Suppose triangle $ABC$ is isosceles with base angles $\alpha+\beta$.
  Let $s = a/c$. Then 
 the following Diophantine
equation is satisfied, where $M$ is the coloring number of the 
tiling:
\begin{eqnarray*}
\frac N {M^2} &=& \frac {1+s}{1-s}
\end{eqnarray*}
An algebraically equivalent form of the equation is 
$$ s = \frac {N-M^2}{N+M^2}. $$
The scale factor $\mu$ is given by 
$$ \mu = M(1+s).$$
\end{theorem}

\begin{example} For the tiling in Figure~\ref{figure:2028}, we have 
$N = 2028$,  and $s = a/c = 1/2$.  The coloring number $M$ is 26, as one can 
see by counting $m$ for each quadratic subtiling with $m^2$ tiles, with a sign 
positive or negative according to the color of the top tile of that subtiling.
The red trapezoid contributes 4, since parallelograms contribute 0 and the 
red trapezoid consists of a $4^2$ quadratic tiling plus a parallelogram.
Thus $M=26$.  Then both sides of the equation 
evaluate to 3.
\end{example}

\begin{proof}  The tile is necessarily rational, by 
Theorem~\ref{theorem:rationaltile}. 
The vertex angle $B$ of $ABC$ is $\alpha = \pi - 2(\alpha + \beta)$, since $3\alpha + 2\beta = \pi$.
 As usual let the sides
of the tile be $(a,b,c)$.
Let $X$ be the length of the two equal sides of $ABC$, namely $AB$ and $BC$, and $Y$ the 
length of the base $AC$. Then twice the area of $ABC$ is $X^2 \sin \alpha$,
and on the other hand it is $N$ times twice the area of the tile,
which is $Nbc\sin \alpha$.  Equating this two expressions we have 
the ``area equation''
\begin{eqnarray}
X^2 \sin \alpha &=& Nbc \sin\alpha \nonumber\\
X^2 &=& Nbc  \label{eq:4311}
\end{eqnarray}
On the other hand by Theorem~\ref{theorem:coloring}, with $M$ the coloring 
number (excess of black over white tiles) of the tiling, we have 
\begin{eqnarray}
M(a+b+c) &=& 2X-Y  \label{eq:4316}
\end{eqnarray}
Since $ABC$ is a triangle, the sum of any two sides is greater than the 
third side, so $2X > Y$.  It follows that
\begin{eqnarray}
M &>& 0  \label{eq:4348}
\end{eqnarray}
By Lemma~\ref{lemma:XY}, we have 
\begin{eqnarray*}
X &=& \mu c \\
Y &=& \mu a
\end{eqnarray*}
Putting these into the area equation (\ref{eq:4311}) and
 the coloring equation (\ref{eq:4316}), respectively, 
 we have
\begin{eqnarray}
\mu^2 c^2 &=& Nbc \label{eq:4328}\\
M(a+b+c) &=& \mu (2c-a) \nonumber
\end{eqnarray}
Dividing by $c$ we have 
\begin{eqnarray*}
M\bigg( \frac a c + \frac b c + 1 \bigg) &=& \mu\big(2- \frac a c\big)
\end{eqnarray*}
Expressing everything in terms of $s$ by Lemma~\ref{lemma:singamma}, we have 
\begin{eqnarray*}
M(s + (1-s^2) + 1) &=& \mu(2-s) \\
 M(2+s-s^2) &=& \mu(2-s)\\
M(2-s)(1+s) &=& \mu(2-s) \\
\end{eqnarray*}
Cancelling $(2-s)$ we have 
\begin{eqnarray*}
\mu &=& M(1+s),
\end{eqnarray*}
which is the last assertion in the theorem.
Putting this expression into (\ref{eq:4328}) we find
\begin{eqnarray*}
c^2 M^2(1+s)^2 &=& Nbc 
\end{eqnarray*}
Dividing by $c^2$ and using Lemma~\ref{lemma:singamma}, we have
\begin{eqnarray*}
M^2(1+s)^2 &=& N(1-s^2) \\
           &=& N(1-s)(1+s) \\
\end{eqnarray*}
Therefore
\begin{eqnarray*}
\frac N {M^2} &=& \frac {1+s}{1-s}.
\end{eqnarray*}
\end{proof}

\begin{corollary}[Shape determines the tile, $M$ and $N$ determine each other]
\label{corollary:shapedeterminestile}
Let $N$ and $ABC$ be given, and suppose $ABC$ is isosceles with 
vertex angle $\alpha$ at $B$ and base angles $\alpha+\beta$,
and $ABC$ is $N$-tiled by a tile with angles $(\alpha,\beta,\gamma)$. 
Let $Y$ be the length of the base $AC$ and
let  $X$ be the length of the equal sides
 $AB$ and $BC$.  
Then the only possible value of $s=a/c$ for a tile that could be 
 used to tile $ABC$ is $s = Y/X$, and the number of tiles $N$ 
 and the coloring number $M$ must satisfy
$$  M^2 = N \frac{1-s}{1+s}  = N \frac {X-Y}{X+Y} $$
The number on the right must be a square.
\end{corollary}

\begin{proof}  
By Lemma~\ref{lemma:XY} we have $s = Y/X$.  By 
Theorem~\ref{theorem:isosceles1} we have the first equation of 
the theorem.  By Lemma~\ref{lemma:XY} we have 
\begin{eqnarray*}
\frac {X-Y}{X+Y} &=& \frac {c-a}{c+a}
\end{eqnarray*}
Dividing numerator and denominator by $c$ and using $s = a/c$ we 
have 
$$ \frac{c-a}{c+a} = \frac {1-s}{1+s}.$$
\end{proof}

\subsection{$N$ is not prime when $ABC$ is isosceles with base angles $\alpha+\beta$}
 
\begin{theorem} \label{theorem:notprime3}
Let $N$ and $ABC$ be given, and suppose $ABC$ is isosceles with 
 base angles $\alpha+\beta$.  Suppose
$ABC$ is $N$-tiled.  Then $N$ is not a prime number.
\end{theorem}

\begin{proof}  Let $X$ be the length of the equal sides 
$AB$ and $BC$, and $Y$ the length of the base $AC$.
  By Theorem~\ref{theorem:rationaltile},
the tile is rational. Let $(a,b,c)$ be the lengths of the 
edges of the tiles; without loss of generality, we may assume that 
$a$, $b$, and $c$ are integers.  
   By Lemma~\ref{lemma:anglesOK}, $b = c-a^2/c$,
so $bc$ is an integer, even if $b$ is not necessarily an integer.
The area equation is 
\begin{equation}
X^2 = Nbc. \label{eq:4681}
\end{equation}

By Theorem~\ref{theorem:isosceles1}, 
$$ s =  a/c = \frac {N-M^2}{N+M^2}. $$
Let $m$ and $k$ be such that $s = m/k$ in lowest terms.
Then by Lemma~\ref{lemma:gsquare},
\begin{equation}
(a,b,c) = (mk,k^2-m^2,k^2). \label{eq:4795}
\end{equation}
Putting that into the previous equation, we have
\begin{eqnarray}
 (N+M^2)a &=& (N-M^2)c \nonumber\\
 (N+M^2)mk &=& (N-M^2)k^2 \nonumber\\
 (N+M^2)m &=& (N-M^2)k \nonumber\\
 N(k-m)&=& M^2(k+m)  \label{eq:8661}
\end{eqnarray}
Then the area equation \eqref{eq:4681} can be written
\begin{eqnarray*}
 X^2 &=& N(k^2-m^2)k^2  \\
 &=& N(k-m)(k+m)k^2 \\
 &=& M^2(k+m)^2 k^2  \mbox{\qquad by (\ref{eq:8661})}
\end{eqnarray*}
Taking the square root we have
\begin{equation*}
X = M(k+m)k
\end{equation*}
By Lemma~\ref{lemma:XY}, we have $Y = (a/c) X = (m/k)X$.
That leads to
\begin{equation}
Y = M(k+m)m  \label{eq:8676}
\end{equation}

Let $d = \gcd(N-M^2,N+M^2)$.  Then $d$ divides the sum
and difference of $N-M^2$ and $N+M^2$, so $d \mid 2N$ and $d \mid 2M^2$.
Since $s > 0$ we have $0 < M^2 < N$, so $\gcd(N,M^2) = 1$ as $N$ is prime;
hence $d \mid 2\gcd(N,M^2) = 2$.  Since $m/k$ is $a/c$ in lowest terms,
$(m,k) = \big((N-M^2)/d,\ (N+M^2)/d\big)$, so $k - m = 2M^2/d \ge M^2$.
Therefore
\begin{equation}
k - m - M \ \ge \ M^2 - M \ \ge \ 0, \label{eq:4829}
\end{equation}
and in particular $0 < M < k$, since $m \ge 1$.

Let $Y = pa + qb + rc$ be the edge decomposition of the base $AC$.
We have $r \ge 1$ by  Lemma~\ref{lemma:pigeonhole},  
since there are no $\gamma$ angles at $A$ or $C$.
 Thus
\begin{equation}
Y  = pa + qb + rc \ \ge \ (k-M)b + c   \label{eq:8692}
\end{equation}
On the other hand, by (\ref{eq:8676}),
\begin{eqnarray*}
(k-M)b - Y &=&  (k-M)(k^2-m^2) - Mm(k+m) \\
 &=& (k+m)\left[(k-M)(k-m) - Mm\right] \\
 &=& k(k+m)(k-m-M),
\end{eqnarray*}
which is nonnegative by (\ref{eq:4829}).  Hence $Y \le (k-M)b $.
Since $c > 0$ we have  $$Y < c+(k-M)b,$$
contradicting (\ref{eq:8692}).
The assumption that $N$ is prime has led to a contradiction.
\end{proof}

\subsection{Existence of tilings of isosceles $ABC$ with base angles $\alpha+\beta$}
We first remark that some of the diagrams in \cite{laczkovich1995} 
imply the existence of tilings.  Those diagrams imply that,
in case the diagrams as shown actually exist, and the figures in them 
are rational, then for large enough $N$, there are $N$-tilings of 
some triangle with the specified shape. We tried to construct tilings
from Laczkovich's diagram for the case of isosceles $ABC$ with base
angles $\alpha + \beta$, but satisfying the many least-common-multiple
conditions on the interior edges requires extremely large values of $N$.
These diagrams did not result in any tilings small enough to 
draw.  We found the examples given above independently, more or less
by trial and error, and initially were hopeful that a similar construction
could be proved to work in general.  In September, 2026, we noticed that
Laczkovich's construction has a similar problem:  there are angles $\alpha$
for which it fails, as will be shown below. In this section we will 
precisely state the conditions under which our 2018 method works,  and 
the conditions under which Laczkovich's 1995 method works, and 
wonder of wonders,  they are exactly complementary,  so for every rational
tile, one of the two methods will produce a tiling of some isosceles
triangle with base angles $\alpha+\beta$.

\begin{theorem}[Laczkovich's construction] \label{theorem:Lconstruction}
Suppose $3\alpha + 2\beta = \pi$, and let $ABC$ be an
isosceles triangle with base angles $\alpha+\beta$.
Suppose the tile $(\alpha,\beta,\gamma)$ is rational, and that
$b/c < x_0$, where $x_0 \approx 0.5698$ is the real root of
$x^3 = (1-x)^2$; equivalently $s > \sqrt{1-x_0} \approx 0.656$.
Then for
some $N$, there is an $N$-tiling of a triangle similar to $ABC$
by the tile $(\alpha,\beta,\gamma)$, where $\gamma = 2\alpha+\beta$.
\end{theorem}

\begin{remark}  The theorem is true without the condition on $b/c$, as we
  prove in Corollary~\ref{corollary:isoscelesall}; but
it does not follow from Laczkovich's construction.  The condition
that the tile is rational cannot be dropped: by
Theorem~\ref{theorem:rationaltile}, a tile that tiles $ABC$ is rational.
\end{remark}

\begin{figure}[ht]
\centering
\includegraphics[width=0.4\textwidth]{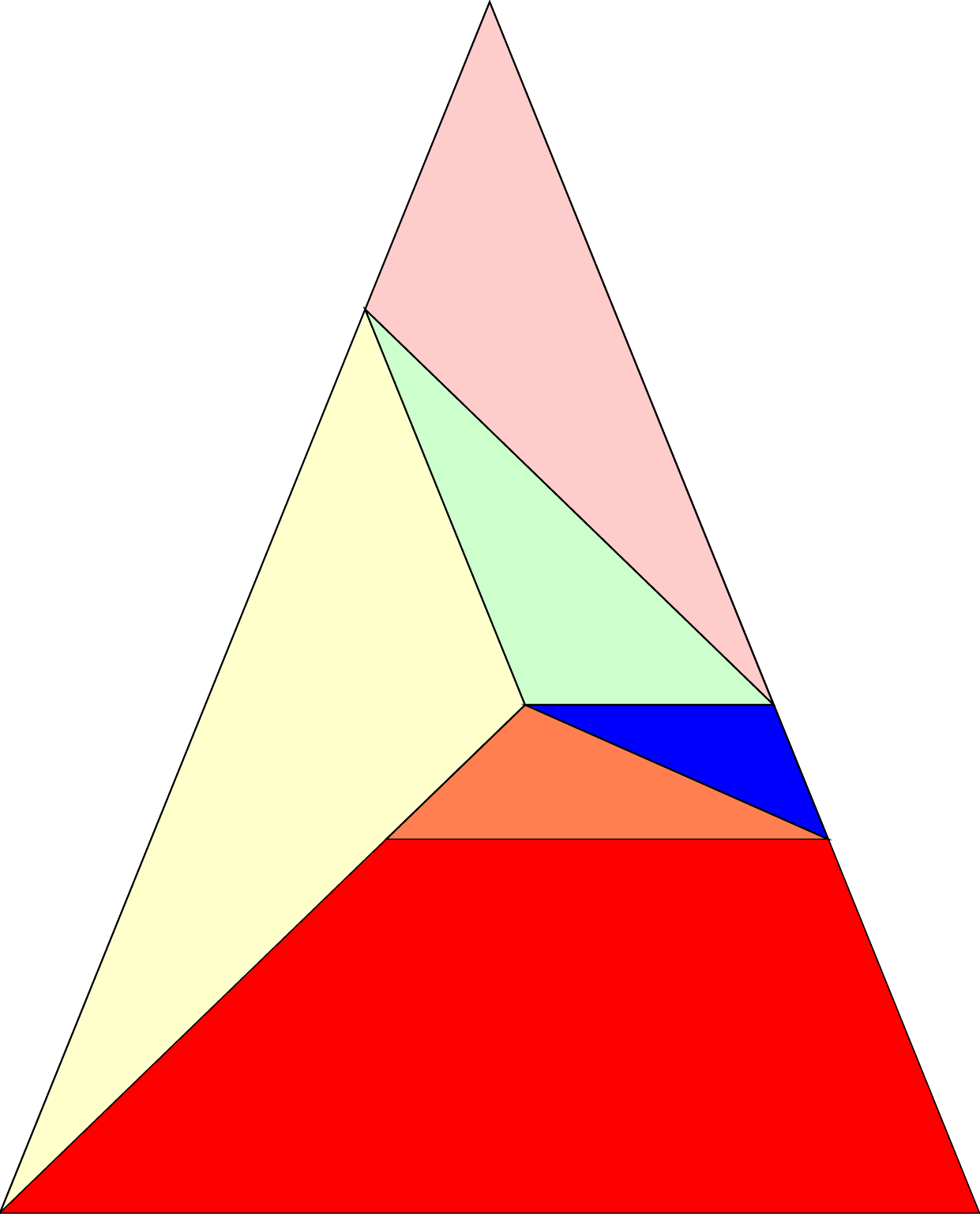}  
\caption{Laczkovich's method for isosceles $\alpha+\beta$, with tile $(12,7,16)$.
}
\label{figure:L1995isosceles}
\end{figure}

\begin{figure}[ht]
\centering
\includegraphics[width=0.4\textwidth]{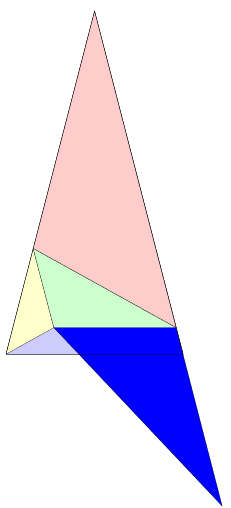}
\caption{Laczcovich's method for isosceles $\alpha+\beta$, with tile $(2,3,4)$.  The method fails
because the blue triangle will not fit in $ABC$.
}
\label{figure:L1995isosceles234}
\end{figure}

\begin{proof}
Laczkovich said \cite{laczkovich1995} that an isosceles triangle $ABC$ with
base angles $\alpha+\beta$, where $3\alpha+2\beta = \pi$,  can be dissected into
similar triangles and a trapezoid, as shown in Figure~\ref{figure:L1995isosceles}.
Since the tile is rational, we can find a triangle with integer sides
$(a,b,c)$ with no common factor having angles $(\alpha,\beta,\gamma)$.  Then
the similar triangles in Figure~\ref{figure:L1995isosceles}  can be scaled (changing the
size of $ABC$ but not its angles) so that triangles in the figure can all be
quadratically tiled.  How to tile the red trapezoid is explained in \cite{laczkovich1995}.
 Then if necessary a further rescaling permits the trapezoid
also to be tiled.  In general the value of $N$ resulting will be large.

The reason that we discuss this construction here instead of just citing it,  is that
it only works when $b/c < x_0$, the root of $x^3 = (1-x)^2$; that is,
when $b^3 < c(c-b)^2$.  In terms of $s$ (since $b/c = 1-s^2$), it is when
$s > \sqrt{1-x_0} \approx 0.656$.   Outside that range, the
blue triangle in the figure does not fall inside triangle $ABC$.  See Figure~\ref{figure:L1995isosceles234}
for what happens with tile $(2,3,4)$. Computation of the coordinates needed to draw the figure
precisely yields the stated condition.
\end{proof}

\FloatBarrier

Now we discuss the methods used to produce tilings like those in Figure~\ref{figure:192}.

Our method starts with $N$ and $M$. 
Define the scaling factor 
$$ \mu = M(1+s).$$   
We generalize the construction of the 
48-tiling in Figure~\ref{figure:48},  which will require 
certain number-theoretic relations between $a$, $b$, and $c$ 
in order that the boundaries of the colored regions in the 
figure line up properly.

The tiling is made of quadratic tilings and parallelograms; adding up
their sizes, the total is $M^2(1+s)/(1-s)$, which is $N$ by
Theorem~\ref{theorem:isosceles1}.  Here is how the figure is 
constructed.  We start with the green quadratic tiling, which will have $p^2$ tiles.
We just pick a $p$ and try it.
Then draw the adjacent light blue tiling; the red tiling 
east of the green; then the yellow quadratic tiling at the top.  At each
color boundary there is an equation to be satisfied relating the number of tiles
in the two adjacent quadratic tilings.  When we come to the yellow quadratic
tiling, sometimes the method fails:
 the yellow quadratic tiling
can  overlap the red,  ruining the tiling. 
\medskip

\begin{example}  The least example of
this is $N = 180$ with the tile $(6,5,9)$: $M = 6$, $\mu = 10 > c$,
and the only admissible $p$ is $6$, for which the yellow tiling overlaps
the red.  With the tile $(12,7,16)$ the least example is $N = 1008$,
shown in Figure~\ref{figure:1008}; there $\mu = 21$.
The method fails to produce a tiling of (any triangle similar to) the isosceles triangle with sides $336$
and base $252$,  although the tiling equation checks:
$$ \frac N {M^2} = \frac {1008} {144} = 7, \quad \frac {1+s}{1-s} = \frac {7/4}{1/4} = 7.$$
\end{example}

\begin{figure}[ht]
\centering
\includegraphics[width=0.4\textwidth]{Figures/tiling1008}  
\raisebox{2.6cm}{%
\includegraphics[width=0.4\textwidth]{Figures/tiling1008B}
}
\caption{$N=1008$.  $ABC$ is isosceles. The tile is $(12,7,16)$.
The tiling method used for previous tilings fails; $p = 12$ is 
too small and $p=24$ is too large. 
}
\label{figure:1008}
\end{figure}

\begin{lemma} \label{lemma:complementary}
Let $(a,b,c)$ be a tile with  $b = c - a^2/c$.
Then
$$ b(a^2+b^2) - a^2c \ = \ -\,\bigl(c^3 - 2bc^2 + b^2c - b^3\bigr), $$
and if $b/c$ is rational, neither side is zero.  Hence for every rational
tile exactly one of the following holds:
\begin{itemize}
\item[(i)] $c + b^2/c < a^2/b + b$, the condition for Laczkovich's
construction (equivalently $b/c < x_0$);
\item[(ii)] $a^2c < b(a^2+b^2)$ (equivalently $b/c > x_0$),
\end{itemize}
where $x_0 \approx 0.5698$ is the real root of $x^3 = (1-x)^2$.
\end{lemma}

\begin{proof}
Since $a^2 = c(c-b)$, we have
$b(a^2+b^2) - a^2c = a^2(b-c) + b^3 = b^3 - c(c-b)^2$, while
$c^3 - 2bc^2 + b^2c - b^3 = c(c-b)^2 - b^3$.  Dividing by $c^3$, the
second expression vanishes if and only if $x = b/c$ satisfies
$x^3 - x^2 + 2x - 1 = 0$, i.e.\ $x^3 = (1-x)^2$.  By the rational root
test the only possible rational roots are $\pm 1$, and neither is a root.
Condition (i), multiplied by $bc$, says $c(c-b)^2 - b^3 > 0$, and (ii)
says the negative of that is positive.
\end{proof}

\begin{theorem} \label{theorem:isoscelesconstruction}
Let the tile $(a,b,c)$ have angles $(\alpha,\beta,\gamma)$ with
$3\alpha+2\beta= \pi$, let $a$, $b$, $c$ be integers with no common factor
(so $c \mid a^2$), and suppose
$$ a^2 c \ < \ b(a^2+b^2). $$
Then for some $N$, some isosceles triangle with base angles $\alpha+\beta$
can be $N$-tiled by the tile.
\end{theorem}
\begin{remark}
This theorem is based on the method used to construct the tilings
 exhibited in Figure~\ref{figure:192}.  As noted above, for certain 
choices of $\alpha$ that method fails.  And also as noted above, 
for certain choices of $\alpha$, Laczkovich's method fails.  As it 
happens, they never {\em both} fail.
\end{remark}

\begin{proof}[Proof of Theorem~\ref{theorem:isoscelesconstruction}]
Let $M$ and $p$ be positive integers, to be chosen below, let
$\mu = M(1+s)$, and let $ABC$ have vertex angle $\alpha$ at $B$, base
$AC$ of length $Y = \mu a$ and legs of length $X = \mu c$.  Place $A$ at
the origin and $C$ on the positive $x$-axis.  Define
$$ q = \frac{pb}{a}, \quad t = \frac{pa}{c}, \quad r = \mu - \frac{pc}{a},
   \quad u = \frac{(\mu-q)c}{b}, \quad v = \frac{(\mu-q)a}{b} = \frac{ua}{c},
   \quad h = qb - va. $$
Write $T(n)$ for a triangle similar to the tile with ratio $n$.  We
assemble $ABC$ from the following pieces.
\begin{enumerate}
\item {\em Green}: $T(p)$ with $\alpha$ at $A$, its $c$ side $AP$ on the
base, $|AP| = pc$, and $\gamma$ at its apex $G$; so $|AG| = pb$ and
$|PG| = pa$.
\item {\em Blue}: $T(q)$ with $\beta$ at $A$, its $a$ side on $AG$
(note $qa = pb$), $\gamma$ at $G$, and its $c$ side $AK$ on $AB$,
$|AK| = qc$.
\item {\em Yellow}: $T(u)$ with $\alpha$ at $B$, its $b$ side $KB$ on
$AB$ and its $c$ side $BL$ on $BC$, $|BL| = uc$.  This fits since
$qc + ub = qc + (\mu - q)c = \mu c = X$.
\item {\em Red}: $T(t)$ with $\alpha$ at $P$, its $c$ side on $PG$
(note $tc = pa$), $\beta$ at $G$, and its $b$ side $PQ$ parallel to
$CB$; together with the parallelogram $PCRQ$ with $|PC| = ra$ and
$|PQ| = |CR| = tb$, tiled by $t$ rows of $r$ cells, each cell two tiles
forming a parallelogram with sides $a$ and $b$.  This fits on the base
since $pc + ra = \mu a = Y$.
\item {\em Pink}: $T(v)$ with $\beta$ at $K$, its $c$ side on $KL$
(note $vc = ua$), $\alpha$ at $L$, and its $a$ side $KV$ on $KG$,
$|KV| = va$.
\item {\em Orange and violet}: the parallelogram $GRLV$.
\end{enumerate}
The angles fit at every vertex: at $A$, $\alpha+\beta$; at $P$,
$\beta+\alpha+\gamma$; at $Q$, $\gamma + (\alpha+\beta)$; at $G$,
$\gamma+\gamma+\beta+\gamma = 2\pi$ (a center); at $K$,
$\alpha+\gamma+\beta$; at $L$, $\beta+\alpha+\gamma$; at $V$,
$\gamma+(\alpha+\beta)$; at $R$, $\gamma+(\alpha+\beta)$.
Since the angle at $Q$ is $\pi$, the points $G$, $Q$, $R$ are collinear,
and $|GR| = ta + ra$.  The line $GR$ is parallel to $AC$, since
$G$ and $R$ have heights $pb\sin\alpha$ and $tb\sin(\alpha+\beta) =
tb\sin\gamma$, which are equal because $pa = tc$ and
$a/\sin\alpha = c/\sin\gamma$.  Next, $t + r = \mu - q$ (since
$q = p/s - ps$), so $|GR| = (\mu-q)a = vb = |LV|$.  The angles at $R$ and
$L$ show that $LV \parallel RG$, so $GRLV$ is a parallelogram with angle
$\gamma$ at $G$ and sides $vb$ and $|GV| = |KG| - |KV| = qb - va = h$.
As a check on $BC$:
$|BL| + |LR| + |RC| = uc + h + tb = \mu c$, because
$uc - va = (\mu-q)(c^2-a^2)/b = (\mu-q)c$ and $q(c-b) = tb$.
Since $GR$ consists of $v$ $b$ edges or of $t+r$ $a$ edges, the
parallelogram $GRLV$ can be tiled by $a \times b$ cells in strips, provided
$h = xa + yb$ with $x, y \ge 0$ integers: $x$ strips of width $a$ with $v$
cells each, and $y$ strips of width $b$ with $t+r$ cells each.

It remains to choose $M$ and $p$ so that $q$, $t$, $r$, $u$, $v$ are
nonnegative integers, $u > 0$, and $h$ is of the form $xa+yb$.  Take
$M = \lambda abc$ and $p = ac\rho$ with positive integers $\lambda$,
$\rho$.  Then $\mu = \lambda ab(a+c)$ and
$$ q = bc\rho, \quad t = a^2\rho, \quad r = \lambda ab(a+c) - c^2\rho, $$
$$ u = c\,w, \quad v = a\,w, \quad h = \rho c(a^2+b^2) - \lambda a^3(a+c),
\quad \mbox{where } w = \lambda a(a+c) - c\rho, $$
all integers.  We need $r \ge 0$ and $h \ge ab$; then $w > 0$ follows
from $r \ge 0$ and $b < c$.  Since $\gcd(a,b) = 1$ by
Lemma~\ref{lemma:gsquare}, every integer $h \ge ab$ has the form $xa+yb$
with $x,y \ge 0$.  So we need an integer $\rho$ with
$$ \frac{\lambda a^3(a+c) + ab}{c(a^2+b^2)} \ \le \ \rho \ \le \
   \frac{\lambda ab(a+c)}{c^2}. $$
The length of this interval is
$$ \frac{\lambda a(a+c)\bigl(b(a^2+b^2) - a^2c\bigr)}{c^2(a^2+b^2)}
   - \frac{ab}{c(a^2+b^2)}, $$
which tends to infinity with $\lambda$, since the coefficient $b(a^2+b^2) - a^2 c$
is positive.   Choose $\lambda$ so large that this length exceeds $1$,
and choose an integer $\rho$ in the interval.  The resulting
tiling has $N = X^2/(bc) = \mu^2 c/b$ tiles, by the area equation.
\end{proof}

\begin{corollary} \label{corollary:isoscelesall}
Let $3\alpha+2\beta = \pi$ and suppose the tile $(\alpha,\beta,\gamma)$ is
rational.  Then for some $N$, some isosceles triangle with base angles
$\alpha+\beta$ can be $N$-tiled by the tile.
\end{corollary}

\begin{proof}
Let $(a,b,c)$ be the integer multiple of the tile with no common factor.
By Lemma~\ref{lemma:complementary}, either $b/c < x_0$, and
Theorem~\ref{theorem:Lconstruction} applies, or $a^2c < b(a^2+b^2)$, and
Theorem~\ref{theorem:isoscelesconstruction} applies.
\end{proof}

\begin{example} With $N = 112$, $M=4$, $\mu = 7$, and 
$(a,b,c) = (12,7,16)$, the tiling equation holds.
   Our construction method fails on this example, since $s = 3/4 > 0.656$,
   so by Lemma~\ref{lemma:complementary} the condition $h \ge 0$ fails for
   every choice of $p$ and $M$.  But Laczkovich's construction
does work in this case, and it will produce a tiling of the given $ABC$ enlarged by
a factor of $\lambda$,  with $N = 112 \lambda^2$ and $M = 4\lambda$,  which satisfies
the tiling equation with the same $s$.   With Laczkovich's method, $\lambda$ is at least 579,
and maybe larger depending on how clever we are in tiling the red trapezoid.  That leads
to an $N$ of at least $37,546,992 = 112 \cdot 579^2$ and $M$ at least $2316$.    So, some triangle 
the shape of $ABC$ can be tiled,  but this does not tell us
 if the isosceles triangle $ABC$ with base 84 and sides 112 can 
be tiled by this tile. 
\end{example}

\subsection{Application to triangles with angles $(2\alpha,\alpha,2\beta)$}
\begin{figure}[ht]
\begin{center}
\includegraphics[width=0.3\textwidth]{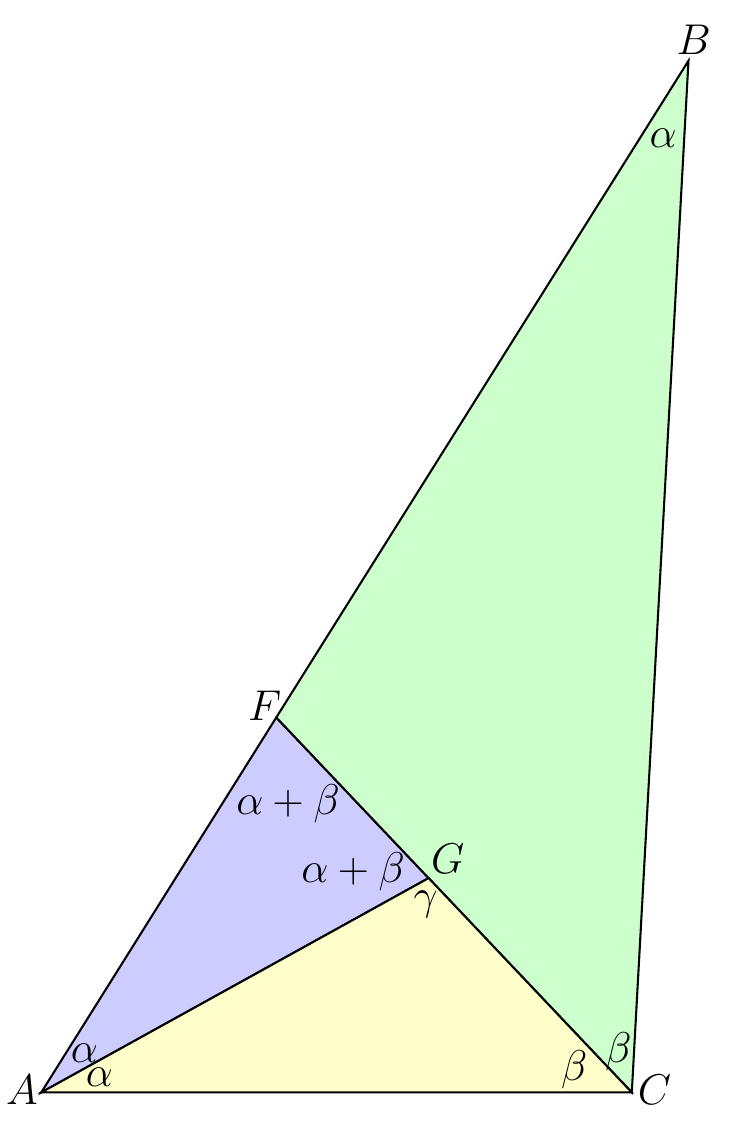}
\end{center}
\caption{Decomposing $(2\alpha,\alpha,2\beta)$ (adapted from \cite[Fig.~5]{laczkovich1995})}
\label{figure:Lfive}
\end{figure}

Corollary~\ref{corollary:isoscelesall} also gives tilings of triangles
 with angles $(2\alpha,\alpha,2\beta)$ for every tile, including  the case
$c > 2b$ left open by Theorem~\ref{theorem:existence2}.  The link is shown
in Figure~\ref{figure:Lfive}.  Bisecting the angle $2\beta$ at $C$ cuts off
a triangle with angles $\beta$ at $C$ and $\alpha$ at $B$, hence similar to
the tile; bisecting the angle $2\alpha$ at $A$ then cuts off a second
triangle with angles $\alpha$ at $A$ and $\beta$ at $C$, also similar to the
tile.  What remains is an isosceles triangle with vertex angle $\alpha$ and
base angles $\alpha+\beta$.  Scale $ABC$ so that $|AC| = c^2$.  Then the law of sines gives these
three triangles the sides $(b+c)(a,b,c)$, $c\,(a,b,c)$, and
$(bc, bc, ab)$ respectively.  In particular the two cuts fall inside
$ABC$ as drawn, since $|BF| = (b+c)b < b(b+2c) = |AB|$ and
$|CG| = ca < (b+c)a = |CF|$.

\begin{corollary} \label{corollary:existence2case2}
Suppose $c$ divides $a^2$, where $a < c$ are integers, and $b = c - a^2/c$.
 Then some triangle with angles $(2\alpha,\alpha,2\beta)$
can be $N$-tiled by the tile $(a,b,c)$, with
$$ N = L^2(b+c)(b+2c), \qquad M = L(b+c-a) $$
for some positive integer $L$.
\end{corollary}

\begin{proof}
By Corollary~\ref{corollary:isoscelesall}, some isosceles triangle with
base angles $\alpha+\beta$, with legs of length $L$ say, can be tiled.
Replacing each tile by a quadratic tiling of scale $bc$ tiles the
isosceles triangle with legs $Lbc$, which is the blue triangle of
Figure~\ref{figure:Lfive} scaled by $L$.  With the same scaling, the other
two triangles are similar to the tile with integer ratios $L(b+c)$ and
$Lc$, so they have quadratic tilings.  The number of tiles is
$L^2\bigl((b+c)^2 + c^2 + bc\bigr) = L^2(b+c)(b+2c)$.
By Theorem~\ref{theorem:coloring}, $M(a+b+c) = X - Y + Z$, which is positive by the
triangle inequality.  Then $M > 0$, and the second tiling equation, in
the form of Theorem~\ref{theorem:secondcasesolution}, gives
$M = L(b+c-a)$.
\end{proof}

\begin{corollary} \label{corollary:secondshapeall}
Let the tile have angles $(\alpha,\beta,\gamma)$ with $3\alpha+2\beta=\pi$.
Some triangle with angles $(2\alpha,\alpha,2\beta)$, not similar to the tile,
 can be tiled by that
tile if and only if $s = 2\sin\frac{\alpha}{2}$ is rational.
\end{corollary}

\begin{proof}
If such a tiling exists, $s$ is rational by
Theorem~\ref{theorem:secondcasesolution}.  Conversely, if $s$ is rational,
let $(a,b,c)$ be the integer multiple of the tile with no common divisor.
Then $b = c - a^2/c$ is an integer, so $c$ divides $a^2$, and
Corollary~\ref{corollary:existence2case2} applies.
\end{proof}

\subsection{The scale factor $\mu$}

\begin{lemma} \label{lemma:mu}
Let $N$ and $ABC$ be given, and suppose $ABC$ is isosceles with 
vertex angle $\alpha$ at $B$ and base angles $\alpha+\beta$.  Suppose
$ABC$ is $N$-tiled by the tile with angles $(\alpha,\beta,\gamma)$ and 
integer sides $(a,b,c)$ with no common factor.  Let $\mu$ be
as in Definition~\ref{definition:scaling}.  Then $\mu = M(1+s)$ and 
$$\mu^2 = \frac {Nb} c.$$  
\end{lemma}

\begin{proof}
We have by Lemma~\ref{lemma:XY}
\begin{eqnarray*}
X &=& \mu c \\
Y &=& \mu a 
\end{eqnarray*}
Since twice the area of $ABC$ is $XY \sin(\alpha+\beta)$, 
and also $N ab \sin \gamma$, we have
\begin{eqnarray*}
XY \sin(\alpha+\beta) &=& N ab \sin \gamma
\end{eqnarray*}
Since $\sin(\alpha+\beta) = \sin \gamma$, we have
\begin{eqnarray*}
Nab &=& XY \\
Nab &=& \mu^2 ac \\
Nb &=& \mu^2 c \\
\mu^2 &=& Nb/c 
\end{eqnarray*}
\end{proof}

\begin{lemma} \label{lemma:halfintegermu}
Let $N$ and $ABC$ be given, and suppose $ABC$ is isosceles with 
vertex angle $\alpha$ at $B$ and base angles $\alpha+\beta$
with $3 \alpha + 2\beta = \pi$.  Suppose $ABC$ 
is $N$-tiled by a tile with angles $(\alpha, \beta, \gamma)$,
and integer sides $(a,b,c)$ with no common factor.  Let $k = gcd(a,c)$.
Then 
\smallskip

(i) The number of $b$ edges of tiles on side $AB$ is divisible by $k$, and
the same for side $BC$.
\smallskip

(ii) $\mu k$ is an integer, and 
$\mu$ is an integer if and only if $k | M$.
\smallskip

(iii) Let $q$ be the number of $b$ edges of tiles on the base $AC$,
and $M$ the coloring number of the tiling. 
Then $M \equiv -q$  mod $k$.
\end{lemma}

\begin{proof}  
To prove (i).  By Lemma~\ref{lemma:XY},  the length $X$ of side $AB$ is equal 
to $\mu c$, and the length $Y$ of side $AC$ is $\mu a$.  Since 
$AB$ and $AC$ are exactly matched by the edges of tiles in the tiling,
and the tile has integer sides, $X$ and $Y$ are integers.  That is,
both $\mu c$ and $\mu a$ are integers.  Hence 
the denominator of $\mu$ divides both $a$ and $c$.  Hence the 
denominator of $\mu$ divides $k = gcd(a,c)$.   Hence $\mu k$ is an integer.

In the tiling, $AB$ is made up of edges of tiles, 
so for some integers $p$, $d$, and $e$, we have 
\begin{eqnarray*}
\mu c &=& pa + db + ec \\
\mu k^2  &=& pa + db + ec \mbox{\qquad since $c = k^2$ by Lemma~\ref{lemma:gsquare}.}
\end{eqnarray*}
Since $\mu k$ is an integer,  $\mu k^2$ is an integer divisible by $k$.
Since $a$ and $c$ are also divisible by $k$, it follows that 
$db$ is divisible by $k$.  Since $b$ and $c$ are relatively prime,
also $b$ and $k$ are relatively prime.  Hence $d$ is divisible by 
$k$.  But $d$ is the number of $b$ edges on side $AB$.  That 
proves the first claim of the theorem.  Since $ABC$ is isosceles,
the same argument applies to side $BC$.   
That proves claim (i) of the theorem.

To prove (ii). 
From $\mu = M(1+s)$ we have $\mu c = M(c+a)$ and hence 
$\mu k^2 = M(c+a) = M(k^2 + a)$.  Dividing by $k^2$ we have 
$$\mu = M + \bigg(\frac M k\bigg) m.$$
Hence $\mu k$ is an integer, namely $M(k+m$), and $\mu$ is an integer
 if and only $k$ divides $M$.

To prove (iii).  The area equation is $X^2 = Nbc$.  Hence $c = k^2$ divides 
$X^2$.  Hence $k | X$.  By the coloring equation (Theorem~\ref{theorem:coloring}),
we have $M(a+b+c) = 2X-Y$.  Mod $k$ we have $X \equiv 0$ and $a \equiv 0$
and $c \equiv 0$.  Let $q$ be the number of $b$
edges on $AC$.  Then $Y \equiv qb$ mod $k$.  Then 
$M(a+b+c) = 2X-Y$ becomes, mod $k$, 
\begin{eqnarray*}
Mb & \equiv & -qb
\end{eqnarray*}
Since $(a,b,c)$ have no common factor, $b$ is relatively prime to $k$.
Therefore
\begin{eqnarray}
M &\equiv& -q \mod k \qedhere \label{eq:8413} 
\end{eqnarray} 
\end{proof}
 
\begin{remark}
Consider the case $N =27$, $M=3$, $\mu = 9/2$, $(a,b,c) =(2, 3, 4)$.
These values of $(N,M)$ satisfy all the equations above. 
In particular $k =2$ and $2M = 6 \equiv 0$ mod 2 but not $M \equiv 0$,
and $\mu$ is not an integer.  However, there does not exist
any tiling corresponding to these values, as will be discussed
below. 
\end{remark}

\begin{theorem} \label{theorem:N45}
Let isosceles triangle $ABC$ with base angles $\alpha+\beta$
 be $N$-tiled by a
triangle with angles $(\alpha, \beta, \gamma)$, with $3\alpha+2\beta = \pi$ 
and $\alpha$ not a rational multiple of $\pi$.
Then  $N \ge 45$.
\end{theorem}

\begin{proof}
Let $(a,b,c)$ be the integer multiple of the tile with no common factor.
Since $b = c - a^2/c$ is an integer, $c$ divides $a^2$, so
$(a,b,c) = (pq,\, q^2-p^2,\, q^2)$ with $0 < p < q$ and $\gcd(p,q) = 1$;
then $s = p/q$ and $k = \gcd(a,c) = q$.  Put $d = q-p$.
No side of $ABC$ is a single $c$ edge, so Lemma~\ref{lemma:twocs} applies
to each side: $X = M(p+q)q > q^2 = c$, and $Y = c$ would mean
$Mp(p+q) = q^2$, so $p = 1$ and $q+1$ divides $q^2$, which is impossible.

By Lemmas~\ref{lemma:XY} and~\ref{lemma:mu}, with the tile at $B$ black
(so that $M = \mu/(1+s) > 0$), the legs $AB$ and $BC$ have
length $X = \mu c = M(p+q)q$ and the base $AC$ has length
$Y = \mu a = Mp(p+q)$, and by Theorem~\ref{theorem:isosceles1},
\begin{eqnarray*}
N &=& \frac{M^2(q+p)}{q-p} \ = \ \frac{M^2(2p+d)}{d}.
\end{eqnarray*}

The tile at $B$ has its $\alpha$ angle at $B$, so it has a $b$ edge on
$AB$ or on $BC$.  By Lemma~\ref{lemma:halfintegermu}, the number of $b$
edges on that side is divisible by $k = q$, so it is at least $q$; by
Lemma~\ref{lemma:twocs}, that side also has at least two $c$ edges.  Hence
\begin{eqnarray*}
M(p+q)q \ = \ X &\ge& qb + 2c \ = \ q(q-p)(q+p) + 2q^2,
\end{eqnarray*}
so $M \ge d + 2q/(p+q)$.  Since $1 < 2q/(p+q) \le 2$, we have
$M \ge d+2$, and therefore
\begin{eqnarray*}
N &\ge& \frac{(d+2)^2(2p+d)}{d} \ \ge\ \frac{(d+2)^3}{d}.
\end{eqnarray*}
The right side is increasing in $d$ and equals $54$ at $d = 4$, so
$N \ge 45$ unless $d \le 3$.  We consider those cases.

If $d = 1$, then $N = M^2(2p+1) \ge 9(2p+1)$, which is at least $45$
unless $p = 1$.  Then the tile is $(2,3,4)$, $N = 3M^2$, and $N < 45$ only
for $M = 3$, $N = 27$.  In that case $Y = 9$.  By Lemma~\ref{lemma:twocs}
the base has at least two $c$ edges, of total length $8$, and the
remaining length $1$ is not a sum of edges of lengths $2$ and $3$.

If $d = 2$, then $p$ is odd and $N = M^2(p+1) \ge 16(p+1)$, which is at
least $45$ unless $p = 1$.  Then the tile is $(3,8,9)$, $N = 2M^2$, and
$N < 45$ only for $M = 4$, $N = 32$.  In that case $Y = 16 < 2c = 18$, so
the base cannot contain two $c$ edges, contradicting Lemma~\ref{lemma:twocs}.

If $d = 3$, then $N = M^2(2p+3)/3 \ge 25(2p+3)/3$, which exceeds $45$
unless $p = 1$.  Then $N = 5M^2/3$, so $3 \mid M$; with $M \ge 5$ this
gives $M \ge 6$ and $N \ge 60$.

Hence $N \ge 45$.
\end{proof}

\begin{remark}
Equality holds in every estimate for $d = 1$, $p = 2$, $M = 3$: the tile
$(6,5,9)$ with $N = 45$, the first entry of
Table~\ref{table:isoscelesalphaplusbeta}.  As in Theorem~11.11 of 
\cite{beeson-isosceles},
isosceles paper, it is not known whether a $45$-tiling exists.
\end{remark}

\subsection{Computational results for isosceles tilings with base angles $\alpha + \beta$}
Corollary~\ref{corollary:isoscelesall} guarantees tilings only for some,
in general very large, $N$.  To decide whether a given $N$ can occur,
we resort to computation.
Luckily, the theoretical considerations above are enough to 
reduce the question of existence of an $N$-tiling of some 
isosceles triangle with base angles $\alpha+\beta$ to computation.
Namely, given $N$, there are only finitely many possibilities for 
the coloring number $M$, and each of these determines a unique
tile $(a,b,c)$  and a particular isosceles triangle $ABC$.
Our theoretical considerations show that the question whether 
some $ABC$ can be $N$-tiled by some tile is equivalent to the 
question whether this {\em particular} $ABC$ can be tiled by 
this particular tile.

  It is a computationally 
decidable question whether a particular tile $(a,b,c)$ can be 
used to $N$-tile a particular triangle $ABC$. So, for each $N$,
we can computationally determine whether there is {\em any} 
$N$-tiling of {\em any} $ABC$.  That program is difficult to write
 and
for large $N$ it is slow.   We therefore develop a computational
``filter'' that checks various necessary conditions.  Then 
only the $N$ that pass these tests are serious candidates for 
further consideration.

Here we summarize the necessary conditions to be checked.

Starting with $N$,  there are finitely many squares $M^2$ less than $N$.
For each of these, the number $s$ is determined by 
$$  s = \frac {N-M^2}{N+M^2}.$$
Putting $s = m/k$ in lowest terms, the tile edge lengths are given by 
$$ (a,b,c) = (mk, k^2-m^2,k^2).$$
Then $\gcd(k,m) = 1$, so $(a,b,c)$ have no common factor.
  By the area equation (\ref{eq:4311}),  if there is a tiling, then
$Nbc = X^2$ is  automatically
  a square; indeed $X = \mu c = M(k+m)k$.  Therefore, if 
$Nbc$ is not a square, there is no $N$-tiling by this tile of an isosceles-$(\alpha+\beta)$ triangle.

  Then we can determine $\mu$ by 
the equation 
$$\mu^2 = Nb/c.$$
Then we can determine the lengths $X$ and $Y$ of the sides of 
$ABC$ by $X = \mu c$ and $Y = \mu a$, according to Lemma~\ref{lemma:XY}.
$X$ is an  integer since $\mu c = \sqrt{Nbc}$.  We have 
$\mu^2 a^2 = Nba^2/c$, which is an integer since $c$ divides $a^2$.
Then $\mu a$ is the rational square root of an integer, so $Y = \mu a$ is 
an integer.  

Next, we recall that in an $N$-tiling of $ABC$,

(i) There are at least two $c$ edges
on each of the three sides of $ABC$, by Lemma~\ref{lemma:twocs}; its
hypothesis~(v) holds, since no side of $ABC$ is a single $c$ edge, as shown
in the proof of Theorem~\ref{theorem:N45}.

(ii) Since the angle at the top vertex $B$ is $\alpha$, either 
side $AB$ or side $BC$ has at least one $b$ edge.

(iii) All the conditions listed in Lemmas~\ref{lemma:halfintegermu}
are satisfied.  In particular $k\mu$ is an integer.

Having determined $X$, $Y$, and $(a,b,c)$, it is now a computational
matter to check whether all the above conditions are met.

\subsection{Computational results for isosceles with base angles $\alpha + \beta$}

\begin{corollary} [to Lemma~\ref{lemma:halfintegermu}]
\label{lemma:july29}
For $N \le 108$, if there are $N$-tilings of an isosceles triangle
with base angles $\alpha+\beta$, then $N$, $M$, $\mu$, and $(a,b,c)$ 
occur in Table~\ref{table:isoscelesalphaplusbeta}. 
\end{corollary}

\begin{table}[ht]
\caption{Possible tilings of isosceles $ABC$ with base angles $\alpha+\beta$}
\label{table:isoscelesalphaplusbeta}
\begin{center}
\begin{tabular}{ccccc}
$N$ & $M$ & $\mu$  & $(a,b,c)$ & tiling exists \smallskip\\
\hline \\
45 & 3 &  5 & (6,5,9) & \nomark \cite{bonfioli2026} \\
48 & 4 &  6 & (2,3,4) & Figure~\ref{figure:48} \\
63 & 3 & 21/4 & (12,7,16) & \nomark \cite{bonfioli2026} \\
64 & 4 & 32/5 &(15,16,25) & \nomark \cite{bonfioli2026} \\
72 & 6 &  8 & (3,8,9) & \nomark \cite{bonfioli2026} \\
75 & 5 & 15/2 & (2,3,4) & Figure~\ref{figure:288} \\
80 & 4 & 20/3 & (6,5,9) & \nomark \cite{bonfioli2026} \\
81 & 3 & 27/5 &(20,9,25) & \nomark \cite{bonfioli2026} \\
90 & 6 & 60/7 &(21,40,49) & \nomark \cite{bonfioli2026} \\
96 & 4 & 48/7 &(35,24,49) & \nomark \cite{bonfioli2026} \\
99 & 3 & 11/2 & (30,11,36) & \nomark \cite{bonfioli2026} \\
100 & 5 &  8 & (15,16,25) & \nomark \cite{bonfioli2026} \\
108 & 6 &  9 & (2,3,4) &  Figure~\ref{figure:48} \\
\end{tabular}
\end{center}
\end{table}%

Bonfioli \cite{bonfioli2026} proves, with Lean proofs, two necessary conditions
for a tiling of the isosceles triangle with base angles $\alpha+\beta$ by the tile
$(mk,\,k^2-m^2,\,k^2)$.  First, $N = t^2 b$ for an integer $t$; then the triangle
has legs $tkb$ and base $tmb$.  Second, $tm \ge 2k$; moreover $t \ge 2k$ if $m = 1$,
and $t \ge k+2$ if $m = k-1$.%
\footnote{Theorems {\tt CorAtomDivFull.noAtomTiling}, {\tt ThetaKill.theta\_te\_ge},
{\tt ThetaKill.theta\_e1} and {\tt ThetaKill.theta\_thick} in the repository
\url{https://github.com/ElVec1o/erdos_634_proof}, commit {\tt 119bec0}.
We have not checked them in Lean.}
The first condition excludes $N = 90$ and $100$; the second excludes the other
rows marked \nomark.  The $48$-tiling, with $t = 4 = 2k$, shows that the bound
is attained.  Bonfolio reports finding a $75$-tiling ($t = 5$), but we could
not find it in his repository.  The state of the art in November 2026 is such that
when Claude failed to find it, Claude just wrote a Python program to search 
for a tiling,  and found the $75$-tiling in Figure~\ref{figure:288} 
and checked it in exact arithmetic (but not in Lean).%
\footnote{The program to find the tiling is 159 lines of Python code and 
ran in under five minutes for $N=75$.}


\begin{remark}Hence
an isosceles triangle with base angles $\alpha+\beta$ can be $N$-tiled only
if $N \ge 45$.  This is the same bound as \cite[Theorem~11.11]{beeson-isosceles},
 for isosceles $ABC$ with base angles $\alpha$ tiled by
$(\alpha,\beta,2\alpha)$, and for the same reason.  Write the tile here
as $(pq, q^2-p^2, q^2)$.  By Lemmas~\ref{lemma:XY} and~\ref{lemma:mu} the
legs have length $X = qt$ with $t^2 = N(q^2-p^2)$.  For the tile
$(k^2, m^2-k^2, km)$ there, $X = kt$ with $t^2 = N(m^2-k^2)$.  In both
cases $N$ times a difference of two coprime squares must be a square, and
the first case that also meets the conditions on the edges is
$3^2-2^2 = 5$, $N = 45$, $t = 15$: the tile $(6,5,9)$ here and $(4,5,6)$
there.  
\end{remark}

\FloatBarrier



\section{The case when $ABC$ is isosceles with base angles $\alpha$}

\begin{figure}[ht]
\includegraphics[width=0.8\textwidth]{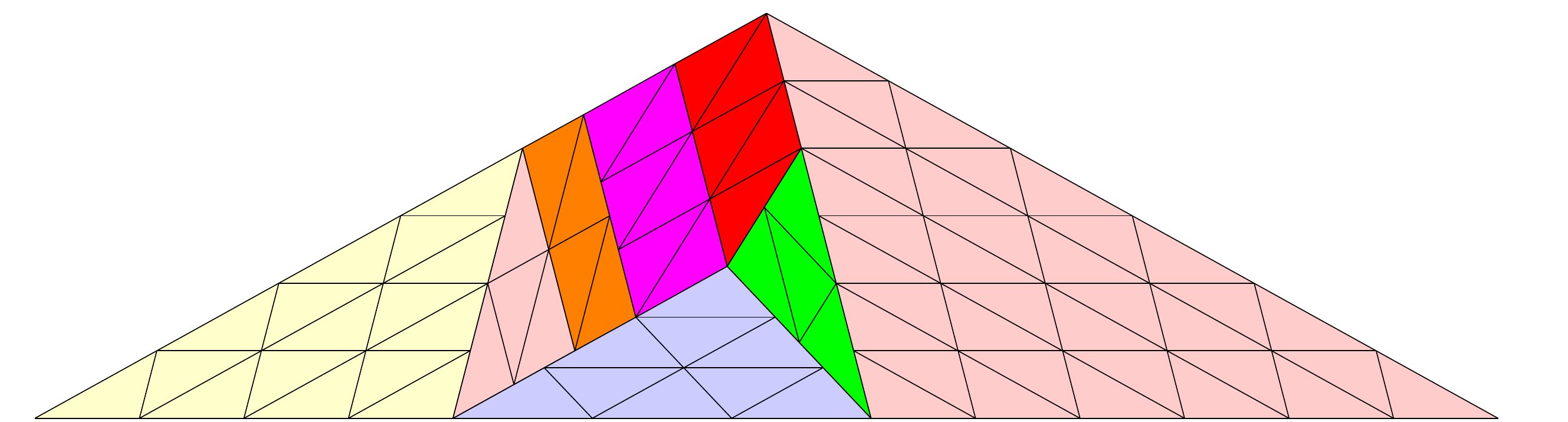} 
\caption{$N = 84$. $ABC$ is isosceles with base angles $\alpha$. $(a,b,c) = (2,3,4)$}
\label{figure:84}
\end{figure}

\begin{figure}[ht]
\includegraphics[width=0.8\textwidth]{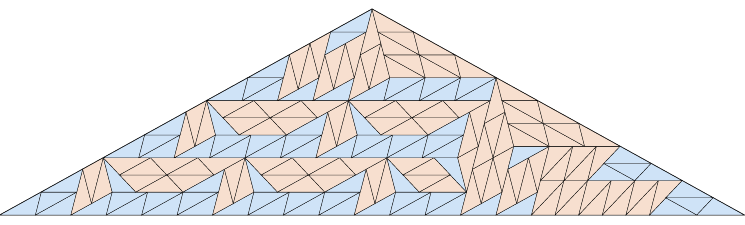} 
\caption{$N = 189$. $(a,b,c) = (2,3,4)$.  This tiling was found by Jan Philip Harries
in July, 2026.}
\label{figure:189}
\end{figure}

\begin{theorem}[Tiling equation for base angles $\alpha$] \label{theorem:tilingequation4}
Suppose $3\alpha + 2\beta = \pi$. 
Let $ABC$ be isosceles with base angles $\alpha$.
Suppose $ABC$ is $N$-tiled by the tile whose 
integer sides $(a,b,c)$ have no common factor, and
whose angles are $(\alpha,\beta,\gamma)$.  
Suppose $ABC$ is not similar to the tile.
Let $s = a/c$, and let $M$ be the coloring number
of the tiling.  

(i)  Then  $M^2 \le 2N$,  and $(M,s)$
is a solution of the ``isosceles-$\alpha$ tiling equation''
\begin{eqnarray*}
 N &=& M^2 \frac {(1+s)(2-s^2)}{(1-s)(2+s)^2}
\end{eqnarray*} 

(ii) If $(a,b,c) = (mk,k^2-m^2,k^2)$ with $k$ and $m$ coprime,
then $2k+m \mid M$.  Setting $j = M/(2k+m)$, the tiling equation can be written
\begin{eqnarray*}
N (k-m) &=& j^2 (m+k)(2k^2-m^2),  
\end{eqnarray*}
in which all the factors are integers. 
\end{theorem}

\begin{proof}  
Consider tilings of a triangle $ABC$ with base angles $\alpha$. 
Then the vertex angle is $\alpha + 2\beta$.  
We have 
\begin{eqnarray*}
\sin(\alpha + 2\beta) &=& \sin (\pi - 2\alpha) \\
&=& \sin 2\alpha \\
&=& (2-s^2) \sin \alpha \mbox{\qquad by Lemma~\ref{lemma:singamma2}}
\end{eqnarray*}

Twice the area of $ABC$ is $X^2 \sin (\alpha + 2\beta)$.  Hence
\begin{eqnarray}
N bc \sin \alpha &=& X^2 \sin(\alpha + 2\beta) \\
   &=& X^2 (2-s^2) \sin \alpha \nonumber \\
 N bc &=& X^2 (2-s^2)  \label{eq:9439}
\end{eqnarray}
Let $\mu = X/b$, so that $X = \mu b$.  Then
\begin{eqnarray}
Nbc &=& \mu^2 b^2 (2-s^2)\nonumber \\
N &=& \mu^2 \frac b c (2-s^2) \nonumber \\
N &=& \mu^2 (1-s^2)(2-s^2) \label{eq:7781}
\end{eqnarray}

Since $X = \mu b$, 
by the law of sines we have 
\begin{eqnarray*}
Y &=&  X \frac {\sin(\alpha + 2\beta)}{ \sin \alpha} \\
  &=&  X \frac {(2-s^2) \sin \alpha}{\sin \alpha} \\
  &=&  X (2-s^2) \\
  &=&  \mu b (2-s^2) 
\end{eqnarray*}

The coloring equation is 
\begin{eqnarray*}
M(a+b+c) &=& 2X+Y  \\    
&=& X (2 + 2-s^2) \\
&=& X (4-s^2) \\
&=& \mu b (4-s^2)
\end{eqnarray*}
Dividing by $c$ and expressing everything in terms of $s$,
\begin{eqnarray*}
M(s+(1-s^2) + 1) &=& \mu(1-s^2)(4-s^2) \\
M(2-s)(1+s) &=& \mu (1-s)(1+s)(2-s)(2+s) \\
M &=& \mu (1-s)(2+s) 
\end{eqnarray*}
Hence
\begin{eqnarray*}
\mu &=&  \frac { M}{(1-s)(2+s)}
\end{eqnarray*}
We can square that equation and substitute for $\mu^2$ in (\ref{eq:7781})
to obtain the area equation for this shape of $ABC$.
\begin{eqnarray}
 N &=& \frac {M^2}{(1-s)^2(2+s)^2} (1-s^2)(2-s^2) \nonumber \\
N &=& M^2 \frac {(1+s)(2-s^2)}{(1-s)(2+s)^2}\label{eq:9484}
\end{eqnarray}
Finally, $M^2 \le 2N$, since the right side of (\ref{eq:9484}) is
$M^2 f(s)$ with $f(s) = (1+s)(2-s^2)/((1-s)(2+s)^2)$ increasing on
$(0,1)$ and $f(0) = 1/2$.
That completes the proof of part~(i).

Multiplying by the denominator and using $s = a/c = m/k$, we have 
\begin{eqnarray}
N &=& M^2 \frac{(m+k)(2k^2-m^2)}{(k-m)(2k+m)^2} \label{eq:9565}
\end{eqnarray}
The leg $X$ of the triangle $ABC$ is given by the area equation,
\begin{eqnarray*}
X^2 \sin (\alpha + 2 \beta) = N bc \sin \alpha   && \\
X^2 (2-s^2) \sin \alpha = N bc \sin \alpha &&  \mbox{\qquad by Lemma~\ref{lemma:singamma2}}\\
X^2 = \frac{ N bc}{2-s^2} &&  \mbox{\qquad by Lemma~\ref{lemma:singamma2}}\\
(a,b,c) = (km,k^2-m^2,k^2) && \mbox{\qquad by Lemmas~\ref{lemma:gsquare} and~\ref{lemma:anglesOK}}\\
X^2 = N\frac {(k^2-m^2)k^4}{2k^2-m^2}  && \mbox{\qquad by the preceding lines}
\end{eqnarray*}
Now substitute for $N$ from (\ref{eq:9565}).
\begin{eqnarray*}
X^2 =  M^2 \frac{(m+k)(2k^2-m^2)}{(k-m)(2k+m)^2}\frac {(k^2-m^2)k^4}{(2k^2-m^2)}  
\end{eqnarray*}
Cancelling $2k^2-m^2$ and using $(k^2-m^2)/(k-m) = k+m$, we have
\begin{eqnarray*}
X^2 &=& M^2 \frac{(m+k)^2 k^4}{(2k+m)^2}   \\
X &=& M \frac {(m+k)k^2}{2k+m}
\end{eqnarray*}
Since $k$ and $m$ are coprime, $2k+m$ is relatively prime to both $m+k$ and $k$.
$X$ is an integer, since $AB$ is composed of tile edges of integral length.  Therefore
$2k+m \mid M$.  Then by (\ref{eq:9565}),
\begin{eqnarray*}
N (k-m) &=& \left( \frac M {2k+m}\right)^2 (m+k)(2k^2-m^2), 
\end{eqnarray*}
in which all the factors are integers. 
\end{proof}

\subsection{$N$ is not prime when $ABC$ is isosceles with base angles $\alpha$.}

\begin{theorem} \label{theorem:notprime5}
 Suppose $3\alpha + 2\beta = \pi$. 
Let $ABC$ be isosceles with base angles $\alpha$.
Suppose $ABC$ is $N$-tiled by a tile 
whose angles are $(\alpha,\beta,\gamma)$.  
Suppose $ABC$ is not similar to the tile.
Then 
$N$ is not a prime number.  
\end{theorem}

\begin{proof}  By Theorem~\ref{theorem:rationaltile},
the tile is rational.  Therefore we may assume without loss of generality
that the tile has
integer sides $(a,b,c)$ with no common factor.  Let $s = a/c$ and let
$a/c = m/k$ in lowest terms, and let $M$ be the coloring number of the
tiling.  According to Theorem~\ref{theorem:tilingequation4},
\begin{eqnarray*}
N (k-m) &=& j^2 (m+k)(2k^2-m^2), \qquad j = \frac M {2k+m},
\end{eqnarray*}
in which all the factors are integers.

We first show that $2k^2-m^2$ is relatively prime to $k-m$.
Since $$2k^2-m^2 = 2(k-m)(k+m) + m^2,$$ any common divisor of $2k^2-m^2$
and $k-m$ divides $m^2$.  But $\gcd(m,k-m) = \gcd(m,k) = 1$, so
$\gcd(m^2,k-m)=1$.  Hence $\gcd(2k^2-m^2,\,k-m) = 1$, and therefore
$2k^2-m^2$ divides $N$.

Moreover $N > 2k^2-m^2$, since $m+k > k-m$ and $j \ge 1$ give
$$ N \ = \ (2k^2-m^2)\,\frac{j^2(m+k)}{k-m} \ > \ 2k^2-m^2 . $$
Finally $2k^2 - m^2 > k^2 \ge 4$, since $m < k$ and $k \ge 2$.  So $N$ has
a divisor strictly between $1$ and $N$, and is not prime.
\end{proof}
 
\subsection{Existence of tilings of isosceles $ABC$ with base angles $\alpha$}

\begin{figure}[ht]
\begin{center}
\includegraphics[width=0.8\textwidth]{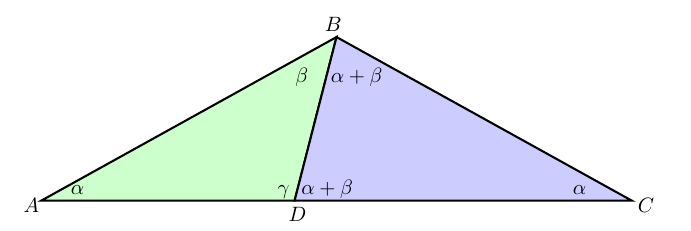}
\end{center}
\caption{Decomposing isosceles $ABC$ with base angles $\alpha$}
\label{figure:alphacut}
\end{figure}

Corollary~\ref{corollary:isoscelesall} also provides tilings of isosceles
triangles with base angles $\alpha$, by means of the dissection shown in
Figure~\ref{figure:alphacut}.  


\begin{corollary} \label{corollary:isoscelesalphaall}
Let the tile have angles $(\alpha,\beta,\gamma)$ with $3\alpha+2\beta=\pi$
and $s$ rational, and let $(a,b,c)$ be the integer multiple of the tile
with no common factor.  Then for some $N$, some isosceles triangle with
base angles $\alpha$ can be $N$-tiled by the tile.  Specifically, if some
isosceles triangle with base angles $\alpha+\beta$ and legs of length $L$
can be $N_1$-tiled, then $N = L^2 + c^2N_1$ tiles suffice.
\end{corollary}

\begin{proof}
By Corollary~\ref{corollary:isoscelesall} there is a tiling of some
isosceles triangle with base angles $\alpha+\beta$; let it have $N_1$
tiles and legs of length $L$.  Replacing each of its tiles by a quadratic
tiling of scale $c$ gives a tiling, with $c^2N_1$ tiles, of the isosceles
triangle with base angles $\alpha+\beta$ and legs $Lc$.  Now take
$X = Lc$ in the dissection of Figure~\ref{figure:alphacut}, whose blue
triangle $DBC$ has legs $|BC| = |DC| = X$.  Then $DBC$ is that triangle,
and $ABD$, whose side $AB$ of length $X$ lies opposite its angle
$\gamma$, is similar to the tile with ratio $X/c = L$, so it has a
quadratic tiling with $L^2$ tiles.  Together these tile $ABC$, using
$L^2 + c^2N_1$ tiles.\end{proof}

\begin{example}
For the tile $(2,3,4)$, the $48$-tiling of Figure~\ref{figure:48} has legs
$L = \mu c = 24$, so the corollary gives $N = 24^2 + 4^2\cdot 48 = 1344$,
with $M = 40$.  Here $1344 = 16 \cdot 84$, so this tiling is much larger
than the $84$-tiling of Figure~\ref{figure:84}.  As with the other shapes,
the construction settles existence but does not produce the smallest $N$.
\end{example}

\subsection{Computational results for isosceles with base angles $\alpha$}

All solutions of the tiling equation in Theorem~\ref{theorem:tilingequation4}
with $N \le 500$,
are listed in the following table:

\begin{center}
\begin{tabular}{lrrr}
$N$ & $M$ & $(a,b,c)$ & tiling exists \\
\hline 
21 & 5 & (2, 3, 4) & \nomark \cite{bonfioli2026} \\
34 & 7 & (3, 8, 9) & \nomark \cite{bonfioli2026} \\
70 & 8 & (6, 5, 9) & \nomark \cite{bonfioli2026} \\
84 & 10 & (2, 3, 4) & \yesmark Figure~\ref{figure:84} \\
$136=2^2\cdot 34$ & 14 & (3, 8, 9) & ? \\
161 & 11 & (12, 7, 16) & \nomark \cite{bonfioli2026} \\
164 & 13 & (15, 16, 25) & \nomark \cite{bonfioli2026} \\
$189 = 3^2\cdot 21$& 15 & (2, 3, 4) & \yesmark Figure~\ref{figure:189}  \\
$280 = 2^2 \cdot 70$ & 16 & (6, 5, 9) & ? \\
294 & 22 & (5, 24, 25) & \nomark \cite{bonfioli2026} \\
$306 = 3^2 \cdot 34$ & 14 & (20, 9, 25) & ? \\
306 & 21 & (3, 8, 9) & \nomark \cite{bonfioli2026} \\
$336 = 2^2\cdot 84$ & 20 & (2, 3, 4) & \yesmark Figure~\ref{figure:84} \\
438 & 19 & (35, 24, 49) & \nomark \cite{bonfioli2026} \\
465 & 27 & (4, 15, 16) & ? \\
\end{tabular}
\end{center}

Bonfioli \cite{bonfioli2026} proved, with a Lean proof, that if an isosceles
triangle with base angles $\alpha$ is $N$-tiled by the tile
$(a,b,c) = (mk,\,k^2-m^2,\,k^2)$, then $N = r^2\, b\,(b+c)$ for some integer $r$.%
\footnote{Theorem {\tt ARA.theoremA\_N\_form\_full}, with the instances
{\tt BADCorFull.noTiling\_34}, {\tt noTiling\_164}, {\tt noTiling\_294},
{\tt noTiling\_438}; also {\tt Refl.MB70.not\_dissection} for $N=70$ and
{\tt N161F.not\_tile\_161} for $N=161$, in the repository
\url{https://github.com/ElVec1o/erdos_634_proof}, commit {\tt 119bec0}.
We have not checked them in Lean.}
This excludes $N = 34$, $164$, $294$, $438$, and $306$ with the tile $(3,8,9)$.
For $N = 70$ he checked this triangle and tile by an exhaustive search
certified in Lean, and $161$ is not a tile count at all.
Bonfioli's condition is equivalent to $(k-m)(2k+m) \mid M$ via
Theorem~\ref{theorem:tilingequation4}, part~(ii).  
As of November 2026, we have no example to refute the conjecture
that there is an $N$-tiling of some isosceles $ABC$ with base angles
$\alpha$ if and only if $N = j^2 b(b+c)$ with $j \ge 2$.  This
is really two conjectures: that $j=1$ is impossible, and that
there are tilings whenever $j \ge 2$.   Both directions are 
open questions.  $N=136$ is a test question with $j=1$.

\section{Searching for tilings}
\label{section:search}
In preceding sections we have given necessary equations, that must 
be solvable if there is an $N$-tiling, one equation for each possible 
shape of $ABC$.   In two of the cases we were able to give matching 
sufficient conditions that gave a complete solution of the question as
to which $ABC$ have an $N$-tiling. In other cases we could not find 
such sufficient conditions, thus leaving the question unsolved for $N$
such that the equations are satisfied but no tiling is known.

Our equations, however, do show that if $N$ is given,  then the possible 
shapes of $ABC$ and the tile are restricted to a finite set, easily 
computable from $N$.  Therefore the question, for a particular $N$,
whether an $N$-tiling of some $ABC$ by some tile exists,  reduces to 
the question whether an $N$-tiling of a particular $ABC$ by a particular 
tile exists.  That question is in turn, in principle,  decidable by 
a trial-and-error search; it amounts to assembling $N$ tiles as a 
jigsaw puzzle into a frame with the shape of $ABC$.

This is a standard search problem, a staple of undergraduate computer 
science, and in 2011 we wrote a standard program to solve it, using C++ for 
efficiency.  We also used some custom touches for efficiency: we searched
by first placing tiles along the boundary (a standard technique in 
human jigsaw-puzzle solving), and if a tile touched a previously-laid
tile, we computed the area of the ``hole'' created, rejecting the tiling
if the hole has an area that is not a multiple of the tile's area.

As an example, we mention some of the results searching for tilings
of isosceles $ABC$ with base angles $\alpha+\beta$. 
For the value $N=12$ it draws 228 pictures, showing each ``blocked''
partial boundary tiling that cannot be continued.   
 $N = 12$ and $N=18$ were ruled out, 
which at the time was new information, since 
Lemma~\ref{lemma:twocs} asserting that there are at least two 
$c$ edges on each side of $ABC$ had not yet been proved, so there were 
four more entries at the beginning of Table~\ref{table:isoscelesalphaplusbeta}.
For isosceles $ABC$ with base angles $\alpha$, there still are two
 small values of $N$ that 
can be eliminated by exhaustive search for a boundary tiling.  The 
largest value of $N$ for which we could complete an exhaustive search
(for any shape of $ABC$) was $N=23$, although the successful search 
in the case $N=28$ led to the discovery of the triquadratic tilings.
For $N=45$, $72$, and $75$, we found thousands of boundary tilings,
but they could not be completed to tilings of $ABC$.  
 We ran this program (in 2011) for $N=99$ for 
more than 48 hours.  Given these disappointing results, we did not 
extend our program to search for a 
completion of a boundary tiling to a full tiling.

Fifteen years later,  in July  2026,  several people asked various AIs  to write such search
programs. 
Specifically 
\begin{itemize}
\item Epameinondas Moustakas reported (personal communication) that there is no $60$-tiling of 
an equilateral $ABC$.%
\footnote{The angles of such a tiling would have $\gamma = \pi/3$ or
$2\pi/3$, and neither is compatible with $3\alpha+2\beta=\pi$: the first
forces $\alpha = -\pi/3$ and the second forces $\beta = 0$.  So the question
is not relevant to this paper.}
\item Vico Bonfioli \cite{bonfioli2026} rediscovered a $44$-tiling by $(2,3,4)$, which 
has been available in earlier versions of this paper since 2018.
But he found it by exhaustive computer search, which is new,  and it is not
exactly the same tiling as was already known. 
His search also found three tilings that were not known before: a
$63$-tiling of the triangle $(2\alpha,\beta,\alpha+\beta)$ with sides $(21,24,18)$, and two
tilings of isosceles triangles with base angles $\beta$, all three by the tile 
$(2,3,4)$.  Those two isosceles tilings were a $99$-tiling of
the isosceles triangle with sides $(24,24,33)$,  and a $176$-tiling of $(32,32,44)$.
These three tilings are shown in Figures~\ref{figure:63},~\ref{figure:99}, and ~\ref{figure:176}.
On 1 October 2026, two days after Paliy's direct construction of a 322-tiling, Bonfioli found a 
completely different 322-tiling by computer search.  It is shown in Figure~\ref{figure:322}.
  Using them as starting cases, he constructs
$7m^2$-tilings of $(2\alpha,\beta,\alpha+\beta)$ and $11m^2$-tilings of the isosceles triangle with base
angles $\beta$ for every $m \ge 2$.  He also reports that there is no $60$-tiling of an equilateral
$ABC$, and he reports exhaustive searches excluding many other values of $N$.

\begin{figure}
\centering
\includegraphics[width=0.5\textwidth]{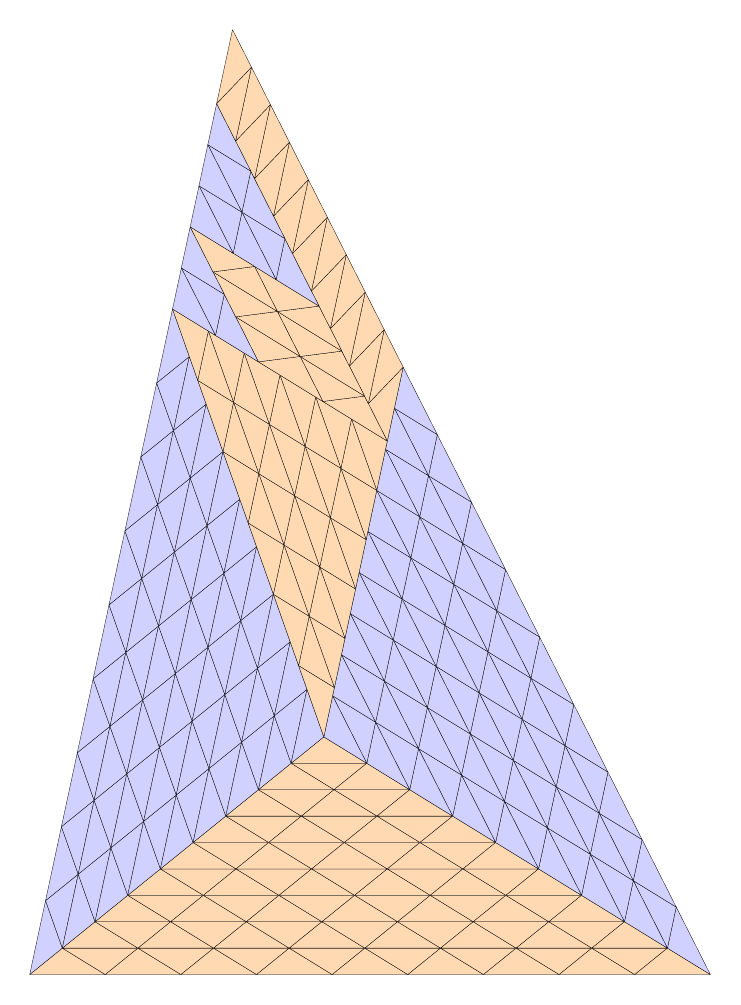}
\caption{$N=322$, $(a,b,c) = (6,5,9)$, $ABC$ of shape $(2\alpha, \alpha, 2\beta)$.
 Vico Bonfioli found this tiling by computer search.}
\label{figure:322}
\end{figure}

\begin{figure}[ht]
\includegraphics[width=0.8\textwidth]{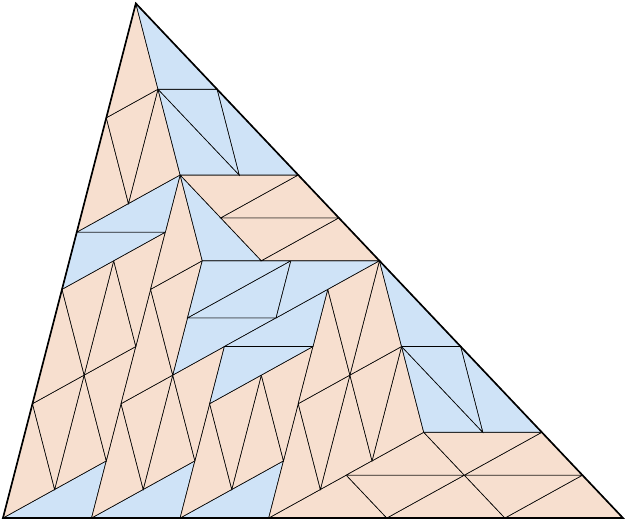}
\caption{$N = 63$. $ABC$ is $(2\alpha,\beta,\alpha+\beta)$ with sides $(21,24,18)$,
$(a,b,c) = (2,3,4)$.   Vico Bonfioli found this tiling by computer search.
}
\label{figure:63}
\end{figure}

If a complicated search program reports a negative result (no tiling found),  one may 
normally be suspicious that the program might have a bug and miss a tiling.  Bonfolio
has solved that problem of trust:  his search program produces not just the word {\em no},
but a complete search tree, containing information describing the entire search.  This
tree can be checked in Lean, producing a Lean-certified proof, for example, ``there is no 
$70$-tiling.''  Anyone can check these proofs;  at that point the details of the 
searching program no longer matter.

\item {Jan Philipp Harries}\cite{harries2026} discovered an $88$-tiling by $(3,5,7)$
and a $189$-tiling by (2,3,4).   His $189$-tiling can be viewed in Figure~\ref{figure:189} 
of this paper.  He also reports the exclusion of $21, 30, 33, 46,
55, 56, 66, 70$,  and $438$. 
\end{itemize}

\section{Conclusions}
In this paper we assumed $3\alpha + 2\beta = \pi$ and studied
$N$-tilings of triangles $ABC$ by a tile with angles $(\alpha,\beta,\gamma)$.
It follows from results of Laczkovich (see Lemma~\ref{lemma:angles}) that
the angles cannot be rational multiples of $\pi$, so there are
only five possible shapes of $ABC$, given $(\alpha,\beta,\gamma)$.
 We showed that the tile must be rational, and by means of the 
area equation and coloring equation, we derived a ``tiling equation''
for each of the five possible shapes of $ABC$.  By means of these
tiling equations, we reduced the existentially quantified question 
whether there is an $N$-tiling of {\em some} $ABC$ by {\em some} tile,
to the specific question of whether there is an $N$-tiling of a 
specific, particular $ABC$ by a specific, particular tile (or tiles)
depending on $N$.  The tiling equations supply necessary conditions 
by which we can rule out a great many values of $N$.  In particular,
we were able to rule out prime values of $N$ for three of the five shapes:

\begin{theorem} 
\label{theorem:notprime}
Let $ABC$ be $N$-tiled by a tile with angles
$(\alpha,\beta,\gamma)$ such that $3\alpha+2\beta = \pi$, and 
$ABC$ is not similar to the tile.  Then $N$ is not prime,
 except possibly for the shapes $(2\alpha, \beta, \alpha+\beta)$
 and $(\beta, 3\alpha,\beta)$.
\end{theorem}

\begin{proof}
By \cite{laczkovich1995} (as detailed in Lemma~\ref{lemma:shapes})
 $ABC$ has one of the five shapes
considered in this paper.  For three of those shapes, we have proved
that $N$ is not prime, in Theorems~\ref{theorem:notprime2}, \ref{theorem:notprime3},
and \ref{theorem:notprime5}.
\end{proof}

\begin{remark} The two shapes $(\beta,3\alpha,\beta)$ and  $(2\alpha,\beta,\alpha+\beta)$ 
that are not settled here have been settled in \cite{beeson-rosu2026}.  They also have
no prime tilings.  
\end{remark}

We list the five shapes and the corresponding necessary tiling
equations in Table~\ref{table:summary}.  In each case, there are 
finitely many possible values of $s$, determined by $N$ and the 
coloring number $M$, which has to be less than $N$, so that for each 
particular $N$ it can be checked by computation, and very efficiently,
whether the equations have a solution or not.  For each of the five shapes
we can give (extra) conditions on the tile that are
sufficient for some triangle of that shape to be tiled, for some $N$;
these are also listed in the table.  In the first row, the condition
$K \mid M^2$ can always be arranged by replacing $(M,K)$ by $(tM,tK)$
with $t = K$, which multiplies $N$ by $K^2$; and in the second row,
$c \mid a^2$ holds for the integer tile with no common factor.  So for
every tile, some triangle of each of the five shapes can be tiled.
These conditions do not in general guarantee
a tiling for every $N$ that satisfies the necessary condition; for
example, $N = 322$ satisfies the second tiling equation with
$s = 2/3$, but whether there is a $322$-tiling is unknown.

\begin{table}[ht]
\caption{Tiling equations when $3\alpha+2\beta = \pi$, with $s = a/c$}
\label{table:summary}
\begin{center}
\begin{tabular}{lrr}
$ABC$ & Necessary &    Sufficient  \smallskip \\
\hline \\
$(2\alpha,\beta,\alpha+\beta)$ & $\displaystyle \frac N{M^2} = \frac{2}{s^2}-1$ &  $s = M/K$ and $K \mid M^2$  \\[10pt]
$(2 \alpha, \alpha, 2\beta)$ & $\displaystyle\frac {N} {M^2 } = \frac{(2-s^2)(3-s^2)}{(1-s)^2(2+s)^2} $ & $c | a^2 $\\[10pt] \\
$(2 \alpha, \alpha, 2\beta)$ & $N = j^2(b+c)(b+2c)$,  &  none\\
& $j = M/(b+c-a)$  & \\
$(\beta,3\alpha,\beta)$ & $\displaystyle\frac {N}{M^2} = \frac{3-s^2}{(1+s)^2}$& $a+c \mid M$ \\[10pt]
$(\alpha+\beta,\alpha+\beta,\alpha)$ & $\displaystyle \frac {N}{M^2} = \frac{1+s}{1-s}$  & $s$ rational  \\[10pt]
$(\alpha,\alpha,\alpha+2\beta)$ 
& $\displaystyle\frac {N}{M^2} = \frac{(1+s)(2-s^2)}{(1-s)(2+s)^2} $&  $s$ rational  
\end{tabular}
\end{center}
\end{table}

Table~\ref{table:numerical_summary} lists the known $N$-tilings, 
for $N \le 100$  and $3\alpha+2\beta = \pi$.  Those tilings, and others for 
larger $N$, are illustrated elsewhere in this paper.  In that
table there is one tiling for each of the five shapes of $ABC$. 
\begin{table}[ht]
\caption{ $N$-tilings with $3\alpha + 2\beta = \pi$ for $N\le 100$}
\label{table:numerical_summary}
\begin{center}
\begin{tabular}{rrrr}
$N$  &       $M$   &    $(a,b,c)$  & $(A,B,C)$  \\
  \hline
28 & 2 & (2, 3, 4) & triquadratic \\
44 & 6 & (2, 3, 4) & isosceles-$\beta$ \\
48 & 4 & (2, 3, 4) & isosceles-$\alpha +\beta$ \\
63 & 3 & (2, 3, 4) & triquadratic \\
77 & 5 & (2, 3, 4) & $(2\alpha,\alpha,2\beta)$ \\ 
84 & 10 & (2, 3, 4) &  isosceles-$\alpha$ \\
\end{tabular}
\end{center}
\end{table}

Using the equations in Table~\ref{table:summary}, we can list, for each
shape, the values of $N$ for which the equations have solutions, and
therefore for which an $N$-tiling is not ruled out. 
 By \cite{beeson-rosu2026},
there are no tilings with $N$ prime and $3\alpha + 2\beta = \pi$.
  For $N \le 100$ the
undecided values are listed in Table~\ref{table:open}.  There are $36$
entries, but only $34$ distinct values, since $71$ and $92$ each occur
for two shapes.

\begin{table}[ht]
\caption{Values of $N \le 100$ for which the existence of a tiling is
undecided}
\label{table:open}
\begin{center}
\begin{tabular}{lp{0.55\textwidth}}
$ABC$ & undecided $N \le 100$ \\
\hline \\
$(2\alpha,\beta,\alpha+\beta)$ & $46$, $49$, $56$, 
   $82$,  $92$, $94$ \\[4pt]
$(2\alpha,\alpha,2\beta)$ & none \\[4pt]
$(\beta,3\alpha,\beta)$ &  $26$, $39$,  $66$,
   $74$,  $92$, $99$ \\[4pt]
$(\alpha+\beta,\alpha,\alpha+\beta)$ & $45$, $64$, $72$, $81$, $90$, $96$,
   $100$ \\[4pt]
$(\alpha,\alpha+2\beta,\alpha)$ & $21$, $34$, $70$
\end{tabular}
\end{center}
\end{table}

The entries for $(2\alpha,\beta,\alpha+\beta)$ use only the conditions
proved here, namely that when $\gcd(M,K) = 1$ there are no $b$ edges on
$AC$ (Corollary~\ref{corollary:nokt}) and that there are at least two $c$
edges there (Lemma~\ref{lemma:twocsapplies}), together with the special
arguments for $N = 7$, $14$, $17$, and $23$.  The entries for
$(\beta,3\alpha,\beta)$ are every solution of the equation for that shape
except $N = 44$, for which a tiling exists.  Most likely a more
systematic search would decide some of these values, as it did for the
four just mentioned.

For $N = 21$, with $ABC$ isosceles with base angles $\alpha$, nonexistence
was reported by Bonfioli \cite{bonfioli2026} as the result of an
exhaustive computer search.  We neither assert nor contradict that result,
and accordingly continue to list $21$ as undecided.

For each $N$, the tiling equations restrict the possible shape of the
tile to a finite number (usually just one) of specific tiles, and for
each tile they determine the triangle $ABC$ up to similarity and the
lengths of its sides in terms of the sides of the tile.  Therefore the
existence or non-existence of a tiling can in principle be decided by
computation: there are only finitely many ways to place the tiles, since
each side of $ABC$ is a sum of tile edges.  In practice such a search is
feasible only for small $N$, which is why the values in
Table~\ref{table:open} remain undecided.

\appendix
\section{Formal Verification in Lean}\label{section:lean}

The results of this paper have been formalized in the Lean~4 proof assistant
(version~4.30.0), using the Mathlib library.  The Lean code is at

\url{https://github.com/mjbeeson/TriangleTiling3Lean}, release v1.0.

Each theorem, lemma and corollary listed in Table~\ref{table:lean} has a Lean
counterpart, named after its label, which states the result as printed here and
whose proof follows the printed proof.  The proofs are checked by the Lean kernel using
only the standard axioms ({\tt propext}, {\tt Classical.choice},
{\tt Quot.sound}), with no {\tt sorry} and no {\tt native\_decide}.
All Lean names are in the namespace {\tt TriangleTiling.ThreeAlpha}; the
isosceles paper's results are in {\tt TriangleTiling}, so equal labels in the two
papers (for example {\tt N45} and {\tt notprime}) name different Lean theorems.
Tilings are formalized as in V.~Bonfioli's Lean library \cite{bonfioli-lean}:
$N$ triangles congruent to a given tile (reflections allowed), with disjoint
interiors, whose union is $ABC$.  The coloring number $M$ is the sum of a proper
coloring of the tiles by $\pm 1$, with the tile at $B$ that has an edge on $AB$
colored $+1$.

\smallskip
\noindent{\em What is not formalized.}
Only Corollary~\ref{lemma:july29}, the table of possible tilings for $N \le 108$,
which rests on a computer search.  Every condition the search uses is
formalized.

\smallskip
\noindent{\em Results taken from Bonfioli's library.}
Laczkovich's classification of the tilings of a triangle by congruent tiles
(Theorems~4.1, 5.1 and 5.3 of \cite{laczkovich1995} and Theorem~3.3 of
\cite{laczkovich2012}), used for Lemmas~\ref{lemma:angles}
and~\ref{lemma:shapes}, is taken from \cite{bonfioli-lean}.  So are the tilings
constructed in Theorems~\ref{theorem:triquadraticexistence},
\ref{theorem:existence2},
\ref{theorem:isoscelesbetaexistence}, \ref{theorem:Lconstruction},
\ref{theorem:isoscelesconstruction} and
Corollary~\ref{corollary:existence2case2}, and the tools used to glue tilings
together in Corollary~\ref{corollary:isoscelesalphaall} and
Theorem~\ref{theorem:existencePaliy}.  The tiling of
Theorem~\ref{theorem:existencePaliy} is assembled in our Lean files, as in the
printed proof, from the tiling of Theorem~\ref{theorem:triquadraticexistence}
and a quadratic tiling.  Our Lean files state these
results as printed here and derive the angles, side lengths and coloring numbers
of the constructed tilings.  The results marked $\dagger$ in
Table~\ref{table:lean} depend on the classification and are in a separate part of
the Lean code, which takes several hours to build; the README of the repository
explains how.

\smallskip

\begin{table}[htbp]
\small
\caption{Results of this paper and their Lean names (all in the namespace
{\tt TriangleTiling.ThreeAlpha}; $\dagger$: in the separately built part).}
\label{table:lean}
\begin{tabular}{ll@{\qquad}ll}
Result & Lean name & Result & Lean name \\
\hline
Lemma~\ref{lemma:angles} & {\tt angles}$^\dagger$ & Corollary~\ref{corollary:nokt} & {\tt nokt}$^\dagger$ \\
Lemma~\ref{lemma:singamma} & {\tt singamma} & Lemma~\ref{lemma:twocsapplies} & {\tt twocsapplies}$^\dagger$ \\
Lemma~\ref{lemma:anglesOK} & {\tt anglesOK} & Lemma~\ref{lemma:No14} & {\tt No14}$^\dagger$ \\
Lemma~\ref{lemma:gsquare} & {\tt gsquare} & Lemma~\ref{lemma:No17} & {\tt No17}$^\dagger$ \\
Lemma~\ref{lemma:gammagreaterpiover2} & {\tt gammagreaterpiover2} & Lemma~\ref{lemma:No7} & {\tt No7}$^\dagger$ \\
Lemma~\ref{lemma:nogamma} & {\tt nogamma}$^\dagger$ & Lemma~\ref{lemma:No23} & {\tt No23}$^\dagger$ \\
Lemma~\ref{lemma:pigeonhole} & {\tt pigeonhole} & Lemma~\ref{lemma:Nolessthan28} & {\tt Nolessthan28}$^\dagger$ \\
Lemma~\ref{lemma:twocs} & {\tt twocs} & Theorem~\ref{theorem:triquadraticexistence} & {\tt triquadraticexistence}$^\dagger$ \\
Lemma~\ref{lemma:cminusb} & {\tt cminusb} & Theorem~\ref{theorem:jshape} & {\tt jshape}$^\dagger$ \\
Lemma~\ref{lemma:standardvertices} & {\tt standardvertices}$^\dagger$ & Lemma~\ref{lemma:6902} & {\tt lemma6902}$^\dagger$ \\
Lemma~\ref{lemma:shapes} & {\tt shapes}$^\dagger$ & Lemma~\ref{lemma:normtheorem} & {\tt normtheorem}$^\dagger$ \\
Lemma~\ref{lemma:signedtiles} & {\tt signedtiles}$^\dagger$ & Lemma~\ref{lemma:tqprime} & {\tt tqprime}$^\dagger$ \\
Theorem~\ref{theorem:coloring} & {\tt coloring}$^\dagger$ & Lemma~\ref{lemma:primeimpliesjone} & {\tt primeimpliesjone} \\
Lemma~\ref{lemma:Mnotzero} & {\tt Mnotzero}$^\dagger$ & Lemma~\ref{lemma:singamma2} & {\tt singamma2} \\
Lemma~\ref{lemma:noboundarylink} & {\tt noboundarylink} & Theorem~\ref{theorem:secondtilingequation} & {\tt secondtilingequation}$^\dagger$ \\
Lemma~\ref{lemma:apropagation} & {\tt apropagation}$^\dagger$ & Lemma~\ref{lemma:Msquared} & {\tt Msquared}$^\dagger$ \\
Lemma~\ref{lemma:incommensurable} & {\tt incommensurable}$^\dagger$ & Lemma~\ref{lemma:scaling} & {\tt scaling}$^\dagger$ \\
Lemma~\ref{lemma:deltacoords} & {\tt deltacoords} & Theorem~\ref{theorem:existence2} & {\tt existence2}$^\dagger$ \\
Lemma~\ref{lemma:deltadirections} & {\tt deltadirections}$^\dagger$ & Theorem~\ref{theorem:existencePaliy} & {\tt existencePaliy}$^\dagger$ \\
Lemma~\ref{lemma:Phitile} & {\tt Phitile}$^\dagger$ & Theorem~\ref{theorem:secondcasesolution} & {\tt secondcasesolution}$^\dagger$ \\
Lemma~\ref{lemma:Phiproperties} & {\tt Phiproperties} & Theorem~\ref{theorem:notprime2} & {\tt notprime2}$^\dagger$ \\
Lemma~\ref{lemma:Phicancellation} & {\tt Phicancellation} & Theorem~\ref{theorem:Paliy} & {\tt Paliy}$^\dagger$ \\
Theorem~\ref{theorem:Mdef} & {\tt Mdef}$^\dagger$ & Corollary~\ref{corollary:everysolution} & {\tt everysolution}$^\dagger$ \\
Lemma~\ref{lemma:Phisigns} & {\tt Phisigns}$^\dagger$ & Theorem~\ref{theorem:solvingsecond} & {\tt solvingsecond}$^\dagger$ \\
Lemma~\ref{lemma:sines} & {\tt sines} & Theorem~\ref{theorem:tilingequation3} & {\tt tilingequation3}$^\dagger$ \\
Theorem~\ref{theorem:M1M2} & {\tt M1M2}$^\dagger$ & Theorem~\ref{theorem:notprime4} & {\tt notprime4}$^\dagger$ \\
Lemma~\ref{lemma:relationlattice} & {\tt relationlattice} & Theorem~\ref{theorem:isoscelesbetaexistence} & {\tt isoscelesbetaexistence}$^\dagger$ \\
Lemma~\ref{lemma:shapeareas} & {\tt shapeareas} & Lemma~\ref{lemma:XY} & {\tt XY}$^\dagger$ \\
Lemma~\ref{lemma:kappa} & {\tt kappa}$^\dagger$ & Theorem~\ref{theorem:isosceles1} & {\tt isosceles1}$^\dagger$ \\
Theorem~\ref{theorem:rationaliso} & {\tt rationaliso}$^\dagger$ & Corollary~\ref{corollary:shapedeterminestile} & {\tt shapedeterminestile}$^\dagger$ \\
Theorem~\ref{theorem:shape45} & {\tt shape45}$^\dagger$ & Theorem~\ref{theorem:notprime3} & {\tt notprime3}$^\dagger$ \\
Theorem~\ref{theorem:arelation} & {\tt arelation}$^\dagger$ & Theorem~\ref{theorem:Lconstruction} & {\tt Lconstruction}$^\dagger$ \\
Theorem~\ref{theorem:tqsquare} & {\tt tqsquare}$^\dagger$ & Lemma~\ref{lemma:complementary} & {\tt complementary} \\
Theorem~\ref{theorem:tqrational} & {\tt tqrational}$^\dagger$ & Theorem~\ref{theorem:isoscelesconstruction} & {\tt isoscelesconstruction}$^\dagger$ \\
Corollary~\ref{corollary:tqparameters} & {\tt tqparameters}$^\dagger$ & Corollary~\ref{corollary:isoscelesall} & {\tt isoscelesall}$^\dagger$ \\
Theorem~\ref{theorem:rationaltile} & {\tt rationaltile}$^\dagger$ & Corollary~\ref{corollary:existence2case2} & {\tt existence2case2}$^\dagger$ \\
Lemma~\ref{lemma:connectedcomponents} & {\tt connectedcomponents} & Corollary~\ref{corollary:secondshapeall} & {\tt secondshapeall}$^\dagger$ \\
Lemma~\ref{lemma:Hcomponents} & {\tt Hcomponents} & Lemma~\ref{lemma:mu} & {\tt mu}$^\dagger$ \\
Lemma~\ref{lemma:dmatrix} & {\tt dmatrix}$^\dagger$ & Lemma~\ref{lemma:halfintegermu} & {\tt halfintegermu}$^\dagger$ \\
Lemma~\ref{lemma:52} & {\tt lemma52}$^\dagger$ & Theorem~\ref{theorem:N45} & {\tt N45}$^\dagger$ \\
Theorem~\ref{theorem:tilingequation} & {\tt tilingequation}$^\dagger$ & Corollary~\ref{lemma:july29} & not formalized (computation) \\
Lemma~\ref{lemma:XYZ} & {\tt XYZ}$^\dagger$ & Theorem~\ref{theorem:tilingequation4} & {\tt tilingequation4}$^\dagger$ \\
Lemma~\ref{lemma:n<K} & {\tt nltK}$^\dagger$ & Theorem~\ref{theorem:notprime5} & {\tt notprime5}$^\dagger$ \\
Lemma~\ref{lemma:X=M(b+c)} & {\tt XeqMbc}$^\dagger$ & Corollary~\ref{corollary:isoscelesalphaall} & {\tt isoscelesalphaall}$^\dagger$ \\
Lemma~\ref{lemma:jmk} & {\tt jmk}$^\dagger$ & Theorem~\ref{theorem:notprime} & {\tt notprime}$^\dagger$ \\
Lemma~\ref{lemma:congruences} & {\tt congruences}$^\dagger$ &  &  \\
\end{tabular}
\end{table}
\FloatBarrier

\vskip 20pt
\noindent {\bf Acknowledgment.} The original version of this paper was
written without AI assistance.  In preparing this revision, the author
used Claude (Anthropic), an AI system, in three ways.  First, to check
the proofs: it found and corrected a gap in the proof of
Lemma~\ref{lemma:twocs}; it found an error in the proof of a former lemma on the
integrality of $\mu$, which led to its removal and to additional entries in
Table~\ref{table:isoscelesalphaplusbeta}, and it corrected the bound in
Theorem~\ref{theorem:Lconstruction}.
Second, it proposed new and shorter proofs of some results, notably
Lemmas~\ref{lemma:No14}--\ref{lemma:Nolessthan28}, and
Theorem~\ref{theorem:N45}.  Third, it wrote the Lean formalization described
in Appendix~\ref{section:lean}, whose proofs are checked by the Lean kernel.
The author checked all proposed changes to the mathematics, rejected or
modified some, checked that the Lean statements say what the theorems of the
paper say, and is responsible for the final text.  The main results are
unchanged.

\begingroup
\raggedright

\begin{thebibliography}{10}

\bibitem{beeson-isosceles}
Michael Beeson.
\newblock Tilings of an isosceles triangle.
\newblock {\em Integers}.
\newblock to appear; \href{https://arxiv.org/abs/1206.1974}{arXiv:1206.1974
  [math.MG]}.

\bibitem{beeson-noseven}
Michael Beeson.
\newblock No triangle can be cut into seven congruent triangles.
\newblock url={https://arxiv.org/abs/1811.09723}, 2019.
\newblock \href {http://arxiv.org/abs/1811.09723} {\path{arXiv:1811.09723}}.

\bibitem{beeson-rosu2026}
Michael Beeson and Grigore Ro\c{s}u.
\newblock Tiling a triangle into a prime number of congruent triangles, 2026.
\newblock \href {http://arxiv.org/abs/2607.23453} {\path{arXiv:2607.23453}}.

\bibitem{beeson-zhang2026}
Michael Beeson and Yan~X Zhang.
\newblock Rationality of certain triangle tilings, 2026.
\newblock \href {http://arxiv.org/abs/2604.01314} {\path{arXiv:2604.01314}}.

\bibitem{bonfioli-lean}
Vico Bonfioli.
\newblock Formalizing triangle tiling in Lean, version 1.0, 22 October 2026.
\newblock Lean~4 formalization, GitHub repository
  \texttt{ElVec1o/formalizing-triangle-tiling-in-lean}.
\newblock URL: \url{https://github.com/ElVec1o/formalizing-triangle-tiling-in-lean}.

\bibitem{bonfioli2026}
Vico Bonfioli.
\newblock A signed-direction invariant for triangle tilings, and the exclusion
  of primes congruent to $3$ modulo $4$.
\newblock GitHub repository \texttt{ElVec1o/erdos\_634\_proof}, 2026.
\newblock Accessed 28 August 2026.
\newblock URL: \url{https://github.com/ElVec1o/erdos_634_proof}.

\bibitem{buehl}
D.~A. Buehl.
\newblock {\em Binary Quadratic Forms}.
\newblock Springer-Verlag, New York, 1989.

\bibitem{golomb}
Solomon~W. Golomb.
\newblock Replicating figures in the plane.
\newblock {\em The Mathematical Gazette}, 48:403--412, 1964.

\bibitem{harries2026}
Jan~Philipp Harries.
\newblock New constructions, obstructions, and multiplier structure for
  {E}rd{\H{o}}s problem 634, 2026.
\newblock Version 0.5, August 28, 2026.
\newblock URL:
  \url{https://github.com/jphme/math-problems/tree/main/progress634}.

\bibitem{laczkovich1995}
M.~Laczkovich.
\newblock {Tilings of triangles}.
\newblock {\em Discrete Mathematics}, 140:79--94, 1995.

\bibitem{laczkovich2012}
Mikl{\'o}s Laczkovich.
\newblock Tilings of convex polygons with congruent triangles.
\newblock {\em Discrete and Computational Geometry}, 38:330--372, 2012.
\newblock \href {http://dx.doi.org/10.1007/s00454-012-9404-x}
  {\path{doi:10.1007/s00454-012-9404-x}}.
  
\bibitem{paliy2026}
Denis Paliy.
\newblock Erd\H{o}s Problem 634: exact arithmetic sectors, constructive spectra,
  and verified finite cases, version 0.3.0.
\newblock GitHub repository \texttt{DenisUkranian/erdos-634-tilings}, 2026.
\newblock See \texttt{docs/two-piece-construction.md} (29 September 2026) for the
  $(2\alpha,\alpha,2\beta)$ construction.

\bibitem{snover1991}
Stephen~L. Snover, Charles Waiveris, and John~K. Williams.
\newblock Rep-tiling for triangles.
\newblock {\em Discrete Mathematics}, 91(2):193--200, 1991.
\newblock \href {http://dx.doi.org/10.1016/0012-365X(91)90110-N}
  {\path{doi:10.1016/0012-365X(91)90110-N}}.

\bibitem{soifer}
Alexander Soifer.
\newblock {\em How Does One Cut a Triangle?}
\newblock Springer, 2009.

\end{thebibliography}
\bibliographystyle{plainurl}

\medskip

\end{document}